\documentclass[11pt,A4paper]{amsbook}
\usepackage{amsfonts,amsmath,amsthm,amssymb}
\usepackage[all]{xy}

\begin{document}

\newtheorem{thm}{Theorem}[chapter]
\newtheorem{lem}[thm]{Lemma}
\newtheorem{lemmadef}[thm]{Lemma-Definition}
\newtheorem{prop}[thm]{Proposition}
\newtheorem{cor}[thm]{Corollary}
\newtheorem{equat}[thm]{Equation}
\newtheorem{form}[thm]{Formula}
\newtheorem{claim}[thm]{Claim}

\theoremstyle{remark}
\newtheorem{rem}[thm]{Remark}
\newtheorem{rems}[thm]{Remarks}

\theoremstyle{definition}
\newtheorem{defn}[thm]{Definition}
\newtheorem{ex}[thm]{Example}
\newtheorem{notation}[thm]{Notation}
\newtheorem{eserci}[thm]{Exercise}

\def\thechapter{\Roman{chapter}} 
\def\thesubsection{\thesection-\Alph{subsection}}

\DeclareRobustCommand{\fineexer}{%
\ifmmode 
\else \leavevmode\unskip\penalty9999 \hbox{}\nobreak\hfill
\fi\quad\hbox{$\triangle$}}
\newenvironment{exer}{\begin{eserci}}{\fineexer\end{eserci}}

\newenvironment{step}[1]{{\sc Step #1:}}{\par\smallskip}
\newcommand{\quot}[2]{
{\lower-.2ex \hbox{$#1$}}{\kern -0.2ex /}
{\kern -0.5ex \lower.6ex\hbox{$#2$}}}

\newcommand{\mapor}[1]{{\stackrel{#1}{\longrightarrow}}}
\newcommand{\mapin}[1]{\smash{\mathop{\hookrightarrow}\limits^{#1}}}
\newcommand{\mapver}[1]{\Big\downarrow\vcenter{\rlap{$\scriptstyle#1$}}}
\newcommand{\liminv}{\smash{\mathop{\lim}\limits_{\leftarrow}\,}}

\newcommand{\pr}{^{\prime}}
\newcommand{\implica}[2]{{$[#1\Rightarrow#2]$}}
\newcommand{\specif}[2]{\left\{#1\,\left|\, #2\right. \,\right\}}
\newcommand{\vale}[1]{$[#1]$}
\newcommand{\h}[1]{{\widehat{#1}}}
\newcommand{\Til}[1]{\widetilde{#1}}
\newcommand{\ds}{\displaystyle}
\renewcommand{\bar}{\overline}
\newcommand{\de}{\partial}
\newcommand{\debar}{{\overline{\partial}}}
\newcommand{\somdir}[2]{\hbox{$\mathrel
{\smash{\mathop{\mathop \bigoplus\limits_{#1}}
\limits^{#2}}}$}}
\newcommand{\tensor}[2]{\hbox{$\mathrel
{\smash{\mathop{\mathop \bigotimes\limits_{#1}}
^{#2}}}$}}
\newcommand{\symm}[2]{\hbox{$\mathrel
{\smash{\mathop{\mathop \bigodot\limits_{#1}}
^{#2}}}$}}
\newcommand{\external}[2]{\hbox{$\mathrel
{\smash{\mathop{\mathop \bigwedge\limits_{#1}}
^{\!#2}}}$}}

\newcommand{\dual}{^{\vee}}
\newcommand{\OP}[2]{\Oh_{\proj^{#1}}(#2)}
\newcommand{\desude}[2]{\dfrac{\de #1}{\de #2}}
\newcommand{\decala}{{\operatorname{dp}}}
\newcommand{\Aut}{{\operatorname{Aut}}}
\newcommand{\Mor}{{\operatorname{Mor}}}
\newcommand{\Def}{{\operatorname{Def}}}
\newcommand{\Hom}{{\operatorname{Hom}}}
\newcommand{\HOM}{\operatorname{\mathcal H}\!\!om}
\newcommand{\DER}{\operatorname{\mathcal D}\!er}
\newcommand{\Spec}{{\operatorname{Spec}}}
\newcommand{\rank}{{\operatorname{rank}}}
\newcommand{\Der}{{\operatorname{Der}}}
\newcommand{\Coder}{{\operatorname{Coder}}}
\newcommand{\coker}{{\operatorname{coker}}}
\newcommand{\KS}{{\operatorname{KS}}}
\newcommand{\End}{{\operatorname{End}}}
\newcommand{\END}{\operatorname{\mathcal E}\!\!nd}
\newcommand{\Image}{\operatorname{Im}}
\newcommand{\Ext}{{\operatorname{Ext}}}

\newcommand{\Set}{\mathbf{Set}}
\newcommand{\Art}{\mathbf{Art}}
\newcommand{\NA}{\mathbf{NA}}

\newcommand{\Q}{{\mathbb Q}}
\newcommand{\F}{{\mathbb F}}
\newcommand{\C}{{\mathbb C}}
\newcommand{\A}{{\mathbb A}}
\newcommand{\N}{{\mathbb N}}
\newcommand{\R}{{\mathbb R}}
\newcommand{\Z}{{\mathbb Z}}
\renewcommand{\H}{{\mathbb H}}
\newcommand{\proj}{{\mathbb P}}
\newcommand{\K}{{\mathbb K}\,}

\newcommand{\ide}[1]{{\mathfrak #1}}
\newcommand{\copa}{{\mathfrak a}}
\newcommand{\copl}{{\mathfrak l}}
\newcommand{\copc}{{\mathfrak c}}

\newcommand{\script}[1]{{\mathcal #1}}
\newcommand{\Oh}{{\mathcal O}}
\newcommand{\sA}{{\mathcal A}}
\newcommand{\sB}{{\mathcal B}}
\newcommand{\sC}{{\mathcal C}}
\newcommand{\sF}{{\mathcal F}}
\newcommand{\sG}{{\mathcal G}}
\newcommand{\sH}{{\mathcal H}}
\newcommand{\sI}{{\mathcal I}}
\newcommand{\sKS}{{\mathcal KS}}
\newcommand{\sL}{{\mathcal L}}
\newcommand{\sM}{{\mathcal M}}
\newcommand{\sP}{{\mathcal P}}
\newcommand{\sT}{{\mathcal T}}
\newcommand{\sU}{{\mathcal U}}
\newcommand{\sV}{{\mathcal V}}
\newcommand{\sX}{{\mathcal X}}
\newcommand{\sY}{{\mathcal Y}}

\newcommand{\piede}{\insert\footins
{\small{\sc Marco Manetti:} \emph{Lectures on 
deformations of complex manifolds}}
}

{\pagestyle{empty}
\frontmatter
\title{\textbf{Lectures on deformations of complex manifolds
\small{(Deformations from differential graded viewpoint)}
}}

\author{Marco Manetti}

\address{Dipartimento di Matematica "G. Castelnuovo",
Universit\`a di Roma ``La Sapienza'', Piazzale Aldo Moro 5,
I-00185 Roma, Italy.}
\urladdr{http://www.mat.uniroma1.it/people/manetti/}
\email{manetti@mat.uniroma1.it}
\date{Spring 2004}

\maketitle

\chapter*{Preface}

This paper is based on a  course  given at the ``Dottorato di
Ricerca in Matematica'' of  the University of Rome ``La Sapienza'' in
the Academic year 2000/2001.\\
The intended aim of the course was to  rapidly introduce,  although not in an
exhaustive way,  the non-expert PhD student
to deformations of compact complex manifolds, from the very beginning to
some recent (i.e. at that time not yet published) results.\\

With the  term "deformation theory", the mathematicians usually intend a 
set of deformation theories, each one of which  
studies small parametric variation of a specific 
mathematical structure, for example: deformation theory of complex manifolds, 
deformation theory of associative algebras,
deformations of schemes, deformations of representations and much more.\\
Every mathematician which tries to
explain and  investigate  deformation theory has to
deal with two opposite features: order and chaos.\\
{\sc Chaos:} the various deformation theories often rely on
theorems which are proved using very different tools, from
families of elliptic differential operators of Kodaira and Spencer
\cite{KS3} to  ringed toposes of Illusie \cite{Illu}.\\
{\sc Order:} all the deformation theories have lots of common
features; for instance they have a vector space of first order
deformations (usually the $H^{1}$ of some complex) and they have
an obstruction space (usually an $H^{2}$).\\
Another unifying aspect of all deformation theories is summarized in the 
slogan ``In
characteristic 0 every deformation problem is governed by a
differential graded Lie algebra'', which underlie some ideas given, mostly in private 
communications, by Quillen, Deligne, Drinfeld and other about 20 years ago.
More recently  (especially in \cite{K3} and 
\cite{K})  these 
ideas have shown a great utility and possibility of development.\\
Nowadays this approach to deformation theory is a very active 
area of research which is usually called \emph{deformation theory via 
DGLA} or \emph{extended deformation theory}.\\

The goal of these notes is to give a soft introduction to extended 
deformation theory. In view of the  aim (and the hope) 
of keeping  this paper selfcontained,  user 
friendly and with a tolerating number of pages, 
we consider only  deformations of compact complex manifolds. Anyhow, 
most part of the formalism and of the results  that we prove here will apply 
to many other deformation problems.\\    
The first part of the paper ({\chaptername}s \ref{CAP:FAMILIES}, \ref{CAP:SEGRE} 
and \ref{CAP:SING}) is a classical introduction to deformations of 
compact complex manifolds; 
the beginners can find here the main definitions, the statements of 
the theorems of Kodaira and Kodaira-Nirenberg-Spencer, an elementary 
description of the semiuniversal deformations of Segre-Hirzebruch 
surfaces and a micro-course in complex analytic singularity theory.\\ 
In the second part (From \chaptername\ \ref{CAP:INFINITESIMAL} to   
\chaptername\ \ref{CAP:TRIVIALK}) we  study deformations in 
the context of dg-objects, where by dg-objects we intend  algebraic 
structures supported on  differential $\Z$-graded vector spaces.\\
Most of this part is devoted to introduce some new objects  which play 
a fundamental role in extended deformation theory, such as for instance: 
deformation functors associated to a differential graded Lie algebra 
and their homotopy invariance, extended deformation functors and 
Gerstenhaber-Batalin-Vilkoviski algebras. The reader of this part 
can also find satisfaction in the  proof of the 
unobstructness of Calabi-Yau manifolds (theorem of 
Bogomolov-Tian-Todorov).\\ 
\chaptername\ \ref{CAP:KAEHLER} is a basic introduction to K\"ahler 
manifold which follows essentially Weil's book \cite{Weil}: some 
modification in the presentation and 
simplification in the proofs are made by using the formalism of 
dg-vector spaces; this partially explain the reason why this 
\chaptername\ is contained in part II of these notes.\\
The third part of the notes ({\chaptername}s \ref{CAPITOLO13} and 
\ref{CAP:TOOLS}) is a basic course in $L_{\infty}$-algebras and their use 
in deformation theory: a nontrivial application of $L_{\infty}$-algebras in made 
in the last section where we give (following \cite{ManettiCCKM})
an algebraic proof of Clemens-Ran theorem ``obstructions to 
deformations annihilate ambient cohomology''.\\

Each \chaptername\ contains: a brief introduction, the main matter, some exercises 
and a survey section. The main matter is organized like a book, while
the survey sections contain bibliographical annotations and 
theorems for which the proof it is not given here.\\

\tableofcontents
\mainmatter}


\chapter[~Smooth families of compact complex manifolds]{Smooth
families of compact
complex manifolds}
\label{CAP:FAMILIES}
\piede

In this chapter we introduce the notion of a family $f\colon\sX\to
B$ of compact complex manifolds as a proper holomorphic submersion
of complex manifolds. Easy examples (\ref{I.1.4}, \ref{I.1.6})
will show that in general the fibres $X_t:=f^{-1}(t)$ are not
biholomorphic each other. Using integration of vector fields we
prove that the family is locally trivial if and only if a certain
morphism $\sKS$ of sheaves over $B$ is trivial, while the
restriction of $\sKS$ at a point $b\in B$ is a linear map
$\KS\colon T_{b,B}\to H^1(X_b,T_{X_b})$, called the 
\emph{Kodaira-Spencer map}, which can interpreted as the first
derivative at the point $b$ of the map
\[B\to\{\hbox{isomorphism classes of complex manifolds}\},\quad t\mapsto
X_t.\]

Then, according to Kodaira, Nirenberg and Spencer we define a 
\emph{deformation} of a complex manifolds $X$ as the data of a family
$\sX\to B$, of a base point $0\in B$ and of an isomorphism
$X\simeq X_0$. The isomorphism class of a deformation involves only
the structure of $f$ in a neighbourhood of $X_0$.\\

In the last section we state, without proof, the principal
pioneer theorems about deformations proved using hard analysis
by Kodaira, Nirenberg and Spencer in the period 1956-58.

\bigskip

\section{Dictionary}

For every complex manifold $M$
we denote by:\begin{itemize}
\item $\Oh_{M}(U)$ the $\C$-algebra of holomorphic functions $f\colon
U\to \C$ defined
on an open subset $U\subset M$.
\item $\Oh_{M}$ the trivial complex line bundle $\C\times M\to M$.
\item $T_{M}$ the holomorphic tangent bundle to $M$.
The fibre of $T_{M}$ at a point $x\in M$, i.e. the complex tangent
space at $x$,  is denoted by
$T_{x,M}$.
\end{itemize}

If $x\in M$ is a point we denote by
$\Oh_{M,x}$  the $\C$-algebra of germs of holomorphic
functions at a point $x\in M$; a choice of local holomorphic coordinates
$z_{1},\ldots,z_{n}$, $z_{i}(x)=0$, gives an isomorphism
$\Oh_{M,x}=\C\{z_{1},\ldots,z_{n}\}$, being $\C\{z_{1},\ldots,z_{n}\}$
the $\C$-algebra of convergent power series.\\

In order to avoid a too heavy  notation we sometimes omit the
subscript $M$, when the underlying complex manifold is clear from the
context.\\

\begin{defn}\label{I.1.1}
A \emph{smooth family of compact complex manifolds} is a proper
holomorphic map $f\colon M\to B$ such that:\begin{enumerate}
\item $M,B$ are nonempty complex manifolds and $B$ is connected.
\item The differential of $f$, $f_{*}\colon T_{p,M}\to T_{f(p),B}$ is
surjective at every point $p\in M$.
\end{enumerate}
Two families $f\colon M\to B$, $g\colon N\to B$ over the same base
are isomorphic if there exists a holomorphic isomorphism $N\to M$
commuting with $f$ and $g$.
\end{defn}

From now on, when there is no risk of confusion, we shall simply say
\emph{smooth family} instead of
smooth family of compact complex manifolds.

Note that if $f\colon M\to B$ is a smooth family then $f$ is open,
closed and surjective. If $V\subset B$ is an open subset then
$f\colon f^{-1}(V)\to V$ is a smooth family; more generally for every
holomorphic map of connected complex manifolds $C\to B$, the
pull-back $M\times_{B}C\to C$ is a smooth family.

For every $b\in B$ we denote  $M_{b}=f^{-1}(b)$: $M_{b}$ is a regular
submanifold of $M$.

\begin{defn}\label{I.1.2}
A smooth family $f\colon M\to B$ is called \emph{trivial} if it
is isomorphic to the product $M_{b}\times B\to B$ for some (and hence
all) $b\in B$.  It is called \emph{locally trivial} if there exists an
open covering $B=\cup U_{a}$ such that every restriction $f\colon
f^{-1}(U_{a})\to U_{a}$ is trivial.\end{defn}

\begin{lem}\label{I.1.3}
Let $f\colon M\to B$ be a smooth family, $b\in B$.
The normal bundle $N_{M_{b}/M}$ of $M_{b}$ in $M$ is
trivial.\end{lem}

\begin{proof} Let $E=T_{b,B}\times M_{b}\to M_{b}$ be the trivial bundle
with fibre $T_{b,B}$. The differential $f_{*}\colon T_{x,M}\to
T_{b,B}$, $x\in M_{b}$ induces a surjective morphism of vector bundles
$(T_{M})_{|M_{b}}\to E$ whose kernel is exactly $T_{M_{b}}$.\\
By definition $N_{M_{b}/M}=(T_{M})_{|M_{b}}/T_{M_{b}}$ and then
$N_{M_{b}/M}=T_{b,B}\times M_{b}$.\end{proof}

By a classical result (Ehresmann's theorem, \cite[Thm. 2.4]{Kobook}), if $f\colon
M\to B$ is a family, then for every $b\in B$ there exists an open
neighbourhood $b\in U\subset B$ and a diffeomorphism $\phi\colon
f^{-1}(U)\to M_{b}\times U$ making the following diagram commutative
\[\xymatrix{&M_{b}\ar[dl]_{i}\ar[dr]^{Id\times\{b\}}&\\
f^{-1}(U)\ar[rr]^{\phi}\ar[dr]_{f}&&M_{b}\times U\ar[dl]^{p_{2}}\\
&U&}\]
being $i\colon M_{b}\to M$ the inclusion.
In particular the diffeomorphism type of the fibre $M_{b}$ is
independent from $b$. Later on (Theorem~\ref{VI.2.1}) we will prove
a result that implies Ehresmann's theorem.\\

The following examples of families show that, in general, if $a,b\in
B$, $a\not=b$, then $M_{a}$ is not biholomorphic to $M_{b}$.

\begin{ex}\label{I.1.4}
Consider $B=\C-\{0,1\}$,
\[ M=\{([x_{0},x_{1},x_{2}],\lambda)\in\proj^{2}\times B\mid
x_{2}^{2}x_{0}=x_{1}(x_{1}-x_{0})(x_{1}-\lambda x_{0})\},\]
and $f\colon M\to B$ the projection. Then $f$ is a family and the
fibre $M_{\lambda}$ is a smooth plane cubic with $j$-invariant
\[ j(M_{\lambda})=2^{8}\frac{(\lambda^{2}-\lambda+1)^{3}}
{\lambda^{2}(\lambda-1)^{2}}.\]
(Recall that two elliptic curves are biholomorphic if and only if they
have the same $j$-invariant.)
\end{ex}

\begin{ex}\label{I.1.5} (Universal family of hypersurfaces)\\
For fixed integers $n,d>0$, consider the projective space $\proj^{N}$,
$N=\ds\left(\begin{array}{c}d+n\\ n\end{array}\right)-1$, with
homogeneous coordinates $a_{i_{0},\ldots,i_{n}}$, $i_{j}\ge 0$,
$\sum_{j} i_{j}=d$, and denote
\[ X=\left\{([x],[a])\in\proj^{n}\times\proj^{N}\,\left|\,\,
\sum_{i_{0}+\ldots+i_{n}=d}a_{i_{0},\ldots,i_{n}}x_{0}^{i_{0}}\ldots
x_{n}^{i_{n}}=0\right.\right\}.\]
$X$ is a smooth hypersurface of $\proj^{n}\times\proj^{N}$, the
differential of the projection $X\to \proj^{N}$ is not surjective at
a point $([x],[a])$ if and only if $[x]$ is a singular point of
$X_{a}$.\\
Let $B=\{[a]\in\proj^{N}\mid X_{a} \hbox{ is smooth }\}$,
$M=f^{-1}(B)$: then $B$ is open (exercise), $f\colon M\to B$ is a family and every smooth
hypersurface of degree $d$ of $\proj^{n}$ is isomorphic to a fibre of
$f$.\end{ex}

\begin{ex}\label{I.1.6} (Hopf surfaces)\\
Let $A\in GL(2,\C)$ be a matrix with eigenvalues of norm
$>1$ and let $\langle A\rangle\simeq\Z\subset GL(2,\C)$ be
the subgroup generated by $A$.
The action of $\langle A\rangle$ on $X=\C^{2}-\{0\}$ is free and properly
discontinuous: in fact
a linear change of coordinates $C\colon \C^{2}\to \C^{2}$ changes the
action of $\langle A\rangle$ into the action of
$\langle C^{-1}AC\rangle$ and  therefore it is
not restrictive to assume $A$ is a lower triangular matrix.\\
Therefore the quotient $S_{A}=X/\langle A\rangle$
is a compact complex manifold
called \emph{Hopf surface}: the holomorphic map $X\to S_{A}$ is the
universal cover and then for every point $x\in S_{A}$ there exists a
natural isomorphism $\pi_{1}(S_{A},x)\simeq \langle A \rangle$.
We have already seen that if $A,B$ are conjugated matrix then $S_{A}$
is biholomorphic to $S_{B}$. Conversely if $f\colon S_{A}\to S_{B}$
is a biholomorphism then $f$ lifts to a biholomorphism
$g\colon X\to X$ such that $gA=B^{k}g$; since $f$ induces an
isomorphism of fundamental groups $k=\pm 1$.\\
By Hartogs' theorem $g$ extends to a biholomorphism
$g\colon\C^{2}\to\C^{2}$ such that $g(0)=0$;
since for every $x\not=0$
$\ds\lim_{n\to\infty}A^{n}(x)=+\infty$ and $\ds\lim_{n\to\infty}B^{-n}(x)=0$
it must be $gA=Bg$.
Taking the differential at $0$ of $gA=Bg$ we get that $A$ is
conjugated to $B$.
\end{ex}

\begin{exer} If $A=e^{2\pi i\tau}I\in GL(2,\C)$, $\tau=a+ib$, $b<0$,
then the Hopf surface
$S_{A}$ is the total space of a holomorphic $G$-principal bundle
$S_{A}\to\proj^{1}$, where $G=\C/(\Z+\tau\Z)$.\end{exer}

\begin{ex}\label{I.1.7} (Complete family of Hopf surfaces)\\
Denote $B=\{(a,b,c)\in \C^{3}\mid |a|>1,\, |c|>1\}$, $X=B\times
(\C^{2}-\{0\})$ and let $\Z\simeq G \subset \Aut(X)$ be the subgroup
generated by
\[ (a,b,c,z_{1},z_{2})\mapsto (a,b,c,az_{1}, bz_{1}+cz_{2})\]
The action of $G$ on $X$ is free and properly discontinuous, let
$M=X/G$ be its quotient and $f\colon M\to B$ the projection on the first
coordinates: $f$ is a family whose fibres are Hopf surfaces.
Every Hopf surface is isomorphic to a fibre of $f$, this motivate the
adjective ``complete''.\\
In particular all the Hopf surfaces are diffeomorphic to $S^{1}\times
S^{3}$ (to see this look at the fibre over $(2,0,2)$).
\end{ex}

\begin{notation}
For every pair of pointed manifolds $(M,x)$, $(N,y)$ we denote by
$\Mor_{\mathbf{Ger}}((M,x),(N,y))$ the set of germs of holomorphic maps $f\colon
(M,x)\to (N,y)$. Every element of $\Mor_{\mathbf{Ger}}((M,x),(N,y))$
is an equivalence class
of pairs $(U,f)$, where $x\in U\subset M$ is an open neighbourhood of
$x$, $f\colon U\to N$ is a holomorphic map such that $f(x)=y$ and
$(U,f)\sim (V,g)$ if and only if there exists an open subset $x\in
W\subset U\cap V$ such that $f_{|W}=g_{|W}$.\\
The category $\mathbf{Ger^{sm}}$ of germs of complex manifolds is the category whose
object are the pointed complex manifold $(M,x)$ and the morphisms are
the $\Mor_{\mathbf{Ger}}((M,x),(N,y))$ defined above.
A germ of complex manifold is
nothing else that an object of $\mathbf{Ger^{sm}}$.
\end{notation}

In \chaptername~\ref{CAP:SING}
we will consider $\mathbf{Ger^{sm}}$ as a full subcategory of the
category of analytic singularities $\mathbf{Ger}$.\\

\begin{exer} $\mathbf{Ger^{sm}}$ is equivalent to its full subcategory whose
objects are $(\C^{n},0)$, $n\in\N$.\end{exer}

Roughly speaking a deformation is a ``framed germ'' of family;
more precisely

\begin{defn}\label{I.1.8} Let $(B,b_{0})$ be a pointed
manifold, a \emph{deformation} $M_{0}\mapor{i}M\mapor{f}(B,b_{0})$
of a compact complex manifold $M_{0}$
over $(B,b_{0})$ is a
pair of holomorphic maps
\[ M_{0}\mapor{i}M\mapor{f}B\]
such that:\begin{enumerate}
\item $fi(M_{0})=b_{0}$.
\item There exists an open neighbourhood $b_{0}\in U\subset B$ such
that $f\colon f^{-1}(U)\to U$ is a proper smooth family.
\item $i\colon M_{0}\to f^{-1}(b_{0})$ is an isomorphism of complex
manifolds.\end{enumerate}
$M$ is called the total space of the deformation and $(B,b_{0})$
the base germ space.\end{defn}

\begin{defn}\label{I.1.9}
Two deformations of $M_{0}$ over the same base
\[ M_{0}\mapor{i}M\mapor{f}(B,b_{0}),\qquad
M_{0}\mapor{j}N\mapor{g}(B,b_{0})\]
are isomorphic if there exists an open neighbourhood $b_{0}\in U\subset
B$,  and a commutative diagram of holomorphic maps
\[\xymatrix{M_{0}\ar[r]^{i}\ar[d]_{j}&f^{-1}(U)\ar[d]^{f}\ar[dl]\\
g^{-1}(U)\ar[r]_{g}&U}\]
with the diagonal arrow  a holomorphic  isomorphism.\end{defn}

For every pointed complex manifold $(B,b_{0})$ we denote by
$\Def_{M_{0}}(B,b_{0})$ the set of isomorphism classes of deformations
of $M_{0}$ with base  $(B,b_{0})$.
It is clear from the definition that
if $b_{0}\in U\subset B$ is open, then
$\Def_{M_{0}}(B,b_{0})=\Def_{M_{0}}(U,b_{0})$.

\begin{exer}
There exists an action of the group $\Aut(M_{0})$ of holomorphic
isomorphisms of $M_{0}$ on the set $\Def_{M_{0}}(B,b_{0})$: if $g\in
\Aut(M_{0})$ and $\xi:\, M_{0}\mapor{i}M\mapor{f}(B,b_{0})$ is a deformation
we define
\[\xi^{g}:\, M_{0}\mapor{ig^{-1}}M\mapor{f}(B,b_{0}).\]
Prove that
$\xi^{g}=\xi$ if and only if $g\colon f^{-1}(b_{0})\to f^{-1}(b_{0})$ can
be extended to an isomorphism $\hat{g}\colon f^{-1}(V)\to
f^{-1}(V)$, $b_{0}\in V$ open neighbourhood, such that $f\hat{g}=f$.
\end{exer}

If $\xi:\, M_{0}\mapor{i}M\mapor{f}(B,b_{0})$ is a deformation
and $g\colon (C,c_{0})\to (B,b_{0})$ is a holomorphic map of pointed
complex manifolds then
\[g^{*}\xi:\, M_{0}\mapor{(i,c_{0})}M\times_{B}C\mapor{pr}(C,c_{0})\]
is a deformation with base point $c_{0}$. It is clear that the
isomorphism class of $g^{*}\xi$ depends only by the class of $g$ in
$\Mor_{\mathbf{Ger}}((C,c_{0}),(B,b_{0}))$.\\
Therefore
every  $g\in \Mor_{\mathbf{Ger}}((C,c_{0}),(B,b_{0}))$ induces a well
defined pull-back morphism
\[g^{*}\colon \Def_{M_{0}}(B,b_{0})\to \Def_{M_{0}}(C,c_{0}).\]

\bigskip

\section{Dolbeault cohomology}

If $M$ is a complex manifold and $E$ is a holomorphic vector bundle
on $M$,  we denote:
\begin{itemize}
\item $E\dual$ the dual bundle of $E$.
\item $\Gamma(U,E)$ the space of holomorphic sections $s\colon U\to E$
on an open subset $U\subset M$.
\item $\Omega^{1}_{M}=T_{M}\dual$ the holomorphic cotangent bundle of $M$.
\item $\Omega^{p}_{M}=\bigwedge^{p} T_{M}\dual$ the bundle of
holomorphic differential $p$-forms.
\end{itemize}

For every open subset $U\subset M$ we denote by
$\Gamma(U,\sA^{p,q}_{M})$ the $\C$-vector space of differential
$(p,q)$-forms on $U$.
If $z_{1},\ldots,z_{n}$ are local holomorphic
coordinates,
then $\phi\in \Gamma(U,\sA^{p,q}_{M})$ is written locally as
$\phi=\sum\phi_{I,J}dz_{I}\wedge d\bar{z}_{J}$, where
$I=(i_{1},\ldots,i_{p})$, ${J}=(j_{1},\ldots,j_{q})$,
$dz_{I}=dz_{i_{1}}\wedge\ldots \wedge dz_{i_{p}}$,
$d\bar{z}_{{J}}=d\bar{z}_{j_{1}}\wedge\ldots \wedge d\bar{z}_{j_{q}}$
and the $\phi_{I,{J}}$ are $C^{\infty}$ functions.\\
Similarly, if $E\to M$ is a holomorphic vector bundle we denote by
$\Gamma(U,\sA^{p,q}(E))$
the space of differential $(p,q)$-forms on $U$ with value in
$E$; locally, if $e_{1},\ldots,e_{r}$ is a local frame for $E$, an
element of $\Gamma(U,\sA^{p,q}(E))$ is written as
$\sum_{i=1}^{r}\phi_{i}e_{i}$, with $\phi_{i}\in
\Gamma(U,\sA^{p,q})$.
Note that there exist natural isomorphisms
$\Gamma(U,\sA^{p,q}(E))\simeq\Gamma(U,\sA^{0,q}(\Omega^{p}_{M}\otimes E))$.

We begin recalling the well known
\begin{lem}[Dolbeault's lemma]\label{I.2.1} Let
\[\Delta^{n}_{R}=\{(z_{1},\ldots,z_{n})\in\C^{n}\mid |z_{1}|<R,\ldots,
|z_{n}|<R\}\]
be a polydisk of radius $R\le+\infty$ ($\Delta^{n}_{+\infty}=\C^{n}$) and
let $\phi\in \Gamma(\Delta^{n}_{R},\sA^{p,q})$,
$q>0$, such that $\bar{\de}\phi=0$.
Then there exists $\psi\in \Gamma(\Delta^{n}_{R},\sA^{p,q-1})$ such that
$\bar{\de}\psi=\phi$.\end{lem}
\begin{proof} \cite[Thm. 3.3]{Kobook}, \cite[pag. 25]{G-H}.\end{proof}

If $E$ is a holomorphic vector bundle,
the $\bar{\de}$ operator extends naturally to the ~ \emph{Dolbeault
operator} $\bar{\de}\colon \Gamma(U,\sA^{p,q}(E))\to
\Gamma(U,\sA^{p,q+1}(E))$ by the rule
$\bar{\de}(\sum_{i}\phi_{i}e_{i})=\sum_{i}(\bar{\de}\phi_{i})e_{i}$.
If $h_{1},\ldots,h_{r}$ is another local frame of $E$ then there
exists a matrix $(a_{ij})$ of holomorphic functions such that
$h_{i}=\sum_{j}a_{ij}e_{j}$ and then
\[\bar{\de}\left(\sum_{i}\phi_{i}h_{i}\right)=
\bar{\de}\left(\sum_{i,j}\phi_{i}a_{ij}e_{j}\right)=
\sum_{i,j}\bar{\de}(\phi_{i}a_{ij})e_{j}=
\sum_{i}(\bar{\de}\phi_{i})a_{ij}e_{j}=\sum_{i}(\bar{\de}\phi_{i})h_{i}.\]
It is obvious that $\bar{\de}^{2}=0$.

\begin{defn}\label{I.2.2}
The Dolbeault's cohomology of $E$,
$H^{p,*}_{\bar{\de}}(U,E)$  is the cohomology  of the complex
\[0\mapor{}\Gamma(U,\sA^{p,0}(E))\mapor{\bar{\de}}
\Gamma(U,\sA^{p,1}(E))\mapor{\debar}\ldots
\mapor{\debar}\Gamma(U,\sA^{p,q}(E))\mapor{\debar}\ldots\]
\end{defn}
Note that $H^{p,0}_{\bar{\de}}(U,E)=\Gamma(U,\Omega^{p}_{M}\otimes E)$
is the space of holomorphic sections.\\
The Dolbeault cohomology has several functorial properties; the most
relevant are:\begin{enumerate}
\item Every holomorphic morphism of holomorphic vector bundles $E\to F$ induces a
morphism of complexes $\Gamma(U,\sA^{p,*}(E))\to \Gamma(U,\sA^{p,*}(F))$ and then
morphisms of cohomology groups $H^{p,*}_{\bar{\de}}(U,E)\to
H^{p,*}_{\bar{\de}}(U,F)$.
\item The wedge product
\[\Gamma(U,\sA^{p,q}(E))\otimes \Gamma(U,\sA^{r,s}(F))\mapor{\wedge}
\Gamma(U,\sA^{p+r,q+s}(E\otimes F)),\]
\[\left(\sum \phi_{i}e_{i}\right)\otimes
\left(\sum \psi_{j}f_{j}\right)\to
\sum \phi_{i}\wedge \psi_{j}e_{i}\otimes e_{j}.\]
commutes with Dolbeault differentials and then induces a \emph{cup}
product
\[\cup\colon H^{p,q}_{\debar}(U,E)\otimes
H^{r,s}_{\debar}(U,F)\to
H^{p+r,q+s}_{\debar}(U,E\otimes F).\]
\item The composition of the wedge product with the trace map $E\otimes
E\dual\to \Oh_{M}$ gives
bilinear morphisms of cohomology groups
\[\cup\colon H^{p,q}_{\debar}(U,E)\times
H^{r,s}_{\debar}(U,E\dual)\to
H^{p+r,q+s}_{\debar}(U,\Oh_{M}).\]\end{enumerate}

\begin{thm}\label{I.2.3}
If $M$ is a compact complex manifold of dimension $n$
and $E\to M$ is a
holomorphic vector bundle then for every $p,q\ge 0$:\begin{enumerate}
\item $\dim_{\C}H^{p,q}_{\debar}(M,E)<\infty$.
\item \emph{(Serre's duality)}
The bilinear map $\Gamma(M,\sA^{p,q}(E))\times
\Gamma(M,\sA^{n-p,n-q}(E\dual))\to \C$,
\[(\phi,\psi)\mapsto \int_{M}\phi\wedge \psi\]
induces a perfect pairing $H^{p,q}_{\debar}(M,E)\times
H^{n-p,n-q}_{\debar}(M,E\dual)\to \C$ and then an isomorphism
$H^{p,q}_{\debar}(M,E)\dual\simeq
H^{n-p,n-q}_{\debar}(M,E\dual)$.\end{enumerate}
\end{thm}

\begin{proof}\cite{Kobook}.\end{proof}

From now on we denote for simplicity
$H^{q}(M,E)=H^{0,q}_{\debar}(M,E)$, $h^{q}(M,E)=\dim_{\C}H^{q}(M,E)$,
$H^{q}(M,\Omega^{p}(E))=H^{p,q}_{\debar}(M,E)$.\\

\begin{defn}\label{I.2.4} If $M$ is a complex manifold of dimension $n$, the
holomorphic line bundle $K_{M}=\bigwedge^{n}T_{M}\dual
=\Omega^{n}_{M}$ is
called the \emph{canonical bundle} of $M$.\end{defn}
Since $\Omega^{p}_{M}=
K_{M}\otimes (\Omega^{n-p}_{M})\dual$,
an equivalent statement of the Serre's
duality is $H^{p}(M,E)\dual\simeq H^{n-p}(M, K_{M}\otimes E\dual)$
for every holomorphic vector bundle $E$ and every $p=0,\ldots,n$.\\

The \emph{Hodge numbers} of a fixed compact complex manifold $M$ are by
definition
\[ h^{p,q}=\dim_{\C}H^{p,q}_{\debar}(M,\Oh)=
\dim_{\C}H^{0,q}_{\debar}(M,\Omega^{p}).\]
The \emph{Betti numbers} of $M$ are the dimensions of the spaces of the
De Rham cohomology of $M$, i.e.
\[ b_{p}=\dim_{\C}H^{p}_{d}(M,\C),\quad H^{p}_{d}(M,\C)=
\frac{\hbox{$d$-closed $p$-forms}}{\hbox{$d$-exact $p$-forms}}.\]

\begin{exer}\label{Hodgeinequality}
Let $p\ge 0$ be a fixed integer and, for every $0\le q\le
p$,  denote by  $F_{q}\subset \quad H^{p}_{d}(M,\C)$ the subspace of
cohomology classes represented by a $d$-closed form
$\eta\in\oplus_{i\le q}\Gamma(M,\sA^{p-i,i})$.
Prove that there exist injective linear morphisms
$F_{q}/F_{q-1}\to H^{p-q,q}_{\debar}(M,\Oh)$. Deduce that $b_{p}\le
\sum_{q}h^{p-q,q}$.
\end{exer}

\begin{exer}\label{ese.affinedolbeault}
Let $f\colon \C^{n}\to \C$ be a holomorphic function and
assume that $X=f^{-1}(0)$ is a regular smooth submanifold; denote $i\colon
X\to \C^{n}$ the embedding.\\
Let $\phi\in \Gamma(\C^{n},\sA^{p,q})$, $q>0$, be a differential form such that
$\debar\phi=0$ in an open  neighbourhood of $X$. Prove that
$i^{*}\phi$ is $\debar$-exact in $X$. (Hint: prove that there exists
$\psi\in \Gamma(\C^{n},\sA^{p,q})$ such that
$\debar\phi=\debar(f\psi)$.)\end{exer}

\begin{exer}
Let $h\colon \C^{n}\to \C$ be holomorphic and let
$U=\{z\in \C^{n}\mid h(z)\not=0\}$.
Prove that $H^{q}(U,\Oh_{U})=0$
for every $q>0$. (Hint: consider the open disk
$\Delta=\{t\in\C\mid |t|<1\}$ and the holomorphic
maps $\phi\colon U\times\Delta\to \C^{n+1}$, $(z,t)\mapsto (z,(1+t)h^{-1}(z))$,
$f\colon \C^{n+1}\to\C$, $f(z,u)=h(z)u-1$; $\phi$ is a biholomorphism
onto the open set $\{(z,u)\in\C^{n+1}\,|, |uh(z)-1|<1\}$;
use Exercise~\ref{ese.affinedolbeault}.)\end{exer}

\begin{exer} Prove that the following facts are
equivalent:\begin{enumerate}
\item For every holomorphic function $f\colon \C\to \C$ there exists
a holomorphic function $h\colon \C\to \C$ such that
$f(z)=h(z+1)-h(z)$ for every $z$.
\item $H^{1}(\C-\{0\},\Oh_{\C})=0$.
\end{enumerate}
(Hint: Denote $p\colon \C\to \C-\{0\}$ the universal covering
$p(z)=e^{2\pi iz}$. Given $f$, use a partition of unity to find a
$C^{\infty}$ function $g$ such that $f(z)=g(z+1)-g(z)$; then $\debar g$ is
the pull back of a $\debar$-closed form on $\C-\{0\}$.)
\end{exer}

\bigskip
\section{\v{C}ech cohomology}

Let $E$ be a holomorphic vector bundle on a complex manifold $M$.
Let $\sU=\{U_{a}\}$, $a\in\sI$,
$M=\cup_{a}U_{a}$ be an open covering.
For every
$k\ge 0$ let $C^{k}(\sU,E)$ be the set of skewsymmetric sequences
$\{f_{a_{0},a_{1},\ldots,a_{k}}\}$, $a_{0},\ldots,a_{k}\in\sI$, where
$f_{a_{0},a_{1},\ldots,a_{k}}\colon U_{a_{0}}\cap\ldots\cap
U_{a_{k}}\to E$ is a holomorphic section. skewsymmetric means that
for every permutation $\sigma\in\Sigma_{k+1}$,
$f_{a_{\sigma(0)},a_{\sigma(1)},\ldots,a_{\sigma(k)}}=(-1)^{\sigma}
f_{a_{0},a_{1},\ldots,a_{k}}$.\\

The \v{C}ech differential $d\colon C^{k}(\sU,E)\to C^{k+1}(\sU,E)$ is
defined as
\[ (df)_{a_{0},\ldots,a_{k+1}}=\sum_{i=0}^{k+1}(-1)^{i}
f_{a_{0},\ldots,\widehat{a_{i}},\ldots,a_{k+1}}.\]
Since $d^{2}=0$ (exercise) we may define cocycles $Z^{k}(\sU,E)=\ker
d\subset C^{k}(\sU,E)$, coboundaries $B^{k}(\sU,E)=\Image d\subset Z^{k}(\sU,E)$
and cohomology groups
$H^{k}(\sU,E)=Z^{k}(\sU,E)/B^{k}(\sU,E)$.\\

\begin{prop}\label{I.2.5} For every holomorphic vector bundle $E$ and every
locally finite covering $\sU=\{U_{a}\}$, $a\in\sI$,
there exists a natural morphism of
$\C$-vector spaces $\theta\colon H^{k}(\sU,E)\to
H^{0,k}_{\debar}(M,E)$.\end{prop}

\begin{proof} Let $t_{a}\colon M\to \C$, $a\in \sI$, be a partition of
unity subordinate to the covering $\{U_{a}\}$:
$supp(t_{a})\subset U_{a}$, $\sum_{a} t_{a}=1$, $\sum\debar
t_{a}=0$.\\
Given $f\in C^{k}(\sU,E)$ and $a\in\sI$ we consider
\[ \phi_{a}(f)=\sum_{c_{1},\ldots,c_{k}}f_{a,c_{1},\ldots,c_{k}}\debar
t_{c_{1}}\wedge \ldots\wedge\debar t_{c_{k}}\in
\Gamma(U_{a},\sA^{0,k}(E)),\]
\[\phi(f)=\sum_{a}t_{a}\phi_{a}(f)\in\Gamma(M,\sA^{0,k}(E)).\]
Since every $f_{a,c_{1},\ldots,c_{k}}$ is holomorphic, it is clear
that $\debar \phi_{a}=0$
and then
\[\debar\phi(f)=\sum_{a}\debar t_{a}\wedge \phi_{a}(f)=
\sum_{c_{0},\ldots,c_{k}}f_{c_{0},\ldots,c_{k}}\debar
t_{c_{0}}\wedge \ldots\wedge\debar t_{c_{k}}.\]
We claim that $\phi$ is a morphism of complexes; in fact
\begin{multline*} 
\phi(df)=\sum_{a}t_{a}
\sum_{c_{0},\ldots,c_{k}}df_{a,c_{0},\ldots,c_{k}}\debar
t_{c_{0}}\wedge \ldots\wedge\debar t_{c_{k}}\\
=\sum_{a}t_{a}\left(\debar \phi(f)-\sum_{i=0}^{k}
\sum_{c_{i}}\debar t_{c_{i}}\wedge
\!\!\!\!\ds\sum_{c_{0},\ldots,\widehat{c_i},\ldots,c_{k}}\!\!\!\!\!
f_{a,c_{0},\ldots,\widehat{c}_{i},\ldots,c_{k}}
\debar t_{c_{0}}\wedge \ldots\wedge
\widehat{\debar t_{c_{i}}}\wedge\ldots\wedge\debar t_{c_{k}}
\right)\\
=\sum_{a}t_{a}\debar\phi(f)=\debar\phi(f).\end{multline*}
Setting $\theta$ as the morphism induced by
$\phi$ in cohomology, we need to prove that $\theta$ is independent
from the choice of the partition of unity.
We first note that, if $df=0$ then,
over $U_{a}\cap U_{b}$, we have
\[{\renewcommand\arraystretch{1.8}
\begin{array}{rl}
\!\phi_{a}(f)-\phi_{b}(f)\!\!\!&=\sum_{c_{1},\ldots,c_{k}}(f_{a,c_{1},\ldots,c_{k}}
-f_{b,c_{1},\ldots,c_{k}})
\debar t_{c_{1}}\wedge \cdots\wedge\debar t_{c_{k}}\\
&=\sum_{c_{1},\ldots,c_{k}}\sum_{i=1}^{k} (-1)^{i-1}
f_{a,b,c_{1},\ldots,\hat{c}_{i},\ldots,c_{k}}
\debar t_{c_{1}}\wedge \cdots\wedge\debar t_{c_{k}}\\
&=\sum_{i=1}^{k} (-1)^{i-1}\sum_{c_{1},\ldots,c_{k}}
f_{a,b,c_{1},\ldots,\hat{c}_{i},\ldots,c_{k}}
\debar t_{c_{1}}\wedge \cdots\wedge\debar t_{c_{k}}\\
&=\ds\sum_{i=1}^{k}\sum_{c_{i}}\debar t_{c_{i}}\wedge
\!\!\!\!\!\!\ds\sum_{c_{1},\ldots,\widehat{c_i},\ldots,c_{k}}
\!\!\!\!\!\! f_{a,b,c_{1},\ldots,\widehat{c}_{i},\ldots,c_{k}}
\debar t_{c_{1}}\wedge \cdots\wedge\widehat{\debar t_{c_{i}}}
\wedge\cdots\wedge
\debar t_{c_{k}}\\
&=0.\end{array}}\]

Let $v_{a}$ be another partition of 1, $\eta_{a}=t_{a}-v_{a}$, and
denote, for $f\in Z^{k}(\sU,E)$,
\[\tilde{\phi}_{a}=\sum_{c_{1},\ldots,c_{k}}f_{a,c_{1},\ldots,c_{k}}\debar
v_{c_{1}}\wedge \ldots\wedge\debar v_{c_{k}},\]
\[\psi_{a}^{j}=
\sum_{c_{1},\ldots,c_{k}}f_{a,c_{1},\ldots,c_{k}}\debar
t_{c_{1}}\wedge \ldots\wedge\debar t_{c_{j-1}}\wedge v_{c_{j}}
\debar v_{c_{j+1}}\wedge \ldots\wedge\debar v_{c_{k}},\quad j=1,\ldots,k.\]
The same argument as above shows that $\tilde{\phi}_{a}=\tilde{\phi}_{b}$
and $\psi_{a}^{j}=\psi_{b}^{j}$ for every $a,b,j$.
Therefore all the $\psi_{a}^{j}$
come from a global section $\psi^{j}\in
\Gamma(M,\sA^{0,k-1}(E))$; moreover
$\phi-\tilde{\phi}=\sum_{j}(-1)^{j-1}\debar \psi^{j}$
and then $\phi$, $\tilde{\phi}$ determine the same
cohomology class.\end{proof}

\begin{exer} In the same situation of Proposition~\ref{I.2.5}
define, for every $k\ge 0$,
$D^{k}(\sU,E)$ as the set of  sequences
$\{f_{a_{0},a_{1},\ldots,a_{k}}\}$, $a_{0},\ldots,a_{k}\in\sI$, where
$f_{a_{0},a_{1},\ldots,a_{k}}\colon U_{a_{0}}\cap\ldots\cap
U_{a_{k}}\to E$ is a holomorphic section.
Denote by $i\colon C^{k}(\sU,E)\to D^{k}(\sU,E)$ the natural
inclusion.
The same definition of the \v{C}ech differential
gives a differential
$d\colon D^{k}(\sU,E)\to D^{k+1}(\sU,E)$ making $i$ a morphism of
complexes. Moreover, it is possible to prove (see e.g. \cite[p. 214]{FAC})
that $i$ induce isomorphisms between cohomology groups.
Prove:\begin{enumerate}
\item Given two holomorphic vector bundles $E,F$  consider the linear
maps
\[ D^{k}(\sU,E)\otimes D^{p-k}(\sU,F)\mapor{\cup}
D^{p}(\sU,E\otimes F),\qquad
(f\cup g)_{a_{0},\ldots,a_{p}}=f_{a_{0},\ldots,a_{k}}\otimes
g_{a_{k},\ldots,a_{p}}.\]
Prove that $\cup$ is associative and $d(f\cup g)=df\cup
g+(-1)^{k}f\cup dg$, where  $f\in D^{k}(\sU,E)$.

\item The antisymmetrizer $p\colon D^{k}(\sU,E)\to C^{k}(\sU,E)$,
\[(pf)_{a_{0},\ldots,a_{n}}=\frac{1}{(n+1)!}\sum_{\sigma}
(-1)^{\sigma}f_{a_{\sigma(0)},\ldots,a_{\sigma(n)}},\qquad
\sigma\in\Sigma_{n+1},\]
is a morphism of complexes and then induce a morphism
$p\colon H^{k}(D^{*}(\sU,E))\to H^{k}(\sU,E)$ such that $pi=Id$
(Hint: the readers who are frightened by combinatorics may use
linearity and compatibility with restriction to open subsets $N\subset
M$ of $d,p$ to reduce the verification of $dp(f)=pd(f)$ in
the case $\sU=\{U_{a}\}$, $a=1,\ldots,m$
finite cover and $f_{a_{1},\ldots,a_{k}}\not=0$  only if
$a_{i}=i$).

\item The same definition of $\phi$ given in the proof of
\ref{I.2.5} gives a morphism of complexes $\phi_{E}\colon
D^{*}(\sU,E)\to \Gamma(M,\sA^{0,*}(E))$
which is equal to the composition of $\phi$ and $p$. In particular
$\phi_{E}$ induces
$\tilde{\theta}\colon H^{k}(D^{*}(\sU,E))\to H^{k}(M,E)$ such that
$\theta p=\tilde{\theta}$.

\item Prove that, if $dg=0$ then
$\phi_{E\otimes F}(f\cup g)=\phi_{E}(f)\wedge \phi_{F}(g)$. (Hint:
write $0=\sum_{b}t_{b}dg_{b,a_{k},\ldots,a_{p}}$.)

\item If $E,F$ are holomorphic vector bundles on $M$
then there exists a functorial cup product
\[\cup\colon H^{p}(\sU,E)\otimes H^{q}(\sU,F)\to
H^{p+q}(\sU,E\otimes F)\]
commuting with
$\theta$ and the wedge product in Dolbeault cohomology.
\end{enumerate}\end{exer}

\begin{thm}[Leray]\label{I.2.6}
Let $\sU=\{U_{a}\}$ be a locally finite covering of a complex
manifold $M$, $E$ a holomorphic vector bundle on $M$: if
$H^{k-q}_{\debar}(U_{a_{0}}\cap\ldots \cap U_{a_{q}},E)=0$ for every
$q<k$ and $a_{0},\ldots,a_{q}$, then $\theta\colon
H^{k}(\sU,E)\to H^{k}_{\debar}(M,E)$ is an isomorphism.\end{thm}

\begin{proof}
The complete proof requires sheaf theory and
spectral sequences; here we prove ``by hand''
only the cases $k=0,1$: this will be sufficient for our
applications.\\
For $k=0$ the theorem is trivial, in fact $H^{0}_{\debar}(M,E)$ and
$H^{0}(\sU,E)$ are both isomorphic to the space of holomorphic
sections of $E$ over $M$. Consider thus the case $k=1$; by assumption
$H^{1}_{\debar}(U_{a},E)=0$ for every $a$.\\
Let $\phi\in \Gamma(M,\sA^{0,1}(E))$ be a $\debar$-closed form, then for every
$a$ there exists $\psi_{a}\in \Gamma(U_{a},\sA^{0,0}(E))$
such that $\debar \psi_{a}=\phi$.
The section $f_{a,b}=\psi_{a}-\psi_{b}\colon U_{a}\cap U_{b}\to E$ is
holomorphic and then $f=\{f_{a,b}\}\in C^{1}(\sU,E)$; since
$f_{a,b}-f_{c,b}+f_{c,a}=0$ for every $a,b,c$ we have $f\in
Z^{1}(\sU,E)$; define $\sigma(\phi)\in H^{1}(\sU,E)$ as the cohomology
class of $f$. It is easy to see that $\sigma(\phi)$ is independent
from the choice of the sections $\psi_{a}$;
we want to prove that $\sigma=\theta^{-1}$. Let $t_{a}$ be a fixed
partition of unity.\\
Let $f\in Z^{1}(\sU,E)$, then $\theta(f)=[\phi]$,
$\phi=\sum_{b}f_{a,b}\debar t_{b}$; we can choose
$\psi_{a}=\sum_{b}f_{a,b}t_{b}$ and then
\[\sigma(\phi)_{a,c}=\sum_{b}(f_{a,b}-f_{c,b})t_{b}=f_{a,c},\qquad
\Rightarrow \sigma\theta=Id.\]
Conversely, if $\phi_{|U_{a}}=\debar \psi_{a}$ then
$\theta\sigma([\phi])$ is the cohomology class of
\[ \debar\sum_{b}(\psi_{a}-\psi_{b})t_{b}=
\debar\sum_{b}\psi_{a}t_{b}-\debar \sum_{b}\psi_{b}t_b  =\phi-\debar
\sum_{b}\psi_{b}t_b.\]
\end{proof}

\begin{rem}\label{I.2.7} The theory of Stein manifolds (see e.g. \cite{Gu-Ro})
says that the hypotheses of
Theorem~\ref{I.2.6} are satisfied for every $k$ whenever every
$U_{a}$ is biholomorphic to an open convex subset of $\C^{n}$.\end{rem}

\begin{ex}\label{I.2.8} Let $T\to \proj^{1}$ be the holomorphic tangent
bundle, $x_{0}, x_{1}$ homogeneous coordinates on $\proj^{1}$,
$U_{i}=\{x_{i}\not=0\}$. Since the tangent bundle of $U_{i}=\C$ is
trivial, by Dolbeault's lemma, $H^{1}(U_{i},T)=0$ and by Leray's theorem
$H^{i}(\proj^{1},T)=H^{i}(\{U_{0},U_{1}\},T)$, $i=0,1$.\\
Consider the affine coordinates $s=x_{1}/x_{0}$, $t=x_{0}/x_{1}$,
then the holomorphic sections of $T$ over $U_{0}, U_{1}$ and
$U_{0,1}=U_{0}\cap U_{1}$ are given respectively by convergent power
series
\[ \smash{\sum_{i=0}^{+\infty}}a_{i}s^{i}{\frac{\de~}{\de s}},\qquad
\sum_{i=0}^{+\infty}b_{i}t^{i}\frac{\de~}{\de t},\qquad
\smash{\sum_{i=-\infty}^{+\infty}}c_{i}s^{i}\frac{\de~}{\de s}.\]
Since, over $U_{0,1}$, $t=s^{-1}$ and
$\ds\frac{\de~}{\de t}=-s^{2}\ds\frac{\de~}{\de s}$, the Cech
differential is given by
\[ d\left(\sum_{i=0}^{+\infty}a_{i}s^{i}\frac{\de~}{\de s},
\sum_{i=0}^{+\infty}b_{i}t^{i}\frac{\de~}{\de t}\right)=
\sum_{i=0}^{+\infty}a_{i}s^{i}\frac{\de~}{\de s}+
\sum_{i=-\infty}^{2}b_{2-i}s^{i}\frac{\de~}{\de s},\]
and then $H^{1}(\{U_{0},U_{1}\},T)=0$ and
\[
H^{0}(\{U_{0},U_{1}\},T)=
\left\langle\left(\frac{\de~}{\de s},-t^2\frac{\de~}{\de t}\right),
\left(s\frac{\de~}{\de s},-t\frac{\de~}{\de t}\right),
\left(s^2\frac{\de~}{\de s},-\frac{\de~}{\de t}\right)\right\rangle.\]
\end{ex}

\begin{ex}\label{I.2.9}
If $X=\proj^{1}\times\C^{n}_{t}$ then $H^{1}(X,T_{X})=0$.
If $\C\subset \proj^{1}$ is an affine open subset
with affine coordinate $s$, then $H^{0}(X,T_{X})$ is the free
$\Oh(\C^{n})$-module generated by
\[ \desude{~}{t_{1}},\,\ldots,\,\desude{~}{t_{n}},\,\desude{~}{s},\,
s\desude{~}{s},\,s^{2}\desude{~}{s}.\]
The proof is essentially the same (replacing the constant terms
$a_{i},b_{i},c_{i}$ with holomorphic functions over $\C^{n}$) of
Example~\ref{I.2.8}.\end{ex}

\bigskip

\section{The Kodaira-Spencer map}

\begin{notation}\label{I.3.1}
Given a holomorphic map $f\colon X\to Y$ of complex manifolds and
complexified vector fields $\eta\in
\Gamma(X,\sA^{0,0}(T_{X}))$, $\gamma\in
\Gamma(Y,\sA^{0,0}(T_{Y}))$ we write
$\gamma=f_{*}\eta$ if for every $x\in X$ we have
$f_{*}\eta(x)=\gamma(f(x))$,
where $f_{*}\colon T_{x,X}\to T_{f(x),Y}$
is the differential of $f$.\end{notation}

Let $f\colon M\to B$ be a fixed smooth family of compact complex
manifolds, $\dim B=n$, $\dim M=m+n$; for every $b\in B$ we let
$M_{b}=f^{-1}(b)$.\\

\begin{defn}\label{I.3.2} A holomorphic coordinate chart
$(z_{1},\ldots,z_{m},t_{1},\ldots,t_{n})\colon U\hookrightarrow\C^{m+n}$,
$U\subset M$ open,
is called \emph{admissible} if $f(U)$ is contained in a coordinate chart
$(v_{1},\ldots,v_{n})\colon V\hookrightarrow \C^{n}$,
$V\subset B$, such that $t_i=v_i\circ f$ for every $i=1,\ldots, n$.
\end{defn}

Since the differential of $f$ has everywhere maximal rank, by the
implicit function theorem,
$M$ admits a locally finite covering of admissible coordinate charts.

\begin{lem}\label{I.3.3} Let $f\colon M\to B$ be a smooth family of compact
complex manifolds. For every $\gamma\in \Gamma(B,\sA^{0,0}(T_{B}))$ there
exists $\eta\in \Gamma(M,\sA^{0,0}(T_{M}))$ such that $f_{*}\eta=\gamma$.
\end{lem}

\begin{proof} Let $M=\cup U_{a}$ be a locally finite covering of
admissible charts; on every $U_{a}$ there exists $\eta_{a}\in
\Gamma(U_{a},\sA^{0,0}(T_{M}))$ such that $f_{*}\eta_{a}=\gamma$.\\
It is then sufficient to take $\eta=\sum_{a}\rho_{a}\eta_{a}$, being
$\rho_{a}\colon U_{a}\to\C$ a partition of unity subordinate to the
covering $\{U_{a}\}$.\end{proof}

Let $T_{f}\subset T_{M}$ be the holomorphic vector subbundle of
tangent vectors $v$ such that $f_{*}v=0$.
If $z_{1},\ldots,z_{m},t_{1},\ldots,t_{n}$
is an admissible system of local coordinates then
$\ds\frac{\de~}{\de z_{1}},\ldots,\ds\frac{\de~}{\de z_{m}}$
is a local frame of $T_{f}$. Note that the restriction of $T_{f}$
to $M_{b}$ is equal to $T_{M_{b}}$.\\

For every open subset $V\subset B$ let $
\Gamma(V,T_{B})$ be the space of
holomorphic vector fields on $V$.\\
For every $\gamma\in\Gamma(V,T_{B})$ take
$\eta\in \Gamma(f^{-1}(V),\sA^{0,0}(T_{M}))$ such that $f_{*}\eta=\gamma$.
In an admissible system of local coordinates $z_{i},t_{j}$ we have
$\eta=\sum_{i}\eta_{i}(z,t)\ds\frac{\de~}{\de
z_{i}}+\sum_{j}\gamma_{i}(t)\ds\frac{\de~}{\de t_{j}}$, with
$\gamma_{i}(t)$ holomorphic,
$\debar\eta=\sum_{i}\debar\eta_{i}(z,t)\ds\frac{\de~}{\de z_{i}}$
and then $\debar\eta\in \Gamma(f^{-1}(V),\sA^{0,1}(T_{f}))$.\\
Obviously $\debar\eta$ is $\debar$-closed and then we
can define the \emph{Kodaira-Spencer} map
\[ \sKS(V)_{f}\colon \Gamma(V,T_{B})\to H^{1}(f^{-1}(V),T_{f}),
\qquad \sKS(V)_{f}(\gamma)=[\debar\eta].\]

\begin{lem}\label{I.3.4}
The map $\sKS(V)_{f}$ is a well-defined homomorphism of
$\Oh(V)$-modules.\end{lem}

\begin{proof} If $\tilde{\eta}\in\Gamma(f^{-1}(V),\sA^{0,0}(T_{M}))$,
$f_{*}\tilde{\eta}=\gamma$, then $\eta-\tilde{\eta}\in
(f^{-1}(V),\sA^{0,0}(T_{f}))$ and $[\debar\tilde{\eta}]=[\debar\eta]\in
H^{1}(f^{-1}(V),T_{f})$.\\
If $g\in\Oh(V)$ then $f_{*}(f^{*}g)\eta=g\gamma$,
$\debar (f^{*}g)\eta=(f^{*}g)\debar\eta$.\end{proof}
If $V_{1}\subset V_{2}\subset B$ then the Kodaira-Spencer maps
$\sKS(V_{i})_{f}\colon \Gamma(V_{i},T_{B})\to H^{1}(f^{-1}(V_{i}),T_{f})$, $i=1,2$,
commute with the restriction maps
$\Gamma(V_{2},T_{B})\to \Gamma(V_{1},T_{B})$,
$H^{1}(f^{-1}(V_{2}),T_{f})\to H^{1}(f^{-1}(V_{1}),T_{f})$.
Therefore we get a well defined  $\Oh_{B,b}$-linear map
\[ \sKS_{f}\colon \Theta_{B,b}\to (R^{1}f_{*}T_{f})_{b},\]
where $\Theta_{B,b}$ and $(R^{1}f_{*}T_{f})_{b}$ are by definition
the direct limits, over the set of open neighbourhood $V$ of $b$, of
$\Gamma(V,T_{B})$ and $H^{1}(f^{-1}(V),T_{f})$ respectively.

If $b\in  B$, then  there exists a
linear map $\KS_{f}\colon T_{b,B}\to H^{1}(M_{b},T_{M_{b}})$ such that
for every open subset $b\in V\subset B$ there exists a
commutative diagram
\[\begin{array}{ccc}
\Gamma(V,T_{B})&\mapor{\sKS(V)_{f}}&H^{1}(f^{-1}(V),T_{f})\\
\mapver{}&&\mapver{r}\\
T_{b,B}&\mapor{\KS_{f}}&H^{1}(M_{b},T_{M_{b}})\end{array}\]
where the vertical arrows are the natural restriction maps.\\
In fact, if $V$ is a polydisk  then $T_{b,B}$ is the quotient
of the complex vector space $\Gamma(V,T_{B})$ by the subspace
$I=\{\gamma\in \Gamma(V,T_{B})\mid \gamma(b)=0\}$;
by $\Oh(V)$-linearity $I$ is contained in the kernel of $r\circ
\sKS(V)_{f}$.

The Kodaira-Spencer map has at least two geometric interpretations:
obstruction to the holomorphic lifting of vector fields and
first-order variation of complex structures (this is a concrete
feature of the general philosophy that deformations are a derived
construction of automorphisms).

\begin{prop}\label{I.3.5} Let $f\colon M\to B$ be a family of compact
complex manifolds and
 $\gamma\in\Gamma(V,T_{B})$, then $\sKS(V)_{f}(\gamma)=0$ if and
only if there exists $\eta\in\Gamma(f^{-1}(V),T_{M})$ such that
$f_{*}\eta=\gamma$.\end{prop}

\begin{proof} One implication is trivial; conversely let
$\eta\in \Gamma(f^{-1}(V),\sA^{0,0}(T_{M}))$ such that $f_{*}\eta=\gamma$.
If $[\debar\eta]=0$ then there exists $\tau\in
\Gamma(f^{-1}(V),\sA^{0,0}(T_{f}))$ such that $\debar(\eta-\tau)=0$,
$\eta-\tau\in\Gamma(f^{-1}(V),T_{M})$ and $f_*(\eta-\tau)=\gamma$.
\end{proof}

To compute the Kodaira-Spencer map in terms of Cech cocycles we
assume that $V$ is a polydisk with coordinates $t_{1},\ldots,t_{n}$ and we
fix a
locally finite covering $\sU=\{U_{a}\}$ of admissible holomorphic
coordinates
$z_{1}^{a},\ldots,z_{m}^{a},t_{1}^{a},\ldots,t_{n}^{a}\colon U_{a}\to
\C$, $t_{i}^{a}=f^{*}t_{i}$.\\
On $U_{a}\cap U_{b}$ we have the transition functions
{\renewcommand{\arraystretch}{1.7}
\[ \left\{\begin{array}{ll} z_{i}^{b}=g_{i,a}^{b}(z^{a},t^{a}),&\quad
i=1,\ldots,m\\
t_{i}^{b}=t_{i}^{a},&\quad i=1,\ldots,n\end{array}\right.\]}
Consider a fixed integer
$h=1,\ldots,n$ and  $\eta\in \Gamma(f^{-1}(V),\sA^{0,0}(T_{M}))$
such that $f_{*}\eta=\dfrac{\de~}{\de
t_{h}}$; in local coordinates we have
\[\eta=\sum_{i}\eta_{i}^{a}(z^{a},t^{a})\dfrac{\de~}{\de
z_{i}^{a}}+\dfrac{\de~}{\de t_{h}^{a}},\qquad
\eta=\sum_{i}\eta_{i}^{b}(z^{b},t^{b})\dfrac{\de~}{\de
z_{i}^{b}}+\dfrac{\de~}{\de t_{h}^{b}}.\]
Since, for every $a$,
$\eta-\dfrac{\de~}{\de t_{h}^{a}}\in \Gamma(U_{a},\sA^{0,0}(T_{f}))$ and
$\debar\left(\eta-\dfrac{\de~}{\de t_{h}^{a}}\right)=\debar\eta$,
$\sKS(V)_{f}\left(\dfrac{\de~}{\de t_{h}}\right)\in H^{1}(\sU,T_{f})$ is represented
by the cocycle
\begin{form}\label{KScociclo}
\[\sKS(V)_{f}\left(\frac{\de~}{\de
t_{h}}\right)_{b,a}=\left(\eta-\frac{\de~}{\de t_{h}^{b}}\right)
-\left(\eta-\frac{\de~}{\de
t_{h}^{a}}\right)=\frac{\de~}{\de t_{h}^{a}}-\frac{\de~}{\de t_{h}^{b}}=
\sum_{i}\frac{\de g_{i,a}^{b}}{\de t_{h}^{a}}\frac{\de~}{\de
z_{i}^{b}}.\]
\end{form}

The above formula allows to prove easily the invariance of the
Kodaira-Spencer maps under base change; more precisely if $f\colon
M\to B$ is a smooth family, $\phi\colon C\to B$ a holomorphic map,
$\hat{\phi}$, $\hat{f}$ the pullbacks of $\phi$ and $f$,
\[
\begin{array}{ccc}M\times_{B}C&\mapor{\hat{\phi}}&M\\
\mapver{\hat{f}}&&\mapver{f}\\
C&\mapor{\phi}&B\end{array}\]
$c\in C$,
$b=f(c)$.
\begin{thm}\label{I.3.7} In the above notation, via the natural isomorphism
$M_{b}=\hat{f}^{-1}(c)$, we have
\[ \KS_{\hat{f}}=\KS_{f}\phi_{*}\colon T_{c,C}\to
H^{1}(M_{b},T_{M_{b}}).\]\end{thm}

\begin{proof}
It is not restrictive to assume $B\subset\C^{n}_{t}$,
$C\subset\C^{s}_{u}$ polydisks,  $c=\{u_{i}=0\}$ and
$b=\{t_{i}=0\}$, $t_{i}=\phi_{i}(u)$.\\
If $z^{a},t^{a}\colon U_{a}\to \C$, $z^{b},t^{b}\colon U_{b}\to \C$
are admissible local coordinate sets with transition functions
$z_{i}^{b}=g_{i,a}^{b}(z^{a},t^{a})$, then
$z^{a},u^{a}\colon U_{a}\times_{B}C\to \C$,
$z^{b},t^{b}\colon U_{b}\times_{B}C\to \C$ are admissible with
transition functions $z_{i}^{b}=g_{i,a}^{b}(z^{a},\phi(u^{a}))$.\\
Therefore
\[\KS_{\hat{f}}\left(\frac{\de~}{\de u_{h}}\right)_{b,a}=
\sum_{i}\frac{\de g_{i,a}^{b}}{\de u_{h}^{a}}\frac{\de~}{\de
z_{i}^{b}}=
\sum_{i,j}\frac{\de g_{i,a}^{b}}{\de t_{j}^{a}}
\frac{\de \phi_{j}}{\de u_{h}^{a}}
\frac{\de~}{\de
z_{i}^{b}}=\KS_{f}\left(\phi_{*}\frac{\de~}{\de u_{h}}\right)_{b,a}.\]
\end{proof}

It is clear that the Kodaira-Spencer map
$\KS_{f}\colon T_{b_{0},B}\to H^{1}(M_{0},T_{M_{0}})$
is defined  for every isomorphism class of
deformation $M_{0}\to M\mapor{f}(B,b_{0})$:
The map $\sKS_{f}\colon \Theta_{B,b_{0}}\to (R^{1}f_{*}T_{f})_{b_{0}}$
is defined  up to isomorphisms of the $\Oh_{B,b_0}$ module
$(R^{1}f_{*}T_{f})_{b_{0}}$.
\bigskip

\begin{defn}\label{I.3.8} Consider a deformation
$\xi:\, M_{0}\mapor{i}M\mapor{f}(B,b_{0})$,
$fi(M_{0})=b_{0}$, with Kodaira-Spencer map $\KS_{\xi}\colon
T_{b_{0},B}\to H^{1}(M_{0},T_{M_{0}})$.
$\xi$ is called:\begin{enumerate}
\item \emph{Versal} if $\KS_{\xi}$ is surjective and
for every germ of complex manifold $(C,c_{0})$
the morphism
\[\Mor_{\mathbf{Ger}}((C,c_{0}),(B,b_{0}))\to \Def_{M_{0}}(C,c_{0}),\qquad
g\mapsto g^{*}\xi\]
is surjective.

\item \emph{Semiuniversal} if it is versal and
$\KS_{\xi}$ is bijective.

\item \emph{Universal} if $\KS_{\xi}$ is bijective and
for every pointed complex manifolds $(C,c_{0})$
the morphism
\[\Mor_{\mathbf{Ger}}((C,c_{0}),(B,b_{0}))\to \Def_{M_{0}}(C,c_{0}),\qquad
g\mapsto g^{*}\xi\]
is bijective.
\end{enumerate}\end{defn}

Versal deformations are also called \emph{complete};  semiuniversal
deformations are also called \emph{miniversal}
or \emph{Kuranishi} deformations.

Note that if $\xi$ is semiuniversal, $g_{1}, g_{2}\in
\Mor_{\mathbf{Ger}}((C,c_{0}),(B,b_{0}))$ and
$g_{1}^{*}\xi=g_{2}^{*}\xi$ then, according to
Theorem~\ref{I.3.7},
$dg_{1}=dg_{2}\colon
T_{c_{0},C}\to T_{b_{0},B}$.

\begin{exer} A universal deformation $\xi:\,
M_{0}\mapor{i}M\mapor{f}(B,b_{0})$
induces a representation (i.e. a homomorphism of
groups)
\[\rho\colon\Aut(M_{0})\to \Aut_{\mathbf{Ger}}((B,b_{0})),\qquad
\rho(g)^{*}\xi=\xi^{g},\quad g\in\Aut(M_{0}).\]
Every other universal deformation over the germ $(B,b_{0})$ gives a
conjugate representation.
\end{exer}

\bigskip

\section{Rigid varieties}

\begin{defn}\label{I.4.1}
A deformation $M_{0}\to M\to (B,b_{0})$ is called \emph{trivial}
if it is isomorphic to
\[M_{0}\mapor{Id\times \{b_{0}\}}M_{0}\times
B\mapor{pr}(B,b_{0}).\]
\end{defn}

\begin{lem}\label{I.4.2}
Let $f\colon M\to \Delta^{n}_{R}$ be a smooth family of
compact complex manifolds, $t_{1},\ldots,t_{n}$  coordinates in the
polydisk $\Delta^{n}_{R}$. If there exist holomorphic vector fields
$\chi_{1},\ldots,\chi_{n}$ on $M$ such that
$f_{*}\chi_{h}=\ds\frac{\de~}{\de t_{h}}$ then there exists $0<r\le R$
such that $f\colon f^{-1}(\Delta^{n}_{r})\to\Delta^{n}_{r}$ is the
trivial family.\end{lem}

\begin{proof}
For every $r\le R$, $h\le n$ denote
\[\Delta^{h}_{r}=\{(z_{1},\ldots,z_{n})\in\C^{n}\mid |z_{1}|<r,\ldots,
|z_{h}|<r, z_{h+1}=0,\ldots, z_{n}=0\}\subset \Delta^{n}_{R}.\]
We prove by induction on $h$ that there exists
$R\ge r_{h}>0$ such that the restriction of the
family $f$ over $\Delta^{h}_{r_{h}}$ is trivial.
Taking $r_{0}=R$ the
statement is obvious for $h=0$. Assume that the family is trivial
over $\Delta^{h}_{r_{h}}$, $h<n$; shrinking $\Delta^{n}_{R}$ if
necessary it is not restrictive to assume $R=r_{h}$
and the family trivial over $\Delta^{h}_{R}$.\\
The integration of the vector field $\chi_{h+1}$ gives an open neighbourhood
$M\times \{0\}\subset U\subset M\times\C$ and a holomorphic map
$H\colon U\to M$ with the following properties (see e.g. \cite[Ch. 
VII]{cartan}):
\begin{enumerate}
\item For every $x\in M$, $\{x\}\times\C\cap U=\{x\}\times\Delta(x)$
with $\Delta(x)$ a disk.
\item For every $x\in M$ the map $H_{x}=H(x,-)\colon \Delta(x)\to
M$ is the solution of the Cauchy problem
\renewcommand\arraystretch{1.8}
\[\left\{\begin{array}{l}
\ds\frac{dH_{x}}{dt}(t)=\chi_{h+1}(H_{x}(t))\\
H_{x}(0)=x\end{array}\right.\]
\renewcommand\arraystretch{1}
In particular if
$H(x,t)$ is defined then
$f(H(x,t))=f(x)+(0,\ldots,t,\ldots,0)$ ($t$ in the $(h+1)$-th coordinate).
\item If $V\subset M$ is open and $V\times \Delta\subset U$ then for
every $t\in \Delta$ the map $H(-,t)\colon V\to M$ is an open embedding.
\end{enumerate}

Since $f$ is proper there exists $r\le R$ such that
$f^{-1}(\Delta^{h}_{r})\times \Delta_{r}\subset U$; then the
holomorphic map
$H\colon f^{-1}(\Delta^{h}_{r})\times \Delta_{r}\to
f^{-1}(\Delta^{h+1}_{r})$ is a biholomorphism (exercise) giving a
trivialization of the family over $\Delta^{h+1}_{r}$.
\end{proof}

\begin{ex}\label{I.4.2,5}
Lemma~\ref{I.4.2} is generally false if $f$ is not proper
(cf. the exercise in Lecture 1 of \cite{Konts}).\\
Consider for instance an irreducible polynomial
$F\in\C[x_{1},\ldots,x_{n},t]$; denote by $f\colon
\C^{n}_{x}\times\C_{t}\to\C_{t}$ the projection on the second
factor  and
\[V=\left\{(x,t)\,\left|\, F(x,t)=\desude{F}{x_{i}}(x,t)=0,\,
i=1,\ldots,n \right. \right\}.\]
Assume that $f(V)$ is a finite set of points and set $B=\C-f(V)$,
$X=\{(x,t)\in\C^{n}\times B\mid F(x,t)=0\}$. Then $X$ is a regular
hypersurface, the restriction
$f\colon X\to B$ is surjective and its differential is surjective
everywhere.\\
$X$ is closed in the affine variety $\C^{n}\times B$,
by Hilbert's Nullstellensatz there exist regular functions
$g_{1},\ldots,g_{n}\in\Oh(\C^{n}\times B)$ such that
\[g:=\sum_{i=1}^{n}g_{i}\desude{F}{x_{i}}\equiv 1\quad\pmod{F}.\]
On the open subset $U=\{g\not=0\}$ the algebraic vector field
\[\chi=\sum_{i=1}^{n}\frac{g_{i}}{g}\left(
\desude{F}{x_{i}}\desude{~}{t}-\desude{F}{t}\desude{~}{x_{i}}\right)
=\desude{~}{t}-\sum_{i=1}^{n}\frac{g_{i}}{g}\desude{F}{t}\desude{~}{x_{i}}
\]
is tangent to $X$ and lifts $\desude{~}{t}$.\\
In general the fibres of $f\colon X\to B$ are not biholomorphic:
consider for example the case $F(x,y,\lambda)=
y^{2}-x(x-1)(x-\lambda)$. Then $B=\C-\{0,1\}$ and $f\colon X\to B$ is
the restriction to the affine subspace $x_{0}\not=0$ of the family
$M\to B$ of
Example~\ref{I.1.4}.\\
The fibre $X_{\lambda}=f^{-1}(\lambda)$ is
$M_{\lambda}-\{\hbox{point}\}$, where $M_{\lambda}$ is
an elliptic curve with $j$-invariant
$j({\lambda})=2^{8}(\lambda^{2}-\lambda+1)^{3}
\lambda^{-2}(\lambda-1)^{-2}$. If $X_{a}$ is biholomorphic to $X_{b}$
then, by Riemann's extension theorem, also $M_{a}$ is biholomorphic to
$M_{b}$ and then $j(a)=j(b)$.\end{ex}

\begin{exer} Find a holomorphic vector field $\chi$ lifting
$\desude{~}{\lambda}$ and tangent to $\{F=0\}\subset \C^{2}\times\C$,
where
$F(x,y,\lambda)=y^{2}-x(x-1)(x-\lambda)$ (Hint: use the Euclidean algorithm
to find $a,b\in\C[x]$ such that
$ay\desude{F}{y}+b\desude{F}{x}=1+2aF$).
\end{exer}

\begin{thm}\label{I.4.4}
A deformation $M_{0}\to M\mapor{f}(B,b_{0})$ of a compact manifold is trivial if and
only if $\sKS_{f}\colon \Theta_{B,b_{0}}\to (R^{1}f_{*}T_{f})_{b_{0}}$ is
trivial.\end{thm}

\begin{proof} One implication is clear; conversely assume $\sKS_{f}=0$,
it is not restrictive to assume $B$ a polydisk
with coordinates $t_{1},\ldots,t_{n}$
and $f$ a smooth family.
After a possible shrinking of $B$ we have
$\sKS(B)_{f}\left(\desude{~}{t_{i}}\right)=0$
for every $i=1,\ldots,n$. According to \ref{I.3.5} there exist
holomorphic vector fields $\xi_{i}$ such that
$f_{*}\xi_{i}=\desude{~}{t_{i}}$; by \ref{I.4.2} the family is
trivial over a smaller polydisk $\Delta\subset B$.\end{proof}

Note that if a smooth family $f\colon M\to B$ is locally trivial, then
for every $b\in B$ the Kodaira-Spencer map
$\KS_{f}\colon T_{b,B}\to H^{1}(M_{b},T_{M_{b}})$
is trivial for every $b\in B$.

\begin{thm}\label{I.4.5} \emph{(Semicontinuity and base change)}\\
Let $E\to M$ be a holomorphic vector bundle on the total
space of a smooth family $f\colon M\to B$.
Then, for every $i\ge 0$:\begin{enumerate}
\item $b\mapsto h^{i}(M_{b},E)$ is upper semicontinuous.
\item If $b\mapsto h^{i}(M_{b},E)$ is constant, then for every $b\in
B$ there exists an open  neighbourhood $b\in U$ and elements
$e_{1},\ldots,e_{r}\in H^{i}(f^{-1}(U),E)$ such that:\begin{enumerate}
\item $H^{i}(f^{-1}(U),E)$ is the free $\Oh(U)$-module generated by
$e_{1},\ldots,e_{n}$.
\item $e_{1},\ldots,e_{r}$ induce a basis of $H^{i}(M_{c},E)$ for
every $c\in U$.
\end{enumerate}
\item If $\,b\mapsto h^{i-1}(M_{b},E)$ and $\,b\mapsto h^{i+1}(M_{b},E)$ are constant
then also $\,b\mapsto h^{i}(M_{b},E)$ is constant.
\end{enumerate}
\end{thm}
\begin{proof} \cite[Ch. 3, Thm. 4.12]{Ba-Sta}, \cite[I, Thm. 2.2]{KS3},
\cite{Kobook}.\end{proof}

\begin{cor}\label{I.4.6}
Let $X$ be a compact complex manifold. If $H^1(X,T_X)=0$ then
every deformation of $X$ is trivial.
\end{cor}

\begin{defn}\label{I.4.7} A compact complex manifold $X$ is called 
\emph{rigid}
if $H^1(X,T_X)=0$.\end{defn}

\begin{cor}\label{I.4.8} Let $f\colon M\to B$ a smooth family of compact complex manifolds. If
$b\mapsto h^1(M_b,T_{M_b})$ is constant and $\KS_f=0$ at every point $b\in B$ then the
family is locally trivial.\end{cor}

\begin{proof} (cf. Example \ref{nontriv}) 
Easy consequence of Theorems~\ref{I.4.4} and \ref{I.4.5}.
\end{proof}

\begin{ex}\label{I.4.9}
Consider the following family of Hopf surfaces $f\colon M\to \C$,
$M=X/G$ where $X=B\times (\C^{2}-\{0\})$ and $G\simeq\Z$ is generated
by $(b,z_{1},z_{2})\mapsto (b,2z_{1},b^{2}z_{1}+2z_{2})$.\\
The fibre $M_{b}$ is the Hopf surface $S_{A(b)}$, where
$A(b)=\left(\begin{array}{cc}2&0\\ b^{2}&2\end{array}\right)$
and then $M_{0}$ is not biholomorphic to $M_{b}$ for every $b\not=0$.\\
This family is isomorphic to $N\times_{\C}B$,
where $B\to \C$ is the map $b\mapsto b^{2}$ and $N$ is the quotient of
$\C\times (\C^{2}-\{0\})$ by the group generated by
$(s,z_{1},z_{2})\mapsto (s,2z_{1},sz_{1}+2z_{2})$. By base-change
property, the Kodaira-Spencer map $\KS_{f}\colon
T_{0,B}\to H^{1}(M_{0},T_{M_{0}})$ is trivial.\\
On the other hand the family is  trivial over $B-\{0\}$, in
fact the map
\[(B-\{0\})\times (\C^{2}-\{0\})\to (B-\{0\})\times (\C^{2}-\{0\}),\qquad
(b,z_{1},z_{2})\mapsto (b,b^{2}z_{1},z_{2})\]
induces to the quotient an
isomorphism $(B-\{0\})\times M_{1}\simeq (M-f^{-1}(0))$. Therefore
the Kodaira-Spencer map $\KS_{f}\colon
T_{b,B}\to H^{1}(M_{b},T_{M_{b}})$ is trivial for every $b$.\\
According to the base-change theorem the dimension of
$H^{1}(M_{b},T_{M_{b}})$ cannot be constant: in fact it is proved in
\cite{KS3} that $h^{1}(M_{0},T_{M_{0}})=4$ and
$h^{1}(M_{b},T_{M_{b}})=2$ for $b\not=0$.
\end{ex}

\begin{ex}\label{I.4.10} Let 
$M\subset\C_{b}\times\proj^{3}_{x}\times\proj^{1}_{u}$ be the subset
defined by the equations
\[ u_{0}x_{1}=u_{1}(x_{2}-bx_{0}),\qquad u_{0}x_{2}=u_{1}x_{3},\]
$f\colon M\to \C$ the projection onto the first factor
and $f^{*}\colon M^{*}=(M-f^{-1}(0))\to (\C-\{0\})$ its restriction.\\
Assume already proved
that $f$ is a family (this will be done in the next chapter); we
want to prove that:\begin{enumerate}
\item  $f^{*}$  is a trivial
family.
\item  $f$ is not locally trivial at $b=0$.\\
\end{enumerate}
\begin{proof}[Proof of 1]
After the linear change of coordinates $x_{2}-bx_{0}\mapsto x_{0}$ the
equations of $M^{*}\subset \C-\{0\}\times\proj^{3}\times\proj^{1}$
become \[ u_{0}x_{1}=u_{1}x_{0},\qquad u_{0}x_{2}=u_{1}x_{3}\]
and there exists an isomorphism of families
$\C-\{0\}\times\proj^{1}_{s}\times\proj^{1}_{u}\to M^{*}$, given by
\[ (b,[t_{0},t_{1}],[u_{0},u_{1}])\mapsto
(b,[t_{0}u_{1},t_{0}u_{0},t_{1}u_{1},t_{1}u_{0}], [u_{0},u_{1}]).\]
\end{proof}
\begin{proof}[Proof of 2]
Let $Y\simeq\proj^{1}\subset M_{0}$ be the subvariety of equation
$b=x_{1}=x_{2}=x_{3}=0$. Assume $f$ locally trivial, then there
exist an open neighbourhood  $0\in U\subset\C$ and a commutative
diagram of holomorphic maps
\renewcommand\arraystretch{1.2}
\[\begin{array}{ccc}
Y\times U&\mapor{j}&M\\
\mapver{pr}&&\mapver{f}\\
U&\mapin{i}&\C\end{array}\]
\renewcommand\arraystretch{1}
where $i$ is the inclusion, $j$ is injective and extends the
identity $Y\times\{0\}\to Y\subset M_{0}$.\\
Possibly shrinking $U$ it is not restrictive to assume that the
image of $j$ is contained in the open subset
$V_{0}=\{x_{0}\not=0\}$. For $b\not=0$ the holomorphic map
$\delta\colon V_{0}\cap M_{b}\to\C^{3}$,
\[\delta(b,[x_{0},x_{1},x_{2},x_{3}], [u_{0},u_{1}])=
\left(\frac{x_{1}}{x_{0}},\frac{x_{2}}{x_{0}},\frac{x_{3}}{x_{0}}\right),\]
is injective; therefore for $b\in U$, $b\not=0$, the holomorphic map
$\delta j(-,b)\colon Y\simeq\proj^{1}\to\C^{3}$ is injective. This
contradicts the maximum principle of holomorphic functions.
\end{proof}
\end{ex}

\begin{ex}\label{nununiv}
In the notation of Example~\ref{I.4.10}, the deformation $M_{0}\to 
M\mapor{b}(\C,0)$ is not universal: in order to see this it is sufficient 
to prove that  M is isomorphic to the 
deformation  $g^{*}M$, where $g\colon (\C,0)\to (\C,0)$ is the 
holomorphic map $g(b)=b+b^{2}$.\\
The equation of $g^{*}M$ is 
\[ u_{0}x_{1}=u_{1}(x_{2}-(b+b^{2})x_{0}),\qquad u_{0}x_{2}=u_{1}x_{3},\]
and the isomorphism of deformations  $g^{*}M\to M$ is given by 
\[ (b,[x_{0},x_{1},x_{2},x_{3}], 
[u_{0},u_{1}])=(b,[(1+b)x_{0},x_{1},x_{2},x_{3}], [u_{0},u_{1}]).\]
\end{ex}

\begin{ex}\label{nontriv}
Applying the
base change $\C\to\C$, $b\mapsto b^{2}$, 
to the family $M\to\C$ of Example~\ref{I.4.10} we 
get a family 
with trivial $\KS$ at every point of the base but not locally trivial 
at $0$.\\ 
We will prove in
\ref{II.1.4} that $H^{1}(M_{b},T_{M_{b}})=0$ for $b\not=0$ and
$H^{1}(M_{0},T_{M_{0}})=\C$.
\end{ex}

\bigskip

\section{Historical survey,~\ref{CAP:FAMILIES}}

The  deformation theory of complex manifolds
began in the years 1957-1960 by a series of papers of
Kodaira-Spencer \cite{KS1}, \cite{KS2}, \cite{KS3}
and Kodaira-Nirenberg-Spencer \cite{KNS}.

The main results of these papers were the completeness and existence
theorem for versal deformations.

\begin{thm}\label{I.5.1} \emph{(Completeness theorem, \cite{KS2})}\\
A deformation $\xi$ over a smooth germ $(B,0)$ of  a compact complex manifold
$M_{0}$ is versal if and only if the Kodaira-Spencer map
$\KS_{\xi}\colon T_{0,B}\to H^{1}(M_{0},T_{M_{0}})$ is
surjective.\end{thm}

Note that if a deformation $M_{0}\mapor{}M\mapor{f}(B,0)$
is versal then
we can take a linear subspace $0\in C\subset B$ making the
Kodaira-Spencer map $T_{0,C}\to H^{1}(M_{0},T_{M_{0}})$ bijective; by
completeness theorem $M_{0}\to M\times_{B}C\to (C,0)$ is semiuniversal.\\

In general, a compact complex
manifold does not have a versal deformation over a smooth germ. The
problem of determining when such a deformation exists is one of the
most difficult in deformation theory.\\
A partial answer is given by

\begin{thm}\label{I.5.2} \emph{(Existence theorem, \cite{KNS})}\\
Let $M_{0}$ be a compact complex manifold. If
$H^{2}(M_{0},T_{M_{0}})=0$ then $M_{0}$ admits a semiuniversal
deformation over a smooth base.\end{thm}

The condition $H^{2}(M_{0},T_{M_{0}})=0$ is sufficient but it is
quite far from being necessary. The ``majority'' of manifolds having
a versal deformation over a smooth germ has the above cohomology group
different from $0$.\\

The next problem is to determine when a semiuniversal deformation is
universal: a sufficient (and almost necessary) condition is given by
the following theorem.

\begin{thm}\label{I.5.3} \emph{(\cite{Sch}, \cite{Wawrik1})}
Let $\xi:\, M_{0}\mapor{}M\mapor{}(B,0)$ be a semiuniversal
deformation of
a compact complex manifold $M_{0}$. If $b\mapsto h^{0}(M_{b},T_{M_{b}})$
is constant (e.g. if $H^{0}(M_{0},T_{M_{0}})=0$)
then $\xi$ is universal.\\
\end{thm}

\begin{rem} If a compact complex manifold $M$ has finite holomorphic
automorphisms then $H^{0}(M,T_{M})=0$, while the converse is generally
false (take as an example the Fermat quartic surface in $\proj^{3}$,
cf. \cite{segreB}).
\end{rem}

\begin{ex}\label{I.5.4} Let $M\to B$ be a smooth family of
compact complex tori of
dimension $n$, then $T_{M_{b}}=\oplus_{i=1}^{n}\Oh_{M_{b}}$ and then
$h^{0}(M_{b},T_{M_{b}})=n$ for every $b$.
\end{ex}

\begin{ex}\label{I.5.5} If $K_{M_{0}}$ is ample then, by a theorem of
Matsumura \cite{Matsu}, $H^{0}(M_{0},T_{M_{0}})=0$.\end{ex}

\begin{exer} The deformation $M_{0}\mapor{}M\mapor{f}\C$, where $f$
is the family of Example~\ref{I.4.10}, is not universal.\end{exer}


\chapter[Deformations of Segre-Hirzebruch surfaces]{Deformations of
Segre-Hirzebruch surfaces}
\label{CAP:SEGRE}
\piede

In this chapter we compute the Kodaira-Spencer map of some
particular deformations and, using the completeness theorem
\ref{I.5.1}, we give a concrete description of the semiuniversal
deformations of the Segre-Hirzebruch surfaces $\F_k$ (Theorem~\ref{semipersegre}).\\
As a by-product we get examples of deformation-unstable
submanifolds (Definition~\ref{II.4.1}). A sufficient condition for
stability of submanifolds is the well known \emph{Kodaira stability
theorem} (Thm.~\ref{II.4.2}) which is stated without proof in the
last section.

\bigskip

\section{Segre-Hirzebruch surfaces}
\label{sec:1}

We consider the following description of the Segre-Hirzebruch surface
$\F_q$, $q\ge 0$.\\
\[\F_q=(\C^2-\{0\})\times(\C^2-\{0\})/\sim,\]
where the equivalence relation $\sim$ is given by the
$(\C^{*})^{2}$-action
\[(l_0,l_1,t_0,t_1)\mapsto (\lambda l_0,\lambda l_1, \lambda^q\mu t_0,
\mu t_1),\qquad \lambda,\mu\in\C^*.\]

The projection $\F_{q}\to \proj^{1}$, 
$[l_{0},l_{1},t_{0},t_{1}]\mapsto [l_{0},l_{1}]$ is well defined and it is 
a $\proj^{1}$-bundle (cf. Example~\ref{II.2.4}).\\
Note that $\F_{0}=\proj^{1}\times\proj^{1}$;
$\F_q$ is covered by four affine planes $\C^2\simeq U_{i,j}=
\{l_it_j\not=0\}$.
In this affine covering we define local coordinates according to
the following table

\begin{center}
{\renewcommand{\arraystretch}{4}
\begin{tabular}{c|c}
$\quad\ds U_{0,0}:\quad z=\ds\frac{l_1}{l_0},\quad
s=\ds\frac{t_1l_0^q}{t_0}\quad$&
$\quad U_{0,1}:\quad z=\ds\frac{l_1}{l_0},\quad
s\pr=\ds\frac{t_0}{t_1l_0^q}\quad$\\
\hline
$\quad U_{1,0}:\quad w=\ds\frac{l_0}{l_1},\quad y\pr=\ds\frac{t_1l_1^q}
{t_0}\quad $&
$\quad U_{1,1}:\quad w=\ds\frac{l_0}{l_1},\quad y=\ds\frac{t_0}{t_1l_1^q}\quad$
\end{tabular}}
\end{center}
\bigskip

We also denote
\[ V_{0}=\{l_{0}\not=0\}=U_{0,0}\cup U_{0,1},\qquad
V_{1}=\{l_{1}\not=0\}=U_{1,0}\cup U_{1,1}.\]

We shall call $z,s$ principal affine coordinates and $U_{0,0}$ principal
affine subset.
Since the changes of coordinates are holomorphic, the above affine
covering gives a structure of complex manifold of dimension 2 on
$\F_{k}$.\\

\begin{exer} If we consider the analogous construction of $\F_{q}$
with $\R$ instead of $\C$ we get $\F_{q}$=\emph{torus} for $q$ even and
$\F_{q}$=\emph{Klein bottle} for $q$ odd.\end{exer}

\begin{defn}\label{II.1.1}
For $q>0$ we set $\sigma_{\infty}=\{t_{1}=0\}$. Clearly
$\sigma_{\infty}$ is isomorphic to $\proj^{1}$.\end{defn}

\begin{prop}\label{II.1.2}
 $ \F_{0}$ is not homeomorphic to $\F_{1}$.\end{prop}

\begin{proof} Topologically $\F_{0}=S^{2}\times S^{2}$ and therefore
$H_{2}(\F_{0},\Z)=\Z[S^{2}\times\{p\}]\oplus\Z[\{p\}\times S^{2}]$,
where $p\in S^{2}$ and $[V]\in H_{2}$ denotes the homology class of
a closed subvariety $V\subset S^{2}\times S^{2}$ of real dimension 2.\\
The matrix of the intersection form $q\colon H_{2}\times H_{2}\to
H_{0}=\Z$ is
\[\left(\begin{array}{cc}0&1\\1&0\end{array}\right)\]
and therefore $q(a,a)$ is even for every $a\in H_{2}(\F_{0},\Z)$.\\
Consider the following subvarieties of $\F_{1}$:
\[ \sigma=\{t_{0}=0\},\qquad \sigma\pr=\{t_{0}=l_{0}t_{1}\}.\]
$\sigma$ and $\sigma\pr$ intersect transversely at the point
$t_{0}=l_{0}=0$ and therefore their intersection product is equal to
$q([\sigma],[\sigma\pr])=\pm 1$. On the other hand the continuous map
\[ r:(\F_{1}-\sigma_{\infty})\times [0,1]\to (\F_{1}-\sigma_{\infty}),\qquad
r((l_{0},l_{1},t_{0},t_{1}),a)=(l_{0},l_{1},at_{0},t_{1})\]
shows that $\sigma$ is a deformation retract of
$(\F_{1}-\sigma_{\infty})$. Since $r_{1}\colon \sigma\pr\to \sigma$ is an
isomorphism we have $[\sigma]=[\sigma\pr]\in
H_{2}(\F_{1}-\sigma_{\infty},\Z)$ and then a fortiori
$[\sigma]=[\sigma\pr]\in
H_{2}(\F_{1},\Z)$. Therefore $q([\sigma],[\sigma])=\pm 1$ is not even
and $\F_{0}$ cannot be homeomorphic to $\F_{1}$.\end{proof}

It is easy to find projective embeddings of the surfaces $\F_{q}$;

\begin{ex}\label{II.1.3} The Segre-Hirzebruch surface $\F_{q}$
is isomorphic to the subvariety $X\subset \proj^{q+1}\times\proj^{1}$ of
equation
\[ u_{0}(x_{1},x_{2},\ldots,x_{q})=u_{1}(x_{2},x_{3},\ldots,x_{q+1}),\]
where $x_{0},\ldots,x_{q+1}$ and $u_{0},u_{1}$ are homogeneous
coordinates in $\proj^{q+1}$ and $\proj^{1}$ respectively.\\
An isomorphism $\F_{q}\to X$ is given by:
\[u_{0}=l_{0},\quad u_{1}=l_{1},\quad  x_{0}=t_{0},\quad
x_{i}=t_{1}l_{0}^{i-1}l_{1}^{q+1-i},\, i=1,\ldots q+1.\]
\end{ex}

Denote by $T\to \F_{q}$ the holomorphic tangent bundle, in order
to compute the spaces $H^{0}(\F_{q},T)$ and $H^{1}(\F_{q},T)$ we
first notice that the open subsets $V_{0},V_{1}$ are isomorphic to
$\C\times\proj^{1}$. Explicit isomorphisms are given by
\[ V_{0}\to \C_{z}\times\proj^{1},\qquad
(l_{0},l_{1},t_{0},t_{1})\mapsto
\left(z=\frac{l_{1}}{l_{0}},[t_{0},t_{1}]\right),\]
\[ V_{1}\to \C_{w}\times\proj^{1},\qquad (l_{0},l_{1},t_{0},t_{1})\mapsto
\left(w=\frac{l_{0}}{l_{1}},[t_{0},t_{1}]\right).\] According to
Example~\ref{I.2.9} $H^{1}(V_{i},T)=0$, $i=0,1$, and then
$H^{0}(\F_{q},T)$ and $H^{1}(\F_{q},T)$ are isomorphic,
respectively, to the kernel and the cokernel of the \v{C}ech
differential
\[ H^{0}(V_{0},T)\oplus H^{0}(V_{1},T)\mapor{d} H^{0}(V_{0}\cap
V_{1},T),\qquad d(\chi,\eta)=\chi-\eta.\]
In the  affine coordinates $(z,s)$, $(w,y)$ we have that:
\begin{enumerate}
\item $H^{0}(V_{0},T)$ is the free $\Oh(\C_{z})$-module generated
by $\desude{~}{z}$, $\desude{~}{s}$, $s\desude{~}{s}$, $s^{2}\desude{~}{s}$.
\item $H^{0}(V_{1},T)$ is the free $\Oh(\C_{w})$-module generated
by $\desude{~}{w}$, $\desude{~}{y}$, $y\desude{~}{y}$, $y^{2}\desude{~}{y}$.
\item $H^{0}(V_{0}\cap V_{1},T)$ is the free $\Oh(\C_{z}-\{0\})$-module generated
by $\desude{~}{z}$, $\desude{~}{s}$, $s\desude{~}{s}$, $s^{2}\desude{~}{s}$.
\end{enumerate}
The change of coordinates is given by
\renewcommand\arraystretch{2}
\[ \left\{\begin{array}{l} z=w^{-1}\\
s=y^{-1}w^{q}\end{array}\right.\qquad\qquad
\left\{\begin{array}{l} w=z^{-1}\\ y=s^{-1}z^{-q}\end{array}\right.\]
and then
\[\left\{\begin{array}{l}
\desude{~}{w}=-z^{2}\desude{~}{z}+qy^{-1}w^{-q-1}\desude{~}{s}=
-z^{2}\desude{~}{z}+qzs\desude{~}{s}\\
\desude{~}{y}=-y^{-2}w^{q}\desude{~}{s}=
-z^{q}s^{2}\desude{~}{s}
\end{array}\right.\]
\renewcommand\arraystretch{1}
\[ d\left(\sum_{i\ge 0}z^{i}\left(a_{i}\desude{~}{z}+
(b_{i}+
c_{i}s+d_{i}s^{2})\desude{~}{s}\right),
\sum_{i\ge 0}w^{i}\left(\alpha_{i}\desude{~}{w}+
(\beta_{i}+
\gamma_{i}y+\delta_{i}y^{2})\desude{~}{y}\right)\right)=\]
\begin{multline*}
=\sum_{i\ge 0}z^{i}\left(a_{i}\desude{~}{z}+b_{i}\desude{~}{s}+
c_{i}s\desude{~}{s}+d_{i}s^{2}\desude{~}{s}\right)\\
+
\sum_{i\ge 0}z^{-i}\left(\alpha_{i}
\left(z^{2}\desude{~}{z}-qzs\desude{~}{s}\right)
+\beta_{i}s^{2}z^{q}\desude{~}{s}+
\gamma_{i}s\desude{~}{s}+\delta_{i}z^{-q}\desude{~}{s}\right)
\end{multline*}
An
easy computation gives the following

\begin{lem}\label{II.1.4}
\[\sum_{i\in\Z}z^{i}\left(a_{i}\desude{~}{z}+b_{i}\desude{~}{s}+
c_{i}s\desude{~}{s}+d_{i}s^{2}\desude{~}{s}\right)\in
H^{0}(V_{0}\cap V_{1}, T)\] belongs to the image of the \v{C}ech
differential if and only if $b_{-1}=b_{-2}=\ldots=b_{-q+1}=0$. In
particular the vector fields
\[ z^{-h}\desude{~}{s}\in H^{0}(V_{0}\cap
V_{1}, T),\qquad h=1,\ldots,q-1\]
represent a basis of $H^{1}(\F_{q},T)$ and
then $h^{1}(\F_{q},T)=\max(0,q-1)$.\end{lem}

\begin{exer} Prove that $h^{0}(\F_{q},T)=\max(6,q+5)$.\end{exer}

\begin{thm}\label{II.1.5} If $a\not=b$ then $\F_{a}$ is not
biholomorphic to $\F_{b}$.\end{thm}

\begin{proof} Assume $a>b$. If $a\ge 2$ then the dimension of
$H^{1}(\F_{a},T_{\F_{a}})$ is bigger than the dimension of
$H^{1}(\F_{b},T_{\F_{b}})$. If $a=1$, $b=0$ we apply
Proposition~\ref{II.1.2}.\end{proof}

We will show in \ref{II.3.3}
that $\F_{a}$ is diffeomorphic to $\F_{b}$ if and
only if $a-b$ is even.

\bigskip

\section{Decomposable bundles on projective spaces}

For $n>0$, $a\in\Z$  we define
\[ \OP{n}{a}=(\C^{n+1}-0)\times\C/\C^{*},\]
where the action of the multiplicative group $\C^{*}=\C-0$ is
\[ \lambda(l_{0},\ldots,l_{n},t)=
(\lambda l_{0},\ldots,\lambda l_{n},\lambda^{a}t),\qquad
\lambda\in\C^{*}.\]

The projection $\OP{n}{a}\to\proj^{n}$, $[l_{0},\ldots,l_{n},t]\mapsto
[l_{0},\ldots,l_{n}]$, is a holomorphic line bundle.
Notice that
$\Oh_{\proj^{n}}=\OP{n}{0}\to\proj^{n}$
is the trivial vector bundle of rank 1.\\

The obvious projection maps give a commutative diagram
\renewcommand\arraystretch{1.3}
\[\begin{array}{ccc}
(\C^{n+1}-0)\times\C&\mapor{}&\OP{n}{a}\\
\mapver{}&&\mapver{p}\\
(\C^{n+1}-0)&\mapor{\pi}&\proj^{n}\end{array}\]
\renewcommand\arraystretch{1}
inducing an isomorphism between $(\C^{n+1}-0)\times\C$ and the fibred
product of $p$ and $\pi$; in particular for every open subset
$U\subset\proj^{n}$ the space $H^{0}(U,\OP{n}{a})$ is naturally
isomorphic to the space of holomorphic maps $f\colon
\pi^{-1}(U)\to\C$ such that $f(\lambda x)=\lambda^{a}f(x)$ for every
$x\in \pi^{-1}(U)$, $\lambda\in\C^{*}$.\\
If $U=\proj^{n}$ then, by Hartogs' theorem, every holomorphic map
$f\colon \pi^{-1}(U)\to \C$ can be extended to a function
$f\colon\C^{n+1}\to\C$. Considering the power series expansion of
$f$ we get a natural isomorphism between
$H^{0}(\proj^{n},\OP{n}{a})$ and the space of homogeneous
polynomials of degree $a$ in the homogeneous
coordinates $l_{0},\ldots,l_{n}$.

\begin{exer} Prove that
$h^{0}(\proj^{n},\OP{n}{a})=\binom{n+a}{n}$.\end{exer}

\begin{exer} Under the isomorphism $\sigma_{\infty}=\proj^{1}$ we have
$N_{\sigma_{\infty}/\F_{q}}=\OP{1}{-q}$.\end{exer}

On the open set $U_{i}=\{l_{i}\not=0\}$ the section $l_{i}^{a}\in
H^{0}(U_{i},\OP{n}{a})$ is nowhere 0 and then gives a trivialization
of $\OP{n}{a}$ over $U_{i}$. The multiplication maps
\[ H^{0}(U_{i},\OP{n}{a})\otimes H^{0}(U_{i},\OP{n}{b})\to
H^{0}(U_{i},\OP{n}{a+b}),\qquad f\otimes g\mapsto fg,\]
give natural isomorphisms of line bundles
\[ \OP{n}{a}\otimes \OP{n}{b}=\OP{n}{a+b},\quad
\HOM(\OP{n}{a},\OP{n}{b})=\OP{n}{b-a}\]
(In particular $\OP{n}{a}\dual=\OP{n}{-a}$.)

\begin{defn}\label{II.2.1} A holomorphic vector bundle
$E\to\proj^{n}$ is called \emph{decomposable} if it is isomorphic
to a direct sum of line bundles of the form $\OP{n}{a}$.\\
Equivalently a vector bundle is decomposable if it is isomorphic
to
\[ (\C^{n+1}-0)\times\C^{r}/\C^{*}\to (\C^{n+1}-0)/\C^{*}=\proj^{n},\]
where the action is $\lambda(l_{0},\ldots,l_{n},t_{1},\ldots,t_{r})=
(\lambda l_{0},\ldots,\lambda l_{n},\lambda^{a_{1}}t_{1},\ldots,
\lambda^{a_{r}}t_{r})$.
\end{defn}

\begin{lem}\label{II.2.2} Two decomposable  bundles of rank $r$,
$E=\oplus_{i=1}^{r}\OP{n}{a_{i}}$, $F=\oplus_{i=1}^{r}\OP{n}{b_{i}}$,
$a_{1}\le a_{2}\cdots\le a_{r}$, $b_{1}\le b_{2}\cdots\le b_{r}$,
are isomorphic if and only if $a_{i}=b_{i}$ for every
$i=1,\ldots,r$.\end{lem}

\begin{proof} Immediate from the formula
\[ h^{0}(\proj^{n},(\oplus_{i}\OP{n}{a_{i}})\otimes\OP{n}{s})=
\sum_i h^{0}(\proj^{n},\OP{n}{a_{i}+s})=\!\!
\sum_{\{i\mid a_{i}+s\ge 0\}}\binom{a_{i}+s+n}{n}.\]
\end{proof}

\begin{ex}\label{II.2.3}
If $n\ge 2$ not every holomorphic vector bundle is
decomposable. Consider for example the surjective  morphism
\[\phi\colon\oplus_{i=0}^{n}\OP{n}{1}e_{i}\to \OP{n}{2},\qquad
\sum f_ie_{i}\mapsto \sum f_{i}l_{i}.\]
We leave it as an exercise to show that the kernel of $\phi$ is not
decomposable (Hint: first prove that $\ker\phi$ is generated by the global
sections $l_{i}e_{j}-l_{j}e_{i}$).
\end{ex}

For every holomorphic vector bundle $E\to X$ on a complex manifold $X$
we denote by $\proj(E)\to X$ the projective bundle whose fibre over
$x\in X$ is $\proj(E)_{x}=\proj(E_{x})$. If $E\to X$ is trivial over
an open subset $U\subset X$ then also $\proj(E)$ is trivial over $U$;
this proves that $\proj(E)$ is a complex manifold and the projection
$\proj(E)\to X$ is proper.\\

\begin{ex}\label{II.2.4} For every $a,b\in\Z$,
$\proj(\OP{1}{a}\oplus\OP{1}{b})=\F_{|a-b|}$.\\
To see this it is not restrictive to assume $a\ge b$; we have
\[
\proj(\OP{1}{a}\oplus\OP{1}{b})=(\C^{2}-0)\times(\C^{2}-0)/\C^{*}\times\C^{*},\]
where the action is $(\lambda,\eta)(l_{0},l_{1},t_{0},t_{1})=
(\lambda l_{0},\lambda l_{1},\lambda^{a}\eta t_{0},\lambda^{b}\eta
t_{1})$. Setting $\mu=\lambda^{b}\eta$ we recover the definition
of
$\F_{a-b}$.\\
More generally if $E\to X$ is a vector bundle and $L\to X$ is a line
bundle then $\proj(E\otimes L)=\proj(E)$.
\end{ex}

\begin{ex}\label{II.2.5}
The tangent bundle $T_{\proj^{1}}$ is isomorphic to $\OP{1}{2}$.
Let $l_{0},l_{1}$ be homogeneous coordinates on $\proj^{1}$;
$s=\ds\frac{l_{1}}{l_{0}}$, $t=\ds\frac{l_{0}}{l_{1}}$ are
coordinates on $U_{0}=\{l_{0}\not=0\}$, $U_{1}=\{l_{1}\not=0\}$
respectively. The sections of $T_{\proj^{1}}$ over an open set $U$
correspond to pairs $\left(f_{0}(s)\ds\frac{\de~}{\de s},
f_{1}(t)\ds\frac{\de~}{\de t}\right)$, $f_{i}\in\Oh(U\cap U_{i})$,
such that
$f_{1}(t)=-t^{2}f_{0}(t^{-1})$.\\
The isomorphism $\phi\colon \OP{1}{2}\to T_{\proj^{1}}$ is given by
$\phi(l_{0}^{a}l_{1}^{2-a})=
\left(s^{2-a}\ds\frac{\de~}{\de s},-t^{a}\ds\frac{\de~}{\de
t}\right)$.
\end{ex}

\begin{thm}[Euler exact sequence]\label{II.2.6}
On the projective space $\proj^{n}$ there exists an exact sequence
of vector bundles
\[ 0\mapor{}\Oh_{\proj^{n}}
\mapor{\sum l_{i}\frac{\de~}{\de l_{i}}}\oplus_{i=0}^{n}\OP{n}{1}
\desude{~}{l_{i}}~\mapor{\phi}~T_{\proj^{n}}\mapor{}0,\]
where on the affine open subset $l_{h}\not=0$, with coordinates
$s_{i}=\ds\frac{l_{i}}{l_{h}}$, $i\not=h$,
\renewcommand\arraystretch{3}
\[\left\{\begin{array}{ll}
\phi\left(l_{i}\desude{~}{l_{j}}\right)=s_{i}\desude{~}{s_{j}}\quad&i,j\not=h\\
\phi\left(l_{h}\desude{~}{l_{j}}\right)=\desude{~}{s_{j}}&j\not=h
\end{array}\right.,\quad
\left\{\begin{array}{ll}
\phi\left(l_{i}\desude{~}{l_{h}}\right)=-\ds\sum_{j\not=h}s_{i}s_{j}\desude{~}{s_{j}}
\quad&i\not=h\\
\phi\left(l_{h}\desude{~}{l_{h}}\right)=
-\ds\sum_{j\not=h}s_{j}\desude{~}{s_{j}}\end{array}
\right.\]
\renewcommand\arraystretch{1}
\end{thm}

\begin{proof} The surjectivity of $\phi$ is clear. Assume
$\phi\left(\sum_{i,j}a_{ij}l_{i}\desude{~}{l_{j}}\right)=0$, looking at the
quadratic terms in the set $l_{h}\not=0$ we get $a_{ih}=0$ for
every $i\not=h$. In the open set $l_{0}\not=0$ we have
\[\phi\left(\sum_{i}a_{ii}l_{i}\desude{~}{l_{i}}\right)=
\sum_{i=1}^{n}a_{ii}s_{i}\desude{~}{s_{i}}-
\sum_{i=1}^{n}a_{00}s_{i}\desude{~}{s_{i}}=0\]
and then the matrix $a_{ij}$ is a multiple of the identity.
\end{proof}

\begin{rem}\label{II.2.7} It is possible to prove that the map
$\phi$ in the Euler exact sequence is surjective at the level of
global sections, this gives an isomorphism
\[ H^{0}(\proj^{n},T_{\proj^{n}})=\quot{gl(n+1,\C)}{\C Id}
=pgl(n+1,\C)=T_{Id}PGL(n+1,\C).\]
Moreover it is possible to prove that every biholomorphism of $\proj^n$ is
a projectivity and the integration of holomorphic vector fields
corresponds to the exponential map in the complex Lie group $PGL(n+1,\C)$.
\end{rem}

\begin{exer} Use the Euler exact sequence and the
surjectivity of $\phi$ on global sections
to prove that for every
$n\ge 2$ the tangent bundle of $\proj^{n}$ is not
decomposable.\end{exer}

\begin{cor}\label{II.2.8} The canonical bundle of $\proj^{n}$ is
$K_{\proj^{n}}=\OP{n}{-n-1}$.\end{cor}

\begin{proof} From the Euler exact sequence we have
\[ \external{}{n}T_{\proj^{n}}\otimes \Oh_{\proj^{n}}=
\external{}{n+1}\left(\oplus_{i=0}^{n}\OP{n}{1}\right)=\OP{n}{n+1}\]
and then
$K_{\proj^{n}}=(\external{}{n}T_{\proj^{n}})\dual=\OP{n}{-n-1}$.
\end{proof}

\begin{exer} Prove that
$h^{n}(\proj^{n},\OP{n}{-a})=\dbinom{a-1}{n}$.\end{exer}

\begin{lem}\label{II.2.9}
Let $E\to\proj^{1}$ be a holomorphic vector bundle of rank $r$.
If:\begin{enumerate}

\item $H^{0}(\proj^{1},E(s))=0$ for $s<<0$, and

\item There exists a constant $c\in\N$ such that
$h^{0}(\proj^{1},E(s))\ge rs-c$ for $s>>0$.
\end{enumerate}
Then
$E$ is decomposable.\end{lem}

\begin{proof}
Using the assumptions 1 and 2 we may construct recursively a sequence
$a_{1},\ldots,a_{r}\in\Z$ and sections $\alpha_{i}\in
H^{0}(\proj^{1},E(a_{i}))$ such that:\begin{enumerate}
\item $a_{h+1}$ is the minimum integer $s$ such that the map
\[ \oplus_{i=1}^{h}\alpha_{i}\colon \somdir{i=1}{h}
H^{0}(\proj^{1},\OP{1}{s-a_{i}})\to H^{0}(\proj^{1},E(s))\]
is not surjective.
\item $\alpha_{h+1}$ does not belong to the image of
\[ \oplus_{i=1}^{h}\alpha_{i}\colon \somdir{i=1}{h}
H^{0}(\proj^{1},\OP{1}{a_{h+1}-a_{i}})\to H^{0}(\proj^{1},E(a_{h+1})).\]
\end{enumerate}
Notice that $a_{1}\le a_{2}\le\ldots\le a_{r}$.\\
We prove now by induction on $h$ that the morphism of vector bundles
\[ \oplus_{i=1}^{h}\alpha_i\colon \somdir{i=1}{h}
\OP{1}{-a_{i}}\to E\]
is injective on every fibre; this implies that
$\oplus_{i=1}^{r}\alpha_{i}\colon \bigoplus_{i=1}^{r}
\OP{1}{-a_{i}}\to E$ is an isomorphism.\\
For $h=0$ it is trivial.
Assume
$\oplus_{i=1}^{h}\alpha_{i}$ injective on fibres and let
$p\in\proj^{1}$. Choose homogeneous coordinates $l_{0},l_{1}$ such
that $p=\{l_{1}=0\}$ and set $s=l_{1}/l_{0}$.\\
Assume that  there exist $c_{1},\ldots,c_{h}\in\C$ such that
$\alpha_{h+1}(p)=\sum c_{i}(l_{0}^{a_{h+1}-a_{i}}\alpha_{i})(p)\in
E(a_{h+1})_{p}$.
If $e_{1},\ldots,e_{r}$ is a local frame for $E$ at $p$ we have locally
\[\alpha_{h+1}-\sum_{i=1}^{h} c_{i}l_{0}^{a_{h+1}-a_{i}}\alpha_{i}=
\sum_{j=1}^{r}f_{j}(s)l^{a_{h+1}}_{0}e_{j}\]
with $f_{j}(s)$ holomorphic functions such that $f_{j}(0)=0$.\\
Therefore $f_{j}(s)/s$ is still holomorphic and
$l_{0}^{-1}(\alpha_{h+1}-\sum
c_{i}l_{0}^{a_{h+1}-a_{i}}\alpha_{i})\in
H^{0}(\proj^{1},E(a_{h+1}-1))$, in contradiction with the minimality of
$a_{h+1}$.\end{proof}

\begin{thm}\label{II.2.10}
Let $0\mapor{}E\mapor{}F\mapor{}G\mapor{}0$ be an exact sequence
of holomorphic vector bundles on $\proj^{1}$.\begin{enumerate}
\item If $F,G$ are decomposable then also $E$ is decomposable.
\item If $E=\oplus\OP{1}{-a_{i}}$ then $\min(a_{i})$ is the
minimum integer $s$ such that $H^{0}(\proj^{1},F(s))\to
H^{0}(\proj^{1},G(s))$ is not injective.
\end{enumerate}\end{thm}

\begin{proof}
The kernel of $H^{0}(\proj^{1},F(s))\to H^{0}(\proj^{1},G(s))$ is
exactly $H^{0}(\proj^{1},E(s))$.\\
If $F=\oplus_{i=1}^{r}\OP{1}{b_{i}}$,
$G=\oplus_{i=1}^{p}\OP{1}{c_{i}}$ then for $s>>0$
$h^{0}(\proj^{1},F(s))=r(s+1)+\sum b_{i}$,
$h^{0}(\proj^{1},G(s))=p(s+1)+\sum c_{i}$ and then the rank of $E$
is $r-p$ and $h^{0}(\proj^{1},E(s))\ge (r-p)(s+1)+\sum b_{i}-\sum
c_{i}$. According to Lemma~\ref{II.2.9}, the vector bundle $E$ is
decomposable.
\end{proof}

We also state, without proof, the following

\begin{thm}\label{II.2.11}
\begin{enumerate}
\item Every holomorphic line bundle on $\proj^{n}$ is decomposable.
\item (Serre) Let $E$ be a holomorphic vector bundle on $\proj^{n}$,
then:
\begin{enumerate}
\item $H^{0}(\proj^{n},E(s))=0$ for $s<<0$.
\item $E(s)$ is generated by
global sections and $H^{p}(\proj^{n},E(s))=0$
for $p>0$, $s>>0$.\end{enumerate}
\item (Bott vanishing theorem) For every $0<p<n$:
\renewcommand\arraystretch{2}
\[ H^{p}(\proj^{n},\Omega^{q}(a))=
\left\{\begin{array}{ll}
\C\qquad&\hbox{ if  ~} p=q,\quad a=0\\
0&\hbox{otherwise}\end{array}\right.\]
\renewcommand\arraystretch{1}
Moreover
$H^{0}(\proj^{n},\Omega^{q}(a))=H^{n}(\proj^{n},\Omega^{n-q}(-a))\dual=0$
whenever $a<q$.
\end{enumerate}\end{thm}
\begin{proof} \cite{Kobook}\end{proof}

\bigskip

\section{Semiuniversal families of Segre-Hirzebruch surfaces}

Let $q>0$ be a fixed integer, define $M\subset
\C_{t}^{q-1}\times\proj^{1}_{l}\times\proj_{x}^{q+1}$ as the set
of points of homogeneous coordinates
$(t_{2},\ldots,t_{q},[l_{0},l_{1}],[x_{0},\ldots,x_{q+1}])$
satisfying the vectorial equation
\begin{equation}\label{II.3.1}
l_{0}(x_{1},x_{2},\ldots,x_{q})=
l_{1}(x_{2}-t_{2}x_{0},\ldots,x_{q}-t_{q}x_{0},x_{q+1}).
\end{equation}
We denote by $f\colon M\to\C^{q-1}$, $p\colon M\to
\C^{q-1}\times\proj^{1}_{l}$ the projections.

\begin{lem}\label{II.3.2} There exists a holomorphic vector bundle of
rank 2, $E\to \C^{q-1}\times\proj^{1}_{l}$ such that the map
$p\colon M\to \C^{q-1}\times\proj^{1}_{l}$ is
a smooth family  isomorphic  to
$\proj(E)\to \C^{q-1}\times\proj^{1}_{l}$.\end{lem}

\begin{proof} Let $\pi\colon \C^{q-1}\times\proj^{1}_{l}\to
\proj^{1}_{l}$ be the projection; define $E$ as the kernel of the
morphism of vector bundles over $\C^{q-1}\times\proj^{1}_{l}$
\[ \somdir{i=0}{q+1}\pi^{*}\Oh_{\proj^{1}}\mapor{A}
\somdir{i=1}{q}\pi^{*}\OP{1}{1},\]
\medskip
\[ A(t_{2},\ldots,t_{q},[l_{0},l_{1}])\left(\begin{array}{c}
x_{0}\\ x_{1}\\ \vdots\\ x_{q+1}\end{array}\right)=
\left(\begin{array}{c}
l_{0}x_{1}-l_{1}(x_{2}-t_{2}x_{0})\\ l_{0}x_{2}-l_{1}(x_{3}-t_{3}x_{0})\\
\vdots\\ l_{0}x_{q}-l_{1}x_{q+1}\end{array}\right).\]
\medskip
We first note that $A$ is surjective on every fibre, in fact for
fixed $t_{2},\ldots,t_{q},l_{0},l_{1}\in\C$, $A(t_{i},l_{j})$ is
represented by the matrix
\[\left(\begin{array}{cccccc}
t_{2}l_{1}&l_{0}&-l_{1}&\ldots&0&0\\
t_{3}l_{1}&0&l_{0}&\ldots&0&0\\
\vdots&\vdots&\vdots&\ddots&\vdots&\vdots\\
0&0&0&\ldots&l_{0}&-l_{1}\end{array}\right).\]
Since either $l_{0}\not=0$ or $l_{1}\not=0$ the above matrix has
maximal rank.\\
By definition we have that $M$ is the set of points of
$x\in \proj(\oplus_{i=0}^{q+1}\pi^{*}\Oh_{\proj^{1}})$
such that $A(x)=0$ and then $M=\proj(E)$.\end{proof}

For every $k\ge 0$ denote by $T_{k}\subset\C^{q-1}_{t}$ the subset
of points of coordinates $(t_{2},\ldots,t_{q})$ such that there
exists a nonzero $(q+2)$-uple of homogeneous polynomials of degree
$k$
\[(x_{0}(l_{0},l_{1}),\ldots,x_{q+1}(l_{0},l_{1}))\]
which satisfy
identically ($t$ being fixed) the Equation~\ref{II.3.1}.
Note that $t\in T_{k}$ if and only if there exists a nontrivial
morphism $\OP{1}{-k}\to E_{t}$ and then $t\in T_{k}$ if and only if
$-k\le -a$. Therefore $t\in T_{k}-T_{k-1}$ if and only if $a=k$.

\begin{lem}\label{II.3.3} In the notation above:\begin{enumerate}
\item $T_{0}=\{0\}$.
\item $T_{k}\subset T_{k+1}$.
\item If $2k+1\ge q$ then $T_{k}=\C^{q-1}$.
\item If $2k\le q$ and
$t\in T_{k}-T_{k-1}$ then $M_{t}=\F_{q-2k}$.
\end{enumerate}\end{lem}

\begin{proof} 1 and 2 are trivial.\\
Denoting by $S_{k}\subset\C[l_{0},l_{1}]$ the space of homogeneous
polynomials of degree $k$, $\dim_{\C}S_{k}=k+1$; interpreting
Equation~\ref{II.3.1} as a linear map (depending on the parameter
$t$) $A_{k}(t)\colon S_{k}^{q+2}\to S_{k+1}^{q}$, we have that
$t\in T_{k}$ if and only if $\ker
A_{k}(t)\not=0$.\\
Since $(q+2)(k+1)> q(k+2)$ whenever $2k>q-2$, item 3 follows
immediately.\\
Let $E_{t}$ be the restriction of the vector bundle $E$ to
$\{t\}\times\proj^{1}$, $E_{t}$ is the kernel of the surjective
morphism $A(t)\colon
\oplus_{i=0}^{q+1}\Oh_{\proj^{1}}\to\oplus_{i=1}^{q}\OP{1}{1}$.
According to Theorem~\ref{II.2.10}, $E_{t}$ is decomposable. Since
$\bigwedge^{2}E_{t}=\OP{1}{-q}$ we have $E_{t}=\OP{1}{-a}\oplus
\OP{1}{a-q}$ with $-a\le a-q$
and $M_{t}=\proj(E_{t})=\F_{q-2a}$.\\
\end{proof}

\begin{lem}\label{II.3.4} In the notation above
$(t_{2},\ldots,t_{q})\in T_{k}$ if and only if there exists a nonzero
triple $(x_{0},x_{1},x_{q+1})\in\oplus \C[s]$ of polynomials of
degree $\le k$ such that
\[ x_{q+1}=s^{q}x_{1}+x_{0}\left(\sum_{i=2}^{q}t_{i}s^{q+1-i}\right).\]
\end{lem}

\begin{proof} Setting $s=l_{0}/l_{1}$ we have by definition that
$(t_{2},\ldots,t_{q})\in T_{k}$ if and only if there exists a
nontrivial sequence $x_{0},\ldots,x_{q+1}\in \C[s]$ of polynomials
of degree $\le k$ such that
$x_{i+1}=sx_{i}+t_{i+1}x_{0}$ for
every $i=1,\ldots,q$ ($t_{q+1}=0$ by convention).
Clearly this set of equation is equivalent
to $x_{i+1}=s^{i}x_{1}+x_{0}\sum_{j=1}^{i}t_{j+1}s^{i-j}$.\\
Given $x_{0},x_{1},x_{q+1}$ as in the statement, we can define
recursively $x_{i}=s^{-1}(x_{i+1}-t_{i+1}x_{0})$ and the sequence
$x_{0},\ldots,x_{q-1}$ satisfies the defining equation of
$T_{k}$.\end{proof}

\begin{cor}\label{II.3.5}
$(t_{2},\ldots,t_{q})\in T_{k}$ if and only if  the
$(q-k-1)\times (k+1)$ matrix
$B_{k}(t)_{ij}=(t_{q-k-i+j})$
has rank $\le k$.\end{cor}

\begin{proof} If $2k+1\le q$ then $T_{k}=\C^{q-1}$, $q-k-1\le k$ and
the result is trivial: thus it is not
restrictive to assume $k+1\le q-k-1$ and then $\rank B_{k}(t)\le
k$ if and only if $\ker B_{k}(t)\not=0$.\\
We note that if $x_{0},x_{1},x_{q+1}$ satisfy the equation
$x_{q+1}=s^{q}x_{1}+x_{0}(\sum_{i=2}^{q}t_{i}s^{q+1-i})$ then
$x_{1},x_{q+1}$ are uniquely determined by $x_{0}$; conversely
a polynomial $x_{0}(s)$ of degree $\le k$  can be extended to a
solution of the equation if and only if all the coefficients of
$s^{k+1},s^{k+2},\ldots,s^{q-1}$ in the polynomial
$x_{0}(\sum_{i=2}^{q}t_{i}s^{q+1-i})$ vanish.
Writing $x_{0}=a_{0}+a_{1}s+\ldots+a_{k}s^{k}$, this last
condition is equivalent to $(a_{0},\ldots,a_{k})\in\ker B_{k}(t)$.
\end{proof}

Therefore $T_{k}$ is defined by the vanishing of the
$\dbinom{q-k-1}{k+1}$
minors of $B_{k}(t)$, each one of which is a homogeneous polynomial of
degree $k+1$ in $t_{2},\ldots,t_{q}$. In particular $T_{k}$ is an
algebraic cone.\\

As an immediate consequence of Corollary~\ref{II.3.5} we have that
for $q\ge 2$, $0<2k\le q$, the subset $\{t_{k+1}\not=0,
t_{k+2}=t_{k+3}=\ldots=t_{q}=0\}$ is contained in $T_{k}-T_{k-1}$.
In particular $\F_{q}$ is diffeomorphic to $\F_{q-2k}$ for every
$k\le q/2$.

\begin{prop}\label{II.3.6} If $2k<q$ then
$T_{k}$ is an irreducible affine variety
of dimension $2k$.\end{prop}

\begin{proof} Denote
\[ Z_{k}=\{([v],t)\in \proj^{k}\times\C^{q-1}\,| v\in\C^{k+1}-0,\,
B_{k}(t)v=0\}\]
and by  $p\colon Z_{k}\to T_{k}$ the projection on the second factor.
$p$ is surjective and if $t_{k+1}=1$, $t_{i}=0$ for $i\not=k+1$, then
$B_{k}(t)$ has rank $k$ and $p^{-1}(t)$ is one point. Therefore it is
sufficient to prove that $Z_{k}$ is an irreducible variety of
dimension $2k$.\\
Let $\pi\colon Z_{k}\to \proj^{k}$ be the projection.
We have $([a_{0},\ldots,a_{k}], (t_{2},\ldots,t_{q}))\in Z_{k}$ if and
only if for every $i=1,\ldots,q-k-1$
\[0=\sum_{j=0}^{k}t_{i+1+j}a_{j}=\sum_{l=2}^{q}t_{l}a_{l-i-1},\]
where $a_{l}=0$ for $l<0$, $l>k$ and then the fibre over
$[a_{0},\ldots,a_{k}]$ is the kernel of the matrix $A_{ij}=(a_{j-i-1})$
$i=1,\ldots,q-k-1$, $j=2,\ldots,q$. Since at least one $a_{i}$ is
$\not=0$ the rank of $A_{ij}$ is exactly $q-k-1$ and then the fibre is
a vector subspace of dimension $k$. By a general result in algebraic
geometry \cite{Shafarevich},\cite{IntroGA} $Z$ is an irreducible variety
of dimension $2k$.\end{proof}

\begin{thm}\label{semipersegre} 
In the above notation 
 the Kodaira-Spencer map $\KS_{f}\colon
T_{0,\C^{q-1}}\to H^{1}(M_{0},T_{M_{0}})$ is
bijective for every $q\ge 1$ and therefore, by completeness theorem \ref{I.5.1}, 
deformation $\F_{q}\to M\to (\C^{q-1},0)$ is semiuniversal.\end{thm}

\begin{proof}
We have seen that $M_{0}=\F_{q}$.
Let $V_{0},V_{1}\subset
\F_{q}$ be the open subset defined in Section~\ref{sec:1}.
Denote $M_{i}\subset
M$ the open subset $\{l_{i}\not=0\}$, $i=0,1$.\\
We have an isomorphism $\phi_{0}\colon \C^{q-1}\times V_{0}\to M_{0}$,
commuting with the projections onto $\C^{q-1}$, given in the affine
coordinates $(z,s)$ by:
\[ l_{0}=1,\quad  l_{1}=z,\quad x_{0}=1,\quad
x_{h}=z^{q-h+1}s-\sum_{j=1}^{q-h}t_{h+j}z^{j}=z(x_{h+1}-t_{h+1}x_{0}),\quad
h>0.\]
Similarly there exists an isomorphism $\phi_{1}\colon \C^{q-1}\times
V_{1}\to M_{1}$,
\[ l_{0}=w,\quad  l_{1}=1,\quad x_{0}=y,\quad
x_{h}=w^{h-1}+y\sum_{j=2}^{h}t_{j}w^{h-j}=wx_{h-1}+t_{h}x_{0},\quad
h>0.\] In the intersection $M_{0}\cap M_{1}$ we have:
\renewcommand\arraystretch{2}
\[\left\{\begin{array}{l}
z=w^{-1}\\
s=\ds\frac{x_{q+1}}{x_{0}}=y^{-1}w^{q}+\ds\sum_{j=2}^{q}t_{j}w^{q+1-j}.
\end{array}\right.\]
\renewcommand\arraystretch{1}
According to Formula~\ref{KScociclo}, for every $h=2,\ldots,q$
\[ \KS_{f}\left(\desude{~}{t_{h}}\right)=
\desude{w^{-1}}{t_{h}}\desude{~}{z}+
\desude{(y^{-1}w^{q}+\sum_{j=2}^{q}t_{j}w^{q+1-j})}{t_{h}}\desude{~}{s}=
z^{h-q-1}\desude{~}{s}.\]
\end{proof}

\bigskip

\section{Historical survey,~\ref{CAP:SEGRE}}

One of the most famous theorems in deformation theory
(at least in algebraic geometry) is the
stability theorem of submanifolds proved by Kodaira in 1963.

\begin{defn}\label{II.4.1} Let $Y$ be a closed submanifold of a compact complex
manifold $X$. $Y$ is called \emph{stable} if for every deformation
$X\mapor{i}\sX\mapor{f}(B,0)$ there exists a deformation
$Y\mapor{j}\sY\mapor{g}(B,0)$ and a commutative diagram of holomorphic
maps
\[\xymatrix{Y\ar[r]^{j}\ar[d]_{i_{|X}}&\sY\ar[dl]\ar[d]^{g}\\
\sX\ar[r]^{f}&B}\]
\end{defn}

The same argument used in Example~\ref{I.4.9} shows that
$\sigma_{\infty}\subset \F_{q}$ is not stable for every $q\ge 2$, while
$\sigma_{\infty}\subset \F_{1}$ is stable because $\F_{1}$ is rigid.\\

\begin{thm}\label{II.4.2} \emph{(Kodaira stability theorem for
submanifolds, \cite{Kod1})}\\
Let $Y$ be a closed submanifold of a compact complex
manifold $X$. If $H^{1}(Y,N_{Y/X})=0$ then $Y$ is stable.
\end{thm}

Just to check Theorem~\ref{II.4.2} in a concrete case, note that
$h^{1}(\sigma_{\infty},N_{\sigma_{\infty}/\F_{q}})=\max(0,q-1)$.\\

Theorem~\ref{II.4.2} has been generalized to arbitrary holomorphic
maps of compact complex manifolds in a series of papers by Horikawa
\cite{Horidhm}.

\begin{defn}\label{II.4.3} Let $\alpha\colon Y\to X$ be a holomorphic
    map of compact
complex manifolds. A \emph{deformation of $\alpha$} over a germ $(B,0)$
is a commutative diagram of holomorphic maps
\renewcommand\arraystretch{1.2}
\[\begin{array}{ccccc}
Y&\mapor{i}&\sY&\mapor{f}&B\\
\mapver{\alpha}&&\mapver{}&&\mapver{Id}\\
X&\mapor{j}&\sX&\mapor{g}&B\end{array}\]
\renewcommand\arraystretch{1}
where $Y\mapor{i}\sY\mapor{f}(B,0)$ and $X\mapor{j}\sX\mapor{g}(B,0)$
are deformations of $Y$ and $X$ respectively.\end{defn}

\begin{defn}\label{II.4.4} In the notation of \ref{II.4.3}, the map $\alpha$
is called:
\begin{enumerate}
\item  \emph{Stable} if every deformation of $X$ can be extended
to a deformation of $\alpha$.

\item  \emph{Costable} if every deformation of $Y$ can be extended
to a deformation of $\alpha$.
\end{enumerate}\end{defn}

Consider two locally finite coverings $\sU=\{U_{a}\}$,
$\sV=\{V_{a}\}$, $a\in \sI$, $Y=\cup U_{a}$, $X=\cup V_{a}$ such that
$U_{a}$, $V_{a}$ are biholomorphic to polydisks and
$\alpha(U_{a})\subset V_{a}$ for every $a$ ($U_{a}$ is allowed to be
the empty set).

Given $a\in \sI$ and local coordinate systems
$(z_{1},\ldots,z_{m})\colon U_{a}\to\C^{m}$,
$(u_{1},\ldots,u_{n})\colon V_{a}\to\C^{n}$ we have
linear morphisms of vector spaces
\renewcommand\arraystretch{3.5}
\[\begin{array}{ll}
\alpha^{*}\colon \Gamma(V_{a},T_{X})\to
\Gamma(U_{a},\alpha^{*}T_{X}),\qquad&
\alpha^{*}\left(\sum_{i}g_{i}\desude{~}{u_{i}}\right)=
\sum_{i}\alpha^{*}(g_{i})\desude{~}{u_{i}}\\
\alpha_{*}\colon \Gamma(U_{a},T_{Y})\to
\Gamma(U_{a},\alpha^{*}T_{X}),\qquad&
\alpha_{*}\left(\sum_{i}h_{i}\desude{~}{z_{i}}\right)=
\sum_{i,j}h_{i}\desude{u_{j}}{z_{i}}\desude{~}{u_{j}}\end{array}\]
\renewcommand\arraystretch{1}

Define $\H^{*}(\alpha_{*})$ as the cohomology of the complex
\[ 0\mapor{}C^{0}(\sU,T_{Y})\mapor{d_{0}}
C^{1}(\sU,T_{Y})\oplus C^{0}(\sU,\alpha^{*}T_{X})\mapor{d_{1}}\ldots\]
where $d_{i}(f,g)=(df,dg+(-1)^{i}\alpha_{*}f)$, being $d$ the usual
\v{C}ech differential.\\

Similarly define $\H^{*}(\alpha^{*})$ as the cohomology of the complex
\[ 0\mapor{}C^{0}(\sV,T_{X})\mapor{d_{0}}
C^{1}(\sV,T_{X})\oplus C^{0}(\sU,\alpha^{*}T_{X})\mapor{d_{1}}\ldots\]
where $d_{i}(f,g)=(df,dg+(-1)^{i}\alpha^{*}f)$.

\begin{thm}[Horikawa] The groups $\H^{k}(\alpha_{*})$ and
$\H^{k}(\alpha^{*})$ do not depend on the choice of the coverings
$\sU,\sV$. Moreover:
\begin{enumerate}
\item  If $\H^{2}(\alpha_{*})=0$ then $\alpha$ is stable.

\item  If $\H^{2}(\alpha^{*})=0$ then $\alpha$ is costable.
\end{enumerate}
\end{thm}

\begin{exer} Give a Dolbeault-type definition of the groups
$\H^{k}(\alpha_{*})$, $\H^{k}(\alpha^{*})$.\end{exer}

\begin{exer} If $\alpha\colon Y\to X$ is a regular embedding then
 $\H^{k}(\alpha_{*})=H^{k-1}(Y,N_{Y/X})$. (Hint:
take $U_{a}=V_{a}\cap Y$, and local systems of coordinates
$u_{1},\ldots,u_{n}$ such that $Y=\{u_{m+1}=\ldots=u_{n}=0\}$. Then
prove that the projection maps
$C^{k+1}(\sU,T_{Y})\oplus C^{k}(\sU,\alpha^{*}T_{X})
\to C^{k}(\sU,N_{Y/X})$ give a quasiisomorphism of complexes.\end{exer}

The following (non trivial) exercise is reserved to experts in
algebraic geometry:
\begin{exer} Let $\alpha\colon Y\to Alb(Y)$
be the Albanese map of a complex projective manifold $Y$. If
$X=\alpha(Y)$ is a curve then $\alpha\colon Y\to X$ is
costable.\end{exer}


\chapter[~Analytic singularities]{Analytic singularities}
\label{CAP:SING}
\piede

Historically, a major step in deformation theory has been
the introduction of  deformations of complex manifolds over (possibly
non reduced) analytic singularities.\\
This chapter is a short introductory course on analytic algebras
and analytic singularities; moreover  we give an elementary proof 
of the Nullstellenstaz for the ring
$\C\{z_1,\ldots,z_n\}$ of convergent complex power series.\\
Quite important in deformation theory are the smoothness criterion
\ref{V.1.12} and the two dimension bounds \ref{V.4.26} and
\ref{V.4.28}.

\bigskip
\section[~~Analytic algebras]{Analytic algebras}

Let $\C\{z_{1},\ldots,z_{n}\}$ be the ring of convergent power series
with complex coefficient.  Every $f\in \C\{z_{1},\ldots,z_{n}\}$
defines a holomorphic function in a nonempty open neighbourhood $U$
of $0\in\C^{n}$; for notational simplicity we still
denote by $f\colon U\to \C$ this function.\\

If $f$ is a holomorphic function in a neighbourhood of $0$ and
$f(0)\not=0$ then ${1}/{f}$ is holomorphic in a (possibly smaller)
neighbourhood of $0$. This implies that $f$ is invertible in
$\C\{z_{1},\ldots,z_{n}\}$ if and only if $f(0)\not=0$ and
therefore $\C\{z_{1},\ldots,z_{n}\}$ is a local ring with maximal
ideal $\ide{m}=\{f\mid f(0)=0\}$. The ideal $\ide{m}$ is
generated by $z_{1},\ldots,z_{n}$.

\begin{defn}
The \emph{multiplicity} of a power series $f\in \C\{z_{1},\ldots,z_{n}\}$
is defined as
\[ \mu(f)=\sup\{s\in\N\mid f\in \ide{m}^{s}\}\in\N\cup \{+\infty\}.\]
The \emph{valuation} $\nu(S)$ of a nonempty subset $S\subset
\C\{z_{1},\ldots,z_{n}\}$ is
\[ \nu(S)=\sup\{s\in\N\mid S\subset \ide{m}^{s}\}=
\inf\{\mu(f)\mid  f\in S\}\in\N\cup \{+\infty\}.\]
\end{defn}
We note that $\nu(S)=+\infty$ if and only if $S=\{0\}$ and
$\mu(f)$ is the smallest integer $d$ such that the power series
expansion of $f$ contains a
nontrivial homogeneous part of degree $d$.

The local ring $\C\{z_{1},\ldots,z_{n}\}$ has the following important properties:
\begin{itemize}
\item $\C\{z_{1},\ldots,z_{n}\}$ is Noetherian
(\cite[II.B.9]{Gu-Ro}, \cite{CAS}).

\item $\C\{z_{1},\ldots,z_{n}\}$ is a unique factorization domain
(\cite[II.B.7]{Gu-Ro}, \cite{CAS}).

\item $\C\{z_{1},\ldots,z_{n}\}$ is a Henselian ring
(\cite{IntroGA}, \cite{Ana-Ste}, \cite{CAS}). \item
$\C\{z_{1},\ldots,z_{n}\}$ is a regular local ring of dimension
$n$  (see e.g. \cite{A-M}, \cite{CAS}, \cite{Matsubook} for the
basics about dimension theory of local Noetherian ring).
\end{itemize}

We recall, for the reader's convenience, that the dimension of a
local Noetherian ring $A$ with maximal ideal $\ide{m}$ is the
minimum integer $d$ such that there exist $f_{1},\ldots,f_{d}\in
\ide{m}$ with the property $\sqrt{(f_{1},\ldots,f_{d})}=\ide{m}$.
In particular $\dim A=0$ if and only if $\sqrt{0}=\ide{m}$, i.e.
if and only if $\ide{m}$ is nilpotent.

We also recall that a  morphism of local rings $f\colon (A,\ide{m})\to
(B,\ide{n})$ is
called local if $f(\ide{m})\subset \ide{n}$.

\begin{defn}\label{V.1.2}
A local $\C$-algebra  is called an \emph{analytic algebra} if it is
isomorphic to $\C\{z_{1},\ldots,z_{n}\}/I$, for some $n\ge 0$
and some ideal $I\subset (z_{1},\ldots,z_{n})$.\\
We denote by ${\mathbf{An}}$ the category with objects the  analytic algebras
and morphisms  the local morphisms of $\C$-algebras.
\end{defn}

Every analytic algebra is a local Noetherian ring. Every local
Artinian $\C$-algebra with residue field $\C$ is an analytic algebra.\\

The ring $\C\{z_{1},\ldots,z_{n}\}$ is, in some sense, a free object
in the category ${\mathbf{An}}$ as explained in the following lemma

\begin{lem}\label{V.1.4}
Let $(R,\ide{m})$ be an analytic algebra.
Then the map
\[ \Mor_{{\mathbf{An}}}(\C\{z_{1},\ldots,z_{n}\},R)\to
\underbrace{\ide{m}\times\ldots\times\ide{m}}_{\hbox{n factors}},\qquad
f\mapsto (f(z_{1}),\ldots,f(z_{n}))\]
is bijective.\end{lem}
\begin{proof}
We first note that, by the lemma of Artin-Rees (\cite[10.19]{A-M}),
$\cap_{n}\ide{m}^{n}=0$ and then
every local homomorphism
$f\colon \C\{z_{1},\ldots,z_{n}\}\to R$ is uniquely determined by
its factorizations
\[f_{s}\colon \C\{z_{1},\ldots,z_{n}\}/(z_{1},\ldots,z_{n})^{s}\to
R/\ide{m}^{s}.\]
Since $\C\{z_{1},\ldots,z_{n}\}/(z_{1},\ldots,z_{n})^{s}$
is a $\C$-algebra generated by $z_{1},\ldots,z_{n}$,
every $f_{s}$ is uniquely determined by $f(z_{i})$; this proves the
injectivity.\\
For the surjectivity it is not restrictive to assume
$R=\C\{u_{1},\ldots,u_{m}\}$; given
$\phi=(\phi_{1},\ldots,\phi_{n})$, $\phi_{i}\in
\ide{m}$, let $U$ be an open subset
$0\in U\subset \C^{m}_{u}$ where the
$\phi_{i}=\phi_{i}(u_{1},\ldots,u_{m})$
are convergent power series.
The map $\phi=(\phi_{1},\ldots,\phi_{n})\colon
U\to\C^{n}$ is holomorphic, $\phi(0)=0$ and $\phi^{*}(z_{i})=\phi_{i}$.
\end{proof}

Another important and useful tool is the following

\begin{thm}[R\"{u}ckert's nullstellensatz]\label{V.1.6}
Let $I,J\subset\C\{z_{1},\ldots,z_{n}\}$ be proper ideals, then
\[\Mor_{{\mathbf{An}}}(\quot{\C\{z_{1},\ldots,z_{n}\}}{I},\C\{t\})=
\Mor_{{\mathbf{An}}}(\quot{\C\{z_{1},\ldots,z_{n}\}}{J},\C\{t\})\quad\iff\quad
\sqrt{I}=\sqrt{J},\]
where  the left equality is intended as equality of subsets of
$\Mor_{{\mathbf{An}}}(\C\{z_{1},\ldots,z_{n}\},\C\{t\})$
\end{thm}

A proof of Theorem~\ref{V.1.6} will be given in
Section~\ref{sec:ruckert}.

\begin{lem}\label{V.1.8}
Every analytic algebra is isomorphic to
$\quot{\C\{z_{1},\ldots,z_{k}\}}{I}$ for some $k\ge 0$ and some
ideal $I\subset (z_{1},\ldots,z_{k})^{2}$.\end{lem}
\begin{proof} Let $A=\C\{z_{1},\ldots,z_{n}\}/I$ be an analytic algebra
such that $I$ is not contained in $(z_{1},\ldots,z_{n})^{2}$; then
there exists $u\in I$ and an index  $i$ such that
$\desude{u}{z_{i}}(0)\not=0$. Up to permutation of indices we may
suppose $i=n$ and then, by inverse function theorem
$z_{1},\ldots,z_{n-1},u$ is a system of local holomorphic
coordinates. Therefore $A$ is isomorphic to
$\C\{z_{1},\ldots,z_{n-1}\}/I^{c}$, where $I^{c}$ is the kernel of the
surjective morphism
\[ \C\{z_{1},\ldots,z_{n-1}\}\to \C\{z_{1},\ldots,z_{n-1},u\}/I=A.\]
The conclusion follows by induction on $n$.\end{proof}

\begin{defn}\label{V.1.10}
An analytic algebra is called \emph{smooth} if it is isomorphic to
the power series algebra $\C\{z_{1},\ldots,z_{k}\}$ for some $k\ge
0$.\end{defn}

\begin{prop}\label{V.1.12}
Let $R=\C\{z_{1},\ldots,z_{k}\}/I$, $I\subset
(z_{1},\ldots,z_{k})^{2}$, be an analytic algebra.\\
The following conditions are equivalent:\begin{enumerate}
\item $I=0$.
\item $R$ is smooth.
\item for every surjective morphism of  analytic algebras
$B\to A$, the morphism
\[\Mor_{{\mathbf{An}}}(R,B)\to \Mor_{{\mathbf{An}}}(R,A)\]
is surjective.
\item for every $n\ge 2$ the morphism
\[\Mor_{{\mathbf{An}}}(R,\quot{\C\{t\}}{(t^{n})})\to
\Mor_{{\mathbf{An}}}(R,\quot{\C\{t\}}{(t^{2})})\]
is surjective.\end{enumerate}\end{prop}

\begin{proof} \implica{1}{2} and \implica{3}{4} are trivial,
while \implica{2}{3} is an immediate consequence of the Lemma~\ref{V.1.4}.\\
To prove \implica{4}{1}, assume $I\not=0$ and let $s=\nu(I)\ge 2$ be the
valuation of $I$, i.e. the greatest integer $s$ such that $I\subset
(z_{1},\ldots,z_{k})^{s}$: we claim that
$\Mor_{{\mathbf{An}}}(R,\C[t]/(t^{s+1}))\to \Mor_{{\mathbf{An}}}(R,\C[t]/(t^{2}))$ is not
surjective.\\
Choosing $f\in I-(z_{1},\ldots,z_{k})^{s+1}$, after a possible
generic linear  change of coordinates of the form $z_{i}\mapsto
z_{i}+a_{i}z_{1}$, $a_{2},\ldots,a_{k}\in\C$,
we may assume that $f$ contains the monomial $z^{s}_{1}$ with a
nonzero coefficient, say $f=cz_{1}^{s}+\ldots$;
let $\alpha\colon R\to \C[t]/(t^{2})$ be the
morphism defined by $\alpha(z_{1})=t$, $\alpha(z_{i})=0$ for $i>1$.\\
Assume that there exists $\beta\colon R\to \C[t]/(t^{s+1})$ that
lifts  $\alpha$, then
$\beta(z_{1})-t,\beta(z_{2}),\ldots,$ $\beta(z_{k})\in (t^{2})$ and
therefore $\beta(f)\equiv ct^{s}\pmod{t^{s+1}}$.\end{proof}

\begin{lem}\label{V.1.14}
For every analytic algebra $R$ with maximal ideal $\ide{m}$
there exist natural isomorphisms
\[ \Hom_{\C}(\quot{\ide{m}}{\ide{m}^{2}},\C)=\Der_{\C}(R,\C)=
\Mor_{{\mathbf{An}}}(R,\quot{\C[t]}{(t^{2})}).\]
\end{lem}

\begin{proof} Exercise.\end{proof}

\begin{exer} The ring of entire holomorphic functions $f\colon \C\to
\C$ is an integral domain but it is not factorial (Hint: consider
the
sine function $\sin(z)$).\\
For every connected open subset $U\subset\C^{n}$, the ring
$\Oh(U)$ is integrally closed in its field of fractions (Hint:
Riemann extension theorem).
\end{exer}

\bigskip
\section[~~Analytic singularities and fat points]{Analytic
singularities and fat points}
\label{sec:fatpoints}

Let $M$ be a complex manifold, as in \chaptername~\ref{CAP:FAMILIES} we
denote by $\Oh_{M,x}$  the ring of  germs of holomorphic functions
at a point $x\in M$. The elements of $\Oh_{M,x}$ are the
equivalence classes of pairs $(U,g)$, where $U$ is open, $x\in
U\subset M$, $g\colon U\to \C$ is holomorphic and $(U,g)\sim
(V,h)$ if there exists an open subset
$W$, $x\in W\subset U\cap V$ such that $g_{|W}=h_{|W}$.\\
By definition of holomorphic function and the identity principle we have
that
$\Oh_{\C^{n},0}$ is isomorphic to the ring of convergent
power series $\C\{z_{1},\ldots,z_{n}\}$.\\

Let $f\colon M\to N$ be a holomorphic map of complex manifolds, for
every open subset $V\subset N$ we have a homomorphism of $\C$-algebras
\[f^{*}\colon \Gamma(V,\Oh_{N})\to \Gamma(f^{-1}(V),\Oh_{M}),\qquad
f^{*}g=g\circ f\]
If $x\in M$ then the limit above maps $f^{*}$, for $V$ varying over
all the open neighbourhood of
$y=f(x)$,
gives a local homomorphism of local $\C$-algebras
$f^{*}\colon \Oh_{N,y}\to \Oh_{M,x}$.\\
It is clear that $f^{*}\colon \Oh_{N,y}\to \Oh_{M,x}$ depends only
on the behavior of $f$ in a neighbourhood of $x$ and then depends
only
on the class of $f$ in the space $\Mor_{\mathbf{Ger}}((M,x),(N,y))$.\\
A choice of local holomorphic
coordinates $z_{1},\ldots,z_{n}$ on $M$ such that
$z_{i}(x)=0$, gives an invertible morphism in
$\Mor_{\mathbf{Ger}}((M,x),(\C^{n},0))$
and then an isomorphism
$\Oh_{M,x}=\C\{z_{1},\ldots,z_{n}\}$.

\begin{exer} Given $f,g \in \Mor_{\mathbf{Ger}}((M,x),(N,y))$, prove that
$f=g$ if and only if $f^{*}=g^{*}$.\end{exer}

\begin{defn}\label{V.2.2}
An \emph{analytic singularity} is a triple $(M,x,I)$ where $M$
is a complex manifold, $x\in M$ is a point and $I\subset \Oh_{M,x}$ is
a proper ideal.\\
The germ morphisms $\Mor_{\mathbf{Ger}}((M,x,I),(N,y,J))$ are the
equivalence classes of morphisms $f\in \Mor_{\mathbf{Ger}}((M,x),(N,y))$
such that $f^{*}(J)\subset I$ and $f\sim g$ if and only if
$f^{*}=g^{*}\colon \Oh_{N,y}/J\to \Oh_{M,x}/I$.\\
We denote by $\mathbf{Ger}$ the category of analytic singularities (also
called germs of complex spaces).\end{defn}

\begin{lem}\label{V.2.4}
The contravariant functor $\mathbf{Ger}\to {\mathbf{An}}$,
\renewcommand\arraystretch{2}
\[ \begin{array}{rl}
Ob(\mathbf{Ger})\to Ob({\mathbf{An}}),&\qquad (M,x,I)\mapsto \Oh_{M,x}/I;\\
\Mor_{\mathbf{Ger}}((M,x,I),(N,y,J))\to
\Mor_{{\mathbf{An}}}\left(\dfrac{\Oh_{N,y}}{J},\dfrac{\Oh_{M,x}}{I}\right),&\qquad
f\mapsto f^{*};
\end{array}\]
\renewcommand\arraystretch{1}
is an equivalence of categories.
Its ``inverse'' ${\mathbf{An}}\to\mathbf{Ger}$
(cf. \cite[1.4]{MacLane}) is called $\Spec$ (sometimes
$\operatorname{Specan}$).
\end{lem}

\begin{proof} Since $\C\{z_{1},\ldots,z_{n}\}/I$ is isomorphic to
$\Oh_{\C^{n},0}/I$ the above functor is surjective on isomorphism
classes.\\
We only need to prove that
$\Mor_{\mathbf{Ger}}((M,x,I),(N,y,J))\to
\Mor_{{\mathbf{An}}}(\Oh_{N,y}/J,\Oh_{M,x}/I)$ is surjective, being injective
by definition of $\Mor_{\mathbf{Ger}}$. To see this it is not restrictive
to assume $(M,x)=(\C^{m}_{u},0)$, $(N,y)=(\C^{n}_{z},0)$.\\
Let $g^{*}\colon \C\{z_{1},\ldots,z_{n}\}/J\to
\C\{u_{1},\ldots,u_{m}\}/I$ be a local homomorphism and choose,
for every $i=1,\ldots,n$, a convergent power series
$f_{i}\in \C\{u_{1},\ldots,u_{m}\}$ such that $f_{i}\equiv
g^{*}(z_{i})$ $\pmod{I}$. Note that $f_{i}(0)=0$.\\
If $U$ is an open set, $0\in U\subset \C^{m}$, such that $f_{i}$ are
convergent in $U$, then we may define  a holomorphic map
$f=(f_{1},\ldots,f_{n})\colon U\to \C^{n}$.
By construction $f^{*}(z_{i})=g^{*}(z_{i})\in
\quot{\C\{u_{1},\ldots,u_{m}\}}{I}$ and then by Lemma~\ref{V.1.4} $f^{*}=g^{*}$.
\end{proof}

\begin{defn}\label{V.2.6}
Given an analytic singularity $(X,x)=(M,x,I)$, the analytic algebra
$\Oh_{X,x}:=\Oh_{M,x}/I$ is called the \emph{algebra of germs of
analytic functions} of $(X,x)$.\\
The \emph{dimension} of  $(X,x)$ is by definition the
dimension of the analytic algebra $\Oh_{X,x}$.
\end{defn}

\begin{defn}
A \emph{fat point} is an analytic singularity of dimension 0.
\end{defn}

\begin{lem}\label{V.2.8}
Let $X=(M,x,I)$ be an analytic singularity; the following
conditions are equivalent.\begin{enumerate}
\item The maximal ideal of $\Oh_{X,x}$ is nilpotent.
\item $X$ is a fat point.
\item The ideal $I$ contains a power of the maximal ideal of
$\Oh_{M,x}$.
\item If $V$ is open, $x\in
V\subset M$, and
$f_{1},\ldots,f_{h}\colon V\to \C$  are holomorphic functions
generating the ideal $I$, then
there exists
an open neighbourhood $U\subset V$ of $x$ such that
\[U\cap \{f_{1}=\ldots=f_{h}=0\}=\{x\}.\]
\item $\Mor_{{\mathbf{An}}}(\Oh_{X,x},\C\{t\})$ contains only the trivial
morphism $f\mapsto f(0)\in \C\subset \C\{t\}$.
\end{enumerate}\end{lem}

\begin{proof} $[1\Leftrightarrow 2\Leftrightarrow 3]$
are trivial.\\
\implica{3}{4} It is not restrictive to assume that $V$ is contained
in a coordinate chart; let $z_{1},\ldots,z_{n}\colon V\to \C$ be
holomorphic coordinates with $z_{i}(x)=0$ for every $i$.
If 3 holds then there exists $s>0$ such that $z_{i}^{s}\in I$ and
then there exists an open subset $x\in U\subset V$ and holomorphic
functions $a_{ij}\colon U\to \C$ such that
$z_{i}^{s}=\sum_{j}a_{ij}f_{j}$. Therefore
$U\cap V\cap \{f_{1}=\ldots=f_{h}=0\}\subset
U\cap\{z_{1}^{s}=\ldots=z_{n}^{s}=0\}=\{x\}$.\\
\implica{4}{5}
Let $\phi\colon
(\C,0)\to (M,x)$ be a germ of holomorphic map such that
$\phi^{*}(I)=0$. If $\phi$ is defined in an open subset $W\subset \C$
and $\phi(W)\subset U$ then $\phi^{*}(I)=0$ implies
$\phi(W)\subset U\cap \{f_{1}=\ldots=f_{h}=0\}$ and therefore
$\Mor_{\mathbf{Ger}}((\C,0,0),(M,x,I))$ contains only the constant
morphism.\\
\implica{5}{1} is a consequence of Theorem~\ref{V.1.6} (with
$J=\ide{m}_{M,x}$).
\end{proof}

\begin{exer} If $f\in\Mor_{\mathbf{Ger}}((M,x,I),(N,y,J))$ we define the
schematic fibre $f^{-1}(y)$ as the singularity
$(M,x,I+f^{*}\ide{m}_{N,y})$.\\
Prove that the dimension of a singularity $(M,x,I)$ is the minimum
integer $d$ such that there exists a morphism
$f\in\Mor_{\mathbf{Ger}}((M,x,I),(\C^{d},0,0))$ such that $f^{-1}(0)$ is a
fat point.\end{exer}

\begin{defn}\label{V.2.10}
The \emph{Zariski tangent space} $T_{x,X}$ of an analytic
singularity $(X,x)$
is the $\C$-vector space $\Der_{\C}(\Oh_{X,x},\C)$.\end{defn}

Note that every morphism of singularities $(X,x)\to (Y,y)$ induces
a linear morphism of Zariski tangent spaces $T_{x,X}\to T_{y,Y}$.

\begin{exer} (Cartan's Lemma)\\
Let $(R,\ide{m})$ be an analytic algebra and $G\subset \Aut(R)$ a finite
group of automorphisms. Denote $n=\dim_{\C}\quot{\ide{m}}{\ide{m}^{2}}$.\\
Prove that there exists an injective homomorphism of groups
$G\to GL(\C^{n})$ and a $G$-isomorphism of analytic algebras
$R\simeq\Oh_{\C^{n},0}/I$ for some $G$-stable ideal $I\subset
\Oh_{\C^{n},0}$. (Hint: there exists a direct sum decomposition
$\ide{m}=V\oplus \ide{m}^{2}$ such that $gV\subset V$ for every $g\in G$.)
\end{exer}

\bigskip
\section[~~The resultant]{The resultant}
\label{sec:resultant}

Let $A$ be a commutative unitary ring and
$p\in A[t]$  a monic polynomial of degree $d$. It is easy to see that
$A[t]/(p)$ is a
free $A$-module of rank $d$ with basis
$1,t,\ldots,t^{d-1}$.\\
For every $f\in A[t]$ we  denote by
$R(p,f)\in A$ the determinant of the multiplication map
$f\colon A[t]/(p)\to A[t]/(p)$.\\

\begin{defn}\label{V.3.2}
In the notation above, the element
$R(p,f)$ is called the \emph{resultant} of $p$ and $f$.
\end{defn}

If $\phi\colon A\to B$ is a morphism of unitary rings then we
can extend it to a morphism $\phi\colon A[t]\to B[t]$, $\phi(t)=t$, and
it is clear from the definition that
$R(\phi(p),\phi(f))=\phi(R(p,f))$.\\
By Binet's theorem $R(p,fg)=R(p,f)R(p,g)$.

\begin{lem}\label{V.3.4}
In the notation above there exist $\alpha,\beta\in A[t]$
with $\deg\alpha<\deg f$, $\deg\beta<\deg p$ such that
$R(p,f)=\beta f-\alpha p$. In particular $R(p,f)$ belongs to the ideal
generated by $p$ and $f$.\end{lem}

\begin{proof} For every $i,j=0,\ldots,d-1$ there exist $h_{i}\in A[t]$
and $c_{ij}\in A$ such that
\[ t^{i}f=h_{i}p+\sum_{j=0}^{d-1}c_{ij}t^{j},\qquad \deg h_{i}<\deg f.\]
By definition $R(p,f)=\det(c_{ij})$; if $(C^{ij})$ is the adjoint
matrix of $(c_{ij})$ we have, by Laplace formula, for every
$j=0,\ldots,d-1$
\[ \sum_{i}C^{0i}c_{ij}=\delta_{0j}R(p,f)\]
and then
\[ R(p,f)=\sum_{i=0}^{d-1}C^{0i}(t^{i}f-h_{i}p)=\beta f-\alpha p.\]
\end{proof}

\begin{lem}\label{V.3.6}
In the notation above, if $A$ is an integral domain and $p,f$ have a
common factor of positive degree then $R(p,f)=0$. The converse hold
if $A$ is a unique factorization domain.\end{lem}

\begin{proof}
Since $A$ injects into its fraction field, the multiplication
$f\colon A[t]/(p)\to A[t]/(p)$ is  injective if and only if
$R(p,f)\not=0$.\\
If $p=qr$ with $\deg r<\deg p$, then the multiplication
$q\colon A[t]/(p)\to A[t]/(p)$ is not injective and then its
determinant is trivial. If $q$ also divides $f$ then, by the theorem of
Binet also $R(p,f)=0$.\\
Assume now that $A$ is a unique factorization domain and
$R(p,f)=0$. There exists $q\not\in (p)$ such that $fq\in (p)$;
by Gauss' lemma $A[t]$ is a UFD and then there exists a irreducible
factor $p_{1}$ of $p$ dividing $f$. Since $p$ is a monic polynomial
the degree of $p_{1}$ is positive.
\end{proof}

\begin{lem}\label{V.3.8}
Let $A$ be an integral domain and $0\not=\ide{p}\subset
A[t]$ a prime ideal such that $\ide{p}\cap A=0$.
Denote by $ K$ the fraction field of $A$
and by $\ide{p}^e\subset K[x]$ the ideal generated by $\ide{p}$.\\
Then:\begin{enumerate}
\item $\ide{p}^e$ is a prime  ideal.
\item $\ide{p}^e\cap A[x]=\ide{p}$.
\item There exists $f\in \ide{p}$ such that for every monic polynomial
$p\not\in\ide{p}$ we have $R(p,f)\not=0$.
\end{enumerate}
\end{lem}

\begin{proof} {[1]} We have
$\ide{p}^e=\left\{\left.\ds\frac{p}{a}\,
\right|\, p\in \ide{p},\, a\in A-\{0\}\right\}$.
If $\ds\frac{p_1}{a_1}\ds\frac{p_2}{a_2}\in \ide{p}^e$ with
$p_i\in A[x]$, $a_i\in A$; then there exists $a\in A-\{0\}$ such that
$ap_1p_2\in \ide{p}$. Since $\ide{p}\cap A=0$ it must be $p_1\in \ide{p}$ or
$p_2\in \ide{p}$. This shows that $\ide{p}^e$ is prime.\\
{[2]} If $q\in \ide{p}^e\cap A[x]$, then there exists $a\in A$, $a\not=0$
such that $aq\in \ide{p}$ and therefore $q\in \ide{p}$.\\
{[3]} Let $f\in\ide{p}-\{0\}$ be of minimal degree, since $K[t]$
is an Euclidean ring, $\ide{p}^{e}=fK[t]$ and, since $\ide{p}^{e}$
is prime, $f$ is irreducible in $K[t]$. If $p\in A[t]\setminus
\ide{p}$ is a monic polynomial then $p\not\in \ide{p}^{e}=fK[t]$
and then, according to Lemma~\ref{V.3.6}, $R(p,f)\not=0$.
\end{proof}

\begin{thm}\label{V.3.10}
Let $A$ be a unitary ring, $\ide{p}\subset A[t]$ a prime ideal,
$\ide{q}=A\cap \ide{p}$.\\
If $\ide{p}\not=\ide{q}[t]$ (e.g. if $\ide{p}$ is proper and contains
a monic polynomial) then
there exists $f\in \ide{p}$ such that for every monic polynomial
$p\not\in\ide{p}$ we have $R(p,f)\not\in\ide{q}$.\\
If moreover $A$ is a unique factorization domain we can choose $f$
irreducible.
\end{thm}

\begin{proof} $\ide{q}$ is prime and $\ide{q}[t]\subset \ide{p}$,
therefore the image of $\ide{p}$ in
$(A/\ide{q})[t]=A[t]/\ide{q}[t]$
is still a prime ideal satisfying the hypothesis of  Lemma~\ref{V.3.8}.\\
It is therefore sufficient to take $f$ as any lifting of the element
described in Lemma~\ref{V.3.8} and use the functorial properties of the
resultant.
If $A$ is UFD and $f$ is not irreducible we can write
$f=hg$ with $g\in\ide{p}$ irreducible;
but $R(p,f)=R(p,h)R(p,g)$ and then also $R(p,g)\not\in \ide{q}$.
\end{proof}

\begin{exer} If $p,q\in A[t]$ are monic polynomials of degrees $d,l>0$
then for every $f\in A[t]$ we have $R(pq,f)=R(p,f)R(q,f)$. (Hint:
write the matrix of the multiplication $f\colon A[t]/(pq)\to
A[t]/(pq)$ in the basis $1,t,\ldots,t^{d-1},p, tp,\ldots,t^{l-1}p$.)
\end{exer}

\bigskip
\section[~~R{\"u}ckert's Nullstellensatz]{R\"{u}ckert's Nullstellensatz}
\label{sec:ruckert}

The aim of this section is to prove the following theorem, also called
\emph{Curve selection lemma}, which is
easily seen to be equivalent to Theorem~\ref{V.1.6}. The proof given 
here is a particular case of the one sketched in \cite{IntroGA}.

\begin{thm}\label{V.4.2}
Let $\ide{p}\subset \C\{z_{1},\ldots,z_{n}\}$ be a proper prime
ideal and $h\not\in\ide{p}$. Then there exists a homomorphism of
local $\C$-algebras $\phi\colon \C\{z_{1},\ldots,z_{n}\}\to \C\{t\}$
such that $\phi(\ide{p})=0$ and $\phi(h)\not=0$.\end{thm}

\begin{cor}\label{V.4.3}
Let $I\subset \C\{z_{1},\ldots,z_{n}\}$ be a proper
ideal and $h\not\in\sqrt{I}$. Then there exists a homomorphism of
local $\C$-algebras $\phi\colon \C\{z_{1},\ldots,z_{n}\}\to \C\{t\}$
such that $\phi(I)=0$ and $\phi(h)\not=0$.\end{cor}

\begin{proof} If $h\not\in\sqrt{I}$ there exists (cf. \cite{A-M}) a
prime ideal $\ide{p}$ such that $I\subset \ide{p}$ and
$h\not\in\ide{p}$.\end{proof}

Before proving Theorem~\ref{V.4.2} we need a series of results that are
of independent interest. We recall the following

\begin{defn}\label{V.4.4}
A power series $p\in\C\{z_{1},\ldots,z_{n},t\}$ is called a 
\emph{Weierstrass polynomial} in $t$ of degree $d\ge 0$ if
\[ p=t^{d}+\sum_{i=0}^{d-1}p_{i}(z_{1},\ldots,z_{n})t^{i},\qquad
p_{i}(0)=0.\]
\end{defn}

In particular if $p(z_{1},\ldots,z_{n},t)$
is a Weierstrass polynomial in $t$ of degree $d$ then
$p(0,\ldots,0,t)=t^{d}$.

\begin{thm}[Preparation theorem]\label{V.4.6}
Let $f\in \C\{z_{1},\ldots,z_{n},t\}$ be a power series such that
$f(0,\ldots,0,t)\not=0$. Then there exists a unique  $e\in
\C\{z_{1},\ldots,z_{n},t\}$ such that $e(0)\not=0$ and $ef$ is a
Weierstrass polynomial in $t$.\end{thm}

\begin{proof} For the proof we refer to \cite{Ana-Ste},
\cite{CAS}, \cite{G-H}, \cite{Kobook},
\cite{Gu-Ro}, \cite{IntroGA}.
We note that the condition that the power
series $\mu(t)=f(0,\ldots,0,t)$ is not trivial is also necessary and
that the degree of $ef$ in $t$ is equal to the multiplicity at $0$
of $\mu$.\end{proof}

\begin{cor}\label{V.4.8}
Let $f\in \C\{z_{1},\ldots,z_{n}\}$ be a power series of multiplicity
$d$. Then, after a possible generic linear change of coordinates there exists
$e\in \C\{z_{1},\ldots,z_{n}\}$ such that $e(0)\not=0$ and $ef$ is a
Weierstrass polynomial of degree $d$ in $z_{n}$.\end{cor}

\begin{proof} After a generic change of coordinates of the form
$z_{i}\mapsto z_{i}+a_{i}z_{n}$, $a_{i}\in \C$, the series
$f(0,\ldots,0,z_{n})$ has multiplicity $d$.\end{proof}

\begin{lem}\label{V.4.10}
Let $f,g\in \C\{x_1,\ldots ,x_n\}[t]$ be polynomials in $t$ with
$g$ in Weierstrass' form. if $f=hg$ for some
$h\in
\C\{x_1,\ldots ,x_n,t\}$ then $h\in \C\{x_1,\ldots
,x_n\}[t]$.\end{lem}

We note that if  $g$ is not a Weierstrass polynomial
then the above result is false; consider for instance
the case $n=0, f=t^3, g=t+t^2$.

\begin{proof}
Write $g=t^s+\sum g_i(x)t^{s-i}$, $g_i(0)=0$, $f=\sum_{i=0}^r
f_i(x)t^{r-i}$
$h=\sum_ih_i(x)t^i$, we need to prove that
$h_i=0$ for every $i>r-s$.\\
Assume the contrary
and choose an index $j>r-s$
such that the multiplicity
of $h_j$ takes the minimum among all the multiplicities
of the power series  $h_i$, $i>r-s$.\\
From the equality
$0=h_j+\sum_{i>0} g_ih_{j+i}$ we get a contradiction.\end{proof}

\begin{lem}\label{V.4.12}
Let $f\in \C\{x_1,\ldots ,x_n\}[t]$ be an
irreducible monic polynomial of degree $d$.
Then the polynomial
$f_{0}(t)=f(0,\ldots,0,t)\in\C[t]$
has a root of multiplicity $d$.\end{lem}

\begin{proof} Let $c\in \C$ be a root of $f_{0}(t)$. If the multiplicity
of $c$ is $l< d$ then the multiplicity of the power series
$f_{0}(t+c)\in\C\{t\}$
is exactly $l$ and therefore $f(x_{1},\ldots,x_{n},t+c)$ is divided
in $\C\{x_1,\ldots ,x_n\}[t]$  by a Weierstrass polynomial of degree
$l$.\end{proof}

\begin{lem}\label{V.4.13}
Let $p\in\C\{x\}[y]$ be a monic polynomial of positive degree $d$ in $y$.
Then there exists a  homomorphism $\phi\colon \C\{x\}[y]\to\C\{t\}$
such that $\phi(p)=0$ and $0\not=\phi(x)\in (t)$.\\
\end{lem}

\begin{proof}
If $d=1$ then $p(x,y)=y-p_{1}(x)$ and we can consider the morphism
$\phi$ given by $\phi(x)=t$, $\phi(y)=p_{1}(t)$. By induction
we can assume the theorem true  for monic polynomials of
degree $<d$.\\
If $p$ is reducible we have done, otherwise, writing
$p=y^{d}+p_{1}(x)y^{d-1}+\ldots+p_{d}(x)$,
after the coordinate change
$x\mapsto x$, $y\mapsto y-p_{1}(x)/d$ we can assume $p_{1}=0$.\\
For every $i\ge 2$ denote by
$\mu(p_{i})=\alpha_i>0$ the multiplicity of $p_i$ (we set
$\alpha_i=+\infty$ if $p_i=0$).\\
Let $j\ge 2$ be a fixed index such that
$\ds\frac{\alpha_j}{j}\le
\frac{\alpha_i}{i}$ for every $i$.
Setting $m=\alpha_{j}$,
we want to prove that the monic polynomial
$p(\xi^j,y)$ is not irreducible.\\
In fact $p(\xi^j,y)=y^{d}+\sum_{i\ge 2} h_i(\xi)y^{d-i}$,
where  $h_i(\xi)=g_i(\xi^{j})$.\\
For every $i$ the multiplicity of $h_i$ is
$j\alpha_i\ge im$ and then
\[q(\xi,y)=p(\xi^j,\xi^m y)\xi^{-dm}=
t^d+\sum\frac{h_i(\xi)}{\xi^{mi}}y^{d-i}=
y^{d}+\sum\eta_i(\xi)y^{d-i}\]
is a well defined element of $\C\{\xi,y\}$.
Since $\eta_{1}=0$ and $\eta_j(0)\not=0$
the polynomial $q$ is not irreducible and then, by induction
there exists a nontrivial morphism
$\psi\colon \C\{\xi\}[y]\to \C\{t\}$ such that
$\psi(q)=0$, $0\not=\psi(\xi)\in (t)$
and we can take $\phi(x)=\psi(\xi^{j})$ and
$\phi(y)=\psi(\xi^{m}y)$.\end{proof}

\begin{thm}[Division theorem]\label{V.4.16}
Let $p\in \C\{z_{1},\ldots,z_{n},t\}$, $p\not=0$,   be a
Weierstrass polynomial of degree $d\ge 0$ in $t$. Then
for every $f\in \C\{z_{1},\ldots,z_{n},t\}$
there exist a unique  $h\in
\C\{z_{1},\ldots,z_{n},t\}$ such that
$f-hp\in \C\{z_{1},\ldots,z_{n}\}[t]$ is a polynomial
of degree $<d$ in $t$.\end{thm}

\begin{proof} For the proof we refer to \cite{Ana-Ste},
\cite{CAS}, \cite{G-H}, \cite{Kobook},
\cite{Gu-Ro}, \cite{IntroGA}.\end{proof}

We note that an equivalent statement for the division theorem is the
following:

\begin{cor}\label{V.4.18}
If $p\in \C\{z_{1},\ldots,z_{n},t\}$, $p\not=0$, is  a
Weierstrass polynomial of degree $d\ge 0$ in $t$, then
$\C\{z_{1},\ldots,z_{n},t\}/(p)$ is a free
$\C\{z_{1},\ldots,z_{n}\}$-module with basis
$1,t,\ldots,t^{d-1}$.\end{cor}

\begin{proof} Clear.\end{proof}

\begin{thm}[Newton-Puiseux]\label{V.4.14}
Let $f\in\C\{x,y\}$ be a power series of positive multiplicity.
Then there exists a nontrivial
local homomorphism $\phi\colon \C\{x,y\}\to\C\{t\}$
such that $\phi(f)=0$.\\
Moreover if $f$ is irreducible then $\ker\phi=(f)$.
\end{thm}

In the above statement nontrivial means that $\phi(x)\not=0$ or
$\phi(y)\not=0$.

\begin{proof}
After a linear change of coordinates we can assume $f(0,y)$ a non zero
power series of multiplicity $d>0$;
by Preparation theorem there
exists an invertible power series $e$ such that
$p=ef$ is a Weierstrass polynomial of degree $d$ in $y$.\\
According to Lemma~\ref{V.4.13}
there exists a homomorphism $\phi\colon \C\{x\}[y]\to \C\{t\}$
such that $\phi(p)=0$ and $0\not=\phi(x)\in (t)$. Therefore
$\phi(p(0,y))\in (t)$ and, being $p$ a Weierstrass polynomial we have
$\phi(y)\in (t)$ and then $\phi$ extends to a local morphism
$\phi\colon \C\{x,y\}\to \C\{t\}$.\\
Assume now $f$ irreducible, up to a possible change of coordinates
and multiplication for an invertible element we may assume that
$f\in \C\{x\}[y]$ is an irreducible Weierstrass polynomial of degree
$d>0$.\\
Let $\phi\colon \C\{x,y\}\to \C\{t\}$ be a nontrivial morphism such
that $\phi(f)=0$, then $\phi(x)\not=0$ (otherwise
$\phi(y)^{d}=\phi(f)=0$) and therefore the restricted morphism
$\phi\colon \C\{x\}\to\C\{t\}$ is injective.\\
Let $g\in\ker(\phi)$, by division theorem there exists $r\in
\C\{x\}[y]$ such that $g=hf+r$ and then $r\in\ker(\phi)$, $R(f,r)\in
\ker(\phi)\cap \C\{x\}=0$. This implies that $f$ divides $r$.
\end{proof}

The division theorem allows to extend the definition of the resultant
to power series. In fact if
$p\in \C\{z_{1},\ldots,z_{n}\}[t]$ is a  Weierstrass
polynomial in $t$  of degree $d$,
for every $f\in \C\{z_{1},\ldots,z_{n},t\}$ we can define the
resultant
$R(p,f)\in \C\{z_{1},\ldots,z_{n}\}$ as the determinant of the
morphism of free  $\C\{z_{1},\ldots,z_{n}\}$-module
\[f\colon \frac{\C\{z_{1},\ldots,z_{n},t\}}{(p)}\to
\frac{\C\{z_{1},\ldots,z_{n},t\}}{(p)}\]
induced by the multiplication with $f$.\\
It is clear that $R(p,f)=R(p,r)$ whenever $f-r\in (p)$.

\begin{lem}\label{V.4.20}
Let $p\in \C\{z_{1},\ldots,z_{n},t\}$ be a  Weierstrass
polynomial of positive degree in $t$
and $V\subset \C\{z_{1},\ldots,z_{n},t\}$  a
$\C$-vector subspace.\\
Then $R(p,f)=0$ for every $f\in V$ if and only if
there exists a Weierstrass
polynomial $q$ of positive degree such that:
\begin{enumerate}
    \item  $q$ divides $p$ in $\C\{z_{1},\ldots,z_{n}\}[t]$

    \item  $V\subset q\C\{z_{1},\ldots,z_{n},t\}$
\end{enumerate}
\end{lem}

\begin{proof} One implication is clear, in fact if $p=qr$ then the
multiplication by $q$ in not injective in
$\C\{z_{1},\ldots,z_{n},t\}/{(p)}$; therefore $R(p,q)=0$ and by
Binet's theorem $R(p,f)=0$ for every $f\in (q)$.\\
For the converse let $p=p_{1}p_{2}\ldots p_{s}$ be the irreducible
decomposition of $p$ in the UFD $\C\{z_{1},\ldots,z_{n}\}[t]$.
If $R(p,f)=0$ and $r=f-hp\in \C\{z_{1},\ldots,z_{n}\}[t]$ is the rest
of the division then $R(p,r)=0$ and by Lemma~\ref{V.3.6} there exists
a factor $p_{i}$  dividing $r$ and therefore also dividing  $f$.\\
In particular, setting $V_{i}=V\cap (p_{i})$, we have $V=\cup_{i}
V_{i}$ and therefore $V=V_{i}$ for at least one index $i$ and we can
take $q=p_{i}$.
\end{proof}

\begin{proof}[Proof of \ref{V.4.2}]
We first consider the easy cases $n=1$ and $\ide{p}=0$.
If $\ide{p}=0$ then, after a possible change of coordinates, we may
assume $h(0,\ldots,0,t)\not=0$ and therefore we can take
$\phi(z_{i})=0$ for $i=1,\ldots,n-1$ and $\phi(z_{n})=t$.\\
If $n=1$ the only prime nontrivial ideal is  $(z_{1})$ and therefore
the trivial morphism $\phi\colon \C\{z_{1}\}\to \C\subset \C\{t\}$
satisfies the statement of the theorem.\\
Assume then $n>1$, $\ide{p}\not=0$ and  fix a nonzero element $g\in
\ide{p}$.  After a possible linear change of coordinates and
multiplication  by invertible elements we may assume both $h$ and $g$
Weierstrass polynomials in the variable $z_{n}$.
Denoting
\[ \ide{r}=\ide{p}\cap \C\{z_{1},\ldots,z_{n-1}\}[z_{n}],\qquad
\ide{q}=\ide{p}\cap \C\{z_{1},\ldots,z_{n-1}\}=\ide{r}\cap
\C\{z_{1},\ldots,z_{n-1}\},\]
according to Theorem~\ref{V.3.10}, there exists $\hat{f}\in \ide{r}$ such
that $R(h,\hat{f})\not\in\ide{q}$. On the other hand, by Lemma~\ref{V.3.4},
$R(g,f)\in \ide{q}$ for every $f\in \ide{p}$.\\
By induction on $n$ there exists a morphism
$\tilde{\psi}\colon \C\{z_{1},\ldots,z_{n-1}\}\to \C\{x\}$ such that
$\tilde{\psi}(\ide{q})=0$ and $\tilde{\psi}(R(h,\hat{f}))\not=0$.
Denoting by $\psi\colon \C\{z_{1},\ldots,z_{n}\}\to \C\{x,z_{n}\}$ the
natural extension of $\tilde{\psi}$ we have
$R(\psi(h), \psi(\hat{f}))\not=0$ and
$R(\psi(g), \psi(f))=0$ for every $f\in \ide{p}$.
Applying Lemma~\ref{V.4.20} to the Weierstrass polynomial $\psi(g)$
and the vector space $V=\psi(\ide{p})$ we prove the existence of an
irreducible factor $p$ of $\psi(g)$ such that $\psi(\ide{p})\subset
p\C\{x,z_{n}\}$.\\
In particular $p$ divides $\psi(\hat{f})$,
therefore $R(\psi(h),p)\not=0$ and $\psi(h)\not\in p\C\{x,z_{n}\}$.\\
By Newton-Puiseux' theorem there exists $\eta\colon \C\{x,z_{n}\}\to
\C\{t\}$ such that $\eta(p)=0$ and $\eta(\psi(h))\not=0$.
It is therefore sufficient to take $\phi$ as the composition of $\psi$ and
$\eta$.\end{proof}
\bigskip

\begin{exer} Prove that $f,g\in \C\{x,y\}$
have a common factor  of positive multiplicity  if and only if
the $\C$-vector space $\C\{x,y\}/(f,g)$ is infinite dimensional.
\end{exer}

\section[~~Dimension bounds]{Dimension bounds}

As an application of Theorem~\ref{V.4.2} we give some  bounds
for the dimension of an analytic algebra; this bounds will be very
useful in deformation and moduli theory. The first bound (Lemma~\ref{V.4.26})
is completely standard and the proof is reproduced here for
completeness; the second bound (Theorem~\ref{V.4.28}, communicated to the
author by H. Flenner) finds application in the ``$T^{1}$-lifting''
approach to deformation problems.\\

We need the following two results of commutative algebra.

\begin{lem}\label{V.4.22}
Let $(A,\ide{m})$ be a local Noetherian ring and
$J\subset I\subset A$ two ideals. If $J+\ide{m}I=I$ then $J=I$.
\end{lem}

\begin{proof} This a special case of Nakayama's lemma \cite{A-M},
\cite{IntroGA}.\end{proof}

\begin{lem}\label{V.4.24}
Let $(A,\ide{m})$ be a local Noetherian ring and
$f\in \ide{m}$, then $\dim A/(f)\ge \dim A-1$.\\
Moreover, if $f$ is nilpotent then $\dim A/(f)=\dim A$, while if $f$
is not a zerodivisor then $\dim A/(f)=\dim A-1$.
\end{lem}

\begin{proof} \cite{A-M}.\end{proof}

\begin{lem}\label{V.4.26}
Let $R$ be an analytic algebra with maximal ideal
$\ide{m}$, then $\dim R\le \dim_{\C}\dfrac{\ide{m}~}{\ide{m}^{2}}$ and
equality holds if and only if $R$ is smooth.\end{lem}

\begin{proof} Let $n=\dim_{\C}\dfrac{\ide{m}~}{\ide{m}^{2}}$ and
$f_{1},\ldots,f_{n}\in \ide{m}$ inducing a basis of
$\dfrac{\ide{m}~}{\ide{m}^{2}}$. If $J=(f_{1},\ldots,f_{n})$ by
assumption $J+\ide{m}^{2}=\ide{m}$ and then by Lemma~\ref{V.4.22}
$J=\ide{m}$, $R/J=\C$ and $0=\dim R/J\ge \dim R-n$.\\
According to Lemma~\ref{V.1.8} we can write
$R=\C\{z_{1},\ldots,z_{n}\}/I$ for some ideal contained in
$(z_{1},\ldots,z_{n})^{2}$. Since $\C\{z_{1},\ldots,z_{n}\}$ is an
integral domain, according to Lemma~\ref{V.4.24} $\dim R=n$ if and
only if $I=0$.\end{proof}

\begin{thm}\label{V.4.28}
Let $R=P/I$ be an analytic algebra, where
$P=\C\{z_{1},\ldots,z_{n}\}$,
$n>0$ is a fixed integer, and $I\subset P$ is a proper
ideal.\\
Denoting by
$\ide{m}=(z_{1},\ldots,z_{n})$ the maximal ideal of $P$
and by $J\subset P$ the ideal
\[ J=\left\{f\in I\,\left|\, \desude{f}{z_{i}}\in I,\, \forall
i=1,\ldots,n\right.\right\}\]
we have
$\dim R\ge n-\dim_{\C}\dfrac{I}{J+\ide{m}I}$.
\end{thm}

\begin{proof} (taken from \cite{FaMa2})
We first introduce the \emph{curvilinear obstruction map}
\[ \gamma_{I}\colon \Mor_{{\mathbf{An}}}(P,\C\{t\})\to
\Hom_{\C}\left(\dfrac{I}{J+\ide{m}I},\C\right).\]
Given $\phi\colon P\to \C\{t\}$, if $\phi(I)=0$ we define
$\gamma_{I}(\phi)=0$; if $\phi(I)\not=0$ and $s$ is the biggest
integer such that $\phi(I)\subset (t^{s})$ we define, for every
$f\in I$,
$\gamma_{I}(\phi)f$ as the coefficient of $t^{s}$ in the power series
expansion of $\phi(f)$.\\
It is clear that $\gamma_{I}(\phi)(\ide{m}I)=0$, while if
$\phi(I)\subset (t^{s})$ and
$f\in J$ we
have $\phi(f)=f(\phi(z_{1}),\ldots,\phi(z_{n}))$,
\[\frac{d\phi(f)}{dt}=
\sum_{i=1}^{n}\desude{f}{z_{i}}(\phi(z_{1}),\ldots,\phi(z_{n}))
\frac{d\phi(z_{i})}{dt}\in (t^{s})\]
and therefore $\phi(f)\in (t^{s+1})$ (this is the point where the
characteristic of the field plays an essential role).\\
The ideal $I$ is finitely generated, say $I=(f_1,\dots,f_d)$, according to
Nakayama's lemma we can assume $f_{1},\ldots,f_{d}$ a basis of
$I/\ide{m}I$.\\
By repeated application of Corollary~\ref{V.4.3}
(and possibly reordering the $f_i$'s) we
can assume that there exists an $h\le d$ such that the following
holds:\begin{enumerate}
\item $f_i\notin\sqrt{(f_1,\dots,f_{i-1})}$ for $i\le h$;
\item for every $i\le h$ there exists a morphism of analytic algebras
$\phi_i:P\to \C\{t\}$
such that $\phi_i(f_i)\ne 0$, $\phi_i(f_j)=0$ if $j<i$
and the multiplicity of
$\phi_i(f_j))$ is bigger than or equal to the
multiplicity of $\phi_i(f_i))$ for every  $j>i$.
\item $I\subset \sqrt{(f_1,\dots,f_h)}$.
\end{enumerate}
Condition 3) implies that $\dim R=\dim P/(f_1,\dots,f_h)\ge n-h$,
hence it is
enough to prove
that $\gamma_{I}(\phi_{1}),\ldots, \gamma_{I}(\phi_{h})$ are linearly
independent in
$\Hom_{\C}\left(\dfrac{I}{J+\ide{m}I},\C\right)$
and this follows immediately from the fact that the matrix
$a_{ij}=\gamma_{I}(\phi_{i}){f_{j}}$, $i,j=1,\ldots,h$, has rank
$h$, being triangular with nonzero elements on the diagonal.
\end{proof}

\begin{exer}
In the notation of Theorem~\ref{V.4.28} prove that $I^2\subset J$. Prove moreover that
$I=J+\ide{m}I$ if and only if $I=0$.
\end{exer}

\begin{exer}
Let $I\subset \C\{x,y\}$ be the ideal generated by
the polynomial $f=x^{5}+y^{5}+x^{3}y^{3}$ and by its partial
derivatives $f_{x}=5x^{4}+3x^{2}y^{3}$, $f_{y}=5y^{4}+3x^{3}y^{2}$.
Prove that $J$ is not contained in $\ide{m}I$,  compute the
dimension of the analytic algebra $\C\{x,y\}/I$ and of the vector
spaces $\dfrac{I}{J+\ide{m}I}$, $\dfrac{I}{\ide{m}I}$.
\end{exer}

\begin{exer} (easy, but for  experts)
In the notation of \ref{V.4.28},
if $I\subset \ide{m}^{2}$ then
\[\Hom_{\C}\left(\dfrac{I}{J+\ide{m}I},\C\right)=
\Ext^{1}_{R}(\Omega_{R},\C).\]
($\Omega_{R}$ is the $R$-module of separated differentials)
\end{exer}

\begin{exer}
In the notation of Theorem \ref{V.4.28}, prove that for every short
exact sequence $0\to E\to F\to G\to 0$ of $R$-modules of finite
length (i.e. annihilated by some power of the maximal ideal) it is
defined a map
\[ ob\colon \Der_{\C}(R,G)\to \Hom_{R}\left(\frac{I}{J},E\right)\]
with the property that $ob(\phi)=0$ if and only if $\phi$ lifts to a
derivation $R\to F$.\\
Moreover, if $\ide{m}_{R}E=0$ then
$\Hom_{R}\left(\dfrac{I}{J},E\right)=
\Hom_{\C}\left(\dfrac{I}{J+\ide{m}I},E\right)$.
\end{exer}

\begin{rem} ($T^{1}$-lifting for prorepresentable functors.)\\
For every
morphism of analytic algebras $f\colon R\to A$ and every $A$-module of
finite length $M$ there exists a bijection between $\Der_{\C}(R,M)$
and the liftings of $f$ to morphisms $R\to A\oplus M$.\\
In the notation of Theorem \ref{V.4.28}, if $I\subset\ide{m}^{2}$,
then
$\Hom_{\C}\left(\dfrac{I}{J+\ide{m}I},\C\right)$ is
the subspace of  $\Hom_{\C}\left(\dfrac{I}{\ide{m}I},\C\right)$ of
obstructions (see \cite[Section 5]{FaMa1}) of the functor $h_{R}$
arising from all the small extensions of the form
$0\to\C\to A\oplus M\mapor{(Id,p)} A\oplus N\to 0$, where $p\colon
M\to N$ is a morphism of $A$-modules and
$A\oplus M\to A$, $A\oplus N\to A$ are the trivial extensions.
\end{rem}

\bigskip

\section[~~Historical survey,~\ref{CAP:SING}]{Historical
survey,~\ref{CAP:SING}}

According to \cite{CAS}, the preparation theorem was proved by
Weierstrass in 1860, while division theorem was proved by
Stickelberger in 1887.\\
The factoriality of $\C\{z_{1},\ldots,z_{n}\}$
was proved by E. Lasker in a, long time ignored, paper published in
1905. The same result was rediscovered by W. R\"uckert (a student of W.
Krull) together the Noetherianity in 1931. In the same paper of
R\"uckert it is implicitly contained the Nullstellensatz. The ideas of
R\"uckert's proof are essentially the same used in the proof given
in \cite{Gu-Ro}. The proof given here is different.\\

All the algebraic
results of this chapter that make sense also for the ring of
formal power series $\C[[z_{1},\ldots,z_{n}]]$ and their quotients,
remain true. In many cases, especially in deformation theory,
we seek for solutions of systems of analytic equations but we are able to
solve these equation only formally; in this situation a great help
comes from the following theorem, proved by M. Artin in 1968.

\begin{thm}\label{V.5.2}
Consider two arbitrary morphisms of analytic algebras $f\colon S\to R$,
$g\colon S\to \C\{z_{1},\ldots,z_{n}\}$
and a positive integer
$s>0$. The inclusion $\C\{z_{1},\ldots,z_{n}\}\subset
\C[[z_{1},\ldots,z_{n}]]$ and the projection
$\C\{z_{1},\ldots,z_{n}\}\to
\dfrac{\C\{z_{1},\ldots,z_{n}\}}{(z_{1},\ldots,z_{n})^{s}}$ give
structures of $S$-algebras also on
$\C[[z_{1},\ldots,z_{n}]]$ and
$\dfrac{\C\{z_{1},\ldots,z_{n}\}}{(z_{1},\ldots,z_{n})^{s}}$.\\
Assume it is given a morphism of analytic $S$-algebras
\[ \phi\colon R\to
\frac{\C\{z_{1},\ldots,z_{n}\}}{(z_{1},\ldots,z_{n})^{s}}=
\frac{\C[[z_{1},\ldots,z_{n}]]}{(z_{1},\ldots,z_{n})^{s}}.\]
If $\phi$
lifts to a $S$-algebra morphism
$R\to \C[[z_{1},\ldots,z_{n}]]$ then
$\phi$ lifts also to a $S$-algebra morphism $R\to
\C\{z_{1},\ldots,z_{n}\}$.\end{thm}

\textbf{Beware.} Theorem~\ref{V.5.2} does not imply that every
lifting $R\to \C[[z_{1},\ldots,z_{n}]]$ is ``convergent''.\\

\begin{proof} This is an equivalent statement of the main theorem of
\cite{Artin68}.
We leave as as an exercise to the reader to proof of the equivalence of the
two statements.\end{proof}

\begin{exer}\label{V.5.4}
Use Theorem~\ref{V.5.2} to prove:
\begin{enumerate}
    \item  Every irreducible
convergent power series $f\in\C\{z_{1},\ldots,z_{n}\}$
is also irreducible in $\C[[z_{1},\ldots,z_{n}]]$.
    \item  $\C\{z_{1},\ldots,z_{n}\}$
is integrally closed in  $\C[[z_{1},\ldots,z_{n}]]$.
\end{enumerate}
\end{exer}

\begin{rem} It is possible to give also an elementary proof of item 2
of Exercise~\ref{V.5.4} (e.g. \cite{IntroGA}),
while I don't know any proof of item 1 which
does not involve Artin's theorem.\end{rem}


\chapter[~Infinitesimal deformations]{Infinitesimal
deformations of complex manifolds}
\label{CAP:INFINITESIMAL}
\piede

In this chapter we pass from the classical language of deformation theory to the 
formalism  of differential graded objects. 
After a brief introduction 
of dg-vector spaces and dg-algebras, we associate to every deformation 
$X_{0}\hookrightarrow\{X_{t}\}_{t\in T}\to (T,0)$ its \emph{algebraic 
data} (Definition \ref{algdata}), 
which is  a pair of morphisms of  sheaves of dg-algebras on $X_{0}$. 
This algebraic data encodes  the Kodaira-Spencer map and also  
all the ``Taylor coefficients'' of  $t\mapsto X_{t}$.\\
We introduce the notion of infinitesimal deformation   as an 
infinitesimal variation of integrable complex structures; this 
definition will be more useful for our purposes. 
The infinitesimal Newlander-Nirenberg theorem, i.e.  the equivalence of this definition 
with the  more standard definition involving flatness, although not 
difficult to prove, would require a considerable amount of preliminaries in 
commutative and homological algebra and it is not given 
in this notes.\\
In Section~\ref{sec:HistInf} we state without proof the Kuranishi's theorem of 
existence of semiuniversal deformations of compact complex manifolds. 
In order to keep this notes short and selfcontained, we avoid the use of complex analytic
spaces and we state only the  "infinitesimal" version of Kuranishi's theorem. 
This is not a great gap for us since we are mainly interested in 
infinitesimal deformations. The interested reader can find sufficient 
material to filling this gap in the papers 
\cite{Pala1}, \cite{Pala2} and references therein.

\medskip

From now on we assume that the reader is
familiar with the notion of sheaf, sheaf  of algebras, ideal and quotient
sheaves, morphisms of sheaves.\\
If $\sF$ is a sheaf on a topological space $Y$ we denote by
$\sF_{y}$, $y\in Y$, the stalk at the point $y$. If $\sG$ is
another sheaf on $Y$ we denote by $\HOM(\sF,\sG)$ the sheaf of
morphisms from $\sF$ to $\sG$ and by
$\Hom(\sF,\sG)=\Gamma(Y,\HOM(\sF,\sG))$.\\

For every complex manifold $X$ we denote by $\sA_{X}^{p,q}$ the
sheaf of differential forms of type $(p,q)$ and
$\sA_{X}^{*,*}=\oplus_{p,q}\sA_{X}^{p,q}$. The sheaf of
holomorphic functions on $X$ is denoted by $\Oh_{X}$;
$\Omega^{*}_{X}$ (resp.: $\bar{\Omega}^{*}_{X}$) is the sheaf of
holomorphic (resp.: antiholomorphic) differential forms. By
definition $\Omega^{*}_{X}=\ker(\debar\colon \sA^{*,0}\to
\sA^{*,1})$, $\bar{\Omega}^{*}_{X}=\ker(\de\colon \sA^{0,*}\to
\sA^{1,*})$; note that $\phi\in \Omega^{*}_{X}$ if and only if
$\bar{\phi}\in
\bar{\Omega}^{*}_{X}$.\\
If $E\to X$ is a holomorphic vector bundle we denote by $\Oh_{X}(E)$ the
sheaf of holomorphic sections of $E$.\\

\bigskip
\section{Differential graded vector spaces}

This section is purely algebraic and every vector space is
considered over a fixed field $\K$;  unless otherwise specified, by the
symbol $\otimes$ we mean the tensor product $\otimes_{\K}$ over the field
$\K$.

\begin{notation}\label{VI.0.2}
We denote by $\mathbf{G}$  the category of $\Z$-graded $\K$-vector
space. The objects of $\mathbf{G}$ are the $\K$-vector spaces $V$
endowed with a $\Z$-graded direct sum decomposition
$V=\oplus_{i\in\Z}V_{i}$. The elements of $V_i$ are called
homogeneous of degree $i$. The morphisms in $\mathbf{G}$ are the
degree-preserving linear maps.
\end{notation}

If $V=\oplus_{n\in \Z}V_{n}\in \mathbf{G}$ we write $\deg(a;V)=i\in\Z$
if $a\in V_{i}$; if there is no possibility of confusion about $V$
we simply denote $\bar{a}=\deg(a;V)$.\\

Given two graded vector spaces $V,W\in\mathbf{G}$
we denote by $\Hom_{\K}^{n}(V,W)$ the vector space of 
$\K$-linear maps $f\colon V\to W$
such that $f(V_{i})\subset W_{i+n}$ for every $i\in \Z$. Observe that
$\Hom_{\K}^{0}(V,W)=\Hom_{\mathbf{G}}(V,W)$ is the space of morphisms in
the category $\mathbf{G}$.

The tensor product, $\otimes\colon \mathbf{G}\times \mathbf{G}\to \mathbf{G}$,
and the \emph{graded $Hom$}, $\Hom^{*}\colon \mathbf{G}\times \mathbf{G}\to \mathbf{G}$,
are defined in the following way:
given $V,W\in \mathbf{G}$ we set
\[ V\otimes W=\somdir{i\in\Z}{}(V\otimes W)_{i},\hbox{~~~where~~~}
(V\otimes W)_{i}=\somdir{j\in\Z}{}V_{j}\otimes W_{i-j},\]
~
\[ \Hom^{*}(V,W)=\somdir{n}{}\Hom_{\K}^{n}(V,W).\]
We denote by
\[ \langle,\rangle:\Hom^{*}(V,W)\times V\to W,\qquad
\langle f,v\rangle=f(v)\]
the natural pairing.
\bigskip

\begin{defn}\label{twistingmap}
If $V, W\in \mathbf{G}$, the \emph{twisting map} 
$T\colon V\otimes W\to W\otimes V$ is the linear map defined by the
rule $T(v\otimes w)=(-1)^{\bar{v}\,\bar{w}}w\otimes v$, for every pair
of  homogeneous elements $v\in V$, $w\in W$.
\end{defn}

Unless otherwise specified we shall use the \textbf{Koszul signs
convention}.  This means that we choose as \emph{natural isomorphism}
between $V\otimes W$ and $W\otimes V$ the twisting map $T$ and we
make every commutation rule compatible with $T$. More informally, to
``get the signs right'', whenever an ``object of degree $d$ passes on
the other side of an object of degree $h$, a sign $(-1)^{dh}$ must be
inserted''.\\
As an example, the natural map
$\langle,\rangle:V\times\Hom^{*}(V,W)\to W$ must be defined as
$\langle v,f\rangle= (-1)^{\bar{f}\,\bar{v}}f(v)$ for homogeneous $f,v$. 
Similarly, if
$f,g\in \Hom^{*}(V,W)$, their tensor product $f\otimes g\in
\Hom^{*}(V\otimes V,W\otimes W)$ must be defined on bihomogeneous
tensors as $(f\otimes g)(u\otimes v)=
(-1)^{\bar{g}\,\bar{u}}f(u)\otimes g(v)$.
\medskip

\begin{notation}\label{VI.0.4}
We denote by $\mathbf{DG}$ the category of $\Z$-graded differential
$\K$-vector spaces (also called complexes of vector spaces). The
objects of $\mathbf{DG}$ are the pairs $(V,d)$ where $V=\oplus V_i$
is an object of $\mathbf{G}$ and $d\colon V\to V$ is a linear map,
called \emph{differential} such that $d(V_i)\subset V_{i+1}$ and
$d^2=d\circ d=0$. The morphisms in $\mathbf{DG}$ are the
degree-preserving linear maps commuting with the differentials.
\end{notation}

For simplicity we will often consider $\mathbf{G}$ as the full subcategory of
$\mathbf{DG}$ whose objects are the complexes $(V,0)$ with trivial differential.

If $(V,d),(W,\delta)\in \mathbf{DG}$ then also
$(V\otimes W,d\otimes Id+Id\otimes \delta)\in \mathbf{DG}$; according to Koszul
signs convention, since $\delta\in \Hom_{\K}^{1}(W,W)$, we have
$(Id\otimes \delta)(v\otimes w)=(-1)^{\bar{v}}v\otimes \delta(w)$.\\
There exists also a natural differential $\rho$ on $\Hom^{*}(V,W)$ given by
the formula
\[\delta \langle f,v\rangle=\langle \rho f,v\rangle+
(-1)^{\bar{f}}\langle f,dv\rangle.\]

Given $(V,d)$ in $\mathbf{DG}$ we denote as usual by $Z(V)=\ker d$ the
space of cycles, by $B(V)=d(V)$ the space of boundaries and by
$H(V)=Z(V)/B(V)$ the homology of $V$. Note that $Z,B$ and $H$ are all
functors from $\mathbf{DG}$ to $\mathbf{G}$.\\
A morphism in $\mathbf{DG}$ is called a quasiisomorphism if it
induces an isomorphism in homology.\\
A differential graded vector space $(V,d)$ is called \emph{acyclic} if
$H(V)=0$.

\begin{defn}\label{VI.0.6}
Two morphisms $f,g\in \Hom_{\K}^{n}(V,W)$ are said to be
\emph{homotopic} if their difference $f-g$ is a boundary in the
complex $\Hom^{*}(V,W)$.
\end{defn}

\begin{exer}
Let $V,W$ be differential graded vector spaces, then:
\begin{enumerate}
\item $\Hom_{\mathbf{DG}}(V,W)=Z^{0}(\Hom^{*}(V,W))$.
\item If $f\in B^{0}(\Hom^{*}(V,W))\subset \Hom_{\mathbf{DG}}(V,W)$
then the induced morphism $f\colon H(V)\to H(W)$ is trivial.
\item If $f,g\in \Hom_{\mathbf{DG}}(V,W)$ are homotopic then they induce the
same morphism in homology.
\item $V$ is acyclic if and only if the identity $Id\colon V\to V$ is
homotopic to $0$.
(Hint: if $C\subset V$ is a complement of $Z(V)$, i.e.
$V=Z(V)\oplus C$, then $V$ is acyclic if and only if
$d\colon C_{i}\to Z(V)_{i+1}$ is an isomorphism for every $i$.)
\end{enumerate}
\end{exer}

The fiber product of two morphisms $B\mapor{f}D$ and $C\mapor{h}D$ in 
the category $\mathbf{DG}$ is defined as the complex 
\[
C\times_{D}B=\somdir{n}{}(C\times_{D}B)_{n},\qquad 
(C\times_{D}B)_{n}=
\{(c,b)\in C_{n}\times B_{n}\mid
h(c)=f(b)\}\]
with differential $d(c,b)=(dc,db)$.\\

A commutative diagram in
$\mathbf{DG}$
\[\xymatrix{A\ar[r]\ar[d]^{g}&B\ar[d]^{f}\\
C\ar[r]^{h}&D}\] is called cartesian if the induced morphism
$A\to
C\times_{D}B$
is an isomorphism; it is an easy exercise in
homological algebra to prove that if $f$ is a surjective (resp.:
injective) quasiisomorphism, then $g$ is a surjective (resp.:
injective) quasiisomorphism. (Hint: if $f$ is a surjective
quasiisomorphism then $\ker f=\ker g$ is acyclic.)\\

For every integer $n\in\Z$ let's choose a formal symbol $1[n]$ of 
degree $-n$ and denote by $\K[n]$ the graded vector space generated 
by $1[n]$. In other terms, the homogeneous components of 
$\K[n]\in\mathbf{G}\subset\mathbf{DG}$  are
\[
\K[n]_i=\left\{\begin{array}{ll}
\K & \mbox{if $i+n=0$}\\
0 & \mbox{otherwise}\\
\end{array}\right.
\]
For every pair of integers $n,m\in\Z$ there
exists a \textbf{canonical}  linear isomorphism $S_{n}^{m}\in
\Hom_{\K}^{n-m}(\K[n],\K[m])$; it is described by the property $
S_{n}^{m}(1[n])=1[m]$.\\

Given  $n\in\Z$, the shift
functor
$[n]\colon \mathbf{DG}\to \mathbf{DG}$ is defined by setting
$V[n]=\K[n]\otimes V$, $V\in \mathbf{DG}$,
$f[n]=Id_{\K[n]}\otimes f$, $f\in \Mor_{\mathbf{DG}}$.\\
More informally, the complex $V[n]$ is the complex $V$ with degrees
shifted
by $n$, i.e. $V[n]_{i}=V_{i+n}$, and differential multiplied by $(-1)^n$.
The shift functors preserve the subcategory $\mathbf{G}$ and commute
with the homology functor $H$. If $v\in V$ we also write 
$v[n]=1[n]\otimes v\in V[n]$.

\begin{exer} There exist  natural isomorphisms
\[\Hom_{\K}^{n}(V,W)=\Hom_{\mathbf{G}}(V[-n],W)=\Hom_{\mathbf{G}}(V,W[n]).\]
\end{exer}

\begin{ex}\label{omegacomplex}
Among the interesting objects in $\mathbf{DG}$ there
are the acyclic complexes
$\Omega[n]=\K[n]\otimes \Omega$, where
$\Omega=(\Omega_{0}\oplus \Omega_{1},d)$,
$\Omega_{0}=\K$, $\Omega_{1}=\K[-1]$  and $d\colon \Omega_{0}\to\Omega_{1}$ is
the canonical linear isomorphism $d(1[0])=1[-1]$.
The projection $p\colon\Omega\to\Omega_{0}=\K$ and the inclusion 
$\Omega_{1}\to\Omega$
are morphisms in $\mathbf{DG}$.
\end{ex}

\begin{exer}
    Let $V,W$ be differential graded vector spaces, then:
\begin{enumerate}
\item In the notation of Example \ref{omegacomplex},
two morphisms of complexes $f,g\colon V\to W$
are homotopic if and only if there exists $h\in
\Hom_{\mathbf{DG}}(V,\Omega\otimes W)$ such that $f-g=(p\otimes
Id_{|W})\circ h$.
\item If $f,g\colon V\to W$ are homotopic then
$f\otimes h$ is homotopic to $g\otimes h$ for every $h\colon V'\to
W'$.
\item (K\"unneth)
If $V$ is acyclic then $V\otimes U$ is acyclic for every $U\in\mathbf{DG}$.
\end{enumerate}
\end{exer}

\bigskip
\section{Review of terminology about algebras}

Let $R$ be a commutative ring, by a nonassociative (= not necessarily
associative) $R$-algebra we mean a $R$-module $M$ endowed with a
$R$-bilinear map $M\times M\to M$.\\
The nonassociative algebra $M$ is called \emph{unitary} if there exist
a ``unity'' $1\in M$ such that $1m=m1=m$ for every $m\in M$.\\
A \emph{left ideal} (resp.: \emph{right ideal}) of $M$ is a submodule
$I\subset M$ such that $MI\subset I$ (resp.: $IM\subset I$).
A submodule is called an \emph{ideal} if it is both a left and right ideal.\\
A homomorphism of $R$-modules $d\colon M\to M$ is called a \emph{derivation} if
satisfies the Leibnitz rule $d(ab)=d(a)b+ad(b)$. A derivation $d$ is called a 
\emph{differential} if $d^2=d\circ d=0$.\\
A $R$-algebra is \emph{associative} if $(ab)c=a(bc)$ for every $a,b,c\in M$.
Unless otherwise specified, we reserve the simple term \emph{algebra} only to
associative algebra (almost all the algebras considered in these notes are either
associative or Lie).\\
If $M$ is unitary, a \emph{left inverse} of $m\in M$ is an
element $a\in M$ such that $am=1$.
A \emph{right inverse} of $m$ is an
element $b\in M$ such that $mb=1$.\\
If $M$ is unitary and associative, an element
$m$ is called invertible if has
left and right inverses. It is easy to see that if $m$ is invertible
then every left inverse of $m$ is equal to every right inverse, in
particular there exists a unique $m^{-1}\in M$ such that
$mm^{-1}=m^{-1}m=1$.\\

\begin{exer}
Let $g$ be a Riemannian metric on an open connected subset
$U\subset \R^{n}$  and let $\phi\colon U\to \R$ be a differentiable
function (called \emph{potential}).\\
Denote by $R=C^{\infty}(U,\R)$ and by $M$ the (free of rank
$n$) $R$-module of vector fields on $U$. If $x_{1},\ldots,x_{n}$ is a
system of linear coordinates on $\R^n$ denote by:
\begin{enumerate}
    \item  $\de_{i}=\ds\desude{~}{x_{i}}\in M$,
    $\phi_{ijk}=\de_{i}\de_{j}\de_{k}\phi\in R$.\\

    \item  $g_{ij}=g(\de_{i},\de_{j})\in R$ and $g^{ij}$ the coefficients of
    the inverse matrix of $g_{ij}$.\\

    \item  $\de_{i}*\de_{j}=\ds\sum_{k,l}\phi_{ijl}g^{lk}\de_{k}$\\
\end{enumerate}
Prove that the $R$-linear extension $M\times M\to M$
of the product $*$ is independent from the choice of the linear
coordinates and write down the (differential) equation that $\phi$
must satisfy in order to have the product $*$ associative.
This equation is
called WDVV (Witten-Dijkgraaf-Verlinde-Verlinde) equation  and it is very  important in
mathematics since 1990.
\end{exer}

\bigskip
\section{dg-algebras and dg-modules}

\begin{defn}\label{cot.1.1}
A graded (associative, $\Z$-commutative)
algebra is a graded
vector space $A=\oplus A_{i}\in \mathbf{G}$ endowed with a product $A_{i}\times
A_{j}\to A_{i+j}$ satisfying the properties:
\begin{enumerate}
\item $a(bc)=(ab)c$.
\item $a(b+c)=ab+ac$, $(a+b)c=ac+bc$.
\item (Koszul sign convention) $ab=(-1)^{\bar{a}\,\bar{b}}ba$ for
$a,b$ homogeneous.
\end{enumerate}
The algebra $A$ is unitary if there exists $1\in A_{0}$ such that
$1a=a1=a$ for every $a\in A$.\end{defn}

Let $A$ be a graded algebra, then $A_{0}$ is a commutative
$\K$-algebra in the usual sense; conversely every commutative
$\K$-algebra can be considered as a graded algebra concentrated in
degree 0. If $I\subset A$ is a homogeneous left (resp.: right)
ideal then $I$ is also a right (resp.: left) ideal and the
quotient $A/I$ has a natural structure of graded algebra.
\medskip

\begin{ex}\label{VI.0.10}
The  exterior algebra $A=\bigwedge^{*}V$ of a $\K$-vector space $V$,
endowed with wedge product,
is a graded algebra with $A_{i}=\bigwedge^{i}V$.
\end{ex}

\begin{ex}\label{cot.1.4}
\emph{(Polynomial algebras.)}~
Given a set $\{x_{i}\}$, $i\in I$,  of homogeneous indeterminates
of integral degree $\bar{x_{i}}\in\Z$ we can consider the graded
algebra $\K[\{x_{i}\}]$. As a $\K$-vector space $\K[\{x_{i}\}]$ is
generated by monomials in the indeterminates $x_{i}$ subjected to the relations
$x_i    x_j=(-1)^{\bar{x_{i}}\,\bar{x_{j}}}x_j x_i$.\\
In some cases, in order to
avoid confusion about terminology, for a monomial $x_{i_{1}}^{\alpha_{1}}\ldots
x_{i_{n}}^{\alpha_{n}}$ it is defined:\begin{itemize}
\item The \emph{internal degree} $\sum_{h}\bar{x_{i_{h}}}\alpha_{h}$.
\item The \emph{external degree} $\sum_{h}\alpha_{h}$.
\end{itemize}
In a similar way it is defined $A[\{x_{i}\}]$ for every graded
algebra $A$.
\end{ex}

\begin{exer} Let $A$ be a graded algebra:
if every $a\not=0$ is invertible then $A=A_{0}$ and therefore $A$ is a
field.\\
Give an example of graded algebra where every homogeneous $a\not=0$ is
invertible but $A\not=A_{0}$.
\end{exer}

\begin{defn}\label{cot.1.5}
A \emph{dg-algebra} (differential graded algebra)
is the data of a graded algebra $A$
and a $\K$-linear map $s\colon A\to A$, called \emph{differential}, with
the properties:\begin{enumerate}
\item $s(A_{n})\subset A_{n+1}$, $s^{2}=0$.
\item (graded Leibnitz rule) $s(ab)=s(a)b+(-1)^{\bar{a}}as(b)$.
\end{enumerate}
A morphism of dg-algebras is a morphism of graded algebras commuting
with differentials; the category of dg-algebras is denoted by
$\mathbf{DGA}$.\end{defn}

\begin{ex}\label{VI.0.12}
Let $U$ be an open subset of a complex variety and denote by 
$A_{i}=\oplus_{p+q=i}\Gamma(U,\sA_{X}^{p,q})$. 
Then
$\Gamma(U,\sA_{X}^{*,*})=\oplus A_{i}$ admits
infinitely many structures of
differential graded algebras where the differential of each one of is
a linear combination $a\de+b\debar$, $a,b\in \C$.
\end{ex}

\begin{exer} Let $(A,s)$ be a unitary dg-algebra;
prove:\begin{enumerate}
\item $1\in Z(A)$.
\item $Z(A)$ is a graded subalgebra of $A$ and $B(A)$ is a
homogeneous ideal of $Z(A)$. In particular $1\in B(A)$ if and only if
$H(A)=0$.
\end{enumerate}
\end{exer}

A differential ideal of a dg-algebra $(A,s)$ is a homogeneous ideal
$I$ of $A$ such that $s(I)\subset I$; there exists an obvious
bijection between differential ideals and kernels of morphisms of
dg-algebras.

On a polynomial algebra $\K[\{x_{i}\}]$ a differential $s$ is uniquely
determined by the values $s(x_{i})$.
\medskip

\begin{ex}\label{cot.1.6}
Let $t, dt$ be indeterminates of degrees
$\bar{t}=0$,
$\bar{dt}=1$; on the polynomial algebra $\K[t,dt]=\K[t]\oplus\K[t]dt$
there exists an obvious differential $d$ such that $d(t)=dt$, $d(dt)=0$.
Since $\K$ has characteristic 0, we have $H(\K[t,dt])=\K$.
More generally if $(A,s)$ is a dg-algebra then $A[t,dt]$ is a dg-algebra
with differential $s(a\otimes p(t))=s(a)\otimes
p(t)+(-1)^{\bar{a}}a\otimes p'(t)dt$,
$s(a\otimes q(t)dt)=s(a)\otimes q(t)dt$.
\end{ex}

\begin{defn}\label{cot.1.7}
A morphism of dg-algebras $B\to A$ is
called a \emph{quasiisomorphism}
if the induced morphism $H(B)\to H(A)$ is an
isomorphism.\end{defn}

Given a morphism of dg-algebras $B\to A$
the space $\Der^{n}_{B}(A,A)$ of $B$-derivations of degree $n$ is by
definition
\[\Der^{n}_{B}(A,A)=
\{ \phi\in\Hom_{\K}^{n}(A,A)\mid
\phi(ab)\!=\!\phi(a)b+(-1)^{n\bar{a}}a\phi(b),\,  \phi(B)\!=\!0 \}.\]
We also consider the graded vector space
\[ \Der^{*}_{B}(A,A)=\somdir{n\in\Z}{}\Der^{n}_{B}(A,A)\in\mathbf{G}.\]

There exist a natural differential $d$ and a natural bracket $[-,-]$ on
$\Der^{*}_{B}(A,A)$ defined as
\[ d\colon \Der^{n}_{B}(A,A)\to\Der^{n+1}_{B}(A,A),\qquad
d\phi=d_{A}\phi-(-1)^{n}\phi d_{A}\]
and
\[ [f,g]=fg-(-1)^{\bar{f}\,\bar{g}}gf.\]

\begin{exer} Verify that, if $f\in \Der_{B}^{p}(A,A)$
and $g\in \Der_{B}^{q}(A,A)$ then
$[f,g]\in
\Der_{B}^{p+q}(A,A)$ and
$d[f,g]=[df,g]+(-1)^{p}[f,dg]$.
\end{exer}

Let $(A,s)$ be a fixed dg-algebra, by an $A$-dg-module we mean a
differential graded vector space $(M,s)$ together two associative
distributive multiplication
maps $A\times M\to M$, $M\times A\to M$ with the
properties:\begin{enumerate}
\item $A_{i}M_{j}\subset M_{i+j}$,~~ $M_{i}A_{j}\subset M_{i+j}$.
\item (Koszul) $am=(-1)^{\bar{a}\,\bar{m}}ma$, for homogeneous $a\in A$, $m\in
M$.
\item (Leibnitz) $s(am)=s(a)m+(-1)^{\bar{a}}as(m)$.
\end{enumerate}
If $A=A_{0}$ we recover the usual notion of complex of
$A$-modules.

\begin{ex}\label{VI.0.14}
For every morphism of dg-algebras $B\to A$
the space $\Der^{*}_{B}(A,A)=\oplus_{p}\Der_{B}^{p}(A,A)$ has a
natural structure of   $A$-dg-module, with left
multiplication  $(af)(b)=a(f(b))$.
\end{ex}

If $M$ is an $A$-dg-module then $M[n]=\K[n]\otimes_{\K} M$ has a natural
structure of $A$-dg-module with multiplication maps
\[ (e\otimes m)a=e\otimes ma,\qquad a(e\otimes
m)=(-1)^{n\bar{a}}e\otimes am,\qquad e\in\K[n],\,
m\in M,\, a\in A.\]

The tensor product $N\otimes_{A}M$ is defined as the quotient of
$N\otimes_{\K}M$ by the graded submodules generated by all the
elements $na\otimes m-n\otimes am$.

Given two $A$-dg-modules $(M,d_{M}),(N,d_{N})$ we denote by
\[
\Hom^{n}_{A}(M,N)=\left\{\left.
f\in\Hom_{\K}^{n}(M,N)\,\right|\,
f(ma)=f(m)a,\, m\in M, a\in A\right\}\]
\[\Hom^{*}_{A}(M,N)=\somdir{n\in\Z}{}\Hom^{n}_{A}(M,N).\]

The graded vector space $\Hom^{*}_{A}(M,N)$ has a natural structure
of $A$-dg-module with left multiplication
$(af)(m)=af(m)$ and differential
\[ d\colon \Hom^{n}_{A}(M,N)\to\Hom^{n+1}_{A}(M,N),\qquad
df=[d,f]=d_{N}\circ f-(-1)^{n}f\circ d_{M}.\]

Note that $f\in\Hom^{0}_{A}(M,N)$ is a morphism of $A$-dg-modules if
and only if $df=0$. A \emph{homotopy} between two morphism of
dg-modules $f,g\colon M\to N$ is a $h\in \Hom^{-1}_{A}(M,N)$ such
that $f-g=dh=d_{N}h+hd_{M}$. Homotopically equivalent morphisms
induce the same morphism in homology.

Morphisms of $A$-dg-modules $f\colon L\to M$, $h\colon N\to P$
induce, by composition, morphisms $f^{*}\colon
\Hom^{*}_{A}(M,N)\to \Hom_{A}^{*}(L,N)$, $h_{*} \colon
\Hom^{*}_{A}(M,N)\to \Hom_{A}^{*}(M,P)$;

\begin{lem}\label{cot.2.2}
In the above notation if $f$ is homotopic
to $g$ and $h$ is homotopic to $l$ then $f^{*}$ is homotopic to
$g^{*}$ and $l_{*}$ is homotopic to $h_{*}$.\end{lem}

\begin{proof} Let $p\in\Hom^{-1}_{A}(L,M)$ be a homotopy between $f$ and
$g$, It is a straightforward verification to see that the
composition with $p$ is a homotopy between $f^{*}$ and $g^{*}$.
Similarly we prove that $h_{*}$ is homotopic to
$l_{*}$.\end{proof}

\begin{lem}\label{cot.2.3}
\emph{(Base change)}
Let $A\to B$ be a morphism of unitary dg-algebras, $M$ an $A$-dg-module, $N$
a $B$-dg-modules. Then there exists a natural isomorphism of
$B$-dg-modules
\[ \Hom_{A}^{*}(M,N)\simeq\Hom_{B}^{*}(M\otimes_{A}B,N).\]
\end{lem}

\begin{proof} Consider the natural maps:
\[\xymatrix{\Hom_{A}^{*}(M,N)\ar@<.5ex>[r]^{L~~~}&
\Hom_{B}^{*}(M\otimes_{A}B,N)\ar@<.5ex>[l]^{R~~~}},\]
\[ Lf(m\otimes b)=f(m)b,\qquad Rg(m)=g(m\otimes 1).\]
We left as exercise the easy verification that $L,R$ are
isomorphisms of $B$-dg-modules and $R=L^{-1}$.\end{proof}

Given a morphism of dg-algebras $B\to A$ and an $A$-dg-module $M$ we set:
\[\Der^{n}_{B}(A,M)=
\{ \phi\in\Hom_{\K}^{n}(A,M)\mid
\phi(ab)\!=\!\phi(a)b+(-1)^{n\bar{a}}a\phi(b),\, \phi(B)\!=\!0 \}\]
\[ \Der^{*}_{B}(A,M)=\somdir{n\in\Z}{}\Der^{n}_{B}(A,M).\]
As in the case of $\Hom^{*}$, there exists a structure of
$A$-dg-module on $\Der^{*}_{B}(A,M)$ with product $(a\phi)(b)=a\phi(b)$
and differential
\[ d\colon \Der^{n}_{B}(A,M)\to\Der^{n+1}_{B}(A,M),\qquad
d\phi=[d,\phi]=d_{M}\phi-(-1)^{n}\phi d_{A}.\]

Given $\phi\in \Der^{n}_{B}(A,M)$ and $f\in \Hom_{A}^{m}(M,N)$ their
composition $f\phi$ belongs to $\Der^{n+m}_{B}(A,N)$.

\bigskip

\section{Kodaira-Spencer's maps in  dg-land}

In this section, we define on the central fibre of a deformation a
sheaf of differential graded algebras $\sB$ which contains (well
hidden)
the ``Taylor coefficients'' of the variation of the complex
structures given by the deformation
(the first derivative being the Kodaira-Spencer map).

\begin{lem}\label{VI.1.1}
Let $U$ be a differential manifold (not necessarily compact),
$\Delta\subset \C^{n}$ a polydisk with coordinates
$t_{1},\ldots,t_{n}$ and
$f(x,t)\in C^{\infty}(U\times \Delta,\C)$.\\
Then there exist $f_{1},\ldots,f_{n},f_{\bar{1}},\ldots,f_{\bar{n}}\in
C^{\infty}(U\times \Delta,\C)$ such that
\[f_{i}(x,0)=\desude{f}{t_{i}}(x,0),\qquad
f_{\bar{i}}(x,0)=\desude{f}{\bar{t}_{i}}(x,0)\quad\hbox{ and }\]
\[ f(x,t)=f(x,0)+\sum t_{i}f_{i}(x,t)+\sum
\bar{t_{i}}f_{\bar{i}}(x,t).\]
\end{lem}

\begin{proof} The first  2 equalities follow from the third.
Writing $t_{j}=u_{j}+iv_{j}$, $\bar{t_{j}}=u_{j}-iv_{j}$,
with $u_{j},v_{j}$ real coordinates on $\C^{n}=\R^{2n}$ we have
\[ f(x,u,v)=f(x,0,0)+\int_{0}^{1}\frac{d}{ds}f(x,su,sv)ds=\]
\[=f(x,0,0)+\sum_{j}u_{j}\int_{0}^{1}\frac{d}{du_{j}}f(x,su,sv)ds
+\sum_{j}v_{j}\int_{0}^{1}\frac{d}{dv_{j}}f(x,su,sv)ds\]
Rearranging  in the  coordinates $t_{j}$,
$\bar{t_{j}}$ we get the proof.\end{proof}

Let $X$ be a fixed complex manifold;
denote by $\DER^{*}_{\bar{\Omega}^{*}_{X}}(\sA^{0,*}_{X},\sA^{0,*}_{X})
\subset \HOM(\sA^{0,*}_{X},\sA^{0,*}_{X})$ the
sheaf of $\bar{\Omega}^{*}_{X}$-derivations of the sheaf of graded
algebras $\sA^{0,*}_{X}$; we have the following

\begin{prop}\label{VI.1.2} In the notation above there exists a
natural isomorphism of sheaves
\[\theta\colon \sA^{0,*}_{X}(T_{X})\mapor{\sim}
\DER^{*}_{\bar{\Omega}^{*}}(\sA^{0,*}_{X},\sA^{0,*}_{X}).\]
In local holomorphic coordinates $z_{1},\ldots,z_{m}$,
\[\theta\colon \sA^{0,p}_{X}(T_X)\to
\DER^{p}_{\bar{\Omega}^{*}_{X}}(\sA^{0,*}_{X},\sA^{0,*}_{X})
\subset\DER^{p}_{\C}(\sA^{0,*}_{X},\sA^{0,*}_{X})\]
is given by $\theta\left(\phi\desude{~}{z_{i}}\right)(fd\bar{z}_{I})=
\phi\wedge\desude{f}{z_{i}}d\bar{z}_{I}$.\\
The Dolbeault differential in $\sA^{0,*}_{X}(T_{X})$ corresponds, via
the isomorphism $\theta$, to the
restriction to
$\DER^{*}_{\bar{\Omega}^{*}_{X}}(\sA^{0,*}_{X},\sA^{0,*}_{X})$
of the adjoint operator
\[ [\debar,-]\colon \DER^{*}_{\C}(\sA^{0,*}_{X},\sA^{0,*}_{X})
\to\DER^{*+1}_{\C}(\sA^{0,*}_{X},\sA^{0,*}_{X}).\]
\end{prop}

\begin{proof} The morphism $\theta$ is injective and well defined.
Let $U\subset X$ be an open polydisk with coordinates
$z_{1},\ldots,z_{m}$. Take $\xi\in
\Gamma(U,\DER^{p}_{\bar{\Omega}^{*}}(\sA^{0,*}_{X},\sA^{0,*}_{X}))$
and denote $\phi_{i}=\xi(z_{i})\in \Gamma(U,\sA^{0,p}_{X})$. We want to
prove that
$\xi=\theta\left(\sum_{i}\phi_{i}\desude{~}{z_{i}}\right)$.\\
Since, over $U$,  $\sA^{0,*}_{X}$ is generated by $\sA^{0,0}_{X}$ and
$\bar{\Omega}^{*}_{X}$, it is sufficient to prove that for every
open subset $V\subset U$, every point $x\in V$ and every
$C^{\infty}$-function $f\in \Gamma(V,\sA^{0,0}_{X})$ the
equality $\xi(f)(x)=\sum_{i}\phi_{i}\desude{f}{z_{i}}(x)$ holds.\\
If $z_{i}(x)=x_{i}\in\C$, then by Lemma~\ref{VI.1.1} we can write
\[ f(z_{1},\ldots,z_{m})=f(x_{1},\ldots,x_{m})+\sum_{i=1}^{m}
(z_{i}-x_{i})f_{i}(z_{1},\ldots,z_{m})+
\sum_{i=1}^{m}
(\bar{z}_{i}-\bar{x}_{i})f_{\bar{i}}(z_{1},\ldots,z_{m})\]
for suitable $C^{\infty}$ functions $f_{{i}}$, $f_{\bar{i}}$;  therefore
\[\xi(f)(x)=\sum_{i=1}^{m}
\xi(z_{i}-x_{i})f_{i}(x_{1},\ldots,x_{m})=
\sum_{i=1}^{m}\phi_{i}\desude{f}{z_{i}}(x).\]
In particular, for $\xi,\eta\in
\Gamma(U,\DER^{p}_{\bar{\Omega}^{*}_{X}}(\sA^{0,*}_{X},\sA^{0,*}_{X}))$,
we have $\xi=\eta$ if and only if
$\xi(z_{i})=\eta(z_{i})$ for $i=1,\ldots,m$.
Since $\debar\,\bar{\Omega}^{*}_{X}\subset \bar{\Omega}^{*}_{X}$, the
adjoint operator $[\debar,-]$ preserves
$\DER^{p}_{\bar{\Omega}^{*}_{X}}(\sA^{0,*}_{X},\sA^{0,*}_{X})$,
moreover
\[ \theta\left(\debar\phi\desude{~}{z_{i}}\right) z_{j}
=(\debar \phi)\delta_{ij}=\debar(\phi\delta_{ij})-
(-1)^{\bar{\phi}}\left(\phi\desude{~}{z_{i}}\right)(\debar z_{j})=
\left[\debar,\theta\left(\phi\desude{~}{z_{i}}\right)\right]z_{j},\]
and then $\theta\debar=[\debar,-]\theta$.
\end{proof}

According to Proposition~\ref{VI.1.2},
the standard bracket on
$\DER^{*}_{\bar{\Omega}^{*}_{X}}(\sA^{0,*}_{X},\sA^{0,*}_{X})$ induces
a bracket on the sheaf $\sA^{0,*}_{X}(T_{X})$
given in local coordinates by
\[\left[f\desude{~}{z_{i}}d\bar{z}_{I},
g\desude{~}{z_{j}}d\bar{z}_{J}\right]=
\left(f\desude{g}{z_{i}}\desude{~}{z_{j}}-
g\desude{f}{z_{j}}\desude{~}{z_{i}}\right)
d\bar{z}_{I}\wedge d\bar{z}_{J}.\]
Note that for $f,g\in \Gamma(U,\sA^{0,0}_{X}(T_{X}))$, $[f,g]$ is the
usual bracket on vector fields on a differentiable manifolds.

\medskip

Let $B\subset \C^{n}$ be an open subset, $0\in B$, and
let $M_{0}\mapor{i}M\mapor{f}(B,0)$ be a deformation of a compact
complex manifold $M_{0}$; let $t_{1},\ldots,t_{n}$ be a set of 
holomorphic coordinates on $B$.\\
It is not restrictive to assume $M_{0}\subset M$ and $i$ the
inclusion map.

\begin{defn}\label{VI.1.3} In the notation above, denote by
$I_{M}\subset \sA_{M}^{*,*}$ the  graded ideal sheaf  generated
by $\bar{t_{i}},dt_{i}, d\bar{t_{i}}$.
Denote by
$\sB^{*,*}_{M}$ the quotient sheaf
$\sA_{M}^{*,*}/I_{M}$.\footnote{It is
also possible to define $\sB$ as the quotient of $\sA$ by the ideal generated
by $\bar{t_{i}},dt_{i}, d\bar{t_{i}}$ and the $C^{\infty}$ functions on
$B$ with vanishing Taylor series at 0:
the results of this chapter will remain essentially unchanged}
\end{defn}

If $z_{1},\ldots,z_{m},t_{1},\ldots,t_{n}$ are admissible (Defn. 
\ref{I.3.2}) local holomorphic coordinates on an admissible chart $W\subset 
M$,  
$W\simeq (W\cap M_{0})\times \Delta$, $0\in \Delta\subset B$
polydisk, then every $\phi\in \Gamma(W,\sB^{*,*}_{M})$ has a
representative in $\Gamma(W,\sA^{*,*}_{M})$ of the form
\[ \phi_{0}(z)+\sum_{i}t_{i}\phi_{i}(z,t),\qquad
\phi_{0}(z)\in \Gamma(W\cap M_{0},\sA^{*,*}_{M_{0}}),
\quad\phi_{i}\in \Gamma(W,\sA^{*,*}_{M}).\]
By a recursive use of Lemma~\ref{VI.1.1} we have that, for every
$s>0$,  $\phi$ is
represented by
\[ \sum_{|I|<s}t^{I}\phi_{I}(z)+\sum_{|I|=s}t^{I}\phi_{I}(z,t).\]

The ideal sheaf $I_{M}$ is preserved by the differential
operators $d,\de,\debar$ and therefore we have the corresponding
induced operators on the sheaf of graded algebras $\sB^{*,*}_{M}$.
Denoting by $\sB^{0,*}_{M}\subset\sB^{*,*}_{M}$ the image of
$\sA^{0,*}_{M}$ we have that $\sB^{0,*}_{M}$ is a sheaf of dg-algebras
with respect to the differential $\debar$.

\begin{lem}\label{VI.1.4} In the notation above, let $U,V\subset M$ 
be open subsets; if
$U\cap M_{0}=V\cap M_{0}$ then
$\Gamma(U,\sB^{*,*}_{M})=\Gamma(V,\sB^{*,*}_{M})$ and
therefore $\sB^{*,*}_{M}$ is a sheaf of
dg-algebras over $M_{0}$.
\end{lem}

\begin{proof}
It is not restrictive to assume $V\subset U$, then
$U=V\cup (U-M_{0})$ and by the sheaf properties it is sufficient
to show that $\Gamma(U-M_{0},\sB^{*,*}_{M})=
\Gamma(V-M_{0},\sB^{*,*}_{M})=0$.
More generally if $U\subset M$ is open and $U\cap M_{0}=\emptyset$
then $\Gamma(U,\sB^{*,*}_{M})=0$; in fact there exists an open covering
$U=\cup U_{i}$ such that $\bar{t_{i}}$ is invertible in $U_{i}$.\\
If $W\subset M_{0}$ is open we define $\Gamma(W,\sB^{*,*}_{M})=
\Gamma(U,\sB^{*,*}_{M})$,
where $U$ is any open subset of $M$ such that $U\cap
M_{0}=W$.\end{proof}

The pull-back $i^{*}\colon \sA_{M}^{*,*}\to
\sA_{M_{0}}^{*,*}$ factors to a surjective morphism
$i^{*}\colon \sB^{*,*}_{M}\to
\sA^{*,*}_{M_{0}}$ of sheaves of
differential graded algebras over $M_{0}$.\\
Note also that the image in $\sB^{*,*}_{M}$ of the sheaf of
antiholomorphic differential forms $\bar{\Omega}^{*}_{M}$ is
naturally isomorphic to the sheaf $\bar{\Omega}^{*}_{M_{0}}$. In fact
if $z_{1},\ldots,z_{m},t_{1},\ldots,t_{n}$ are local admissible
coordinates at a point $p\in M_{0}$ and
$\psi\in \bar{\Omega}^{q}_{M}$ then
\[ \psi\equiv\sum
\psi_{{j_{1}},\ldots,{j_{q}}}(z)d\bar{z}_{{j_{1}}}\wedge
\ldots\wedge d\bar{z}_{{j_{q}}}\quad
\pmod{\bar{t_{i}},d\bar{t_{i}}},\qquad
\de\psi_{{j_{1}},\ldots,{j_{q}}}=0.\]
Therefore to every deformation
$M_{0}\mapor{i}M\mapor{f}(B,0)$ we can associate an injective  morphism
of sheaves of dg-algebras on $M_{0}$:
\[ \bar{\Omega}^{*}_{M_{0}}\mapor{\hat{f}}\sB^{0,*}_{M}
\subset\sB^{*,*}_{M}.\]

\begin{defn}\label{algdata}
The \emph{algebraic data} of a  deformation
$M_{0}\mapor{i}M\mapor{f}(B,0)$ is the pair of morphisms
of sheaves of dg-algebras on $M_{0}$:
\[ \bar{\Omega}^{*}_{M_{0}}\mapor{\hat{f}}\sB^{*,*}_{M}\mapor{i^{*}}
\sA^{*,*}_{M_{0}}.\]
\end{defn}
We note that $\hat{f}$ injective, $i^{*}$ surjective and $i^{*}\hat{f}$
the natural inclusion. Moreover $\hat{f}$ and $i^{*}$ commute with 
both differentials $\de,\debar$.\\

If $M_{0}\mapor{j}N\mapor{g}(B,0)$ is an isomorphic
deformation then there exists an isomorphism of sheaves of
dg-algebras $\sB^{*,*}_{M}\to\sB^{*,*}_{N}$ which makes commutative the
diagram
\[\xymatrix{
\bar{\Omega}^{*}_{M_{0}}\ar[r]^{\hat{f}}\ar[d]_{\hat{g}}
&\sB^{*,*}_{M}\ar[d]^{i^{*}}\ar[dl]\\
\sB^{*,*}_{N}\ar[r]_{j^{*}}&\sA^{*,*}_{M_{0}}}\]

Similarly if $(C,0)\to (B,0)$ is a germ of holomorphic map,
then the  pull-back of differential forms
induces  a commutative diagram
\[\xymatrix{
\bar{\Omega}^{*}_{M_{0}}\ar[r]\ar[d]
&\sB^{*,*}_{M}\ar[d]\ar[dl]\\
\sB^{*,*}_{M\times_{B}C}\ar[r]&\sA^{*,*}_{M_{0}}}\]

Before going further in the theory, we will show that  
the Kodaira-Spencer map of a deformation
$M_{0}\mapor{i} M\mapor{f} (B,0)$ of a compact connected manifold $M_0$
can be recovered from its algebraic data
$\bar{\Omega}^{*}_{M_{0}}\mapor{\hat{f}}\sB^{*,*}_{M}\mapor{i^{*}}
\sA^{*,*}_{M_{0}}$

\begin{lem}\label{VI.1.5} In the notation above, 
consider $\sA^{0,*}_{M_{0}}$ as a sheaf of
$\sB^{0,*}_{M}$-modules with the structure induced by $i^{*}$ and 
denote for every $j\ge 0$.
\[ \sT^{j}_M=\frac{\DER^{j}_{\bar{\Omega}^{*}}(\sB^{0,*}_{M},\sA^{0,*}_{M_{0}})}
{i^{*}\DER^{j}_{\bar{\Omega}^{*}}(\sA^{0,*}_{M_{0}},\sA^{0,*}_{M_{0}})}.\]
Then there exists a natural linear isomorphism
\[ T_{0,B}=\ker(\Gamma(M_{0},\sT^{0}_M)\to \Gamma(M_{0},\sT^{1}_M),\quad 
h\mapsto \debar_{\sA}h-h\debar_{\sB}).\] 
\end{lem}

\begin{proof}
We consider $T_{0,B}$ as the $\C$-vector space of
$\C$-derivations $\Oh_{B,0}\to\C$.
Let $h\in \Gamma(M_{0},
\DER^{*}_{\bar{\Omega}^{*}}(\sB^{0,*}_{M},\sA^{0,*}_{M_{0}}))$ be
such that 
$\debar_{\sA}h-h\debar_{\sB}\in
i^{*}\DER^{1}_{\bar{\Omega}^{*}}(\sA^{0,*}_{M_{0}},\sA^{0,*}_{M_{0}})$; 
in particular $\debar h(t_i)=0$ for every $i$,  
the function $h(t_i)$ is holomorphic and then
constant.  Moreover,
$h(t_i)=0$ for every $i$ if and only if  $h(\ker i^{*})=0$ if and only if 
$h\in i^{*}\DER^{0}_{\bar{\Omega}^{*}}(\sA^{0,*}_{M_{0}},\sA^{0,*}_{M_{0}})$.\\ 
This gives a linear injective 
morphism
\[\ker(\Gamma(M_{0},\sT^{0}_M)\to \Gamma(M_{0},\sT^{1}_M))\to T_{0,B}.\]
To prove the surjectivity, consider a derivation $\delta\colon \Oh_{B,0}\to\C$ and 
let $M_{0}=\cup U_{a}$, $a\in \sI$,
be a locally finite covering with every $U_{a}$ open polydisk with
coordinate systems $z_{1}^{a},\ldots,z_{m}^{a}\colon U_{a}\to
\C$. Let $t_{1},\ldots,t_{n}$ be coordinates on $B$.\\
Over $U_{a}$, every $\phi\in \sB^{0,*}_{M}$ can be written
as $\phi_{0}(z)+\sum t_{i}\phi_{i}(z)+\sum t_{i}t_{j}\ldots$, with
$\phi_{i}\in \sA^{0,*}_{M_{0}}$.
Setting $h_{a}(\phi)=\sum_{i}\delta(t_{i})\phi_{i}$ we see
immediately that $h_{a}$ is a
$\bar{\Omega}^{*}_{U_{a}}$-derivation lifting $\delta$.
Taking a partition of unity $\rho_{a}$ subordinate to the covering
$\{U_{a}\}$, we can take $h=\sum_{a}\rho_{a}h_{a}$.
\end{proof}

Let $h\in \Gamma(M_{0},
\DER^{*}_{\bar{\Omega}^{*}}(\sB^{0,*}_{M},\sA^{0,*}_{M_{0}}))$ be
such that 
$\psi=\debar_{\sA}h-h\debar_{\sB}\in
i^{*}\DER^{1}_{\bar{\Omega}^{*}}(\sA^{0,*}_{M_{0}},\sA^{0,*}_{M_{0}})$ and let 
$\delta\colon \Oh_{B,0}\to\C$ be the corresponding derivation, $\delta(t_i)=h(t_i)$.\\ 
According to the isomorphism (Proposition~\ref{VI.1.2})
$\DER^{j}_{\bar{\Omega}^{*}}(\sA^{0,*}_{M_{0}},\sA^{0,*}_{M_{0}})=
\sA^{0,j}_{M_{0}}(T_{M_{0}})$ we have
$\psi\in \Gamma(M_{0},\sA^{0,1}(T_{M_{0}}))$.\\
Moreover, being $\psi$ exact in the complex
$\DER^{*}_{\bar{\Omega}^{*}}(\sB^{0,*}_{M},\sA^{0,*}_{M_{0}})$, it is closed
in $\DER^{0}_{\bar{\Omega}^{*}}(\sA^{0,*}_{M_{0}},\sA^{0,*}_{M_{0}})$,
$\psi$ is a $\debar$-closed
form of $\Gamma(M_{0},\sA^{0,1}(T_{M_{0}}))$ and the cohomology class
$[\psi]\in H^{1}(M_{0},T_{M_{0}})$ is depends only on 
the class of $h$ in $\Gamma(M_{0},\sT^{0}_M)$. It is now easy to prove 
the following

\begin{prop}\label{VI.1.6} In the above notation,
$[\psi]=[\debar h-h\debar]=\KS_{f}(\delta)$.
\end{prop}

\begin{proof} (sketch) Let $\eta\in \Gamma(M,\sA^{0,0}_{M}(T_{M}))$ be a
complexified vector field such that $(f_{*}\eta)(0)=\delta$.
We may interpret $\eta$ as a $\bar{\Omega}^{*}_{M}$-derivation
of degree 0 $\eta\colon \sA^{0,*}_{M}\to \sA^{0,*}_{M}$;
passing to the quotient we get a $\bar{\Omega}^{*}_{M_{0}}$-derivation
$h\colon \sB^{0,*}_{M}\to \sA^{0,*}_{M_{0}}$.
The condition $(f_{*}\eta)(0)=\delta$ means that $h$ lifts $\delta$
and therefore $\psi$ corresponds to the restriction of $\debar\eta$ to the
fibre $M_{0}$.\end{proof}

\bigskip

\section{Transversely holomorphic trivializations}

\begin{defn}\label{VI.2.1} A transversely holomorphic trivialization
of a deformation $M_{0}\mapor{i}M\mapor{f}(B,0)$ is a diffeomorphism
$\phi\colon M_{0}\times \Delta\to f^{-1}(\Delta)$ such that:
\begin{enumerate}
\item  $\Delta\subset B$ is an open neighbourhood of the base
point ${0}\in B$

\item  $\phi(x,0)=i(x)$ and $f\phi$ is the projection on the
second factor.

\item  For every $x\in M_{0}$, $\phi\colon \{x\}\times \Delta\to
M$ is a holomorphic function.
\end{enumerate}\end{defn}

\begin{thm}\label{VI.2.2} Every deformation of a compact complex
manifold admits a transversely holomorphic trivialization.\end{thm}

\begin{proof} (cf. also \cite{clemens}, \cite{Voisin})
Let $f\colon M\to B$ be a deformation of $M_{0}$;
it is not restrictive to assume $B\subset \C^{n}$ a polydisk with
coordinates $t_{1},\ldots,t_{n}$ and $0\in B$ the base point of the
deformation. We identify $M_{0}$ with the central fibre $f^{-1}(0)$.\\
After a possible shrinking of $B$ there exist a finite open
covering $M=\cup W_{a}$, $a=1,\ldots,r$, and holomorphic projections
$p_{a}\colon W_{a}\to U_{a}=W_{a}\cap M_{0}$ such that
$(p_{a},f)\colon W_{a}\to U_{a}\times B$ is a biholomorphism for
every $a$ and $U_{a}$ is a local chart with
coordinates $z^{a}_{i}\colon U_{a}\to \C$, $i=1,\ldots,m$.\\
Let $\rho_{a}\colon M_{0}\to [0,1]$ be a $C^{\infty}$ partition of unity
subordinate to the covering $\{U_{a}\}$ and denote
$V_{a}=\rho_{a}^{-1}(]0,1])$; we note that $\{V_{a}\}$ is a covering
of $M_{0}$ and $\bar{V_{a}}\subset U_{a}$. After a possible shrinking
of $B$ we may assume $p_{a}^{-1}(\bar{V_{a}})$ closed in $M$.\\
For every subset $C\subset\{1,\ldots,r\}$ and every
$x\in M_{0}$ we denote
\[ H_{C}=
\left(\bigcap_{a\in C}W_{a}-\bigcup_{a\not\in
C}p_{a}^{-1}(\bar{V_{a}})\right)
\times \left(\bigcap_{a\in C}U_{a}-\bigcup_{a\not\in
C}\bar{V_{a}}\right)
\subset M\times M_{0},\]
\[ C_{x}=\{a\mid x\in \bar{V_{a}}\,\},\qquad   H=\bigcup_{C}H_{C}.\]
Clearly $(x,x)\in H_{C_{x}}$ and then $H$ is an open subset of
$M\times M_{0}$ containing the graph $G$ of the inclusion $M_{0}\to
M$. Since the projection $pr\colon M\times M_{0}\to M$ is open and $M_{0}$ is
compact, after a possible shrinking of $B$ we may assume $pr(H)=M$.\\
Moreover if $(y,x)\in H$ and $x\in \bar{V_{a}}$ then
$(y,x)\in H_{C}$ for some $C$ containing $a$ and therefore $y\in
W_{a}$.\\
For every $a$ consider the  $C^{\infty}$ function $q_{a}\colon H\cap (M\times
U_{a})\to\C^{m}$,
\[ q_{a}(y,x)=\sum_{b}\rho_{b}(x)\desude{z^{a}}{z^{b}}(x)
(z^{b}(p_{b}(y))-z^{b}(x)).\]
By the properties of $H$, $q_{a}$ is well defined and separately
holomorphic in the variable $y$.
If $(y,x)\in H\cap
(M\times (U_{a}\cap U_{c}))$ then
\[ q_{c}(y,x)=\desude{z^{c}}{z^{a}}(x)q_{a}(y,x)\]
and then
\[\Gamma=\{(y,x)\in H\mid q_{a}(y,x)=0\, \hbox{ whenever } x\in
U_{a}\}\]
is  a well defined closed subset of $H$.\\
If $y\in V_{a}\subset M_{0}$ and $x$ is sufficiently near to $y$ then
$x\in (\bigcap_{b\in C_{y}}U_{b}-\bigcup_{b\not\in C}\bar{V_{b}})$ and, for
every $b\in C_{y}$,
\[z^{b}(y)=z^{b}(x)+\desude{z^{b}}{z^{a}}(x)(z^{a}(y)-z^{a}(x))+
o(\|z^{a}(y)-z^{a}(x)\|).\]
Therefore
\[q_{a}(y,x)=z^{a}(y)-z^{a}(x)+o(\|z^{a}(y)-z^{a}(x)\|).\]
In particular the map $x\mapsto q_{a}(y,x)$ is a local diffeomorphism at
$x=y$.\\
Denote $K\subset H$ the open subset of points $(y,x)$ such that, if
$y\in p^{-1}_{a}(V_{a})$ then
$u\mapsto q_{a}(y,u)$ has maximal rank at $u=x$; note that $K$ contains
$G$.\\
Let $\Gamma_{0}$ be the connected component of $\Gamma\cap K$ that
contains $G$; $\Gamma_{0}$ is a $C^{\infty}$-subvariety of $K$ and the
projection $pr\colon \Gamma_{0}\to M$ is a local diffeomorphism.
Possibly shrinking $B$ we may assume that $pr\colon\Gamma_{0}\to M$
is a diffeomorphism.\\
By implicit function theorem $\Gamma_{0}$ is the graph of a
$C^{\infty}$ projection $\gamma\colon M\to M_{0}$.\\
After a possible shrinking of $B$, the map $(\gamma,f)\colon M\to
M_{0}\times B$ is a diffeomorphism, take $\phi=(\gamma,f)^{-1}$.\\
To prove that, for every $x\in M_{0}$, the map
$t\mapsto \phi(x,t)$ is holomorphic we note
that $f\colon \phi(\{x\}\times
B)\to B$ is bijective and therefore
$\phi(x,-)=f^{-1}pr\colon \{x\}\times B\to \phi(\{x\}\times
B)$.\\
The map $f^{-1}\colon B\to \phi(\{x\}\times B)$
is holomorphic if and only if
$\phi(\{x\}\times B)=\gamma^{-1}(x)$ is a holomorphic subvariety and
this is true because for $x$ fixed every
map $y\mapsto q_{a}(y,x)$ is holomorphic.\end{proof}

\medskip

Let $z_{1},\ldots,z_{m},t_{1},\ldots,t_{n}$ be an admissible system
of local coordinates at a point $p\in M_{0}\subset M$.
$z_{1},\ldots,z_{m},t_{1},\ldots,t_{n}$ is also
a system of local coordinates over $M_{0}\times B$.\\
In these systems, a transversely holomorphic
trivialization $\phi\colon M_{0}\times B\to M$ is written as
\[ \phi(z,t)=(\phi_{1}(z,t),\ldots,\phi_{m}(z,t),t_{1},\ldots,t_{n}),\]
where every $\phi_{i}$, being holomorphic in $t_{1},\ldots,t_{n}$, can be
written as
\[ \phi_{i}(z,t)=z_{i}+\sum_{I>0}t^{I}\phi_{i,I}(z),\qquad
I=(i_{1},\ldots,i_{n}),\quad \phi_{i,I}\in C^{\infty}.\]

In a neighbourhood of $p$,
\[ \phi^{*}dz_{i}=dz_{i}+\sum_{I>0}t^{I}\sum_{j=1}^{m}
\left(\desude{\phi_{i,I}}{z_{j}}dz_{j}+
\desude{\phi_{i,I}}{\bar{z_{j}}}d\bar{z_{j}}\right),\quad
\pmod{I_{M_{0}\times B}},\]
\[ \phi^{*}d\bar{z_{i}}=d\bar{z_{i}},\qquad\pmod{I_{M_{0}\times B}}.\]

\begin{lem}\label{VI.2.3} Every transversely holomorphic
trivialization $\phi\colon M_{0}\times B\to M$ induces  isomorphisms
of sheaves of graded algebras over $M_{0}$
\[ \phi^{*}\colon \sB^{*,*}_{M}\to\sB^{*,*}_{M_{0}\times B},\qquad
\phi^{*}\colon \sB^{0,*}_{M}\to\sB^{0,*}_{M_{0}\times B}\]
making commutative the diagrams
\[\xymatrix{
\bar{\Omega}^{*}_{M_{0}}\ar[r]\ar[d]
&\sB^{*,*}_{M}\ar[d]\ar[dl]^{\phi^*}\\
\sB^{*,*}_{M_{0}\times B}\ar[r]&\sA^{*,*}_{M_{0}}}\qquad
\xymatrix{
\bar{\Omega}^{*}_{M_{0}}\ar[r]\ar[d]
&\sB^{0,*}_{M}\ar[d]\ar[dl]^{\phi^*}\\
\sB^{0,*}_{M_{0}\times B}\ar[r]&\sA^{0,*}_{M_{0}}}\]
\end{lem}

\textbf{Beware:} It is not true in general that, for $p>0$,
$\phi^{*}(\sB^{p,q})\subset \sB^{p,q}$.
\begin{proof}
For every open subset $U\subset M$,
the pull-back
\[ \phi^{*}\colon \Gamma(U,\sA^{*,*}_{M})\to
\Gamma(\phi^{-1}(U),\sA^{*,*}_{M_{0}\times B})\]
is an isomorphism preserving the ideals $I_{M}$ and
$I_{M_{0}\times B}$.
Since $U\cap M_{0}=\phi^{-1}(U)\cap M_{0}$, the pull-back
$\phi^{*}$ induces to the quotient an isomorphism of sheaves of
graded algebras $\phi^{*}\colon \sB^{*,*}_{M}\to
\sB^{*,*}_{M_{0}\times B}$.\\
From the above formulas follows that $\phi^{*}(\sB^{p,k-p}_{M})\subset
\oplus_{q\le p}\sB^{q,k-q}_{M_{0}\times B}$ and
$\phi^{*}$ is the identity on
$\bar{\Omega}^{*}_{M_{0}}$. This shows that
$\phi^{*}(\sB^{0,*}_{M})=\sB^{0,*}_{M_{0}\times B}$ and
proves the commutativity of the diagrams.\end{proof}

The $\debar$ operator on $\sA^{*,*}_{M}$ factors to $\sB^{0,*}_{M}$ and
therefore induces operators
\[ \debar\colon
\sB^{0,*}_{M}\to \sB^{0,*+1}_{M},\qquad
\debar_{\phi}=\phi^{*}\debar(\phi^{*})^{-1}\colon
\sB^{0,*}_{M_{0}\times B}\to \sB^{0,*+1}_{M_{0}\times B}.\]

If $z_{1},\ldots,z_{m},t_{1},\ldots,t_{n}$ are admissible local
coordinates at $p\in M_{0}$, we have
\[ (\phi^{*})^{-1}dz_{i}=dz_{i}+\sum_{j=1}^{m}
a_{ij}dz_{j}+
b_{ij}d\bar{z_{j}},\quad
\pmod{I_{M}},\]
where $a_{ij}$, $b_{ij}$ are $C^{\infty}$ functions vanishing on $M_{0}$ and
\[ (\phi^{*})^{-1}d\bar{z_{i}}=d\bar{z_{i}},\qquad\pmod{I_{M}}.\]
Thus we get immediately that $\debar_{\phi}(d\bar{z_{i}})=0$.
Let's now $f$ be a $C^{\infty}$ function in a neighbourhood of $p\in
U\subset M_{0}\times B$ and let $\bar{\pi}\colon \sA^{*,*}_{M}\to
\sA^{0,*}_{M}$ be the projection.
By definition $\debar_{\phi}f$ is the
class in $\sB^{0,*}_{M_{0}\times B}$ of
\[\phi^{*}\bar{\pi}d(\phi^{*})^{-1}f=\phi^{*}\bar{\pi}(\phi^{*})^{-1}df=
\sum_{i=1}^{m}\desude{f}{z_{i}}
\phi^{*}\bar{\pi}(\phi^{*})^{-1}dz_{i}+
\sum_{i=1}^{m}\desude{f}{\bar{z_{i}}}
\phi^{*}\bar{\pi}(\phi^{*})^{-1}d\bar{z_{i}}\]
and then
\[\debar_{\phi}f=\debar
f+\sum_{ij}b_{ij}\desude{f}{z_{i}}d\bar{z_{j}}.\]
If $\psi\colon M_{0}\times B\to M$
is another transversely holomorphic trivialization and
$\theta=\phi^{*}(\psi^{*})^{-1}$ then
$\debar_{\psi}=\theta\debar_{\phi}\theta^{-1}$.

\bigskip

\section{Infinitesimal deformations}

Let $M_{0}\mapor{i}M\mapor{f}(B,0)$ be a deformation of a compact
complex manifold and $J\subset \Oh_{B,0}$ a proper ideal such that
$\sqrt{I}=\ide{m}_{B,0}$; after a possible shrinking of $B$ we can assume
that:
\begin{enumerate}
\item  $B\subset \C^{n}$ is a polydisk with coordinates
$t_{1},\ldots,t_{n}$ and $J$ is generated by a
finite number of holomorphic functions on $B$.

\item  $f\colon M\to B$ is a family admitting a transversely
holomorphic trivialization $\phi\colon M_{0}\times B\to M$.
\end{enumerate}

Denote by $(X,0)$ the fat point $(B,0,J)$ and by
$\Oh_{X,0}=\Oh_{B,0}/J$ its associated analytic algebra.
If $\ide{m}_{B,0}^{s}\subset J$ then the holomorphic functions
$t^{I}$, $I=(i_{1},\ldots,i_{n})$, $|I|<s$, generate $\Oh_{X,0}$ as
a $\C$-vector space.

Denote by $I_{M,J}\subset \sA_{M}^{*,*}$ the graded ideal sheaf generated
by $I_{M}$ and $J$,
$\sB^{*,*}_{M,J}=\sA_{M}^{*,*}/I_{M,J}=\sB_{M}^{*,*}/(J)$,
$\Oh_{M,J}=\Oh_{M}/(J)$. The same argument used in
Lemma~\ref{VI.1.4}
shows that $\sB^{*,*}_{M,J}$ and $\Oh_{M,J}$ are sheaves over $M_{0}$.
In the same manner we define $\sB^{*,*}_{M_{0}\times B,J}$
\begin{lem}\label{VI.4.1}
Let $U\subset M_{0}$ be an open subset, then there exist
isomorphisms
\[ \Gamma(U,\Oh_{M_{0}\times
B,J})=\Gamma(U,\Oh_{M_{0}})\otimes_{\C}\Oh_{X,0},\qquad
\Gamma(U,\sB^{*,*}_{M_{0}\times
B,J})=\Gamma(U,\sA^{*,*}_{M_{0}})\otimes_{\C}\Oh_{X,0}.\]
The same holds for $M$ instead of $M_{0}\times B$ provided that
$U$ is contained in an admissible coordinate chart.\end{lem}

\begin{proof}
We have seen that every $\phi\in \Gamma(U,\sB^{p,q}_{M_{0}\times
B,J})$ is represented by a form
$\sum_{|I|<s}t^{I}\phi_{I}$, with $\phi_{I}\in
\Gamma(U,\sA^{p,q}_{M_{0}})$. Writing every $t^{I}$ as a linear
combination of  the elements of a fixed
basis of $\Oh_{X,0}$ and rearranging the terms we get the desired result.
The same argument applies to $\Oh_{M_{0}\times B,J}$ and, if $U$ is
sufficiently small, to $\sB^{*,*}_{M,J}$, $\Oh_{M,J}$.
\end{proof}

\begin{cor}\label{VI.4.2} $\Oh_{M,J}=\ker(\debar\colon
\sB^{0,0}_{M,J}\to\sB^{0,1}_{M,J})$.\end{cor}

\begin{proof} If $U\subset M_{0}$ is a sufficiently small open subset, we
have
$\Gamma(U,\sB^{*,*}_{M,J})=\Gamma(U,\sA^{*,*}_{M_{0}})\otimes_{\C}\Oh_{X,0}$
and then
\[\ker\left(\debar\colon
\Gamma(U,\sB^{0,0}_{M,J})\to\Gamma(U,\sB^{0,1}_{M,J})\right)=\]
\[=\ker\left(\debar\colon
\Gamma(U,\sA^{0,0}_{M_{0}})\otimes \Oh_{X,0}
\to\Gamma(U,\sA^{0,1}_{M_{0}})\otimes \Oh_{X,0}\right)=
\Gamma(U,\Oh_{M,J}).\]
\end{proof}

The transversely holomorphic trivialization $\phi$ gives a commutative
diagram of morphisms of sheaves of graded algebras
\[\xymatrix{
\bar{\Omega}^{*}_{M_{0}}\otimes \Oh_{X,0}\ar[r]\ar[d]
&\sB^{0,*}_{M,J}\ar[d]\ar[dl]^{\phi^*}\\
\sB^{0,*}_{M_{0}\times B,J}\ar[r]&\sA^{0,*}_{M_{0}}}\]
with $\phi^{*}$ an isomorphism.
The operator $\debar_{\phi}=\phi^{*}\debar(\phi^{*})^{-1}$ is a
$\Oh_{X,0}$-derivation of degree 1 such that
$\debar_{\phi}^{2}=\frac{1}{2}[\debar_{\phi},\debar_{\phi}]=0$
and then
$\eta_{\phi}=\debar_{\phi}-\debar\colon \sB^{0,*}_{M_{0}\times B,J}\to
\sB^{0,*+1}_{M_{0}\times B,J}$ is a $\bar{\Omega}^{*}_{M_{0}}\otimes
\Oh_{X,0}$-derivation.\\

According to Lemma~\ref{VI.4.1} we have $\sB^{0,*}_{M_{0}\times
B,J}=\sA^{0,*}_{M_{0}}\otimes \Oh_{X,0}$; moreover,
if $g_{0}=1,g_{1}(t),\ldots,g_{r}(t)$ is  a basis of $\Oh_{X,0}$ with
$g_{i}\in \ide{m}_{X,0}$ for $i>0$, then we can write
$\eta_{\phi}=\sum_{i>0}g_{i}(t)\eta_{i}$, with every $\eta_{i}$ a
$\bar{\Omega}^{*}_{M_{0}}$-derivation of degree 1 of
$\sA^{0,*}_{M_{0}}$. By Proposition~\ref{VI.1.2}
$\eta_{\phi}\in \Gamma(M_{0},\sA^{0,1}(T_{M_{0}}))\otimes \ide{m}_{X,0}$.\\
In local holomorphic coordinates $z_{1},\ldots,z_{m}$ we have
$\debar_{\phi}(d\bar{z_{i}})=0$ and
\[\debar_{\phi}f=\debar
f+\sum_{i,j,k}g_{i}(t)b^{i}_{j,k}(z)\desude{f}{z_{j}}d\bar{z_{k}}\]
for every $C^{\infty}$ function $f$. The $b^{i}_{j,k}$ are $C^{\infty}$
functions on $M_{0}$.\\

A different choice of transversely holomorphic trivialization
$\psi\colon M_{0}\times B\to M$ gives a conjugate operator
$\debar_{\psi}=\theta\debar_{\phi}\theta^{-1}$, where
$\theta=\phi^{*}(\psi^{*})^{-1}$.

\smallskip

This discussion leads naturally to the definition of deformations of a
compact complex manifolds  over a fat points.

\begin{defn}\label{VI.4.3} A deformation of $M_{0}$ over a fat point $(X,0)$
is a section
\[\eta\in \Gamma(M_{0},\sA^{0,1}(T_{M_0}))\otimes \ide{m}_{X,0}=
\Der^{1}_{\bar{\Omega}^{*}_{M_{0}}}(\sA^{0,*}_{M_{0}},
\sA^{0,*}_{M_{0}})\otimes \ide{m}_{X,0}\]
such that the operator $\debar+\eta\in
\Der^{1}_{\Oh_{X,0}}(\sA^{0,*}_{M_{0}}\otimes\Oh_{X,0},
\sA^{0,*}_{M_{0}}\otimes\Oh_{X,0})$ is a differential,
i.e. $[\debar+\eta,\debar+\eta]=0$.\\
Two deformations $\eta,\mu\in \Gamma(M_{0},\sA^{0,1}(T_{M_{0}}))\otimes
\ide{m}_{X,0}$ are isomorphic if and only if there exists an automorphism of
sheaves of graded algebras $\theta\colon
\sA^{0,*}_{M_{0}}\otimes\Oh_{X,0}\to
\sA^{0,*}_{M_{0}}\otimes\Oh_{X,0}$ commuting with the projection
$\sA^{0,*}_{M_{0}}\otimes\Oh_{X,0}\to \sA^{0,*}_{M_{0}}$ and leaving
point fixed the subsheaf
$\bar{\Omega}^{*}_{M_{0}}\otimes\Oh_{X,0}$ such that
$\debar+\mu=\theta(\debar+\eta)\theta^{-1}$.\end{defn}

According to \ref{VI.1.2} the adjoint operator $[\debar,-]$ corresponds to
the Dolbeault differential in the complex $\sA^{0,*}(T_{M_{0}})$ and
therefore $(\debar+\eta)^{2}=0$ if and only if
$\eta\in \Gamma(M_{0},\sA^{0,1}(T_{M_{0}}))\otimes
\ide{m}_{X,0}$ satisfies the
\textbf{Maurer-Cartan} equation
\[ \debar\eta+\frac{1}{2}[\eta,\eta]=0\in
\Gamma(M_{0},\sA^{0,2}(T_{M_{0}}))\otimes
\ide{m}_{X,0}.\]

We denote with both $\Def_{M_{0}}(X,0)$ and
$\Def_{M_{0}}(\Oh_{X,0})$
the set of isomorphism classes of
deformations of $M_{0}$ over $(X,0)$. By an \emph{infinitesimal
deformation} we mean a deformation over a fat point; by a {\em
first order deformation} we mean a deformation over a fat point
$(X,0)$ such that $\ide{m}_{X,0}\not=0$ and $\ide{m}_{X,0}^{2}=0$.

The Proposition~\ref{VI.1.6} allows to extend naturally the definition of
the Kodaira-Spencer map $\KS\colon T_{0,X}\to H^{1}(M_{0},T_{M_{0}})$
to every infinitesimal deformation over $(X,0)$.\\
Consider in fact $\delta\in \Der_{\C}(\Oh_{X,0},\C)=T_{0,X}$, then
\[h=Id\otimes\delta\colon \sA^{0,*}_{M_{0}}\otimes \Oh_{X,0}
\to\sA^{0,*}_{M_{0}}\]
is a $\bar{\Omega}^{*}_{M_{0}}$-derivation
lifting $\delta$.
Since
\[(\debar h-h(\debar+\eta))(f\otimes
1)=h(-\eta(f))\]
we may define $\KS(\delta)$ as the cohomology class
of the derivation
\[\sA^{0,*}_{M_{0}}\to \sA^{0,*+1}_{M_{0}},\qquad
f\mapsto h(-\eta(f)),\]
which corresponds, via the isomorphism of Proposition~\ref{VI.1.2}, to
\[(Id\otimes\delta)(-\eta)\in
\Gamma(M_{0},\sA^{0,1}(T_{M_{0}})),\]
where $Id\otimes\delta\colon\Gamma(M_{0},\sA^{0,1}(T_{M_{0}}))\otimes
\ide{m}_{X,0}$.
According to the Maurer-Cartan
equation $\debar\eta=-\frac{1}{2}[\eta,\eta]\in
\Gamma(M_{0},\sA^{0,2}(T_{M_{0}}))\otimes
\ide{m}_{X,0}^{2}$ and then
\[ \debar((Id\otimes\delta)(-\eta))=(Id\otimes\delta)(-\debar\eta)=0.\]

A morphism of fat points $(Y,0)\to (X,0)$ is the same of a morphism
of local $\C$-algebras $\alpha\colon \Oh_{X,0}\to \Oh_{Y,0}$;
It is natural to set $Id\otimes\alpha(\eta)\in
\Gamma(M_{0},\sA^{0,1}(T_{M_{0}}))\otimes \ide{m}_{Y,0}$ as the pull-back
of the deformation $\eta$. It is immediate to see that the
Kodaira-Spencer map of $Id\otimes\alpha(\eta)$ is the composition
of the Kodaira-Spencer map of $\eta$ and the differential
$\alpha\colon T_{Y,0}\to T_{X,0}$.

\bigskip

\section{Historical survey,~\ref{CAP:INFINITESIMAL}}
\label{sec:HistInf}

The importance of infinitesimal deformations increased considerably
after the proof (in the period 1965-1975) of several
ineffective existence results of semiuniversal deformations of
manifolds, of maps etc..,
over singular bases.\\
The archetype of these results is the well known theorem of
Kuranishi (1965) \cite{Kura}, asserting the existence of the semiuniversal
deformation of a compact complex manifold over a base which is an
analytic singularity.
An essentially equivalent formulation of Kuranishi theorem is the
following
\begin{thm}\label{VI.5.1} Let $M_{0}$ be a compact complex manifold with
$n=h^{1}(M_{0},T_{M_{0}})$, $r=h^{2}(M_{0},T_{M_{0}})$.\\
Then there exist a polydisk $\Delta\subset \C^{n}$, a section
$\eta\in \Gamma(M,\sA^{0,1}(T_{f}))$, being $M=M_{0}\times \Delta$
and $f\colon M\to \Delta$ the projection, and
$q=(q_{1},\ldots,q_{r})\colon
\Delta\to \C^{r}$ a holomorphic map such that:
\begin{enumerate}
\item $q(0)=0$ and $\desude{q_{i}}{t_{j}}(0)=0$ for every $i,j$,
being $t_{1},\ldots,t_{n}$ holomorphic coordinates on $\Delta$.
\item $\eta$ vanishes on $M_{0}$
and it is holomorphic in $t_{1},\ldots,t_{n}$; this means
that it is possible to write
\[ \eta=\sum_{I>0}t^{I}\eta_{I},\qquad
I=(i_{1},\ldots,i_{n}),\quad \eta_{I}\in
\Gamma(M_{0},\sA^{0,1}(T_{M_{0}})).\]
\item $\eta$ satisfies the Maurer-Cartan equation to modulus
$q_{1},\ldots,q_{s}$, i.e.
\[ \debar\eta+\frac{1}{2}[\eta,\eta]\in \sum q_{i}
\Gamma(M,\sA^{0,2}(T_{f})).\]
\item Given a fat point $(X,0)$ the natural map
\[ \eta\colon \Mor_{{\mathbf{An}}}(\Oh_{\Delta,0}/(q_{1},\ldots,q_{s}),
\Oh_{X,0})\to \Def_{M_{0}}(X,0),\qquad \alpha\mapsto \alpha(\eta)\]
is surjective for every $(X,0)$
and bijective whenever $\Oh_{X,0}=\quot{\C[t]}{(t^{2})}$.
\end{enumerate}
\end{thm}

It is now clear that the study of infinitesimal deformations can be
used to deduce the structure of the holomorphic map $q$ and the
existence of the semiuniversal deformation over a smooth base.
For example we have the following

\begin{cor}\label{VI.5.2}
Let $M_{0}$ be a compact complex manifolds such that for
every $n\ge 2$ the natural map $\Def_{M_{0}}(\C[t]/(t^{n}))\to
\Def_{M_{0}}(\C[t]/(t^{2}))$ is surjective. Then $M_{0}$ has a
semiuniversal deformation
$M_{0}\mapor{}M\mapor{}(H^{1}(M_{0},T_{M_{0}}),0)$.\end{cor}

\begin{proof} (sketch) In the notation of
Theorem~\ref{VI.5.1} we have $(q_{1},\ldots,q_{s})\subset
\ide{m}_{\Delta,0}^{2}$ and then,
according to Proposition~\ref{V.1.12}, $q_{1}=\ldots=q_{s}=0$.
In particular $\eta$ satisfies the Maurer-Cartan equation and by the
Newlander-Nirenberg's theorem (cf. \cite[1.4]{Catacime}, \cite{Voisin})
the small variation of almost complex structure
\cite[2.1, 2.5]{Catacime}, \cite{Voisin}
\[ -\eta\colon \sA^{1,0}_{M}\to \sA^{0,1}_{M},\qquad
-^{t}\!\eta\colon T^{0,1}_{M}\to T^{1,0}_{M}\]
is integrable and gives a  complex structure on $M$
with structure sheaf $\Oh_{M,\eta}=\ker(\debar+\eta\de\colon
\sA^{0,0}_{M}\to \sA^{0,1}_{M})$.\\
The projection map $(M,\Oh_{M,\eta})\to \Delta$ is a family with
bijective Kodaira-Spencer map, by completeness theorem \ref{I.5.1} it is a
semiuniversal deformation.\end{proof}

It is useful to remind here the following result proved by Malgrange \cite{malgrange}

\begin{thm} Let  $q_1,\ldots, q_m\colon (\C^n,0)\to \C$ be germs of holomorphic
functions and $f\colon (\C^n,0)\to \C$ a germ of $C^{\infty}$ function. If
$\debar f\equiv 0$, $\pmod{q_1,\ldots,q_m}$ then there exists a germ of
holomorphic function $g\colon (\C^n,0)\to \C$ such that
$f\equiv g$, $\pmod{q_1,\ldots,q_m}$.\end{thm}


\chapter{Differential graded Lie algebras (DGLA)}
\label{CAP:DGLA}
\piede

The classical formalism  (Grothendieck-Mumford-Schlessinger) of 
infinitesimal deformation theory is described by the procedure 
(see e.g \cite{Artinbook}, \cite{Rim})
\[\hbox{Deformation problem}\quad\leadsto\quad\hbox{Deformation 
functor/groupoid}\]
The above picture is rather easy and suffices for many applications; 
unfortunately in this way we forget information which can be 
useful.  

It has been suggested by several people (Deligne, Drinfeld, Quillen, 
Kontsevich \cite{Konts}, Schlessinger-Stasheff \cite{SS1, SS2}, 
Goldman-Millson \cite{GoMil1,GoMil2} and many others) that 
a possible and useful 
way to preserve information is to consider a factorization  
\[\hbox{Deformation problem}\quad\leadsto\quad
DGLA\quad\leadsto\quad
\hbox{Deformation functor/groupoid}\]
where by $DGLA$ we mean a differential graded Lie algebra 
depending from the data of the deformation problem and
the construction 
\[DGLA\,\leadsto\,
\hbox{Deformation functor},\qquad L\leadsto \Def_{L},\] 
is a well defined, functorial
procedure explained in  this \chaptername.\\ 
More precisely, we introduce (as in \cite{K}) the \emph{deformation functor} associated to a 
differential graded Lie algebra and we prove in particular 
(Corollary~\ref{VIII.4.4})
that \emph{quasiisomorphic 
differential graded Lie algebras give isomorphic deformation 
functors}: this is done in the framework of Schlessinger's theory and 
extended deformation functors.\\
We refer to \cite{GoMil1} for a similar construction which associate to 
every DGLA a  deformation groupoid.\\

Some additional comments on this procedure will be  done in
Section~\ref{sec:histdgla}; for the moment 
we only point out that, for most deformation problems, the correct $DGLA$ is only 
defined up to quasiisomorphism and then the results of this 
\chaptername\  are the 
necessary background for the whole theory.\\

\smallskip

In this chapter $\K$ will be a fixed field of characteristic 0. We 
assume that the reader is familiar with basic concepts about Lie 
algebras and their representations \cite{Hum}, \cite{Ja}; 
unless otherwise stated we allow 
the  Lie algebras to be  infinite dimensional.\\

\bigskip

\section{Exponential and logarithm} 
\label{sec:expo}

For every associative $\K$-algebra $R$ we denote by $R_{L}$ the associated 
Lie algebra with bracket $[a,b]=ad(a)b=ab-ba$; the linear operator  
$ad(a)\in \End(R)$ is called the \emph{adjoint} of $a$, the 
morphism $ad\colon R_{L}\to \End(R)$ is a morphism of Lie algebras.
If $I\subset R$ is an ideal then $I$ is also a Lie ideal of $R_{L}$ 

\begin{exer}\label{eseVIII.1} 
Let $R$ be an associative $\K$-algebra, $a,b\in R$, prove:
\begin{enumerate}
\item  
\[ ad(a)^{n}b=\sum_{i=0}^{n}(-1)^{i}\binom{n}{i}a^{n-i}ba^{i}.\]

\item  If $a$ is nilpotent in $R$ then also $ad(a)$ is nilpotent in 
$\End(R)$ and 
\[ e^{ad(a)}b:=\sum_{n\ge 0}\frac{ad(a)^{n}}{n!}b=e^{a}be^{-a}.\]
\end{enumerate}\end{exer}

Let $V$ be a fixed $\K$-vector space and denote
\[ P(V)=\left\{\left. \sum_{n=0}^{\infty}v_{n}\,\right|\, v_{n}\in 
\tensor{}{n}V \right\}\simeq \prod_{n=0}^{\infty}\tensor{}{n}V.\]
With the natural notion of sum and Cauchy product $P(V)$ becomes an 
associative $\K$-algebra; the vector subspace
\[ \ide{m}(V)=\left\{ \left.\sum_{n=1}^{\infty}v_{n}\,\right|\, v_{n}\in 
\tensor{}{n}V \right\}\subset P(V)\]
is an ideal,  $\ide{m}(V)^{s}=\{ \sum_{n=s}^{\infty}v_{n}\}$ for 
every $s$ and $P(V)$ is complete for the $\ide{m}(V)$-adic 
topology:  this means that a series $\ds\sum_{i=0}^{\infty} x_{i}$ 
is convergent whenever  $x_{i}\in \ide{m}(V)^{i}$ for every $i$.\\ 
In particular, it is well defined the \emph{exponential}
\[ e\colon \ide{m}(V)\to  E(V):=1+\ide{m}(V)=\{ 
1+\sum_{n=1}^{\infty}v_{n}\mid v_{n}\in 
\tensor{}{n}V \}\subset P(V),\quad 
e^{x}=\sum_{n=0}^{\infty}\frac{x^{n}}{n!}\]
and  and the \emph{logarithm}
\[ \log\colon E(V)\to  \ide{m}(V),\qquad 
\log(1+x)=\sum_{n=1}^{\infty}(-1)^{n-1}\frac{x^{n}}{n}.\]
We note that $E(V)$
is a multiplicative subgroup of the set of invertible
elements of $P(V)$ (being $\sum_{n=0}^{\infty}x^{n}$ the 
inverse of $1-x$, $x\in \ide{m}(V)$).
It is well known that exponential and logarithm are  one 
the inverse of the other. Moreover  if $[x,y]=xy-yx=0$ then 
$e^{x+y}=e^{x}e^{y}$ and $\log((1+x)(1+y))=\log(1+x)+\log(1+y)$.
\\

Every linear morphism of $\K$-vector spaces $f_{1}\colon V\to W$ 
induces a natural, homogeneous and continuous homomorphism of $\K$-algebras
$f\colon P(V)\to P(W)$. It is clear that $f(\ide{m}(V))\subset 
\ide{m}(W)$, $f\colon E(V)\to E(W)$ is a group homomorphism  and 
$f$ commutes with the exponential and the logarithm.\\
Consider for instance the three homomorphisms 
\[ \Delta, p,q\colon P(V)\to P(V\oplus V)\]
induced respectively by the diagonal $\Delta_{1}(v)=(v,v)$, by 
$p_{1}(v)=(v,0)$ and by $q_{1}(v)=(0,v)$.\\
We define 
\[ \widehat{l}(V)=\{x\in P(V)\mid \Delta(x)=p(x)+q(x)\},\qquad
\widehat{L}(V)=\{x\in P(V)\mid \Delta(x)=p(x)q(x)\}.\]
It is immediate to observe that $V\subset \widehat{l}(V)\subset 
\ide{m}(V)$ and $1\in \widehat{L}(V)\subset E(V)$.

\begin{thm}
In the above notation we have:
\begin{enumerate} 
\item $\widehat{l}(V)$ is a Lie subalgebra of $P(V)_{L}$.
\item $\widehat{L}(V)$ is a multiplicative subgroup of $E(V)$.
\item Let $f_{1}\colon V\to W$ be a linear map and $f\colon P(V)\to 
P(W)$ the induced algebra homomorphism. 
Then $f(\widehat{l}(V))\subset \widehat{l}(W)$ and 
$f(\widehat{L}(V))\subset \widehat{L}(W)$.  
\item The exponential gives a bijection between $\widehat{l}(V)$ and 
$\widehat{L}(V)$.
\end{enumerate}
\end{thm}    
    
\begin{proof}    
We first note that for every $n\ge 0$ and every 
pair of vector spaces $U,W$ we have a canonical isomorphism 
\[ \tensor{}{n}(U\oplus W)=
\somdir{i=0}{n}(\tensor{}{i}U\oplus 
\tensor{}{n-i}W)\]    
and therefore 
\[P(U\oplus W)=\prod_{i,j=0}^{\infty}\tensor{}{i}U\oplus 
\tensor{}{j}W.\]
In particular for every $x\in P(U)\otimes\K\subset P(U\oplus W)$, 
$y\in \K\otimes P(W)\subset P(U\oplus W)$ we have 
$xy=yx$. In our case, i.e. when $U=W=V$ this implies that 
$p(x)q(y)=q(y)p(x)$ for every $x,y\in P(V)$.\\
Let $x,y\in \widehat{l}(V)$ then 
\[
\begin{split}
\Delta([x,y])&=\Delta(x)\Delta(y)-\Delta(y)\Delta(x)\\
&=(p(x)+q(x))(p(y)+q(y))-
(p(y)+q(y))(p(x)+q(x))\\
&=p([x,y])+q([x,y]).
\end{split}
\]
If $x,y\in \widehat{L}(V)$ then 
\[\Delta(yx^{-1})=\Delta(y)\Delta(x)^{-1}=p(y)q(x)q(x)^{-1}p(x)^{-1}=
p(yx^{-1})q(yx^{-1})\] 
and therefore $yx^{-1}\in \widehat{L}(V)$.\\
If $g\colon P(V\oplus V)\to P(W\oplus W)$ is the algebra homomorphism 
induced by $f_{1}\oplus f_{1}\colon V\oplus V\to W\oplus W$ it is 
clear that $\Delta f=g\Delta$, $pf= gp$ and $qf=gq$. This implies 
immediately item 3.\\ 
If $x\in \widehat{l}(V)$ then the equalities
\[ 
\Delta(e^{x})=e^{\Delta(x)}=e^{p(x)+q(x)}=e^{p(x)}e^{q(x)}=p(e^{x})q(e^{x})\]
prove that $e(\widehat{l}(V))\subset\widehat{L(V)}$. Similarly if $y\in 
\widehat{L}(V)$ then 
\[ 
\Delta(\log(y))=\log(\Delta(y))=
\log(p(y)q(y))=\log(p(y))+\log(q(y))=p(\log(y))+q(\log(y))\]
and therefore $\log(\widehat{L}(V))\subset \widehat{l}(V)$.
\end{proof}

\begin{cor} For every vector space $V$ the binary operation 
\[ *\colon \widehat{l}(V)\times \widehat{l}(V)\to \widehat{l}(V),\qquad
x*y=\log(e^{x}e^{y})\]
induces a group structure on the Lie algebra $\widehat{l}(V)$.\\
Moreover for every linear map $f_{1}\colon V\to W$ the induced 
morphism of Lie algebras $f\colon \widehat{l}(V)\to \widehat{l}(W)$ is also a 
homomorphism of groups.\end{cor}

\begin{proof}Clear.\end{proof}

In the next sections we will give an explicit formula for the 
product $*$ which involves only  the bracket of the Lie algebra 
$\widehat{l}(V)$.

\bigskip

\section{Free Lie algebras and the Baker-Campbell-Hausdorff formula}
\label{sec:freelie}

Let $V$ be a vector space over $\K$, we denote  by 
\[ T(V)=\somdir{n\ge 0}{}\tensor{}{n}V,\qquad
\bar{T(V)}=\somdir{n\ge 1}{}\tensor{}{n}V\subset
T(V).\]
The tensor product induce on $T(V)$ a structure of unital associative 
algebra, the natural embedding $T(V)\subset P(V)$ is a morphism of 
unitary algebras and 
$\bar{T(V)}$ is the ideal  $T(V)\cap\ide{m}(V)$.\\
The algebra $T(V)$ is called \emph{tensor algebra} generated by $V$ 
and $\bar{T(V)}$ is called the \emph{reduced tensor algebra} generated 
by $V$.\\

\begin{lem}
Let $V$ be a $\K$-vector space and $\imath\colon V=\tensor{}{1}V\to 
\bar{T(V)}$ the natural inclusion. For every associative $\K$-algebra 
$R$ and every linear map $f\colon V\to R$ there exists a unique 
homomorphism of $\K$-algebras $\phi\colon \bar{T(V)}\to R$ such that 
$f=\phi\imath$.\end{lem}

\begin{proof} Clear.\end{proof}

\begin{defn}\label{VIII.1.1} 
Let $V$ be a $\K$-vector space;
the \emph{free Lie algebra} generated by $V$ is the smallest Lie subalgebra 
$l(V)\subset \bar{T(V)}_{L}$ which contains $V$.\end{defn}

Equivalently 
$l(V)$ is the intersection of all the Lie subalgebras of $T(V)_{L}$ 
containing $V$.\\ 
For every integer $n>0$ we 
denote by $l(V)_{n}\subset\tensor{}{n}V$ the linear subspace 
generated by all the elements
\[ [v_{1},[v_{2},[\ldots,[v_{n-1},v_{n}]].]],\qquad n>0,\quad 
v_{1},\ldots,v_{n}\in V.\]
By definition $l(V)_{n}=[V,l(V)_{n-1}]$ and therefore 
$\oplus_{n>0}l(V)_{n}\subset l(V)$. On the other hand  
the Jacobi identity 
$[[x,y],z]=[x,[y,z]]-[y,[x,z]]$ 
implies that 
\[ [l(V)_{n},l(V)_{m}]\subset 
[V,[l(V)_{n-1},l(V)_{m}]]+[l(V)_{n-1},[V,l(V)_{m}]]\] 
and therefore, 
by induction on $n$, $[l(V)_{n},l(V)_{m}]\subset l(V)_{n+m}$.\\
As a consequence 
$\oplus_{n>0}l(V)_{n}$ is a Lie subalgebra of $l(V)$ and 
then $\oplus_{n>0}l(V)_{n}=l(V)$, $l(V)_{n}=l(V)\cap \tensor{}{n}V$.\\

Every morphism of vector spaces $V\to W$ induce a morphism of 
algebras $\bar{T(V)}\to \bar{T(W)}$ which restricts to a morphism of 
Lie algebras $l(V)\to l(W)$.

The name \emph{free Lie algebra} of $l(V)$ is motivated by the 
following universal property:\\  
\emph{Let $V$ be a vector space, $H$ a Lie algebra and $f\colon V\to H$ a 
linear map. Then there exists a unique homomorphism of Lie algebras 
$\phi\colon l(V)\to H$ which extends $f$.}\\
We will prove this property in 
Theorem~\ref{VIII.1.2}.

Let $H$ be a Lie algebra with bracket $[,]$ 
and $\sigma_{1}\colon V\to H$ a linear map.\\ 
Define recursively, for every $n\ge 2$,  the linear map 
\[ \sigma_{n}\colon \tensor{}{n}V\to H,\qquad
\sigma_{n}(v_{1}\otimes\ldots\otimes v_{n})=
[\sigma_{1}(v_{1}),\sigma_{n-1}(v_{2}\otimes\ldots\otimes v_{n})].\]
For example, if $V=H$ and $\sigma_{1}$ is the identity
then $\sigma_{n}(v_{1}\otimes\ldots\otimes v_{n})=
[v_{1},[v_{2},[\ldots,[v_{n-1},v_{n}]].]]$.

\begin{thm}[Dynkyn-Sprecht-Wever]\label{VIII.1.2} 
In the notation above, the linear map 
\[\sigma=\sum_{n=1}^{\infty}\frac{\sigma_{n}}{n}\colon l(V)\to H\]
is the unique homomorphism of Lie algebras extending $\sigma_{1}$.
\end{thm}

\begin{proof}
The adjoint representation $\theta\colon V\to \End(H)$, 
$\theta(v)x=[\sigma_{1}(v),x]$ extends to a unique morphism of associative 
algebras $\theta\colon \bar{T(V)}\to \End(H)$ by the composition rule 
\[\theta(v_{1}\otimes\ldots\otimes v_{s})x=
\theta(v_{1})\theta(v_{2})\ldots \theta(v_{s})x.\]
We note that, if $v_{1},\ldots,v_{n},w_{1},\ldots,w_{m}\in V$ then 
\[\sigma_{n+m}(v_{1}\otimes\ldots\otimes v_{n}
\otimes w_{1}\otimes\ldots\otimes w_{m})=
\theta(v_{1}\otimes\ldots\otimes v_{n})\sigma_{m}
(w_{1}\otimes\ldots\otimes w_{m}).\]
Since every element of $l(V)$ is a linear combination of 
homogeneous elements it is sufficient to prove, by induction on $n\ge 
1$, the following properties
\begin{description}
    \item[$A_{n}$] If $m\le n$, $x\in l(V)_{m}$ and $y\in l(V)_{n}$ 
    then $\sigma(xy-yx)=[\sigma(x),\sigma(y)]$.
    \item[$B_{n}$] If $m\le n$,  $y\in l(V)_{m}$ and $h\in H$ then 
    $\theta(y)h=[\sigma(y),h]$.
\end{description}
The initial step $n=1$ is straightforward, assume therefore $n\ge 2$.\\
\implica{A_{n-1}+B_{n-1}}{B_{n}} 
We can consider only the case $m=n$.
The element $y$ is a linear combination of elements of the form
$ab-ba$, $a\in V$, $b\in l(V)_{n-1}$ and, using $B_{n-1}$ we get 
\[ \theta(y)h=[\sigma(a),\theta(b)h]-\theta(b)[\sigma(a),h]=
[\sigma(a),[\sigma(b),h]]-[\sigma(b),[\sigma(a),h].\]
Using $A_{n-1}$ we get therefore 
\[ \theta(y)h=[[\sigma(a),\sigma(b)],h]=[\sigma(y),h].\]
\implica{B_{n}}{A_{n}}  
\[\begin{split} 
\sigma_{n+m}(xy-yx)&=\theta(x)\sigma_{n}(y)-\theta(y)\sigma_{m}(x)
=[\sigma(x),\sigma_{n}(y)]-[\sigma(y),\sigma_{m}(x)]\\
&=n[\sigma(x),\sigma(y)]-m[\sigma(y),\sigma(x)]=
(n+m)[\sigma(x),\sigma(y)].
\end{split}\]
Since $l(V)$ is generated by $V$ as a Lie algebra, the unicity of  
$\sigma$ follows.
\end{proof}

\begin{cor}\label{VIII.1.2bis}
For every vector space $V$     
the linear map 
\[\sigma\colon \bar{T(V)}\to 
l(V),\qquad 
\sigma(v_{1}\otimes\ldots\otimes v_{n})=
\frac{1}{n}[v_{1},[v_{2},[\ldots,[v_{n-1},v_{n}]].]]\]
is a projection.
\end{cor}

\begin{proof} 
The identity on $l(V)$ is the unique Lie homomorphism extending the 
natural inclusion $V\to l(V)$.
\end{proof} 
    
The linear map $\sigma$ defined in Corollary~\ref{VIII.1.2bis}
extends naturally to a projector 
$\sigma\colon P(V)\to P(V)$. We have the following theorem

\begin{thm}[Friedrichs]\label{VIII.1.3} 
In the notation above 
\[\widehat{l}(V)=\{x\in P(V)\mid \sigma(x)=x\}\quad\hbox{ and } 
\quad l(V)=T(V)\cap \widehat{l}(V).\]\end{thm}

\begin{proof} 
The two equalities are  equivalent, we will prove the 
second.    
We have already seen that $T(V)$ and $\widehat{l}(V)$ are Lie subalgebras 
of  $P(V)_{L}$ containing $V$ and then $l(V)\subset T(V)\cap 
\widehat{l}(V)$.\\
Define the linear map 
\[\delta\colon T(V)\to T(V\oplus V),\qquad 
\delta(x)=\Delta(x)-p(x)-q(x).\]
By definition $T(V)\cap \widehat{l}(V)=\ker\delta$ and we need to prove 
that if $\delta(x)=0$ for some homogeneous $x$ then $x\in l(V)$. 
For later computation we point out that, 
under the identification  $T(V\oplus V)=T(V)\otimes T(V)$, 
for every  monomial $\prod_{i}x_{i}$ with 
$x_{i}\in \ker\delta$ we have 
\[ \delta(\prod_{i}x_{i})=\prod_{i}(x_{i}\otimes 1+1\otimes 
x_{i})-(\prod_{i}x_{i})\otimes 1-1\otimes (\prod_{i}x_{i}).\]
In particular if $x\in \bar{T(V)}$ then $\delta(x)$ is the natural 
projection of $\Delta(x)$ onto the subspace $\somdir{i,j\ge 1}{}
\tensor{}{i}V\otimes \tensor{}{j}V$.\\ 

Let $\{y_{i}\mid i\in \sI\,\}$ be a fixed homogeneous 
basis of $l(V)$.
We can find a total 
ordering on the set $\sI$
such that if $y_{i}\in l(V)_{n}$, 
$y_{j}\in l(V)_{m}$ and $n<m$ then
$i<j$. For every index $h\in\sI$ we denote by $J_{h}\subset T(V)$ the ideal 
generated by $y_{h}^{2}$ and the $y_{i}$'s for every $i>h$, 
then $J_{h}$ is a 
homogeneous ideal and $y_{h}\not\in J_{h}$.\\
A \emph{standard monomial}  is a monomial 
of the form $y=y_{i_{1}}y_{i_{2}}\ldots y_{i_{h}}$ with 
$i_{1}\le \ldots\le i_{h}$.
The \emph{external degree}  of the above standard monomial $y$ is by 
definition the positive integer $h$.\\
Since $y_{i}y_{j}=y_{j}y_{i}+\sum_{h} a_{h}y_{h}$, $a_{h}\in \K$, the 
standard monomials 
generate $\bar{T(V)}$ as a vector space and the standard monomials of 
external degree 1 are a basis of $l(V)$. 
\begin{claim}\label{st-mo-claim}
For every $n>0$ the following hold:
\begin{enumerate}
    \item  The image under $\delta$ of the standard monomials of 
    external degree $h$ with $2\le h\le n$ are linearly independent. 

    \item  The standard monomials of external degree $\le n$ are 
    linearly independent.  
\end{enumerate}
\end{claim}
\begin{proof}[Proof of Claim] 
Since the standard monomials of external degree 1 are linearly 
independent and contained in the kernel of $\delta$ it is immediate 
to see the implication \implica{1}{2}.\\ 
We prove \vale{1} by induction on $n$, 
being the statement true  for $n=1$.\\  
Consider a nontrivial, finite linear combination $l.c.$ of standard monomials 
of external degree $\ge 2$ and $\le n$.
There exists an index $h\in \sI$ such that we can write 
$l.c. =z+\sum_{i=1}^{n}y_{h}^{i}w_{i}$, where $z$, $w_{i}$ are linear 
combination of standard monomials in $y_{j}$, $j>h$ and at least one 
of the $w_{i}$ is non trivial.   
If we consider the composition $\phi$ of $\delta\colon T(V)\to T(V\oplus 
V)=T(V)\otimes T(V)$ with the projection $T(V)\otimes T(V)\to 
T(V)/J_{h}\otimes T(V)$ we have 
\[ \phi(l.c.)=\sum_{i=1}^{n}iy_{h}\otimes y_{h}^{i-1}w_{i}=
y_{h}\otimes \sum_{i=1}^{n}iy_{h}^{i-1}w_{i}.\]
Since $\sum_{i=1}^{n}iy_{h}^{i-1}w_{i}$ is a nontrivial linear 
combination of standard monomials of external degrees $\le n-1$, by 
inductive assumption, it is different from 0 on $T(V)$.   
\end{proof}

From the claim follows immediately that the kernel of $\delta$ is 
generated by the standard monomials of degree $1$ and therefore  
$\ker\delta=l(V)$.
\end{proof}

\begin{exer}\label{eseVIII.3} 
Let $x_{1},\ldots,x_{n},y$ be linearly independent
vectors in a vector space $V$. Prove that the $n!$ vectors
\[ \sigma_{n+1}(x_{\tau(1)}\ldots x_{\tau(n)}y),\qquad 
\tau\in\Sigma_{n},\]
are linearly independent in the free Lie algebra $l(V)$.\\
(Hint: Let $W$ be a vector space with basis $e_{0},\ldots,e_{n}$ 
and consider the subalgebra $A\subset \End(W)$ generated by the 
endomorphisms $\phi_{1},\ldots,\phi_{n}$, $\phi_{i}(e_{j})=
\delta_{ij}e_{i-1}$. Take a suitable morphisms of Lie algebras
$l(V)\to A\oplus W$.)\end{exer}

Our main use of the projection $\sigma\colon P(V)\to\widehat{l}(V)$ consists in 
the proof of the  
an explicit description of the product $*\colon \widehat{l}(V)\times 
\widehat{l}(V)\to \widehat{l}(V)$. 

\begin{thm}[Baker-Campbell-Hausdorff formula]\label{VIII.1.5}
For every $a,b\in \widehat{l}(V)$ we have 
\renewcommand\arraystretch{0.6}
\[a*b=\sum_{n>0}
\frac{(-1)^{n-1}}{n}\!\!\!\!\!\!\!\!
\sum_{\begin{array}{c}\scriptstyle
p_{1}+q_{1}>0\\
\cdot\\
\scriptstyle p_{n}+q_{n}>0\end{array}
}\!\!\!\!\!
\frac{\ds\left(\sum_{i=1}^{n}(p_{i}+q_{i})\right)^{\!\!\!-1}}
{p_{1}!q_{1}!\ldots p_{n}!q_{n}!}
ad(a)^{p_{1}}ad(b)^{q_{1}}\ldots ad(b)^{q_{n}-1}b.\]
\renewcommand\arraystretch{1}
In particular 
$a*b-a-b$ belongs to the Lie ideal 
of $\widehat{l}(V)$ generated by $[a,b]$.
\end{thm}

\begin{proof}
Use the formula of the statement to define momentarily
a binary operator  $\bullet$ on $\widehat{l}(V)$; we
want to prove that $\bullet=*$.\\ 
Consider first the case
$a,b\in V$, in this situation  
\renewcommand\arraystretch{0.6}
\[a*b=\sigma\log(e^{a}e^{b})=\sigma\left(\sum_{n>0}
\frac{(-1)^{n-1}}{n}\left(\sum_{p+q>0}
\frac{a^{p}b^{q}}{p!q!}\right)^{\!\!n}\,\,\right)=\]
\[=\sigma\left(\sum_{n>0}
\frac{(-1)^{n-1}}{n}\sum_{\begin{array}{c}
\scriptstyle p_{1}+q_{1}>0\\
\cdot\\
\scriptstyle p_{n}+q_{n}>0\end{array}
}\frac{a^{p_{1}}b^{q_{1}}\ldots a^{p_{n}}b^{q_{n}}}
{p_{1}!q_{1}!\ldots p_{n}!q_{n}!}\right)\]
\[=\sum_{n>0}
\frac{(-1)^{n-1}}{n}\sum_{\begin{array}{c}\scriptstyle
p_{1}+q_{1}>0\\
\cdot\\
\scriptstyle p_{n}+q_{n}>0\end{array}
}\frac{1}{m}
\frac{\sigma_{m}(a^{p_{1}}b^{q_{1}}\ldots a^{p_{n}}b^{q_{n}})}
{p_{1}!q_{1}!\ldots p_{n}!q_{n}!},\qquad 
m:=\sum_{i=1}^{n}(p_{i}+q_{i}).\]
The elimination of the operators $\sigma_{m}$ gives 
\[a*b=\sum_{n>0}
\frac{(-1)^{n-1}}{n}\!\!\!\!\!\!\!\!
\sum_{\begin{array}{c}\scriptstyle
p_{1}+q_{1}>0\\
\cdot\\
\scriptstyle p_{n}+q_{n}>0\end{array}
}\!\!\!\!\!
\frac{\ds\left(\sum_{i=1}^{n}(p_{i}+q_{i})\right)^{\!\!\!-1}}
{p_{1}!q_{1}!\ldots p_{n}!q_{n}!}
ad(a)^{p_{1}}ad(b)^{q_{1}}\ldots ad(b)^{q_{n}-1}b.\]
\renewcommand\arraystretch{1}

Choose a vector space $H$ and a surjective linear map 
$H\to \widehat{l}(V)$, its composition with the inclusion 
$\widehat{l}(V)\subset \ide{m}(V)\subset P(V)$ 
extends to a continuous morphism of  associative algebras $q\colon P(H)\to
P(V)$; since $\widehat{l}(V)$ is a Lie subalgebra of $P(V)$ we have 
$q(l(H)_{n})\subset \widehat{l}(V)$ for every $n$ and then 
$q(\widehat{l}(H))\subset \widehat{l}(V)$.
Being $q\colon \widehat{l}(H)\to \widehat{l}(V)$ a morphism of Lie algebras, we have 
that $q$ commutes with $\bullet$.\\
On the other hand $q$ also commutes with exponential and logarithms and 
therefore $q$ commutes with the product $*$. 
Since $*=\bullet\colon H\times H\to \widehat{l}(H)$ the proof is done.  
\end{proof}

The first terms of the  Baker-Campbell-Hausdorff formula are:
\[ a*b=a+b+\frac{1}{2}[a,b]+\frac{1}{12}[a,[a,b]]
-\frac{1}{12}[b,[b,a]]+\ldots\]

\medskip

\section{Nilpotent Lie algebras}
\label{sec:nilplie}

We recall that every Lie 
algebra $L$ has a \emph{universal enveloping algebra} $U$ characterized by the 
properties \cite[17.2]{Hum}, \cite[Ch. V]{Ja}:\begin{enumerate}
\item $U$ is an associative algebra and there exists an injective 
morphism of Lie algebras $i\colon L\to U_{L}$.
\item For every associative algebra $R$ and every morphism $f\colon 
L\to R_{L}$ of Lie algebras there exists a unique morphism of 
associative algebras $g\colon U\to R$ such that $f=gi$.
\end{enumerate}

A concrete exhibition of the universal enveloping algebra is given by 
$U=\bar{T(L)}/I$, where $I$ is the ideal generated by all the elements 
$a\otimes b-b\otimes a-[a,b]$, $a,b\in L$.
The only non trivial condition to check is the injectivity of the 
natural map $L\to U$. This is usually proved using the well 
known Poincar\'e-Birkhoff-Witt's theorem \cite[Ch. V]{Ja}.\\

\begin{exer} Prove that, for every vector space $V$, 
$\bar{T(V)}$ is the universal enveloping 
algebra of $l(V)$.\end{exer}

\begin{defn}
The \emph{lower central series} 
of a Lie algebra $L$ is defined recursively by 
$L^{1}=L$, $L^{n+1}=[L,L^{n}]$.\\ 
A Lie algebra $L$ is called \emph{nilpotent} 
if $L^{n}=0$ for $n>>0$. 
\end{defn}

It is clear that if $L$ is a nilpotent Lie algebra then the adjoint 
operator $ad(a)=[a,-]\colon L\to L$ is nilpotent for every $a\in L$.
According to Engel's theorem \cite[3.2]{Hum} the converse is true  
if $L$ is finite dimensional.

\begin{ex} It is immediate from the construction that 
the lower central series of the free Lie algebra $l(V)\subset\bar{T(V)}$
is $l(V)^{n}=\somdir{i\ge n}{}l(V)_{i}=
l(V)\cap \somdir{i\ge n}{}\tensor{}{i}V$.
\end{ex}

If $V$ is  a nilpotent Lie algebra,  
then the Baker-Campbell-Hausdorff
formula  defines a product  
$V\times V\mapor{*}V$, 
\renewcommand\arraystretch{0.6}
\[a*b=\sum_{n>0}
\frac{(-1)^{n-1}}{n}\!\!\!\!\!\!\!\!
\sum_{\begin{array}{c}\scriptstyle
p_{1}+q_{1}>0\\
\cdot\\
\scriptstyle p_{n}+q_{n}>0\end{array}
}\!\!\!\!\!
\frac{\ds\left(\sum_{i=1}^{n}(p_{i}+q_{i})\right)^{\!\!\!-1}}
{p_{1}!q_{1}!\ldots p_{n}!q_{n}!}
ad(a)^{p_{1}}ad(b)^{q_{1}}\ldots ad(a)^{p_{n}}ad(b)^{q_{n}-1}b.\]
\renewcommand\arraystretch{1}

It is clear from the definition that the  product $*$ 
commutes with every morphism of nilpotent Lie algebra.
The identity on $V$ induce a morphism of Lie algebras 
$\pi\colon {l}(V)\to V$  such that 
$\pi(l(V)_{n})=0$ for $n>>0$;    
this implies that $\pi$ can be extended to a  
morphism of Lie algebras
$\pi\colon \widehat{l}(V)\to V$.

\begin{prop} The Baker-Campbell-Hausdorff product 
$*$ induces a group structure on every nilpotent Lie algebras 
$V$.\end{prop}

\begin{proof} The morphism 
of Lie algebras $\pi\colon \widehat{l}(V)\to V$ is surjective and commutes 
with $*$.
\end{proof}

It is customary to denote by $exp(V)$ the group $(V,*)$. Equivalently 
it is possible to define 
$exp(V)$ as the set  
of formal symbols $e^{v}$, $v\in V$, endowed with the group 
structure  $e^{v}e^{w}=e^{v*w}$.

\begin{ex} Assume that $V\subset M=M(n,n,\K)$ is the Lie subalgebra of 
strictly upper triangular matrices. Since the product of $n$ matrices 
of $V$ is always equal to 0,  the inclusion $V\to M$ 
extends to a morphism of associative algebras $\phi\colon P(V)\to M$
and the morphism 
\[\phi\colon exp(V)\to GL(n,\K),\qquad 
\phi(e^{A})=\sum_{i=0}^{\infty}\frac{A^{i}}{i!}\in GL(n,\K).\] 
is a homomorphism of groups.
\end{ex}

The above example can be generalized in the following way

\begin{ex}\label{VIII.1.6} 
Let $R$ be an associative unitary  
$\K$-algebra, $R^{*}\subset R$ the multiplicative group of 
invertible elements  
and $N\subset R$ a nilpotent subalgebra  (i.e. $N^{n}=0$ 
for $n>>0$).\\
Let $V$ be a nilpotent Lie algebra and $\xi\colon V\to N
\subset R$ a representation. This means that $\xi\colon V\to N_{L}$  
is a morphism of Lie algebras.\\
Denoting by $\imath\colon V\mapin{}U$ the universal enveloping algebra, we
have   a commutative diagram
\renewcommand\arraystretch{1.3}
\[\begin{array}{ccccc}
{l(V)}&\mapor{\pi}&V&\mapor{\xi}&N_{L}\\
\mapver{}&&\mapver{}&&\mapver{}\\
\bar{T(V)}&\mapor{\eta}&U&\mapor{\psi}&R
\end{array}\]
where $\pi$, $\xi$ are morphisms of Lie algebras and $\eta,\psi$ 
homomorphisms of algebras.
Since the image of the composition $\phi=\psi\eta$ is contained in the 
nilpotent subalgebra $N$   
the above diagram  extends to 
\[\begin{array}{ccc}
\widehat{l(V)}&\mapor{}&P(V)\\
\mapver{\pi}&&\mapver{\phi}\\
V&\mapor{\xi}&R\end{array}\]
with $\phi$ homomorphism of associative algebras. If $f\in N$ it 
makes sense its exponential $e^{f}\in R$. 
For every $v\in V$ we have
$e^{\xi(v)}=\phi(e^{v})$  and for every $x,y\in V$ 
\[ e^{\xi(x)}e^{\xi(y)}=\phi(e^{x})\phi(e^{y})
=\phi(e^{x}e^{y})=\phi(e^{x*y})=e^{\xi(x*y)}.\]
The same assertion can be stated by saying that the exponential map
$e^{\xi}\colon (V,*)=exp(V)\to R^{*}$ is a homomorphism of groups. 
\renewcommand\arraystretch{1}
\end{ex}

\bigskip

\section{Differential graded Lie algebras}
\label{sec:dgla}

\begin{defn}\label{VIII.2.1} A \emph{differential graded Lie algebra} (DGLA ) 
$(L,[,],d)$ is the data of a $\Z$-graded vector space 
$L=\oplus_{i\in \Z}L^i$ together a bilinear bracket
$[,]\colon L\times L\to L$ and a linear map $d\in\Hom^{1}(L,L)$ satisfying 
the following condition:\begin{enumerate}
\item $[~,~]$ is homogeneous skewsymmetric: 
this means $[L^i,L^j]\subset L^{i+j}$ and 
$[a,b]+(-1)^{\bar{a}\bar{b}}[b,a]=0$ for every $a,b$ homogeneous. 

\item Every triple of homogeneous elements 
$a,b,c$  satisfies the (graded) Jacobi identity 
\[ [a,[b,c]]=[[a,b],c]+(-1)^{\bar{a}\bar{b}}[b,[a,c]].\]

\item $d(L^i)\subset L^{i+1}$, $d\circ d=0$ and 
$d[a,b]=[da,b]+(-1)^{\bar{a}}[a,db]$. The map $d$ is called the 
differential of $L$.
\end{enumerate}\end{defn}

\begin{exer}\label{eseVIII.5}
Let $L=\oplus L^i$ be a DGLA and $a\in L^i$:\begin{enumerate}
\item If $i$ is even then $[a,a]=0$.
\item If $i$ is odd then $[a,[a,b]]=\ds\frac{1}{2}[[a,a],b]$ 
for every $b\in L$ and $[[a,a],a]=0$.
\end{enumerate}\end{exer}

\begin{ex}\label{VIII.2.5}
If $L=\oplus L^i$ is a DGLA
then $L^0$ is a Lie algebra in the 
usual sense. Conversely, every Lie algebra can be considered as a 
DGLA concentrated in degree 0.
\end{ex}

\begin{ex}\label{VIII.2.6}
Let $(A,d_{A})$, $A=\oplus A_{i}$, be a dg-algebra over $\K$ 
and $(L,d_{L})$, $L=\oplus L^i$, a DGLA.\\ 
Then $L\otimes_{\K}A$ has a 
natural structure of DGLA by setting:
\[(L\otimes_{\K}A)^n=\oplus_i (L^{i}\otimes_{\K}A_{n-i}),\]
\[d(x\otimes a)=d_{L}x\otimes a+(-1)^{\bar{x}}x\otimes d_{A}a,\quad 
[x\otimes a,y\otimes b]=(-1)^{\bar{a}\,\bar{y}}[x,y]\otimes ab.\]
\end{ex}

\begin{ex}\label{VIII.2.7}
Let $E$ be a holomorphic vector bundle on a complex manifold
$M$. 
We may define  a DGLA $L=\oplus L^{p}$, 
$L^{p}=\Gamma(M, \sA^{0,p}(\END(E)))$ 
with the Dolbeault 
differential and the natural bracket. More precisely if $e,g$ are  local 
holomorphic sections of $\END(E)$ and $\phi,\psi$  differential forms we 
define $d(\phi e)=(\bar{\de}\phi)e$, 
$[\phi e,\psi g]=\phi\wedge\psi[e,g]$.\end{ex}

\begin{ex}\label{VIII.2.8} 
Let $(\sF^{*},d)$ be a sheaf of dg-algebras on a topological space; 
the space $\Der^{*}(\sF^{*},\sF^{*})$ is a DGLA with 
bracket $[f,g]=fg-(-1)^{\bar{f}\,\bar{g}}gf$ and
differential $\delta(f)=[d,f]$.
\end{ex}

\begin{defn} 
We shall say that a DGLA $L$ is \emph{$ad_{0}$-nilpotent} if for every $i$ 
the image of the adjoint action $ad:L^{0}\to \End(L^{i})$ is 
contained in a nilpotent (associative) subalgebra.
\end{defn}

\begin{exer}~\\
1) Every nilpotent DGLA (i.e. a DGLA whose descending 
central series is definitively trivial) is $ad_{0}$-nilpotent.\\ 
2) If $L$ is $ad_{0}$-nilpotent then $L^{0}$ is a nilpotent Lie 
algebra.\\
3) The converses of 1) and 2) are generally false.
\end{exer}

\begin{defn}\label{VIII.2.2}
A linear map $f\colon L\to L$ is called a \emph{derivation of 
degree $n$} if $f(L^i)\subset L^{i+n}$ and satisfies the graded 
Leibnitz rule $f([a,b])=[f(a),b]+(-1)^{n\bar{a}}[a,f(b)]$.
\end{defn}

We note that the Jacobi identity is equivalent to the assertion that, 
if $a\in L^i$ then $ad(a)\colon L\to L$, $ad(a)(b)=[a,b]$, 
is a derivation of degree $i$. The differential 
$d$ is a derivation of degree 1.

By following the standard notation we denote by 
$Z^i(L)=\ker(d\colon L^i\to L^{i+1})$, 
$B^i(L)=\Image(d\colon 
L^{i-1}\to L^{i})$, $H^i(L)=Z^i(L)/B^i(L)$.

\begin{defn}\label{VIII.2.3}
The \emph{Maurer-Cartan equation} 
(also called the deformation equation) of a DGLA $L$ is 
\[da+\frac{1}{2}[a,a]=0,\qquad a\in L^1.\]
The solutions  $MC(L)\subset L^{1}$ of the Maurer-Cartan 
equation are called the Maurer-Cartan elements of the DGLA $L$.
\end{defn}

There is  an obvious notion of morphisms of DGLAs;
we denote by $\mathbf{DGLA}$ the category of differential graded Lie 
algebras.\\
Every morphism of 
DGLAs induces a morphism between cohomology groups. 
It is  moreover clear  that 
morphisms of DGLAs preserve the solutions 
of the Maurer-Cartan equation.

A \emph{quasiisomorphism} of DGLAs
is a morphism inducing isomorphisms in cohomology. 
Two DGLA's are \emph{quasiisomorphic} 
if they are equivalent under the equivalence relation 
generated by quasiisomorphisms.

The cohomology of a DGLA is itself a differential graded Lie algebra with 
the induced bracket and zero differential:  

\begin{defn}\label{VIII.2.4}
A DGLA $L$ is called \emph{Formal} if it is quasiisomorphic 
to its cohomology DGLA $H^*(L)$.
\end{defn}

\begin{exer}\label{eseVIII.4} 
Let $D\colon L\to L$ be a derivation, then the kernel of 
$D$ is a graded Lie subalgebra.\end{exer}

\begin{ex}\label{VIII.2.9} 
Let $(L,d)$ be a DGLA and denote $\Der^{i}(L,L)$ the space of 
derivations $f\colon L\to L$ of degree $i$. 
The space $\Der^{*}(L,L)=\oplus_{i}\Der^{i}(L,L)$ 
is a DGLA with 
bracket $[f,g]=fg-(-1)^{\bar{f}\,\bar{g}}gf$ and
differential $\delta(f)=[d,f]$.
\end{ex}

For a better understanding of some of next topics it is useful to 
consider the following functorial construction. Given  
a DGLA $(L, [,],d)$ we can construct a new DGLA $(L', [,]', d')$ by 
setting 
$(L')^i=L^i$ for every $i\not=1$, $(L')^1=L^1\oplus\K d$ 
(here $d$ is considered as a formal symbol of degree 1)
with the bracket and the differential
\[[a+vd,b+wd]'=[a,b]+vd(b)+(-1)^{\bar{a}}wd(a),\qquad 
d'(a+vd)=[d,a+vd]'=d(a).\]

The natural inclusion $L\subset L'$ is a morphism of DGLA; for a better 
understanding of the 
Maurer-Cartan equation it is convenient to consider the affine 
embedding $\phi\colon L^{1}\to (L')^{1}$, $\phi(a)=a+d$. 
For an element $a\in L^1$ we have 
\[ d(a)+\frac{1}{2}[a,a]=0\quad\iff\quad [\phi(a),\phi(a)]'=0.\]

Let's now introduce the notion of gauge action on the Maurer-Cartan 
elements of an $ad_{0}$-nilpotent  DGLA. 
Note that  $[L^{0},L^{1}\oplus\K d]\subset L^{1}$; in particular 
if $L$ is $ad_{0}$-nilpotent  then also $L'$ 
is $ad_{0}$-nilpotent.

Given an $ad_{0}$-nilpotent DGLA $N$,  
the exponential of the adjoint action 
gives  homomorphisms of groups 
\[ exp(N^{0})=(N^{0},*)\to GL(N^{i}),\quad e^{a}\mapsto e^{ad(a)},\qquad 
i\in\Z\]
where $*$ is the product given by the 
Baker-Campbell-Hausdorff formula.\\
These homomorphisms induce  actions of the group 
$exp(N^{0})$ onto the vector spaces $N^{i}$ given by 
\[ e^{a}b=e^{ad(a)}b=
\sum_{n\ge 0}\frac{1}{n!}ad(a)^{n}(b).\]

\begin{lem}\label{VIII.2.11} In the above notation, if 
$W$ is a linear subspace of $N^{i}$ and 
$[N^{0},N^{i}]\subset W$ then the 
exponential adjoint action preserves the affine 
subspaces $v+W$, $v\in N_{i}$.\end{lem}

\begin{proof} Let $a\in N^{0}$, $v\in N^{i}$, $w\in  W$, then 
\[ e^{a}(v+w)=v+
\sum_{n\ge 1}\frac{1}{n!}ad(a)^{n-1}([a,v])+
\sum_{n\ge 0}\frac{1}{n!}ad(a)^{n}(w).
\]
\end{proof}

\begin{lem}\label{VIII.2.12} In the above notation  the exponential 
adjoint action preserves the quadratic cone 
$Z=\{v\in N^1\mid [v,v]=0\}$.\\
For every $v\in Z$ and $u\in N^{-1}$ the element $exp([u,v])$ belongs 
to the stabilizer of $v$.
\end{lem}

\begin{proof}
By Jacobi identity $2[v,[a,v]]=-2[v,[v,a]]=[a,[v,v]]$
for every $a\in N^{0}$, $v\in N^{1}$.\\
Let $a\in N^{0}$ be a fixed element, for every $u\in N^{1}$ 
define the polynomial function 
$F_u\colon\K\to N^{2}$ by 
\[F_u(t)=e^{-ad(ta)}[e^{ad(ta)}u,e^{ad(ta)}u].\]
For every $s,t\in\K$, if $v=e^{ad(sa)}u$ then 
\[F_u(t+s)=e^{ad(-sa)}F_v(t),\quad 
\desude{F_v}{t}(0)=-[a,[v,v]]+2[v,[a,v]]=0\]
\[\desude{F_u}{t}(s)=e^{ad(-sa)}u\desude{F_v}{t}(0)=0.\] 
Since the field $\K$ has characteristic 0 every function $F_v$ is 
constant, proving the invariance of $Z$.\\
If $u\in N^{-1}$  and $v\in Z$, then by Jacobi identity
$[[u,v],v]=ad([u,v])v=0$ and then $exp([u,v])v=v$.
\end{proof}

If $L$ is an $ad_{0}$-nilpotent DGLA then \ref{VIII.2.11} and \ref{VIII.2.12} 
can be applied to $N=L'$.
Via the affine embedding $\phi\colon L^{1}\to (L')^{1}$, 
the exponential of the adjoint action on $L'$ induces the so called
\emph{Gauge action} of $exp(L^0)$ over 
the set of solution of the Maurer-Cartan equation, given 
explicitly by 
\renewcommand\arraystretch{2.5}
\[\begin{array}{cl} 
exp(a)(w)\!\!\!&=\ds\phi^{-1}\left(e^{ad(a)}\phi(w)\right)=
\sum_{n\ge 0}\frac{1}{n!}ad(a)^{n}(w)-
\sum_{n\ge 1}\frac{1}{n!}ad(a)^{n-1}(da)\\
&=\ds w+\sum_{n\ge 0}\frac{ad(a)^{n}}{(n+1)!}([a,w]-da).
\end{array}\]
\renewcommand\arraystretch{1}

\begin{rem}\label{VIII.2.13} 
If $w$ is a solution of the Maurer-Cartan equation and $u\in 
L^{-1}$ then $[w,u]+du=[w+d,u]\in L^{0}$ belongs to the 
stabilizer of $w$ under the gauge action.\\
For every $a\in L^{0}$, $w\in L^{1}$, the polynomial
$\gamma(t)=exp(ta)(w)\in L^{1}\otimes\K[t]$ is the solution of the ``Cauchy 
problem''
\renewcommand\arraystretch{1.4}
\[ \left\{\begin{array}{l}
\ds\frac{d\gamma(t)}{dt}=[a,\gamma(t)]-da\\
\gamma(0)=w\end{array}\right.\]
\renewcommand\arraystretch{1}
\end{rem}

\bigskip

\section{Functors of Artin rings}
\label{sec:functors}

\subsection{Basic definitions}

We denote by: 

\begin{itemize}
    \item  $\Set$ the category of sets in a fixed universe; we 
    also make the choice of a fixed set $\{0\}\in \Set$ 
    of cardinality 1. 
    \item  $\mathbf{Grp}$  the category of groups.
    
    \item  $\Art_{\K}$ the category of local Artinian
$\K$-algebras with residue field $\K$ (with as morphisms the local
homomorphisms). If $A\in \Art_{\K}$, we will denote by $\ide{m}_A$ its maximal
ideal.
\end{itemize}

A \emph{ small extension} $e$
in $\Art_{\K}$  is an exact sequence of abelian groups
\[e:\quad 0\mapor{}M\mapor{i} B\mapor{p}A\mapor{}0\]
such that  $B\mapor{p} A$ is a  morphism
in $\Art_{\K}$ and 
$\ker p=i(M)$ is annihilated by the
maximal ideal of $B$ (that is, as a $B$-module it is a
$\K$-vector space).\\ 

Given a surjective morphism  $B\to A$  in $\Art_{\K}$ with kernel 
$J$,
there exists a sequence of small extensions  
\[ 0\mapor{}\quot{\ide{m}_{B}^{n}J}{\ide{m}_{B}^{n+1}J}\mapor{} 
\quot{B}{\ide{m}_{B}^{n+1}J}\mapor{}
\quot{B}{\ide{m}_{B}^{n}J}\mapor{}0,\qquad n\ge 0.\]
Since,
by Nakayama's lemma, there exists $n_{0}\in\N$ such that 
$\ide{m}_{B}^{n}J=0$  for every $n\ge n_{0}$ we get that every surjective 
morphism is $\Art_{\K}$ is the composition of a finite number of small 
extensions.\\

\begin{defn}\label{VIII.3.1} 
A \emph{Functor of Artin rings} is a covariant functor 
$F\colon \Art_{\K}\to \Set$
such that $F(\K)\simeq \{0\}$.\end{defn}

\begin{ex}\label{VIII.3.2} 
If $V$ is a $\K$-vector space we may interpret $V$ as 
a functor of Artin rings 
$V\colon \Art_{\K}\to \Set$, $V(A)=V\otimes_{\K}\ide{m}_{A}$. 
If $V=0$ we get the \emph{trivial functor} $0\colon \Art_{\K}\to 
\Set$.\end{ex}   

The functors of Artin rings are the object of a new category 
whose morphisms are the natural transformation of functors.
A natural transformation $\eta\colon F\to G$ is an 
isomorphism if and only if $\eta(A)\colon F(A)\to G(A)$ is bijective
for every $A\in\Art_{\K}$.

\begin{defn}\label{VIII.3.3} 
Let $ F,G\colon \Art_{\K}\to \Set$ be two functors of 
Artin rings and $\eta\colon F\to G$ a natural transformation; 
$\eta$ is called \emph{smooth} if
for every small extension
\[0\mapor{}M\mapor{} B\mapor{p}A\mapor{}0\]
the map 
\[ (\eta,p)\colon F(B)\to G(B)\times_{G(A)}F(A)\]
is surjective.\\
A functor of Artin rings $F$ is called \emph{smooth} if the morphism
$F\to 0$ is smooth.\\
\end{defn}

\begin{exer}\label{eseVIII.6} 
$ F\colon \Art_{\K}\to \Set$ is smooth if and only if 
for every surjective morphism $B\to A$ is $\Art_{\K}$, 
the map $F(B)\to F(A)$ is also surjective.\\
If $V$ is a vector space then $V$ is smooth as a functor of Artin 
rings (cf. Example~\ref{VIII.3.2}).
\end{exer}

\begin{exer}\label{eseVIII.7} 
Let $R$ be an analytic 
algebra and let $h_{R}\colon \Art_{\C}\to\Set$ be the functor of 
Artin rings defined by $h_{R}(A)=\Mor_{\mathbf{An}}(R,A)$.\\
Prove that $h_{R}$ is smooth if and only if $R$ is smooth.
\end{exer}

\begin{ex}\label{VIII.3.4} 
Let $M_{0}$ be a compact complex manifold and define 
for every $A\in\Art_{\C}$ 
\[ \Def_{M_{0}}(A)=\Def_{M_{0}}(\Oh_{X,0})=\Def_{M_{0}}(X,0)\]
where $(X,0)=\Spec(A)$ is a fat point such that $\Oh_{X,0}=A$; 
since it is always possible to write $A$ as a quotient of 
$\C\{z_{1},\ldots,z_{n}\}$ for some $n\ge 0$, such a fat point 
$(X,0)$ always exists. 
According to \ref{V.2.4} the isomorphism class 
of $(X,0)$ depends only on $A$.\\
Every morphism $\Oh_{X,0}\to \Oh_{Y,0}$ in $\Art_{\C}$ is induced by a 
unique morphism $(Y,0)\to (X,0)$. The pull-back of infinitesimal 
deformations gives a morphism 
$\Def_{M_{0}}(X,0)\to \Def_{M_{0}}(Y,0)$. Therefore 
$\Def_{M_{0}}\colon \Art_{\C}\to \Set$ is a functor of Artin rings.
\end{ex}
\medskip

\begin{defn}\label{VIII.3.5} 
The \emph{tangent space} to a functor of Artin rings 
$F\colon\Art_{\K}\to\Set$ is by definition 
\[ t_{F}=F\left(\frac{\K[t]}{(t^{2})}\right)=F(\K\oplus\K\epsilon),\qquad 
\epsilon^{2}=0.\]
\end{defn}

\begin{exer}\label{eseVIII.8}
Prove that, for every analytic algebra $R$ there exists
a natural isomorphism $t_{h_{R}}=\Der_{\C}(R,\C)$ (see 
Exercise~\ref{eseVIII.7}).\end{exer}

\medskip

\subsection{Automorphisms functor}

In this section  every tensor 
product is intended over $\K$, i.e $\otimes=\otimes_{\K}$.
Let $S\mapor{\alpha} R$ be a morphism of graded $\K$-algebras,
for every $A\in \Art_{\K}$ we have natural 
morphisms $S\otimes A\mapor{\alpha}R\otimes A$ and 
$R\otimes_{\K}A\mapor{p}R$, $p(x\otimes a)=x\bar{a}$, where $\bar{a}\in\K$ is 
the class of $a$ in the residue field of $A$.\\

\begin{lem}\label{VIII.3.6} Given $A\in \Art_{\K}$ and a commutative diagram of 
morphisms of graded $\K$-algebras
\[\xymatrix{S\otimes A\ar[r]^{\alpha}\ar[d]^{\alpha}&R\otimes 
A\ar[d]^{p}\\
R\otimes A\ar[ur]^{f}\ar[r]^{p}&R}\]
we have that $f$ is an isomorphism and 
$f(R\otimes J)\subset R\otimes J$ for every ideal $J\subset 
A$.\end{lem}  

\begin{proof} $f$ is a morphism of graded $A$-algebras, 
in particular for every 
ideal $J\subset A$, $f(R\otimes J)\subset Jf(R\otimes A)\subset R\otimes J$.
In particular, if $B=A/J$, then $f$ induces a morphism of graded 
$B$-algebras $\bar{f}\colon R\otimes B\to R\otimes B$.\\
We claim that if $\ide{m}_{A}J=0$ then $f$ is the identity on $R\otimes J$; 
in fact for every $x\in R$, $f(x\otimes 1)-x\otimes 1\in \ker{p}=R\otimes 
\ide{m}_{A}$
and then if $j\in J$, $x\in R$. 
\[ f(x\otimes j)=jf(x\otimes 1)=x\otimes j+j(f(x\otimes 1)-x\otimes 1)=
x\otimes j.\]
Now we prove the lemma by induction on $n=\dim_{\K}A$, being $f$ the 
identity for $n=1$. Let 
\[0\mapor{}J\mapor{} A\mapor{}B\mapor{}0\]
be a small extension with $J\not=0$. Then we have a commutative 
diagram with exact rows
\[\begin{array}{ccccccccc}
0&\mapor{}&R\otimes J&\mapor{}&R\otimes A&\mapor{}&R\otimes 
B&\mapor{}&0\\
&&\mapver{Id}&&\mapver{f}&&\mapver{\bar{f}}&&\\
0&\mapor{}&R\otimes J&\mapor{}&R\otimes A&\mapor{}&R\otimes 
B&\mapor{}&0\end{array}\]
By induction $\bar{f}$ is an isomorphism and by snake lemma also $f$ 
is an isomorphism.\end{proof}

\begin{defn}\label{VIII.3.7} 
For every $A\in\Art_{\K}$ let $\Aut_{R/S}(A)$ be the set 
of commutative diagrams of graded $\K$-algebra morphisms 
\[\xymatrix{S\otimes A\ar[r]\ar[d]&R\otimes A\ar[d]\\
R\otimes A\ar[ur]^{f}\ar[r]&R}\]
\end{defn}

According to Lemma~\ref{VIII.3.6} 
$\Aut_{R/S}$ is a functor from the category $\Art_{\K}$ to the 
category of groups $\mathbf{Grp}$. Here we consider $\Aut_{R/S}$ as a 
functor of Artin rings (just forgetting the group structure).

Let $\Der_{S}^{0}(R,R)$ be the space of $S$-derivations  $R\to R$ 
of degree 0. 
If $A\in\Art_{\K}$ and $J\subset \ide{m}_{A}$ is an ideal then, 
since $\dim_{\K}J<\infty$ there exist natural isomorphisms 
\[ \Der_{S}^{0}(R,R)\otimes J=\Der_{S}^{0}(R,R\otimes J)=
\Der_{S\otimes A}^{0}(R\otimes A,R\otimes J),
\]
where $d=\sum_{i} d_{i}\otimes j_{i}\in \Der_{S}^{0}(R,R)\otimes J$ 
corresponds to the $S\otimes A$-derivation 
\[d\colon R\otimes A\to R\otimes J\subset R\otimes A,\qquad
d(x\otimes a)=\sum_{i} d_{i}(x)\otimes j_{i}a.\]
For every $d\in \Der_{S\otimes A}^{0}(R\otimes A,R\otimes A)$ denote 
$d^{n}=d\circ\ldots\circ d$ the iterated composition of $d$ with 
itself $n$ times. The generalized Leibnitz rule gives 
\[ d^{n}(uv)=\sum_{i=0}^{n}\binom{n}{i}d^{i}(u)d^{n-1}(v), \qquad 
u,v\in R\otimes A.\]
Note in particular that if $d\in \Der_{S}^{0}(R,R)\otimes \ide{m}_{A}$ then 
$d$ is a nilpotent endomorphism of $R\otimes A$ and 
\[ e^{d}=\sum_{n\ge 0 } \frac{d^{n}}{n!}\]
is a morphism of $\K$-algebras belonging to $\Aut_{R/S}(A)$.

\begin{prop}\label{VIII.3.8} For every $A\in \Art_{\K}$ the exponential
\[exp\colon \Der_{S}^{0}(R,R)\otimes \ide{m}_{A}\to \Aut_{R/S}(A)\]
is a bijection.\\
\end{prop}

\begin{proof} This is obvious if $A=\K$; by induction on the dimension 
of $A$ we may assume that there exists a nontrivial small extension
\[0\mapor{}J\mapor{} A\mapor{}B\mapor{}0\]
 such that $exp\colon \Der_{S}^{0}(R,R)\otimes \ide{m}_{B}\to \Aut_{R/S}(B)$
is bijective.\\
We first note that if $d\in \Der_{S}^{0}(R,R)\otimes \ide{m}_{A}$, $h\in 
\Der_{S}^{0}(R,R)\otimes J$ then $d^{i}h^{j}=h^{j}d^{i}=0$ whenever $j>0$, 
$j+i\ge 2$ and then $e^{d+h}=e^{d}+h$; this easily implies that  
$exp$ is injective.\\ 
Conversely take a $f\in 
\Aut_{R/S}(A)$; by the inductive assumption there exists $d\in 
\Der_{S}^{0}(R,R)\otimes \ide{m}_{A}$ such that 
$\bar{f}=\bar{e^{d}}\in  
\Aut_{R/S}(B)$; denote $h=f-e^{d}\colon R\otimes A\to R\otimes J$.
Since 
$h(ab)=f(a)f(b)-e^{d}(a)e^{d}(b)=h(a)f(b)+e^{d}(a)h(b)=h(a)\bar{b}+\bar{a}h(b)$ we have 
that $h\in \Der_{S}^{0}(R,R)\otimes J$ and then 
$f=e^{d+h}$.\end{proof}

The same argument works also if $S\to R$ is a morphism of sheaves of 
graded $\K$-algebras over a topological space and $\Der_{S}^{0}(R,R)$, 
$\Aut_{R/S}(A)$ are respectively the vector space of $S$-derivations 
of degree 0 of $R$  and 
the  $S\otimes A$-algebra automorphisms of $R\otimes A$
lifting the identity on $R$.\\

\begin{ex}\label{VIII.3.9}
Let $M$ be a complex manifold, $R=\sA^{0,*}_{M}$,
$S=\bar{\Omega}_{M}^{*}$. According to Proposition~\ref{VI.1.2}
$\Der_{S}^{0}(R,R)=\Gamma(M,\sA^{0,0}(T_{M}))$ and then the 
exponential  gives isomorphisms  
\[exp\colon \Gamma(M,\sA^{0,0}(T_{M}))\otimes \ide{m}_{A}\to \Aut_{R/S}(A).\]
Since $exp$ is clearly functorial in $A$, interpreting the vector 
space $\Gamma(M,\sA^{0,0}(T_{M}))$ as a functor ( 
Example~\ref{VIII.3.2}), we have an isomorphism of functors 
$exp\colon \Gamma(M,\sA^{0,0}(T_{M}))\to \Aut_{R/S}$.
\end{ex}

\medskip

\subsection{The exponential functor} 

Let $L$ be a Lie algebra over  $\K$, $V$ 
a $\K$-vector space and $\xi\colon L\to \End(V)$ a representation of 
$L$.\\
For every $A\in\Art_{\K}$ the morphism $\xi$ can be extended
naturally to a morphism of Lie algebras 
$\xi\colon L\otimes A\to \End_{A}(V\otimes A)$.
Taking the exponential we get a functorial map 
\[exp(\xi)\colon L\otimes \ide{m}_{A}\to GL_{A}(V\otimes A)
,\qquad
exp(\xi)(x)=e^{\xi(x)}=\sum_{i=0}^{\infty}\frac{\xi^{n}}{n!}x,\]
where $GL_{A}$ denotes the group of $A$-linear invertible morphisms.\\ 
Note that $exp(\xi)(-x)=(exp(\xi)(x))^{-1}$. 
If $\xi$ is injective then also $exp(\xi)$ is injective (easy exercise).

\begin{thm}\label{VIII.3.10} In the notation above the image of 
$exp(\xi)$ is a subgroup. More precisely for every $a,b\in 
L\otimes \ide{m}_{A}$ there exists $c\in L\otimes \ide{m}_{A}$ such that 
$e^{\xi(a)}e^{\xi(b)}=e^{\xi(c)}$ and 
$a+b-c$ belong to the Lie ideal of $L\otimes \ide{m}_{A}$ generated by 
$[a,b]$.\end{thm}

\begin{proof} This is an immediate consequence of the  
Campbell-Baker-Hausdorff formula.\end{proof}

In the above notation 
denote $P=\End(V)$ and let $ad(\xi)\colon L\to \End(P)$ be the
adjoint  representation of $\xi$, 
\[ad(\xi)(x)f=[\xi(x),f]=\xi(x)f-f\xi(x).\]
Then for every $a\in L\otimes \ide{m}_{A}$, 
$f\in \End_{A}(V\otimes A)=P\otimes A$ 
we have (cf. Exercise~\ref{eseVIII.1}, \cite[2.3]{Hum})
\[ e^{ad(\xi)(a)}f=e^{\xi(a)}fe^{-\xi(a)}.\]

\bigskip

\section{Deformation functors associated to a DGLA}
\label{sec:defofun}

Let $L=\oplus L^i$ be a DGLA over  $\K$, we 
can define the following three functors:\begin{enumerate}
\item The Gauge functor 
$G_{L}\colon \Art_{\K}\to \mathbf{Grp}$, defined by 
$G_{L}(A)=exp(L^0\otimes \ide{m}_A)$. It is immediate to see that 
$G_{L}$ is smooth.

\item The Maurer-Cartan functor 
$MC_{L}\colon \Art_{\K}\to \Set$ defined by 
\[MC_L(A)=MC(L\otimes \ide{m}_{A})=\left\{x\in L^1\otimes \ide{m}_A\,\left|\, 
dx+\frac{1}{2}[x,x]=0\right.\right\}.\]

\item The gauge action of the group $exp(L^0\otimes \ide{m}_A)$ on the set 
$MC(L\otimes \ide{m}_{A})$ is functorial in $A$ and gives an 
action of the group 
functor $G_{L}$ over $MC_{L}$. We call $\Def_{L}=MC_L/G_{L}$ the 
corresponding quotient.
By definition $\Def_L(A)=MC_L(A)/G_{L}(A)$ for every $A\in 
\Art_{\K}$.\\
The functor $\Def_L$ is called the \emph{deformation functor} 
associated to the DGLA $L$.
\end{enumerate}

The reader should make attention to the difference  
between the deformation functor $\Def_{L}$  
associated to a DGLA $L$ and the functor of deformations of a DGLA 
$L$.\\

\begin{prop}\label{VIII.4.1} Let $L=\oplus L^{i}$ be a DGLA. If 
$[L^{1},L^{1}]\cap Z^{2}(L)\subset B^{2}(L)$ (e.g. if $H^{2}(L)=0$) 
then $MC_{L}$ and $\Def_{L}$ are smooth functors.\end{prop}

\begin{proof}  It is sufficient to prove that for every small extension 
\[0\mapor{}J\mapor{} A\mapor{\alpha}B\mapor{}0\] 
the map $MC(L\otimes \ide{m}_{A})\mapor{\alpha}MC(L\otimes \ide{m}_{B})$ is 
surjective.\\
Given $y\in L^{1}\otimes \ide{m}_{B}$ such that $dy+\frac{1}{2}[y,y]=0$ 
we first choose 
$x\in L^{1}\otimes \ide{m}_{A}$ such that $\alpha(x)=y$; we need to prove 
that there exists $z\in L^{1}\otimes J$ such that 
$x-z\in MC(L\otimes \ide{m}_{A})$.\\
Denote
$h=dx+\ds{1\over 2}[x,x]\in L^2\otimes J$; we have 
\[dh=d^2x+[dx,x]=[h,x]-\frac{1}{2}[[x,x],x].\]
Since $[L^2\otimes J,L^1\otimes \ide{m}_A]=0$ we have $[h,x]=0$, by 
Jacobi identity $[[x,x],x]=0$ and then $dh=0$, 
$h\in Z^{2}(L)\otimes J$.\\
On the other hand $h\in ([L^{1},L^{1}]+B^{2}(L))\otimes \ide{m}_{A}$, 
using the assumption of the Proposition
$h\in (B^{2}(L)\otimes \ide{m}_{A})\cap L^{2}\otimes J$
and then there exist $z\in L^1\otimes \ide{m}_A$ such that $dz=h$.\\
Since
$Z^{1}(L)\otimes \ide{m}_{A}\to Z^{1}(L)\otimes \ide{m}_{B}$ is surjective
it is possible to take
$z\in L^{1}\otimes J$: it is now immediate 
to observe that $x-z\in MC(L\otimes \ide{m}_{A})$.\end{proof}

\begin{exer}\label{eseVIII.9}
Prove that if $MC_L$ is smooth then $[Z^1,Z^1]\subset B^2$.
\end{exer}

\begin{prop}\label{VIII.4.2} If $L\otimes \ide{m}_{A}$ is abelian then 
$\Def_{L}(A)=H^{1}(L)\otimes \ide{m}_{A}$. In particular 
$t_{\Def_{L}}=H^{1}(L)\otimes \K\epsilon$, $\epsilon^{2}=0$.
\end{prop}

\begin{proof} The Maurer-Cartan equation reduces to $dx=0$ and then 
$MC_{L}(A)=Z^{1}(L)\otimes \ide{m}_{A}$. If $a\in L^{0}\otimes \ide{m}_{A}$ 
and $x\in L^{1}\otimes \ide{m}_{A}$ we have 
\[ exp(a)x=x+\sum_{n\ge 0}\frac{ad(a)^{n}}{(n+1)!}([a,x]-da)=x-da\]
and then $\Def_{L}(A)=\ds\frac{Z^{1}(L)\otimes \ide{m}_{A}}{d(L^{0}\otimes 
\ide{m}_{A})}= H^{1}(L)\otimes \ide{m}_{A}$.\end{proof}

\begin{exer}\label{eseVIII.10} 
If $[Z^{1},Z^{1}]=0$ then $MC_{L}(A)=Z^{1}\otimes \ide{m}_{A}$ 
for every $A$.
\end{exer}

It is clear that every morphism $\alpha\colon L\to N$ of 
DGLA induces morphisms of functors $G_{L}\to G_{N}$, 
$MC_{L}\to MC_{N}$. These morphisms are compatible with the gauge actions 
and therefore induce a morphism between the deformation functors
$\Def_{\alpha}\colon \Def_L\to \Def_N$.

The following Theorem~\ref{VIII.4.3} 
(together its Corollary~\ref{VIII.4.4}) is 
sometimes called the \emph{basic theorem 
of deformation theory}. 
For the clarity of exposition the (nontrivial) proof of \ref{VIII.4.3} 
is postponed at the end of  Section~\ref{sezioneIV.6}.

\begin{thm}\label{VIII.4.3}
Let $\phi\colon L\to N$ be a morphism of differential graded Lie 
algebras  and denote 
by $H^i(\phi)\colon H^i(L)\to H^i(N)$ the induced maps in cohomology.
\begin{enumerate}
\item If $H^1(\phi)$ is surjective and $H^2(\phi)$ injective then the 
morphism $\Def_{\phi}\colon \Def_L\to \Def_N$ is smooth.

\item If $H^0(\phi)$ is surjective, $H^{1}(\phi)$ is bijective 
and $H^2(\phi)$ is injective then
$\Def_{\phi}\colon\Def_L\to \Def_N$ is an isomorphism.\end{enumerate}
\end{thm}

\begin{cor}\label{VIII.4.4}
Let $L\to N$ be a quasiisomorphism of DGLA.  
Then the induced morphism $\Def_L\to \Def_N$ is an isomorphism.
\end{cor}

\begin{exer}\label{eseVIII.11}
Let $L$ be a formal DGLA, then $\Def_L$ is smooth if and 
only if the induced bracket $[~,~]\colon H^1\times H^1\to H^2$ is zero.
\end{exer}

\begin{ex}\label{VIII.4.5} Let $L=\oplus L^{i}$ be a DGLA and choose
a vector space decomposition $N^1\oplus B^1(L)=L^1$.\\ 
Consider the DGLA  $N=\oplus N^i$ where $N^{i}=0$ if $i<1$ and 
$N^{i}=L^{i}$ if $i>1$ with the differential and bracket induced by $L$.
The natural inclusion 
$N\to L$ gives isomorphisms $H^i(N)\to H^i(L)$ for every $i\ge 1$. In 
particular the morphism $\Def_N\to \Def_L$ is 
smooth and induce an isomorphism on tangent spaces
$t_{\Def_{N}}=t_{\Def_{L}}$ .\end{ex}
  
\textbf{Beware.} One of the most frequent wrong interpretations 
of Corollary~\ref{VIII.4.4} 
asserts that if $L\to N$ is a quasiisomorphism of nilpotent DGLA then 
$MC(L)/exp(L^{0})\to MC(N)/exp(N^{0})$ is a bijection. This is false 
in general: consider for instance $L=0$ and $N=\oplus N^{i}$ with 
$N^{i}=\C$ for $i=1,2$, $N^{i}=0$ for $i\not=1,2$, $d\colon N^{1}\to 
N^{2}$ the identity and $[a,b]=ab$ for $a,b\in N^{1}=\C$.\\
\medskip

Let $T_M$ be the holomorphic tangent bundle of a complex manifold $M$. 
The Kodaira-Spencer DGLA is defined as 
\[KS(M)=\oplus KS(M)^{p},\qquad KS(M)^{p}= \Gamma(M,\sA^{0,p}(T_M))\] 
with the Dolbeault differential
and the bracket (cf. Proposition~\ref{VI.1.2}) 
\[ [\phi d\bar{z}_I, \psi d\bar{z}_J]=[\phi,\psi]d\bar{z}_I\wedge 
d\bar{z}_J\] 
for $\phi,\psi\in \sA^{0,0}(T_M)$, $I,J\subset\{1,...,n\}$ 
and  $z_1,...,z_n$ local holomorphic coordinates.

\begin{thm}\label{VIII.4.6} 
Let $L=KS(M_{0})$ be the Kodaira-Spencer differential 
graded Lie algebra of a compact complex manifold $M_{0}$. Then there 
exists an isomorphism of functors
\[ \Def_{M_{0}}=\Def_{L}.\]
\end{thm}

\begin{proof} Fix $A\in\Art_{\C}$, according to  
Propositions~\ref{VI.1.2} and \ref{VIII.3.8} the exponential 
\[ exp\colon L^{0}\otimes 
\ide{m}_{A}=\Gamma(M_{0},\sA^{0,0}(T_{M_{0}}))\otimes \ide{m}_{A}\to
\Aut_{\sA^{0,*}/\bar{\Omega}^{*}}(A)\]
is an isomorphism.\\
Therefore $\Def_{M_{0}}$ is the quotient of 
\[ MC_{L}(A)=\left\{\eta\in 
\Gamma(M_{0},\sA^{0,1}(T_{M_{0}}))\otimes \ide{m}_{A}\,\left|\,\,
\debar\eta+\frac{1}{2}[\eta,\eta]=0\right.\right\},\]
by the equivalence relation $\sim$, given by
$\eta\sim\mu$ if and only if there exists $a\in L^{0}\otimes \ide{m}_{A}$ 
such that
\[ \debar+\mu=e^{a}(\debar +\eta)e^{-a}=e^{ad(a)}(\debar+\eta)\]
or, equivalently, if and only if $\phi(\mu)=e^{ad(a)}\phi(\eta)$, 
where $\phi$ is the affine embedding defined above.\\
Keeping in mind the definition of the gauge action on the 
Maurer-Cartan elements we get immediately that this equivalence 
relation on $MC_{L}(A)$ is exactly the one induced by the gauge action 
of $exp(L^{0}\otimes \ide{m}_{A})$.\end{proof}

\begin{cor}\label{VIII.4.7} 
    Let $M_{0}$ be a compact complex manifold. If either 
$H^{2}(M_{0},T_{M_{0}})=0$ or its Kodaira-Spencer DGLA 
$KS(M_{0})$ is quasiisomorphic to an abelian DGLA, then 
$\Def_{M_{0}}$ is smooth.\end{cor}

\bigskip

\section{Extended deformation functors (EDF)}	
\label{sezioneIV.5}

We will always work over a fixed field $\K$ of characteristic 0.
All vector spaces, linear maps, algebras, tensor products  
etc. are understood of being over $\K$, unless otherwise 
specified.

We denote by:\begin{itemize} 
\item $\NA$  the category of all differential $\Z$-graded 
associative (graded)-commutative nilpotent finite dimensional 
$\K$-algebras.
\item By $\NA\cap\mathbf{DG}$ we denote the full subcategory of 
$A\in \NA$ with trivial multiplication, i.e. $A^{2}=0$.
\end{itemize}

In other words an object in $\NA$ is a finite dimensional complex 
$A=\oplus A_{i}\in \mathbf{DG}$ endowed with a structure of dg-algebra 
such that $A^{n}=AA\ldots A=0$ for $n>>0$. 
Note that if $A=A_{0}$ is concentrated in degree 0, then $A\in \NA$ if 
and only if $A$ is the maximal ideal of a local artinian $\K$-algebra 
with residue field $\K$.\\
If $A\in\NA$ and $I\subset A$ is a differential ideal, then 
also $I\in \NA$ and the inclusion $I\to A$ is a morphism of 
dg-algebras.\\

\begin{defn}\label{VIII.5.1}
A \emph{small extension} in $\NA$ is a short exact sequence in $\mathbf{DG}$ 
\begin{equation*}
\label{smex}
0\mapor{}I\mapor{}A\mapor{\alpha}B\mapor{}0
\end{equation*}
such that $\alpha$ is a morphism in $\NA$ and $I$ is an ideal of 
$A$ such that  $AI=0$; in addition it is called \emph{acyclic} 
if  $I$ is an acyclic complex, or equivalently if
$\alpha$ is a quasiisomorphism.   				
\end{defn}

\begin{exer}\label{eseVIII.12}~~~
\begin{itemize}   
\item  Every surjective morphism $A\mapor{\alpha} B$ in the category 
$\NA$ is the composition of a finite number of small extensions.   

\item  If $A\mapor{\alpha} B$ is a surjective quasiisomorphism in  
$\NA$ and $A_i=0$ for every $i>0$ then $\alpha$ is the 
composition of a finite number of acyclic small extensions. This is 
generally false if $A_{i}\not=0$ for some $i>0$. 
\end{itemize}
\end{exer}					
					
\medskip

\begin{defn}\label{VIII.5.2}
A covariant functor $F\colon \NA \to 
\Set$ is called a \emph{predeformation functor} if the following
conditions are satisfied:
\begin{enumerate}
\item $F(0)=0$ is the one-point set.
\item  For every pair of 
morphisms 
$\alpha\colon A\to C$, $\beta\colon B\to C$ in 
$\NA$ consider the map
\[\eta\colon F(A\times_C B)\to F(A)\times_{F(C)}F(B)\]
Then:
\begin{enumerate}
\item\label{schlessi1} $\eta$ is surjective when $\alpha$ is surjective.
\item\label{schlessi2} $\eta$ is bijective when $\alpha$ is 
surjective and $C\in\NA\cap \mathbf{DG}$ is an acyclic complex.
\end{enumerate}
\item\label{quasismooth} 
For every acyclic small extension
\[ 0\mapor{} I\mapor{} A\mapor{}B\mapor{}0\]
the induced map $\colon F(A)\to F(B)$ is surjective.
\end{enumerate}
\end{defn} 

If we consider the above definition for a functor defined only for 
algebras concentrated in degree 0, then condition 3 is empty, while 
conditions 1 and 2 are essentially the classical Schlessinger's 
conditions \cite{Sch}, \cite{FaMa1}, \cite{ManettiDGLA}.

\begin{lem}\label{VIII.5.3}
For a covariant functor $F\colon\NA\to \mathbf{Set}$ with $F(0)=0$ it 
is sufficient to check condition \ref{schlessi2} of definition 
\ref{VIII.5.2}  
when $C=0$ and when $B=0$ separately.
\end{lem}

\begin{proof} 
Follows immediately from the equality 
\[ A\times_CB=(A\times B)\times_C 0\]
where $A\mapor{\alpha}C$, $B\mapor{\beta}C$ are as in 
\ref{schlessi2} of \ref{VIII.5.2} and 
the fibred product on the right comes from the morphism $A\times 
B\to C$, $(a,b)\mapsto \alpha(a)-\beta(b)$.
\end{proof}

\begin{defn}\label{VIII.5.4} 
A predeformation functor $F\colon \NA \to 
\mathbf{Set}$ is called a 
\emph{deformation functor} 
if $F(I)=0$ for every acyclic complex $I\in \NA\cap \mathbf{DG}$.
\end{defn}

The predeformation functors (resp.: deformation functors) together their 
natural transformations form a category which we denote by 
$\mathbf{PreDef}$ (resp.: $\mathbf{Def}$).

\begin{lem}\label{VIII.5.5} 
Let  $F\colon \NA\to \mathbf{Set}$ be a deformation 
functor. Then:
\begin{enumerate}
    \item  For every acyclic small extension
\[0\mapor{} I\mapor{} A\mapor{}B\mapor{}0\]
the induced map $\colon F(A)\to F(B)$ is bijective.

    \item  For every pair of complexes $I,J\in\NA\cap\mathbf{DG}$ and 
    every pair of homotopic morphisms $f,g\colon I\to J$, we have 
    $F(f)=F(g)\colon F(I)\to F(J)$.
\end{enumerate}
\end{lem}
					
\begin{proof}
We need to prove that for every acyclic small extension
\[0\mapor{} I\mapor{} A\mapor{\rho}B\mapor{}0\]
the diagonal map $F(A)\to F(A)\times_{F(B)}F(A)$ is surjective; in order 
to prove this it is sufficient to prove that the diagonal map 
$A\to A\times_BA$ induces a surjective map $F(A)\to F(A\times_BA)$. 
We have a canonical isomorphism $\theta\colon A\times I\to A\times_BA$, 
$\theta(a,x)=(a,a+x)$ which sends $A\times \{0\}$ onto the diagonal; since 
$F(A\times I)=F(A)\times F(I)=F(A)$ the proof of item 1 is concluded.\\
For item 2, we can write $I=I^{0}\times I^{1}$, 
$J=J^{0}\times J^{1}$,
with $d(I^{0})=d(J^{0})=0$ 
and $I^{1},J^{1}$ acyclic. Then the inclusion $I^{0}\mapor{i}I$ and the 
projection $J\mapor{p} J^{0}$ induce bijections 
$F(I^{0})=F(I)$, $F(J^{0})=F(J)$. It is now sufficient to note that 
$pfi=pgi\colon I^{0}\to J^{0}$.\end{proof}

A standard argument in Schlessinger's theory \cite[2.10]{Sch} 
shows that for every 
predeformation functor $F$ and every $A\in \NA\cap \mathbf{DG}$ there 
exists a natural structure of vector space on $F(A)$, where the sum and 
the scalar multiplication are described by the maps
\[
A\times A\mapor{+}A\quad \Rightarrow\quad F(A\times A)=F(A)\times 
F(A)\mapor{+}F(A)\]
\[ s\in\K,\quad A\mapor{\cdot s}A\qquad 
\Rightarrow\qquad F(A)\mapor{\cdot s}F(A)\]
We left as an exercise to check that the vector space axioms are 
satisfied; if $A\to B$ is a morphism in $\NA\cap \mathbf{DG}$ then  the 
commutativity of the diagrams
\[ \begin{array}{ccc}
A\times A&\mapor{+}&A\\
\mapver{}&&\mapver{}\\
B\times B&\mapor{+}&B\end{array},\qquad
\begin{array}{ccc}
A&\mapor{\cdot s}&A\\
\mapver{}&&\mapver{}\\
B&\mapor{\cdot s}&B\end{array},\quad s\in\K\]
shows that $F(A)\to F(B)$ is $\K$-linear. Similarly if $F\to G$ is a 
natural transformations of predeformation functors, the 
map $F(A)\to G(A)$ is $\K$-linear for every 
$A\in \NA\cap \mathbf{DG}$.\\

In particular, for every predeformation functor $F$ and  for every integer $n$ 
the sets $F(\Omega[n])$ (see Example~\ref{omegacomplex}) and 
$F(\K[n])$ are vector spaces and the projection $p\colon 
\Omega[n]\to\K[n]$ induce a linear map  
$F(\Omega[n])\to F(\K[n])$

\begin{defn}\label{VIII.5.6} Let $F$ be a predeformation functor, the 
\emph{tangent space} of $F$ is the graded vector space $TF[1]$, where
\[ TF=\somdir{n\in\Z}{}T^nF,\qquad 
T^{n+1}F=TF[1]^{n}=\coker(F(\Omega[n])\mapor{p} F(\K[n])),\quad n\in\Z.\]
A natural transformation $F\to G$ of predeformation functors is 
called a \emph{quasiisomorphism} if induces an isomorphism on tangent 
spaces, i.e. if $T^{n}F\simeq T^{n}G$ for every $n$.
\end{defn}

We note that if $F$ is a deformation functor then $F(\Omega[n])=0$ 
for every $n$ and 
therefore $TF[1]^{n}=T^{n+1}F=F(\K\epsilon)$, where 
$\epsilon$ is an indeterminate of degree $-n\in\Z$ such 
that $\epsilon^{2}=0$.\\ 
In particular $T^{1}F=t_{F^{0}}$, 
where $F^{0}\colon \Art_{\K}\to \mathbf{Set}$, $F^{0}(A)=F(\ide{m}_{A})$, 
is the \emph{truncation} of $F$ in degree 0.\\
\medskip

One of the most important 
examples of deformation functors is the deformation 
functor associated to  a differential graded Lie algebra.\\

Given a DGLA $L$ and $A\in \NA$, the tensor product $L\otimes{A}$ has 
a 
natural structure of nilpotent DGLA with 
\[ (L\otimes{A})^i=\somdir{j\in\Z}{}L^j\otimes A_{i-j}\]
\[ d(x\otimes a)=dx\otimes a+(-1)^{\bar{x}}x\otimes da\]
\[ [x\otimes a, y\otimes b]=(-1)^{\bar{a}\,\bar{y}}[x,y]\otimes ab\]
Every morphism of DGLA, $L\to N$ and every morphism $A\to B$ in $\NA$ give 
a natural commutative diagram of morphisms of differential graded Lie 
algebras
\[\begin{array}{ccc}
L\otimes{A}&\mapor{}&N\otimes{A}\\
\mapver{}&&\mapver{}\\
L\otimes{B}&\mapor{}&N\otimes{B}\end{array}\]
The Maurer-Cartan functor $MC_L\colon\NA\to\mathbf{Set}$ of a DGLA $L$ is by 
definition 
\[ MC_L(A)=MC(L\otimes{A})
=\left\{x\in (L\otimes A)^{1}\,\left|\, 
dx+\frac{1}{2}[x,x]=0\,\right.\right\}.\]

\begin{lem}\label{VIII.5.7} 
For every differential graded Lie algebra $L$, 
$MC_L$ is a predeformation functor.
\end{lem}

\begin{proof}

It is evident that $MC_L(0)=0$ and for every pair of morphisms 
$\alpha\colon A\to C$, $\beta\colon B\to C$ in $\NA$ we have
\[ MC_L(A\times_C B)=MC_L(A)\times_{MC_L(C)}MC_L(B)\]
Let $0\mapor{} I\mapor{} A\mapor{\alpha}B\mapor{}0$ be an acyclic small 
extension and $x\in MC_L(B)$. Since $\alpha$ is surjective there 
exists $y\in (L\otimes{A})^1$ such that $\alpha(y)=x$. Setting 
\[ h=dy+\frac{1}{2}[y,y]\in (L\otimes{I})^2\]
we have 
\[dh=\frac{1}{2}d[y,y]=[dy,y]=[h,y]-\frac{1}{2}[[y,y],y].\] 
By 
Jacobi identity $[[y,y],y]=0$ and, since $AI=0$ also $[h,y]=0$; thus 
$dh=0$ and, being $L\otimes{I}$ acyclic by K\"{u}nneth formula, there 
exists 
$s\in (L\otimes{I})^1$ such that $ds=h$. The element $y-s$ lifts $x$ 
and satisfies the Maurer-Cartan equation. 
We have therefore proved that $MC_L$ is a predeformation functor.
\end{proof}

\begin{exer}\label{eseVIII.13}
Prove that  $MC\colon\mathbf{DGLA}\to \mathbf{PreDef}$ 
is a faithful functor and
every differential graded Lie algebra can be
recovered, up to isomorphism, from its Maurer-Cartan functor.
\end{exer}

It is interesting to point out 
that, if $A\to B$ is a surjective quasiisomorphism in $\NA$, then in 
general $MC_L(A)\to MC_L(B)$ is not surjective.
As an example take $L$ a finite-dimensional non-nilpotent complex Lie 
algebra, considered as a DGLA concentrated in degree 0 and fix $a\in L$ 
such that $ad(a)\colon L\to L$ has an eigenvalue $\lambda\not=0$. Up to 
multiplication of $a$ by $-\lambda^{-1}$ we can assume $\lambda=-1$.
Let $V\subset L$ be the image of $ad(a)$, the linear map $Id+ad(a)\colon 
V\to V$ is not surjective and then there exists $b\in L$ such that the 
equation $x+[a,x]+[a,b]=0$ has no solution in $L$.\\
Let $u,v,w$ be indeterminates of degree 1 and consider the dg-algebras 
\[ B=\C u\oplus \C v,\qquad B^2=0,\, d=0\]
\[ A=\C u\oplus \C v\oplus\C w\oplus\C dw,\qquad uv=uw=dw,\, vw=0\]
The projection $A\to B$ is a quasiisomorphism but the element $a\otimes 
u+b\otimes v\in MC_L(B)$ cannot lifted to $MC_L(A)$. In fact if there 
exists 
$\xi=a\otimes u+b\otimes v+x\otimes w\in MC_L(A)$ then 
\[ 0=d\xi+\frac{1}{2}[\xi,\xi]=(x+[a,x]+[a,b])\otimes dw\]
in contradiction with the previous choice of $a,b$.

\medskip

For every DGLA $L$ and every $A\in \NA$ we define $\Def_L(A)$ as the 
quotient of $MC(L\otimes{A})$ by the gauge action of the group 
$exp((L\otimes{A})^0)$. The gauge action commutes with 
morphisms in $\NA$ and with morphisms of differential graded Lie algebras; 
therefore the above definition gives a  functor 
$\Def_L\colon\NA\to 
\mathbf{Set}$.

\begin{thm}\label{VIII.5.8}
For every DGLA $L$,  $\Def_L\colon \NA\to \mathbf{Set}$ 
is a deformation functor with $T^i\Def_L=H^i(L)$.
\end{thm}

\begin{proof} 
If $C\in \NA\cap\mathbf{DG}$ is a complex 
then $L\otimes{C}$ is an abelian DGLA and according to 
Proposition~\ref{VIII.4.2}, 
$MC_{L}(C)=Z^1(L\otimes{C})$ and 
$\Def_L(C)=H^1(L\otimes{C})$.  In particular 
$T^i\Def_L=H^1(L\otimes\K[i-1])=H^i(L)$ and, by K\"unneth formula, 
$\Def_{L}(C)=0$ if $C$ is acyclic.\\

Since $\Def_L$ is the quotient of a predeformation functor, 
conditions 1 and 3 of \ref{VIII.5.2} are trivially verified and then 
it is sufficient to 
verify  condition 2.\\ 
Let $\alpha\colon A\to 
C$, $\beta\colon B\to C$ morphism in $\NA$ with $\alpha$ surjective. 
Assume there are given $a\in MC_L(A)$, $b\in MC_L(B)$ such that 
$\alpha(a)$ and $\beta(b)$ give the same element in $\Def_L(C)$; then there 
exists $u\in (L\otimes{C})^0$ such that 
$\beta(b)=e^u\alpha(a)$. Let $v\in 
(L\otimes{A})^0$ be a lifting of $u$, 
changing if necessary $a$ with its gauge 
equivalent element $e^va$, we may suppose $\alpha(a)=\beta(b)$ and 
then the pair 
$(a,b)$ lifts to $MC_L(A\times_C B)$: this proves that the map 
\[ \Def_L(A\times_C B)\to \Def_L(A)\times_{\Def_L(C)}\Def_L(B)\]
is surjective.

If $C=0$ then the gauge action $exp((L\otimes(A\times B))^0)\times 
MC_L(A\times B)\to MC_L(A\times B)$ is the direct product of the gauge 
actions $exp((L\otimes{A})^0)\times MC_L(A)\to MC_L(A)$,  
$exp((L\otimes{B})^0)\times MC_L(B)\to MC_L(B)$ 
and therefore $\Def_L(A\times B)=\Def_L(A)\times 
\Def_L(B)$.

Finally assume $B=0$, $C$ acyclic complex and denote $D=\ker\alpha\simeq 
A\times_C B$. Let $a_1,a_2\in MC_L(D)$, $u\in (L\otimes{A})^0$ be 
such that $a_2=e^u a_1$; we need to prove that there exists 
$v\in (L\otimes {D})^0$ such that $a_2=e^v a_1$.\\
Since $\alpha(a_1)=\alpha(a_2)=0$ and $L\otimes{C}$ is an abelian DGLA we 
have $0=e^{\alpha(u)}0=0-d\alpha(u)$ and then $d\alpha(u)=0$. 
$L\otimes C$ is acyclic and then 
there exists $h\in 
(L\otimes{A})^{-1}$ such that $d\alpha(h)=-\alpha(u)$ and $u+dh\in 
(L\otimes{D})^0$. 
Setting  $w=[a_1,h]+dh$, then, according to Remark~\ref{VIII.2.13}, 
$e^w a_1=a_1$ and $e^u e^w a_1=e^v a_1=a_2$, 
where $v=u*w$ is determined by Baker-Campbell-Hausdorff formula. 
We claim that $v\in L\otimes{D}$: in fact $v=u*w\equiv u+w\equiv 
u+dh\pmod{[L\otimes 
{A}, L\otimes{A}]}$ and since $A^2\subset D$ we have 
$v=u*w\equiv u+dh\equiv 0 \pmod{L\otimes{D}}$.
\end{proof}

\begin{lem}\label{VIII.5.9} For every DGLA $L$, the projection
$\pi\colon MC_L\to\Def_L$ is a quasiisomorphism.\end{lem} 

\begin{proof}
Let $i\in\Z$ be fixed; in the notation of \ref{VIII.5.6}
we can write $\Omega[i-1]=\K\epsilon\oplus\K d\epsilon$, where 
$\epsilon^{2}=\epsilon d\epsilon=(d\epsilon)^{2}=0$ and 
$\bar{\epsilon}=1-i$, $\bar{d\epsilon}=2-i$. We have
\[ MC_{L}(\K \epsilon)=\{x\epsilon\in (L\otimes\K 
\epsilon)^{1}| d(x\epsilon)=0\}=Z^{i}(L)\otimes \K\epsilon\]
\begin{equation*}\begin{split}
MC_{L}(\K\epsilon\oplus\!\K d\epsilon)&=\{x\epsilon+y d\epsilon\in 
(L\otimes \Omega[i-1])^{1}\mid dx\epsilon
+(-1)^{1-i}x d\epsilon+dy 
d\epsilon=0\}\\
&=\{(-1)^{i}dy\epsilon+y d\epsilon|\, y\in L^{i-1}\}.
\end{split}\end{equation*}
Therefore the image of $p\colon MC_{L}(\K\epsilon\oplus\!\K d\epsilon)
\to MC_{L}(\K\epsilon)$ is exactly $B^{i}(L)\otimes \K\epsilon$ and then 
\[ MC_{L}(\Omega[i-1])\mapor{p}MC_{L}(\K[i-1])
\mapor{\pi}\Def_{L}(\K[i-1])\mapor{}0\]
is exact.\end{proof}

\bigskip

\section[The inverse function theorem]{Obstruction 
theory and the inverse function theorem for 
deformation functors}	
\label{sezioneIV.6}

\begin{lem}\label{VIII.6.1}
Let $F\colon \NA\to \mathbf{Set}$ be a deformation functor; for every 
complex $I\in\NA\cap\mathbf{DG}$ there exists a natural isomorphism 
\[F(I)=\somdir{i\in\Z}{}TF[1]^i\otimes H_{-i}(I)=
\somdir{i\in\Z}{}T^{i+1}F\otimes H_{-i}(I)=H^{1}(TF\otimes{I}).\]
\end{lem}

\begin{proof}
Let $s\colon H_*(I)\to Z_*(I)$ be a linear section of the natural 
projection, then the composition of $s$ with the natural embedding 
$Z_*(I)\to I$ is unique up to homotopy and its cokernel is an 
acyclic complex, therefore it  
gives a well defined isomorphism  $F(H_*(I))\to F(I)$. This says that it 
is not restrictive to prove the lemma for complexes with zero differential.
Moreover since $F$ commutes with direct sum of complexes we can reduce to 
consider the case when $I=\K^s[n]$ is a vector space concentrated 
in degree $-n$.
Every $v\in I$ gives a morphism $TF[1]^{n}=F(\K[n])\mapor{v}F(I)$ and we 
can define a natural map $TF[1]^{n}\otimes I\to F(I)$, $x\otimes 
v\mapsto v(x)$.  It is easy to verify that this map is an isomorphism of vector 
spaces.\end{proof}

\begin{thm}\label{VIII.6.2}
Let $0\mapor{}I\mapor{\iota}A\mapor{\alpha}B\mapor{}0$ be an exact 
sequence of morphisms in $\NA$ and let $F\colon \NA\to \mathbf{Set}$ be 
a deformation functor.
\begin{enumerate}
    \item  If $AI=0$ then there exist  natural
    transitive  actions of the abelian group 
    $F(I)$ on the nonempty fibres of $F(A)\to F(B)$. 

    \item  If $AI=0$ then there exists a natural ``obstruction map'' 
    $F(B)\mapor{ob}F(I[1])$ with the property that 
    $ob(b)=0$ if and only if $b$ belongs to the image of $F(A)\to 
    F(B)$.

    \item  If $B$ is a complex, i.e. $A^{2}\subset I$, then 
    there exist  natural
    transitive  actions of the abelian group 
    $F(B[-1])$ on the nonempty fibres of $F(I)\to F(A)$.
\end{enumerate}
Here natural means in particular that commutes with natural 
transformation of functors.\end{thm}

\begin{proof} \vale{1} 
There exists an isomorphism of dg-algebras
\[
A\times I\mapor{} A\times_B A;\qquad (a,t)\mapsto (a,a+t)
\]
and then  there exists a natural 
surjective map
\[
\vartheta_F: F(A)\times F(I)=F(A\times I)\to F(A)\times_{F(B)}F(A)
\]
The commutativity of the diagram 
\[
\begin{array}{ccc}
A\times I\times I&\mapor{}&A\times I\\
\mapver{}&&\mapver{}\\
A\times I&\mapor{}&A\end{array},\qquad
\begin{array}{ccc}
(a,t,s)&\mapsto&(a,t+s)\\
\mapver{}&&\mapver{}\\
(a+t,s)&\mapsto&(a+t+s)\end{array}\]
implies in particular that the
composition of $\vartheta_{F}$ with the projection in the second factor
give a  natural transitive 
action of the abelian group 
$F(I)$ on the fibres of the map $F(A)\to F(B)$.\\ 
    
\vale{2} We introduce the mapping cone of $\iota$
as the dg-algebra $C=A\oplus I[1]$ with the 
product 
$(a,m)(b,n)=(ab,0)$ (note that, as a graded algebra, $C$ is the 
trivial extension of $A$ by $I[1]$) and differential 
\[
d_C=\left(\begin{array}{cc}
d_A&\iota\\ 
0&d_{I[1]}\\
\end{array}\right)
\colon A\oplus I[1]\to A[1]\oplus I[2]
\]
We left as exercise the easy 
verification that  $C\in \NA$, the inclusion 
$A\to C$  and the projections $C\to I[1]$, $C\to B$ are morphisms in $\NA$.\\
The kernel of $C\to B$ is isomorphic to 
$I\oplus I[1]$ with differential 
\[
\left(\begin{array}{cc}
d_I&Id_{I[1]}\\ 
0&d_{I[1]}\\
\end{array}\right).
\]
Therefore 
$0\mapor{}I\oplus I[1]\mapor{}C\mapor{}B\mapor{}0$
is an acyclic small extension and then $F(C)=F(B)$.\\
On the other hand $A=C\times_{I[1]}0$ and then the map 
\[ F(A)\to F(C)\times_{F(I[1])}0\]
is surjective. It is sufficient to define $ob$ as the composition of 
the inverse of $F(C)\to F(B)$ with $F(C)\to F(I[1])$.\\

3) The 
derived inverse mapping cone is the dg-algebra $D=A\oplus B[-1]$ 
with product $(x,m)(x,n)=(xy,0)$ 
and differential 
\[
d_D=\left(\begin{array}{cc}
d_A&0\\ \alpha&d_{B[-1]}\\
\end{array}\right)
\colon A\oplus B[-1]\to A[1]\oplus B
\]
Here the projection $D\to A$ 
and the inclusions 
inclusion $I\to D$, $B[-1]\to D$ are morphisms in $\NA$.\\
Since $0\mapor{}B[-1]\mapor{}D\mapor{}A\mapor{}0$ is a small extension, 
by Item 1, there exist natural actions of 
$F(B[-1])$ on the nonempty fibres of $F(D)\to F(A)$.  
The quotient of $I\to D$ is the acyclic complex $B\oplus B[-1]$, and 
then, according to 2b of \ref{VIII.5.2}, $F(I)\to F(D)$ is an isomorphism. 
\end{proof}

\begin{exer}\label{eseIV.14}
Prove that the stabilizers of the actions described in 
Theorem~\ref{VIII.6.2} are vector subspaces.\end{exer}

Given two integers $p\le q$ we denote by $\NA_{p}^{q}$ the full 
subcategory of $\NA$ whose objects are algebras $A=\oplus A_{i}$ such 
that $A_{i}\not=0$ only if $p\le i\le q$.\\

\begin{thm}\label{VIII.6.3}
Let $\theta\colon F\to G$ be a morphism of deformation functors. 
Assume that $\theta\colon TF[1]^{i}\to TG[1]^{i}$ is 
surjective for $p-1\le i\le q$ and injective for $p\le i\le q+1$.
Then:\begin{enumerate}

\item for every surjective morphism $\alpha\colon A\to B$ in the 
category $\NA_{p-1}^{q}$  the morphism 
\[ (\alpha,\theta)\colon F(A)\to F(B)\times_{G(B)}G(A)\] 
is surjective. 

\item  $\theta\colon F(A)\to G(A)$ is 
surjective for every $A\in \NA_{p-1}^{q}$.

\item $\theta\colon F(A)\to G(A)$ is 
a bijection for every $A\in \NA_{p}^{q}$.
\end{enumerate}
\end{thm}

\begin{proof}

The proof  uses the natural generalization to the 
differential graded case of some standard techniques in Schlessinger's 
theory, cf. \cite{FaMa1}.

We first note that, according to Lemma~\ref{VIII.6.1}, for every complex
$I\in\NA_{p}^{q}\cap\mathbf{DG}$ we have that 
$\theta\colon F(I)\to G(I)$ is bijective, 
$\theta\colon F(I[1])\to G(I[1])$ is injective
and $\theta\colon F(I[-1])\to G(I[-1])$ is surjective.\\
Moreover, since $F(0)=G(0)=0$, we have $F(0)\times_{G(0)}G(A)=G(A)$ 
and then Item 2 is an immediate consequence of Item 1.

\begin{step}{1}
For every small extension in $\NA_{p-1}^{q}$,
\[ 0\mapor{} I\mapor{} A\mapor{\alpha} B\mapor{}0\]
and every $b\in F(B)$ we have either 
$\alpha^{-1}(b)=\emptyset$ or $\theta(\alpha^{-1}(b))=
\alpha^{-1}(\theta(b))$.

In fact we have a commutative diagram
\[
\begin{array}{ccc}
F(A)&\mapor{\alpha}&F(B)\\
\mapver{\theta}&&\mapver{\theta}\\
G(A)&\mapor{\alpha}&G(B)\\
\end{array}
\]
and compatible transitive actions of the abelian 
groups $F(I)$, $G(I)$ on the fibres of 
the horizontal maps. Since $F(I)\to G(I)$ is surjective 
this proves Step 1.
\end{step}

\begin{step}{2} Let 
\[ 0\mapor{} I\mapor{\iota} A\mapor{\alpha} B\mapor{}0\]
be a small extension 
in $\NA_{p-1}^{q}$
and $b\in F(B)$. Then $b$ lifts to $F(A)$ if and 
only if $\theta(b)$ lifts to $G(A)$.

The \emph{only if} part is trivial, let's prove the \emph{if} part. 
If $\theta(b)$ lifts to $G(A)$ then $ob(\theta(b))=0$ in $G(I[1])$; 
since the obstruction maps commute with natural transformation of 
functors and $F(I[1])\to G(I[1])$ is injective, also 
$ob(b)=0$ in $F(I[1])$ and then $b$ lifts to $F(A)$.
\end{step}

\begin{step}{3}
For every surjective morphism $\beta\colon A\to C$ in the 
category $\NA_{p-1}^{q}$,  the morphism 
\[ (\alpha,\theta)\colon F(A)\to F(C)\times_{G(C)}G(A)\] 
is surjective. 

Let $J$ be the kernel of $\beta$ and consider the sequence of 
homogeneous differential ideals $J=J_{0}\supset J_{1}=AJ_{0}\supset 
J_{2}=AJ_{1}\cdots$. Since $A$ is nilpotent we have $J_{n}\not=0$ and 
$J_{n+1}=0$ for some $n\ge 0$. Denoting by $I=J_{n}$ and $B=A/I$  we have a 
small extension 
\[ 0\mapor{} I\mapor{} A\mapor{\alpha} B\mapor{}0\]
By induction on $\dim_{\K}A$ we can assume that 
the natural morphism
$F(B)\to F(C)\times_{G(C)}G(B)$ is surjective and therefore it is 
sufficient to prove that 
$F(A)\to F(B)\times_{G(B)}G(A)$ is surjective.\\
Let
$\tilde{a}\in G(A)$ be  fixed element 
and let $b\in F(B)$ such that $\theta(b)=\alpha(\tilde{a})$. By Step 2 
$\alpha^{-1}(b)$ is not empty and then by Step 1 $\tilde{a}\in 
\theta(F(A))$.
\end{step}

\begin{step}{4}
For every surjective morphism $f\colon A\to B$ 
in the category $\NA_{p}^{q}$
and every $a\in F(A)$
we define 
\[ S_F(a,f)=\{\xi\in F(A\times_B A)\, | 
\, \xi\mapsto (a,a)\in F(A)\times_{F(B)}F(A)\subset F(A)\times F(A)\}. \]
By definition, if $f$ is a small extension and $I=\ker f$ then $S_F(a,f)$ 
is naturally isomorphic to the stabilizer of $a$ under the action of 
$F(I)$ on the fibre $f^{-1}(f(a))$. 
It is also clear that:\begin{enumerate} 
    \item $\theta(S_F(a,f))\subset S_G(\theta(a),f)$.
    \item If $\alpha\colon B\to C$ is a surjective morphism if $\NA$, then 
$S_F(a,f)=h^{-1}(S_F(a,\alpha f))$, where $h\colon F(A\times_B A)\to
F(A\times_C A)$ is induced by the natural inclusions $A\times_B 
A\subset A\times_C A$.
    \end{enumerate}
\end{step}

\begin{step}{5}
For  every surjective morphism $f\colon A\to 
B$ in $\NA_{p}^{q}$ and 
every $a\in F(A)$ 
the map $\theta\colon S_F(a,f)\to S_G(\theta(a),f)$ is surjective.

This is trivially true if $B=0$, 
we prove the general assertion by induction on $\dim_{\K}B$.
Let 
\[
0\mapor{}I\mapor{}B\mapor{\alpha}C\mapor{}0
\]
be a small extension with $I\not=0$, set  $g=\alpha f$ and denote by 
$h\colon A\times_C A\to I$ the surjective morphism in $\NA_{p}^{q}$ 
defined 
by $h(a_1,a_2)=f(a_1)-f(a_2)$; we have an exact sequence
\[ 0\mapor{}A\times_B A\mapor{\iota}A\times_C A\mapor{h}I\mapor{}0.\]

According to 2a of \ref{VIII.5.2} the maps
\[ F(A\times_B A)\to F(A\times_C A)\cap h^{-1}(0);\qquad
S_F(a,f)\to S_F(a,g)\cap h^{-1}(0)
\]
are surjective.

Let $\tilde{\xi}\in S_G(\theta(a),f)$ and let $\eta\in S_F(a,g)$ such 
that $\theta(\eta)=\iota(\tilde{\xi})$. Since $F(I)=G(I)$ we have 
$h(\eta)=0$ and then $\eta$ lifts to some $\xi_1\in S_F(a,f)$.
According to Theorem~\ref{VIII.6.2} there exist 
surjective maps commuting with $\theta$
\[
F(A\times_B A)\times F(I[-1]) \mapor{\varrho}F(A\times_B 
A)\times_{F(A\times_C A)}F(A\times_B A)
\]
\[
G(A\times_B A)\times G(I[-1]) \mapor{\varrho}G(A\times_B 
A)\times_{G(A\times_C A)}G(A\times_B A)
\]
Since $F(I[-1])\to G(I[-1])$ is surjective 
there exists  
$v\in F(I[-1])$ such that 
$\varrho(\theta(\xi_1),\theta(v))=(\theta(\xi_1),\tilde{\xi})$; defining  
$\xi\in F(A\times_{B}A)$ by the formula 
$\varrho(\xi_1,v)=(\xi_1,\xi)$ we get  $\theta(\xi)=\tilde{\xi}$ and 
then $\xi\in S_F(a,f)$.
\end{step}

\begin{step}{6} 
For every $A\in \NA_{p}^{q}$ the map $\theta\colon F(A)\to G(A)$ is 
injective. 

According to Lemma~\ref{VIII.6.1}
this is true  if $A^2=0$; if $A^2\neq 0$ we can suppose by 
induction that there exists a small extension 
\[ 0\mapor{} I\mapor{\iota} A\mapor{\alpha} B\mapor{}0\]
with $I\neq 0$ and $\theta\colon F(B)\to G(B)$ injective.

Let $a_1,a_2\in F(A)$ be two elements such that 
$\theta(a_1)=\theta(a_2)$; by assumption $f(a_1)=f(a_2)$ and then there 
exists $t\in F(I)$ such that  $\vartheta_F(a_1,t)=(a_1,a_2)$. Since $\vartheta$ is a natural
transformation
$\vartheta_G(\theta(a_1),\theta(t))=(\theta(a_1),\theta(a_2))$ and then 
$\theta(t)\in S_G(\theta(a_1),\alpha)$. By Step 5 there exists $s\in 
S_F(a_1,\alpha)$ such that $\theta(s)=\theta(t)$ and by injectivity of 
$\theta\colon F(I)\to G(I)$ we get $s=t$ and then $a_1=a_2$.
\end{step}
~\end{proof}

As an immediate consequence we have:

\begin{cor}\label{VIII.6.4}
A morphism of deformation functors $\theta\colon F\to G$ is an 
isomorphism if and only if it is a quasiisomorphism.
\end{cor}

\begin{proof}[Proof of Theorem~\ref{VIII.4.3}]
We apply Theorem~\ref{VIII.6.3} to the morphism of deformation functors
$\theta=\Def_{\phi}\colon 
\Def_L\to \Def_N$.\\
According to Theorem~\ref{VIII.5.8}, the first item of \ref{VIII.4.3} 
is exactly the first item of 
\ref{VIII.6.3} for $p=1, q=0$, while the second item of  \ref{VIII.4.3}
is exactly the third item of \ref{VIII.6.3} for $p=q=0$.    
\end{proof}

\bigskip

\section{Historical survey, \ref{CAP:DGLA}}
\label{sec:histdgla}

The material Sections \ref{sec:expo}, \ref{sec:freelie} and 
\ref{sec:nilplie}
is standard and well exposed in 
the literature about Lie algebras; in Sections
\ref{sec:dgla}, \ref{sec:functors} and 
\ref{sec:defofun} we follows the approach of \cite{ManettiDGLA}, while 
the material of Sections \ref{sezioneIV.5} and 
\ref{sezioneIV.6} comes from \cite{ManettiEDF}.

\medskip

Some remarks on the introduction of this \chaptername:\\
~\\
A) Given a deformation problem, in general it is not an easy task to find 
a factorization as in the introduction  and in some cases it is still unknown.\\
~\\
B) Even in the simplest examples,  
the governing DGLA is only 
defined up to (non canonical) 
quasiisomorphism and then the Theorem~\ref{VIII.4.3}
is a necessary background for the whole theory.\\
For example, there are very good reasons to consider, 
for the study of deformations of a compact 
complex manifold  $M$, the DGLA $L=\oplus L^{i}$, where $L^{i}$ is the 
completion of $\Gamma(M,\sA^{0,i}(T_{M}))$ is a suitable Sobolev's 
norm. According to elliptic regularity the inclusion 
$\KS(M)\subset L$ is a quasiisomorphism of DGLA.\\
In general a correct procedure gives, for every deformation problem 
$P$ with associated deformation functor $\Def_{P}$, 
a connected subcategory $\mathbf{P}\subset\mathbf{DGLA}$ with the following 
properties:
\begin{enumerate}
    \item  If $L$ is an object of $\mathbf{P}$ then $\Def_{L}=\Def_{P}$.

    \item Every morphism in $\mathbf{P}$ is a quasiisomorphism of DGLA.

    \item  If $\Mor_{\mathbf{P}}(L,N)\not=\emptyset$ then the induced 
    isomorphism $\Def_{\alpha}\colon \Def_{L}\to \Def_{N}$
    is independent from the choice of $\alpha\in \Mor_{\mathbf{P}}(L,N)$.
\end{enumerate}
~\\
C) It may happen that two people, say Circino and Olibri, starting 
from the same deformation problem, get two non-quasiisomorphic 
DGLA governing the problem. This is possible because the DGLA governs 
an extended (or derived) deformation problem. If Circino and Olibri
have in mind two different extensions of the problem 
then they get different DGLA.
~\\
D) Although the interpretation of deformation problems in terms of solutions 
of Maurer-Cartan equation is  very useful on its own, 
in many situation it is unavoidable to recognize that 
the category of DGLA is too rigid for a ``good'' theory.  
The appropriate way of extending this category will be the 
introduction of  $L_{\infty}$-algebras; these new 
objects will be described in \chaptername~\ref{CAP:TOOLS}.


\chapter[~K\"{a}hler manifolds]{K\"{a}hler manifolds}
\label{CAP:KAEHLER}
\piede

This chapter provides a basic introduction to K\"ahler 
manifolds. We first study the local theory, following essentially 
Weil's book \cite{Weil} and then, assuming harmonic and elliptic 
theory, we give a proof of the $\de\debar$-lemma which is presented 
both  in the 
classical  version 
(Theorem~\ref{X.4.9}, Item 2) and  in the ``homological'' version 
(Theorem~\ref{X.4.9}, Item 1).\\
The material of this \chaptername\  is widely present 
in the literature, with the possible exception of the homological 
version of $\de\debar$-lemma; I only tried to simplify the presentation and some 
proofs. The main references are \cite{Weil}, \cite{Wells} and \cite{DGMS} 

\bigskip
\section{Covectors on complex vector spaces}
\label{sec:covectors}

Given a complex vector space $E$ of dimension $n$ we denote by:
\begin{itemize}
\item  $E\dual=\Hom_{\C}(E,\C)$ its dual.

\item  $E_{\C}=E\otimes_{\R}\C$, with the structure of $\C$-vector 
space induced by the scalar multiplication $a(v\otimes b)=v\otimes 
ab$. 

\item  $\bar{E}$ its complex conjugate.
\end{itemize}
The conjugate $\bar{E}$ is defined as the set of formal symbols 
$\bar{v}$, $v\in E$ with the vector space structure given by 
\[ \bar{v}+\bar{w}=\bar{v+w},\qquad a\bar{v}=\bar{\bar{a}v}.\]
The conjugation $~^{\bar{~}}\colon E\to \bar{E}$, $v\mapsto \bar{v}$ is a 
$\R$-linear isomorphism. 

There exists a list of natural isomorphisms (details left as exercise)
\begin{enumerate}
\item  $(E_{\C})\dual=(E\dual)_{\C}=\Hom_{\R}(E,\C)$

\item  $\bar{E\dual}=\bar{E}\dual$ given by $\bar{f}(\bar{v})=\bar{f(v)}$, 
$f\in E\dual$, $v\in E$.\\

\item  $ E\oplus \bar{E}\to E_{\C}$,  $\quad(v,\bar{w})\mapsto v\otimes 
1-iv\otimes i+w\otimes 1+iw\otimes i$, being $i$ a square root of $-1$.

\item  $E\dual\oplus \bar{E\dual}\to E_{\C}\dual=
\Hom_{\R}(E,\C)$,  $\quad(f,\bar{g})(v)=f(v)+\bar{g(v)}$.
\end{enumerate}

Under these isomorphisms, the image of $E\dual$ (resp.: $\bar{E}\dual$)
inside $E_{\C}\dual$ is the subspace of $f$ such that $f(iv)=if(v)$ 
(resp.: $f(iv)=-if(v)$). Moreover $E\dual=\bar{E}^{\perp}$, 
$\bar{E}\dual=E^{\perp}$.\\

For $0\le p,q\le n$ we set $\sA^{p,q}=
\bigwedge^{p}E\dual\otimes \bigwedge^{q}\bar{E}\dual$: this is called 
the space of $(p,q)$-covectors of $E$.
We also set $\sA^{p}=\oplus_{a+b=p}\sA^{a,b}$ (the space of 
$p$-covectors) and $\sA=\oplus_{a,b}\sA^{a,b}$. 
Denote by $P_{a,b}\colon\sA\to\sA^{a,b}$, $P_{p}\colon \sA\to \sA^{p}$ 
the projections.\\
If $z_{1},\ldots,z_{n}$ is a basis of $E\dual$ then
$\bar{z_{1}},\ldots,\bar{z_{n}}$ is a basis of $\bar{E}\dual$ and 
therefore
\[z_{i_{1}}\wedge\ldots\wedge z_{i_{p}}\wedge 
\bar{z_{j_{1}}}\wedge\ldots\wedge\bar{z_{j_{q}}},\qquad  
i_{1}<\ldots<i_{p},\, j_{1}<\ldots<j_{q}\]
is a basis of $\sA^{p,q}$.
Since $E_{\C}\dual=E\dual\oplus\bar{E}\dual$, 
we have 
$\external{}{*}E_{\C}\dual=\sA$.\\
The complex conjugation is defined in $\sA$ and gives a $\R$-linear isomorphism 
$~^{\bar{~}}\colon \sA\to \sA$. On the above basis, the conjugation 
acts as
\[\overline{z_{i_{1}}\wedge\ldots\wedge z_{i_{p}}\wedge 
\bar{z_{j_{1}}}\wedge\ldots\wedge\bar{z_{j_{q}}}}=
(-1)^{pq}z_{j_{1}}\wedge\ldots\wedge z_{j_{q}}\wedge 
\bar{z_{i_{1}}}\wedge\ldots\wedge\bar{z_{i_{p}}}.\]
Since
$\bar{\sA^{a,b}}=\sA^{b,a}$, we have 
$P_{a,b}(\bar{\eta})=\bar{P_{b,a}(\eta)}$.

\begin{defn}\label{X.1.1} The operator $C\colon 
\sA\to\sA$ is 
defined by the formula 
\[ C=\sum_{a,b}i^{a-b}P_{a,b}.\]
\end{defn}
Note that $\bar{C(u)}=C(\bar{u})$  (i.e. $C$ is a real operator) and 
$C^{2}=\sum_{p}(-1)^{p}P_{p}$.\\

\section{The exterior algebra of a Hermitian space}
\label{sec:exteriorherm}

Let $E$ be a complex vector space of dimension $n$.
A Hermitian form on $E$ is a $\R$-bilinear map $h\colon E\times 
E\to\C$ satisfying the conditions 
\begin{enumerate}
\item  $h(av,w)=ah(v,w)$, ~$h(v,aw)=\bar{a}h(v,w)$,~~~ $a\in\C$, $v,w\in E$.

\item  $h(w,v)=\bar{h(v,w)}$,~~~ $v,w\in E$. 
\end{enumerate}

Note that $h(v,v)\in \R$ for every $v$. $h$ is called positive 
definite if $h(v,v)>0$ for every $v\not=0$. 

\begin{defn}\label{X.1.2} A \emph{Hermitian space} is a pair $(E,h)$ where $h$ is a 
positive definite Hermitian form on $E$.\end{defn}

It is well known that a 
Hermitian form $h$ on a finite dimensional vector space $E$
is positive definite if and only if it 
admits a unitary basis, i.e. a basis $e_{1},\ldots,e_{n}$ of $E$ such 
that $h(e_{i},e_{j})=\delta_{ij}$.

Every Hermitian space $(E,h)$ induces canonical Hermitian structures 
on the associated vector spaces. For example 
\[\bar{h}\colon \bar{E}\times\bar{E}\to \C,\qquad
\bar{h}(\bar{v},\bar{w})=\bar{h(v,w)}\]
and 
\[h^{p}\colon \external{}{p}E\times \external{}{p}E\to \C,\qquad
h^{p}(v_{1}\wedge\ldots\wedge v_{p},w_{1}\wedge\ldots\wedge w_{p})=
\det(h(v_{i},w_{j}))\]
are Hermitian forms. If $e_{1},\ldots,e_{n}$ is a unitary basis of $E$ 
then $\bar{e_{1}},\ldots, \bar{e_{n}}$ is a unitary basis for 
$\bar{h}$ and $e_{i_{1}}\wedge \ldots\wedge e_{i_{p}}$, 
$i_{1}<\ldots<i_{p}$, is a unitary basis for $h^{p}$.\\

Similarly, if $(F,k)$ is another Hermitian space then we have natural 
Hermitian structures on $E\otimes F$ and $\Hom_{\C}(E,F)$ given by 
\[ hk\colon E\otimes F\to \C,\qquad hk(v\otimes f, w\otimes 
g)=h(v,w)k(f,g)\]
\[ h\dual k\colon \Hom_{\C}(E,F)\to \C,\qquad
h\dual k(f,g)=\sum_{i=1}^{n}k(f(e_{i}), g(e_{i}))\]
where $e_{i}$ is a unitary basis of $E$. It is an easy exercise (left 
to the reader) to prove that $h\dual k$ is well defined and positive 
definite.\\
In particular the complex dual $E\dual$ is a Hermitian space and the dual basis 
of a unitary basis for $h$ is 
a unitary basis for $h\dual$.\\ 

Let $e_{1},\ldots,e_{n}$ be a basis of $E$, 
$z_{1},\ldots,z_{n}\in E\dual$ its dual basis; then 
\[ h(v,w)=\sum_{i,j}h_{ij}z_{i}(v)\bar{z_{j}(w)}\]
where $h_{ij}=h(e_{i},e_{j})$. We have $h_{ji}=\bar{h_{ij}}$ and the 
basis is unitary if and only if $h_{ij}=\delta_{ij}$. We then 
write $h=\sum_{ij}h_{ij}z_{i}\otimes\bar{z_{j}}$; in doing this we 
also consider $h$ as an element of $E\dual\otimes 
\bar{E}\dual=(E\otimes\bar{E})\dual$.

Taking the real and the imaginary part of $h$ we have $h=\rho-i\omega$, 
with $\rho,\omega\colon E\times E\to \R$. It is immediate to observe 
that $\rho$ is symmetric, $\omega$ is skewsymmetric and 
\[ \rho(iv,iw)=\rho(v,w),\quad \omega(iv,iw)=\omega(v,w),\quad 
\rho(iv,w)=\omega(v,w).\]
Since $z_{i}\wedge\bar{z_{j}}=z_{i}\otimes\bar{z_{j}}-\bar{z_{j}}\otimes 
z_{i}$, we can write  
\[ \omega=\frac{i}{2}(h-\bar{h})=
\frac{i}{2}\sum_{ij}h_{ij}z_{i}\wedge\bar{z_{j}}\in 
\sA^{1,1}.\]
Note that $\omega$ is real, i.e. $\bar{\omega}=\omega$, and the 
Hermitian form is positive definite if and only if for every 
$v\not=0$, $h(v,v)=\rho(v,v)=\omega(v,iv)>0$. The basis 
$e_{1},\ldots,e_{n}$ is unitary if and only if 
$\omega=\ds\frac{i}{2}\sum_{i}z_{i}\wedge \bar{z_{i}}$.\\

Let now $e_{1},\ldots,e_{n}$ be a fixed unitary basis 
of a Hermitian space $(E,h)$ with dual basis $z_{1},\ldots,z_{n}$ and denote 
$u_{j}=\ds\frac{i}{2}z_{j}\wedge 
\bar{z_{j}}$; if $z_{j}=x_{j}+iy_{j}$ then $u_{j}=x_{j}\wedge y_{j}$ 
and 
\[ \frac{\omega^{\wedge n}}{n!}=u_{1}\wedge \ldots\wedge 
u_{n}=x_{1}\wedge y_{1}\wedge\ldots\wedge x_{n}\wedge y_{n}.\]
Since $x_{1},y_{1},\ldots,x_{n},y_{n}$ is a system of coordinates on 
$E$, considered as a real oriented 
vector space of dimension $2n$ and 
the quadratic form $\rho$ is written in this coordinates 
\[ \rho(v,v)=\sum_{i=1}^{n} (x_{i}(v)^{2}+y_{i}(v)^{2}),\]
we get from the above formula that 
$\omega^{\wedge n}/n!\in \external{\R}{2n}\Hom_{\R}(E,\R)$  
is the volume element associated to the scalar product 
$\rho$ on $E$. 

For notational simplicity, if $A=\{a_{1},\ldots,a_{p}\}\subset
\{1,\ldots,n\}$  and $a_{1}<a_{2}<\ldots<a_{p}$, we denote $|A|=p$ and
\[ z_{A}=z_{a_{1}}\wedge\ldots\wedge z_{a_{p}},\quad
\bar{z}_{A}=\bar{z}_{a_{1}}\wedge\ldots\wedge \bar{z}_{a_{p}},\quad
 u_{A}=u_{a_{1}}\wedge\ldots\wedge u_{a_{p}}.\]
 
For every decomposition of $\{1,\ldots,n\}=A\cup B\cup M\cup N$ into four 
disjoint subsets, we denote 
\[z_{A,B,M,N}=\frac{1}{\sqrt{2^{|A|+|B|}}}
z_{A}\wedge\bar{z}_{B}\wedge u_{M}\in \sA^{|A|+|M|,|B|+|M|}.\] 
These elements give a basis of $\sA$
which we call \emph{standard basis}.\\ 
Note that $\bar{z_{A,B,M,N}}=(-1)^{|A|\,|B|}z_{B,A,M,N}$.
 
\begin{defn}\label{X.1.3} The $\C$-linear 
operator $*\colon \sA^{p,q}\to \sA^{n-q,n-p}$ is 
defined as 
\[ *z_{A,B,M,N}=sgn(A,B)i^{|A|+|B|}z_{A,B,N,M},\]
where $sgn(A,B)=\pm 1$ is  
the sign compatible with the formulas
\begin{equation}
z_{A,B,M,N}\wedge *\bar{z_{A,B,M,N}}=
z_{A,B,M,N}\wedge\bar{*z_{A,B,M,N}}=
u_{1}\wedge\ldots\wedge 
u_{n}.\end{equation}
\begin{equation}
C^{-1}*z_{A,B,M,N}=(-1)^{\frac{(|A|+|B|)(|A|+|B|+1)}{2}}z_{A,B,N,M}=
(-1)^{\frac{(p+q)(p+q+1)}{2}+|M|}z_{A,B,N,M}.
\end{equation}
\end{defn}

\begin{exer}\label{eseX.1} 
Verify that Definition~\ref{X.1.3} is well 
posed.\end{exer} 

In particular   
\[*^{2}z_{A,B,M,N}=(-1)^{|A|+|B|}z_{A,B,M,N}=
(-1)^{|A|+|B|+2|M|}z_{A,B,M,N}\] 
and then 
\[(C^{-1}*)^{2}=Id,\qquad *^{2}=C^{2}=\sum_{p}(-1)^{p}P_{p}.\]

If we denote $vol\colon \C\to \sA^{n,n}$ the multiplication for the 
``volume element'' $\omega^{\wedge n}/n!$, then $vol$ is an 
isomorphism and we can consider the $\R$-bilinear maps 
\[ (\, ,\,)\colon\sA^{a,b}\times \sA^{a,b}\to \C,\qquad 
(v,w)=vol^{-1}(v\wedge\bar{*w})=vol^{-1}(v\wedge*\bar{w}).\]
Clearly $(\, ,\,)$ is $\C$-linear on the first member and 
$\C$-antilinear in the second; since 
\renewcommand\arraystretch{1.8}
\[(z_{A,B,M,N},z_{A',B',M',N'})=\left\{
\begin{array}{ll} 1&\hbox{ if }A=A', B=B', M=M', N=N'\\
0&\hbox{ otherwise }\end{array}\right.\]
\renewcommand\arraystretch{1}
we have that $(\, ,\,)$ 
is a positive definite Hermitian form with the 
$z_{A,B,M,N}$'s, $|A|+|M|=a$, $|B|+|M|=b$, a unitary basis; since $*$ 
sends unitary basis into unitary basis we also get that  
$*\colon \sA^{a,b}\to \sA^{n-b,n-a}$ is an isometry.\\

\begin{lem} The Hermitian form $(\, ,\,)$ is the Hermitian form associated
to the Hermitian space $(E,h/2)$. In particular $(\, ,\,)$ and $*$ are independent 
from the choice of the unitary basis $e_{1},\ldots,e_{n}$.\end{lem}

\begin{proof} The basis $\sqrt{2}e_{1},\ldots, \sqrt{2}e_{n}$ is a 
unitary basis for $h/2$ and then the standard basis is a unitary 
basis for the associated Hermitian structures on $\sA$.\\  
From the formula $(v,w)\omega^{\wedge n}=n!(v\wedge\bar{*w})$ and from 
the fact that the wedge product is nondegenerate follows that $*$ 
depends only by $\omega$ and $(\, ,\,)$.
\end{proof}

Consider now, for every $j=1,\ldots,n$, the linear operators
\[ L_{j}\colon \sA^{p,q}\to \sA^{p+1,q+1},\quad 
L_{j}(\eta)=\eta\wedge u_{j},\]
\[ \Lambda_{j}\colon \sA^{p,q}\to \sA^{p-1,q-1},\quad 
\Lambda_{j}(\eta)=\eta\dashv\left(\frac{2}{i}\bar{e_{j}}\wedge 
e_{j}\right),\]
where $\dashv$ denotes the contraction on the right. 
More concretely, in the standard  basis
\renewcommand\arraystretch{1.6}
\[ L_{i}z_{A,B,M,N}=\left\{\begin{array}{ll}
z_{A,B,M\cup\{i\},N-\{i\}}\quad&\hbox{ if } i\in N\\
0&\hbox{ otherwise}\end{array}\right.\]
\[ \Lambda_{i}z_{A,B,M,N}=\left\{\begin{array}{ll}
z_{A,B,M-\{i\},N\cup\{i\}}\quad&\hbox{ if } i\in M\\
0&\hbox{ otherwise}\end{array}\right.\]
\renewcommand\arraystretch{1}
It is therefore immediate to observe that $L_{i}*=*\Lambda_{i}$ and $*L_{i}=\Lambda_{i}*$. 
Setting $L=\sum_{i}L_{i}$, $\Lambda=\sum_{i}\Lambda_{i}$ we have 
therefore 
\[ L(\eta)=\eta\wedge \omega,\qquad \Lambda=*^{-1}L*=*L*^{-1}.\]

\begin{lem}\label{X.1.4} The operators $L$ and $\Lambda$ do not depend from the 
choice of the unitary basis.\end{lem}

\begin{proof} $\omega$ and $*$ do not depend.\end{proof} 

\begin{prop}\label{X.1.5} The following commuting relations hold:
\[ [L,C]=0,\quad [\Lambda,C]=0,\quad [*,C]=0,\quad 
[\Lambda,L]=\sum_{p=0}^{2n}(n-p)P_{p}.\]
\end{prop}
\begin{proof} Only the last is nontrivial, 
we have:
\[ Lz_{A,B,M,N}=\sum_{i\in N}z_{A,B,M\cup\{i\},N-\{i\}},\qquad 
\Lambda z_{A,B,M,N}=\sum_{i\in M}z_{A,B,M-\{i\},N\cup\{i\}},\]
\[\Lambda Lz_{A,B,M,N}=\sum_{i\in N}z_{A,B,M,N}
+\sum_{j\in M}\sum_{i\in N}
z_{A,B,M\cup\{i\}-\{j\},N\cup\{j\}-\{i\}},\] 
\[L\Lambda
z_{A,B,M,N}=\sum_{i\in M}z_{A,B,M,N} +\sum_{j\in M}\sum_{i\in N}
z_{A,B,M\cup\{i\}-\{j\},N\cup\{j\}-\{i\}}.\] 
Therefore we get 
\[ (\Lambda 
L-L\Lambda)z_{A,B,M,N}=(|N|-|M|)z_{A,B,M,N}=(n-|A|-|B|-2|M|)z_{A,B,M,N}.\]
and then 
\[ [\Lambda,L]=\sum_{p=0}^{2n}(n-p)P_{p}.\]
\end{proof}

\bigskip

\section{The Lefschetz decomposition}

The aim of this section is to study the structure of 
$\bigwedge^{*,*}E\dual$ as a module over the algebra $\Phi$ generated 
by the linear operators $C^{-1}*,L,\Lambda$.\\ 
In the notation of the 
previous section, it is immediate to see that there exists 
a direct sum decomposition of $\Phi$-modules 
$\bigwedge^{*,*}E\dual=\bigoplus V_{A,B}$, where $V_{A,B}$ is the 
subspace generated by the $2^{n-|A|-|B|}$ elements $z_{A,B,M,N}$, 
being $A,B$ fixed.\\
It is also clear that every $V_{A,B}$ is isomorphic to one of the 
$\Phi$-modules $V(h,\tau)$, $h\in\N$, $\tau=\pm 1$, defined in the 
following way:
\begin{enumerate}
\item  $V(h,\tau)$ is 
the $\C$-vector space  
with basis $u_{M}$, $M\subset\{1,\ldots,h\}$. 
 
\item  The linear operators $L,\Lambda$ and $C^{-1}*$ act on 
$V(h,\tau)$ as
\[ Lu_{M}=\sum_{i\not\in M}u_{M\cup\{i\}},\quad
\Lambda u_{M}=\sum_{i\in M}u_{M-\{i\}},\quad
C^{-1}*u_{M}=\tau u_{M^{c}},\]
where $M^{c}=\{1,\ldots,h\}-M$ denotes the complement of $M$.
\end{enumerate}

We have a direct sum decomposition 
\[ V(h,\tau)=\somdir{\alpha\equiv h\!\!\pmod{2}}{}V_{\alpha},\]
~\\
where $V_{\alpha}$ is the subspace generated by the $u_{M}$ with 
$|M^{c}|-|M|=\alpha$. An element of $V_{\alpha}$ is called homogeneous of 
weight $\alpha$. Set $P_{\alpha}\colon V(h,\tau)\to 
V_{\alpha}$ the projection.

Note that $L\colon V_{\alpha}\to V_{\alpha-2}$, $\Lambda\colon 
V_{\alpha}\to V_{\alpha+2}$ and $C^{-1}*\colon V_{\alpha}\to 
V_{-\alpha}$. 

We have already seen that 
\[[\Lambda,L]=\sum_{\alpha\in\Z}\alpha 
P_{\alpha},\qquad LC^{-1}*=C^{-1}*\Lambda,\qquad C^{-1}*L=\Lambda 
C^{-1}*.\]
A simple combinatorial argument shows that for every $r\ge 0$, 
\[L^{r}u_{M}=r!\sum_{M\subset N, |N|=|M|+r}u_{N}.\]

\begin{lem}\label{X.2.1} For every $r\ge 1$ we have
\[ [\Lambda,L^{r}]=\sum_{\alpha} r(\alpha-r+1)L^{r-1}P_{\alpha}.\]
\end{lem}

\begin{proof} This has already done for $r=1$, we prove the general 
statement for induction on $r$. We have 
\[ [\Lambda,L^{r+1}]=[\Lambda,L^{r}]L+L^{r}[\Lambda,L]=
\sum_{\alpha} r(\alpha-r+1)L^{r-1}P_{\alpha}L+\sum_{\alpha}\alpha 
P_{\alpha}.\]
Since $P_{\alpha}L=LP_{\alpha+2}$ we have 
\[ [\Lambda,L^{r+1}]=\sum_{\alpha}
r(\alpha-r+1)L^{r}P_{\alpha+2}+\sum_{\alpha}\alpha  P_{\alpha}=\sum_{\alpha}
(r(\alpha-r-1)+\alpha)L^{r}P_{\alpha}.\]
\end{proof}

\begin{defn}\label{X.2.2} 
A homogeneous vector $v\in V_{\alpha}$ is called \emph{primitive} 
if $\Lambda v=0$.\end{defn}

\begin{prop}\label{X.2.3} Let 
$v\in V_{\alpha}$ be a  primitive element, then:
\begin{enumerate}
    \item  $L^{q}v=0$ for 
every $q\ge\max(\alpha+1,0)$. In particular if $\alpha<0$ then 
$v=L^{0}v=0$.

    \item  If $\alpha\ge 0$, 
    then for every $p>q\ge 0$ 
\[ \Lambda^{p-q}L^{p}v=
\prod_{r=q+1}^{p}r(\alpha-r+1)L^{q}v;\]
in particular
$\Lambda^{\alpha}L^{\alpha}v=\alpha!^{2}v$.
\end{enumerate}
\end{prop}

\begin{proof} We first note that for $s,r\ge 1$ 
\[ \Lambda^{s}L^{r}v=\Lambda^{s-1}[\Lambda,L^{r}]v=
r(\alpha-r+1)\Lambda^{s-1}L^{r-1}v.\]
and then for every $p>q\ge 0$ 
\[ \Lambda^{p-q}L^{p}v=
\prod_{r=q+1}^{p}r(\alpha-r+1)L^{q}v.\]
If $p>q>\alpha$ then $r(\alpha-r+1)\not=0$ for every $r>q$ and then 
$L^{q}v=0$ if and only if $ \Lambda^{p-q}L^{p}v=0$: 
taking $p>>0$ we get the required vanishing.\end{proof}

\begin{lem}\label{X.2.4} Let $\alpha\ge 0$, 
$m=(h-\alpha)/2$ and $v=\sum_{|M|=m}a_{M}u_{M}\in V_\alpha$, 
$a_{M}\in\C$. If $v$ is primitive, then for every $M$
\[ a_{M}=(-1)^{m}\!\!\sum_{N\subset M^{c}, |N|=m}\!\! a_{N}.\]
\end{lem}

\begin{proof} For $m=0$ the above equality becomes 
$a_{\emptyset}=a_{\emptyset}$ and therefore we can assume $m>0$.
Let $M\subset \{1,\ldots,h\}$ be a fixed subset of cardinality $m$, 
since  
\[ 0=\Lambda v=\sum_{|H|=m}a_{H}\sum_{i\in H}u_{H-\{i\}}=
\sum_{|N|=m-1}u_{N}\sum_{i\not\in N}a_{N\cup\{i\}}\]
we get for every $N\subset\{1,\ldots,h\}$ 
of cardinality $m-1$ the equality
\[ R_{N}:\qquad \sum_{i\in M-N}a_{N\cup\{i\}}=-\sum_{i\not\in M\cup 
N}a_{N\cup\{i\}}.\]
For every $0\le r\le m$ denote by 
\[ S_{r}=\sum_{|H|=m, |H\cap M|=r}a_{H}.\]
Fixing an integer $1\le r\le m$ and taking the sum of  the equalities 
$R_{N}$, for all $N$ such that $|N\cap M|=r-1$ we get 
\[ rS_{r}=-(m-r+1)S_{r-1}\]
and then 
\[ a_{M}=S_{m}=-\frac{S_{m-1}}{m}=\frac{2S_{m-2}}{m(m-1)}=\ldots=
(-1)^{m}\frac{m!}{m!}S_{0}=(-1)^{m}\!\!\!\!\sum_{N\subset M^{c}, 
|N|=m}\!\!\!\! a_{N}.\]
\end{proof}

\begin{lem}\label{X.2.5} 
If $v\in V_{\alpha}$, $\alpha\ge 0$, is primitive, then for 
every $0\le r\le \alpha$
\[ C^{-1}*L^{r}v=\tau (-1)^{m} \frac{r!}{(\alpha-r)!}L^{\alpha-r}v,\]
where $m=(h-\alpha)/2$.
\end{lem}

\begin{proof} Consider first the case $r=0$; writing 
$v=\sum a_{N}u_{N}$ with $|N|=m$, $a_{N}\in \C$, we have:
\[ \frac{L^{\alpha}v}{\alpha!}=\sum_{|N|=m}a_{N}
\sum_{\genfrac{}{}{0pt}{1}{N\subset M}{ 
|M|=m+\alpha}}u_{M}=
\sum_{|N|=m}a_{N}\sum_{\genfrac{}{}{0pt}{1}{M\subset N^{c}}{|M|=m}}
u_{M^{c}}=
\sum_{|M|=m}u_{M^{c}}
\sum_{\genfrac{}{}{0pt}{1}{N\subset M^{c}}{ 
|N|=m}}a_{N}.\]
\[C^{-1}*v=\tau\sum_{|M|=m}a_{M}u_{M^{c}}.\]
The equality $C^{-1}*v=\tau (\alpha!)^{-1}L^{\alpha}v$ follows 
immediately from Lemma~\ref{X.2.4}.
If $r\ge 1$ then 
\[ C^{-1}*L^{r}v=\Lambda^{r}C^{-1}*v=
\frac{\tau(-1)^{m}}{\alpha!}
\Lambda^{r}L^{\alpha}v.\]
Using the formula of \ref{X.2.3} we get 
\[ C^{-1}*L^{r}v=\frac{\tau(-1)^{m}}{\alpha!}
\prod_{j=\alpha-r+1}^{\alpha}j(\alpha-j+1)L^{\alpha-r}v=
\tau (-1)^{m} \frac{r!}{(\alpha-r)!}L^{\alpha-r}v.\]
\end{proof}

\begin{thm}\label{X.2.6} \emph{(Lefschetz decomposition)}
\begin{enumerate}
\item  Every $v\in V_{\alpha}$ can be written in a unique way 
as 
\[ v=\sum_{r\ge \max(-\alpha,0)}L^{r}v_{r}\]
with every $v_{r}\in V_{\alpha+2r}$ primitive.

\item  For a fixed $q\ge h$ there exist noncommutative 
polynomials $G^{q}_{\alpha,r}(\Lambda,L)$ with rational 
coefficients  such that $v_{r}=G^{q}_{\alpha,r}(\Lambda,L)v$ 
for every $v\in V_{\alpha}$.
\end{enumerate}
\end{thm}
\begin{proof} Assume first $\alpha\ge 0$, we prove the existence of 
the decomposition $v=\sum_{r\ge 0} L^{r}v_{r}$ as above by induction on 
the minimum $q$ such that $\Lambda^{q}v=0$. If $q=1$ then $v$ is 
already primitive.
If $\Lambda^{q+1}v=0$ then $w=\Lambda^{q}v\in V_{\alpha+2q}$ is 
primitive and then, setting  $\gamma=\prod_{r=1}^{q}r(\alpha+2q-r+1)$,
we have $\gamma>0$ and
\[\Lambda^{q}\left(v-L^{q}\frac{w}{\gamma}\right)=
w-\Lambda^{q}L^{q}\frac{w}{\gamma}=0.\]
This prove the existence when $\alpha\ge 0$. If $\alpha<0$ then 
$C^{-1}*v\in V_{-\alpha}$ and we can write:
\[ C^{-1}*v=\sum_{r\ge 0}L^{r}v_{r},\qquad 
v=\sum_{r\ge 0}C^{-1}*L^{r}v_{r},\quad v_{r}\in V_{-\alpha+2r}.\]
According to Lemma~\ref{X.2.5} 
\[ v=\sum_{r\ge 0}\gamma_{r}L^{-\alpha+r}v_{r}=
\sum_{r\ge -\alpha}\gamma_{r+\alpha}L^{r}v_{r}\]
for suitable rational coefficients $\gamma_{r}$.\\
The unicity of the decomposition and item 2 are proved at the same 
time. If
\[ v=\sum_{r=\max(-\alpha,0)}^{q}L^{r}v_{r}\]
is a decomposition with every $v_{r}\in V_{\alpha+2r}$ primitive,
then $L^{\alpha+q}v=L^{\alpha+2q}v_{q}$ and  
\[ v_{q}=
\frac{1}{(\alpha+2q)!^{2}}\Lambda^{\alpha+2q}L^{\alpha+2q}v_{q}=
\frac{1}{(\alpha+2q)!^{2}}\Lambda^{\alpha+2q}L^{\alpha+q}v.\]
Therefore $v_{q}$ is uniquely determined by $v$ and we can take
$G^{q}_{\alpha,q}=(\alpha+2q)!^{-2}\Lambda^{\alpha+2q}L^{\alpha+q}$.\\
Since $v-L^{q}v_{q}=
(1-L^{q}G^{q}_{\alpha,q})v=
\sum_{r=\max(-\alpha,0)}^{q-1}L^{r}v_{r}$ we can 
proceed by decreasing induction on $q$.\end{proof}

\begin{cor}\label{X.2.7} 
$v\in V_{\alpha}$, $\alpha\ge 0$, is primitive if and
only  if $L^{\alpha+1}v=0$.\end{cor}

\begin{proof} Let $v=\sum_{r\ge 0} L^{r}v_{r}$ be the 
Lefschetz decomposition of $v$,  then 
$\sum_{r>0}L^{\alpha+r+1}v_{r}$ is the Lefschetz decomposition of 
$L^{\alpha+1}v$. Therefore $L^{\alpha+1}v=0$ if and only if 
$v=v_{0}$.\end{proof}

It is clear that Theorem~\ref{X.2.6} and Corollary~\ref{X.2.7} hold also for 
every finite direct sum of $\Phi$-modules of type $V(h,\tau)$.

For later use we reinterpret Lemma~\ref{X.2.5} for the 
$\Phi$-module $\sA$: we have 
\[\sA=\somdir{A,B}{}V_{A,B},\qquad  
V_{A,B}=V\left(n-|A|-|B|,(-1)^{\frac{(|A|+|B|)(|A|+|B|+1)}{2}}\right)\]
where the sum is taken over all pairs of disjoint subsets $A,B$ of 
$\{1,\ldots,n\}$.
The space $\sA_{\alpha}=\bigoplus (V_{A,B})_{\alpha}$ 
is precisely the space $\bigoplus_{a}\sA^{a,n-\alpha-a}$ 
of $(n-\alpha)$-covectors. We then get the 
following 

\begin{lem}\label{X.2.8} 
If $v\in \sA$ is a primitive $p$-covector, $p\le n$, then 
\renewcommand\arraystretch{1}{2}
\[ C^{-1}*L^{r}v=
\left\{\begin{array}{ll}
\ds(-1)^{\frac{p(p+1)}{2}} 
\frac{r!}{(n-p-r)!}L^{n-p-r}v\quad &\hbox{ if }r\le n-p\\
0&\hbox{ if } r>n-p\end{array}\right.\]
\renewcommand\arraystretch{1}
\end{lem}

\bigskip

\section{K\"{a}hler identities}

Let $M$ be a complex manifold of dimension $n$ 
and denote by $\sA^{*,*}$ the sheaf of differential forms on $M$. 
By definition $\sA^{a,b}$ is the sheaf of sections of the complex 
vector bundle 
$\bigwedge^{a}T_{M}\dual\otimes\bigwedge^{b}\bar{T_{M}}\dual$. 
The operators $P_{a,b}$, $P_{p}$ 
and $C$, defined on the fibres of the above  bundles, extend in 
the obvious way to operators in the sheaf $\sA^{*,*}$.   

If $d\colon \sA^{*,*}\to \sA^{*,*}$ 
is the De Rham differential we denote: 
\[ d^{C}=C^{-1}dC,\qquad \de=\frac{d+id^{C}}{2},\qquad 
\debar=\frac{d-id^{C}}{2},\] 
\[d=Cd^{C}C^{-1},\qquad    d=\de+\debar,\qquad d^{C}=i(\debar-\de).\]
If $\eta$ is a $(p,q)$-form then we can write 
$d\eta=\eta'+\eta''$ with $\eta'\in\sA^{p+1,q}$, 
$\eta''\in\sA^{p,q+1}$ and then 
\[ d^{C}(\eta)=C^{-1}d(i^{p-q}\eta)=\frac{i^{p-q}}{i^{p-q+1}}\eta'+
\frac{i^{p-q}}{i^{p-q-1}}\eta''=i^{-1}\eta'+i\eta'',\qquad
\de\eta=\eta',\quad \debar\eta=\eta''.\]
Since $0=d^{2}=\de^{2}+\de\debar+\debar\de+\debar^{2}$ we get 
$0=\de^{2}=\de\debar+\debar\de=\debar^{2}$ and then 
$(d^{C})^{2}=0$, $dd^{C}=2i\de\debar=-d^{C}d$.

Using the structure of graded Lie algebra on the space of 
$\C$-linear operators of the sheaf of graded algebras $\sA^{*,*}$
(with the total degree $\bar{v}=a+b$ if $v\in\sA^{a,b}$), 
the above relation can be rewritten as 
\[ [d,d]=dd+dd=2d^{2}=0,\qquad [d^{C},d^{C}]=[d,d^{C}]=
[\de,\de]=[\debar,\debar]=[\de,\debar]=0.\]
Note finally that $d$ and $C$ are  real operators and then also 
$d^{C}$ is; moreover $\debar{\bar{\eta}}=\bar{\de\eta}$. 

A Hermitian metric on $M$ is a positive definite Hermitian form $h$ 
on the tangent vector bundle $T_{M}$. If $z_{1},\ldots,z_{n}$ are 
local holomorphic coordinates then $h_{ij}=
h\left(\desude{~}{z_{i}},\desude{~}{z_{j}}\right)$ is a
smooth function and the matrix $(h_{ij})$ is 
Hermitian and positive definite. The local expression of $h$ is then 
$h=\sum_{ij}h_{ij}dz_{i}\otimes d\bar{z}_{j}$ and the  differential 
form 
\[\omega=\frac{i}{2}\sum_{i,j}h_{ij}dz_{i}\wedge d\bar{z}_{j}
\in \Gamma(M,\sA^{1,1})\]
is globally definite and gives the imaginary part of $-h$;
$\omega$ is called the (real, $(1,1)$) associated form to $h$.\\

The choice of a Hermitian metric on $M$ induces, for every open 
subset $U\subset M$, linear operators 
\[ L\colon \Gamma(U,\sA^{a,b})\to \Gamma(U,\sA^{a+1,b+1}),\qquad 
Lv=v\wedge \omega,\]
\[ *\colon \Gamma(U,\sA^{a,b})\to \Gamma(U,\sA^{n-b,n-a}),\]
\[ \Lambda\colon \Gamma(U,\sA^{a,b})\to \Gamma(U,\sA^{a-1,b-1}),
\qquad \Lambda=*^{-1}L*=(C^{-1}*)^{-1}LC^{-1}*.\]

The commuting  relations between them  
\[ [L,C]=[\Lambda,C]=[*,C]=[L,*^{2}]=0,\qquad 
[\Lambda,L^{r}]=\sum_{p}r(n-p-r+1){P_{p}}\]
are still valid.\\
A  differential form $v$ is primitive if $\Lambda v=0$; the 
existence of the polynomials $G^{n}_{n-p,r}(\Lambda,L)$ 
(cf. Theorem~\ref{X.2.6}) 
gives the existence and unicity of 
Lefschetz decomposition for every differential $p$-form 
\[ v=\sum_{r\ge\max(p-n,0)}L^{r}v_{r},\qquad \Lambda v_{r}=0.\]

We set:
\[ \delta=-*d*,\qquad \delta^{C}=-*d^{C}*=C^{-1}\delta C,\]
\[ \de^{*}=-*\debar *=\frac{\delta-i\delta^{C}}{2},\qquad 
\debar^{*}=-*\de *=\frac{\delta+i\delta^{C}}{2}.\]

\begin{defn}\label{X.3.1} The Hermitian metric $h$ is called a \emph{K\"{a}hler 
metric} if $d\omega=0$.\end{defn}

Almost all the good properties of K\"ahler metrics come from the 
following 
\begin{thm}\label{X.3.4} \emph{(K\"{a}hler identities)}
Let $h$ be a K\"{a}hler metric on a complex manifold, 
then:
\renewcommand\arraystretch{2}
\[ 
\begin{tabular}{|l|l|l|l|}\hline
$[L,d]=0$&$[L,d^{C}]=0$&$[L,\de]=0$&$[L,\debar]=0$\\ \hline
$[\Lambda,d]=-\delta^{C}$&$[\Lambda,d^{C}]=\delta$&
$[\Lambda,\de]=i\debar^{*}$&$[\Lambda,\debar]=-i\de^{*}$\\ \hline
$[L,\delta]=d^{C}$&$[L,\delta^{C}]=-d$&
$[L,\de^{*}]=i\debar$&$[L,\debar^{*}]=-i\de$\\ \hline
$[\Lambda,\delta]=0$&$[\Lambda,\delta^{C}]=0$&
$[\Lambda,\de^{*}]=0$&$[\Lambda,\debar^{*}]=0$\\ \hline
\end{tabular}\]
\renewcommand\arraystretch{1}
\end{thm}

\begin{proof} It is sufficient to prove that 
$[L,d]=0$ and $[\Lambda,d]=-\delta^{C}$. 
In fact, since $\Lambda=*^{-1}L*=*L*^{-1}$ we have 
$[\Lambda,\delta]+*[L,d]*=0$ and $[L,\delta]+*[\Lambda,d]*=0$: 
this will prove the first column. The second column follows from 
the first using the fact that $C$ commutes 
with $L$ and $\Lambda$. The last two columns are linear
combinations of the first two.\\
If $v$ is a $p$-form then, since $d\omega=0$,  
\[ [L,d]v=dv\wedge \omega-d(v\wedge \omega)=-(-1)^{p}v\wedge 
d\omega=0.\]
According to the Lefschetz decomposition it is sufficient to prove 
that  
$[\Lambda,d]L^{r}u=-\delta^{C}L^{r}u$ for every $r\ge 0$ and every 
primitive $p$-form $u$ ($p\le n$). 
We first note that, being $u$ primitive, 
$L^{n-p+1}u=0$ and then $L^{n-p+1}du=dL^{n-p+1}u=0$. 
This implies that the Lefschetz 
decomposition of $du$ is $du=u_{0}+Lu_{1}$.\\
Setting $\alpha=n-p$, we have $u\in V_{\alpha}$, $u_{0}\in 
V_{\alpha-1}$, $u_{1}\in V_{\alpha+1}$: 
\[ [\Lambda,d]L^{r}u=\Lambda L^{r}du-d\Lambda L^{r}u=
\Lambda L^{r}u_{0}+\Lambda L^{r+1}u_{1}- r(\alpha-r+1)dL^{r-1}u=\]
\[=r(\alpha-r)L^{r-1}u_{0}+(r+1)(\alpha-r+1)L^{r}u_{1}-
r(\alpha-r+1)L^{r-1}u_{0}-r(\alpha-r+1)L^{r}u_{1}=\]
\[=-rL^{r-1}u_{0}+(\alpha-r+1)L^{r}u_{1}.\]
On the other hand we have by \ref{X.2.8}
\[ \begin{split}
-\delta^{C}L^{r}u=C^{-1}*d*CL^{r}u&=
C^{-1}*dC^{2}C^{-1}*L^{r}u\\
&=C^{-1}*dC^{2}
(-1)^{p(p+1)/2}\frac{r!}{(\alpha-r)!}L^{\alpha-r}u
\end{split}\]
and then 
\[ -\delta^{C}L^{r}u=(-1)^{p(p-1)/2}
\frac{r!}{(\alpha-r)!}C^{-1}*L^{\alpha-r}(u_{0}+Lu_{1}).\]
Again by \ref{X.2.8},
\[ C^{-1}*L^{\alpha-r}u_{0}=(-1)^{(p+1)(p+2)/2}
\frac{(\alpha-r)!}{(r-1)!}L^{r-1}u_{0},\]
\[C^{-1}*L^{\alpha-r+1}u_{1}=(-1)^{(p-1)p/2}
\frac{(\alpha-r+1)!}{r!}L^{r}u_{1}.\]
Putting all the terms together we obtain the  result.
\end{proof}

\begin{cor}\label{X.3.5} If $\omega$ is the associated form of a K\"ahler metric 
$h$ then $d\omega^{\wedge p}=\delta\omega^{\wedge p}=0$ 
for every $p\ge 0$ .\end{cor}

\begin{proof} 
The equality $d\omega^{\wedge p}=0$ follows immediately from 
the Leibnitz rule. Since $\omega^{\wedge p}$ is a $(p,p)$ form, we 
have $C\omega^{\wedge p}=\omega^{\wedge p}$ and then also 
$d^{C}\omega^{\wedge p}=0$.\\
We prove $\delta\omega^{\wedge p}=0$ by induction on $p$, being the 
result trivial when $p=0$. If $p>0$ we have 
\[ 0=d^{C}\omega^{\wedge p-1}=L\delta \omega^{\wedge p-1}-
\delta L\omega^{\wedge p-1}=-\delta\omega^{\wedge p}.\]
\end{proof}

The gang of Laplacians is composed by:\begin{enumerate}
\item $\Delta_{d}=\Delta=[d,\delta]=d\delta+\delta d$.
\item $\Delta_{d^{C}}=\Delta^{C}=C^{-1}\Delta C=
[d^{C},\delta^{C}]=d^{C}\delta^{C}+\delta^{C}d^{C}$.
\item $\Delta_{\de}=\square=[\de,\de^{*}]=\de\de^{*}+\de^{*}\de$.
\item $\Delta_{\debar}=\bar{\square}=[\debar,\debar^{*}]=
\debar\,\debar^{*}+\debar^{*}\,\debar$.
\end{enumerate}
A straightforward computation shows that 
$\Delta+\Delta^{C}=2\square+2\bar{\square}$.

\begin{cor}\label{X.3.7} In the above notation, 
if $h$ is a K\"{a}hler metric then:
\[ [d,\delta^{C}]=[d^{C},\delta]=[\de,\debar^{*}]=[\debar,\de^{*}]=0,
\qquad
\frac{1}{2}\Delta=\frac{1}{2}\Delta^{C}=\square=\bar{\square}.\] 
In particular $\Delta$ is bihomogeneous of degree $(0,0)$.
\end{cor}

\begin{proof} 
According to Theorem~\ref{X.3.4} and the Jacobi identity we have 
\[ [d,\delta^{C}]=[d,[d,\Lambda]]=\frac{1}{2}[[d,d],\Lambda]=0.\]
The proof of $[d^{C},\delta]=[\de,\debar^{*}]=[\debar,\de^{*}]=0$ 
is similar and left as exercise.
For the equalities among Laplacians 
it is sufficient to shows that $\Delta=\Delta^{C}$ and
$\square=\bar{\square}$. According to the K\"{a}hler identities 
\[ \Delta=[d,\delta]=[d,[\Lambda,d^{C}]]=[[d,\Lambda],d^{C}]+
[\Lambda,[d,d^{C}]].\]
Since $[d,d^{C}]=dd^{C}+d^{C}d=0$ we have 
\[\Delta=[d,\delta]=[[d,\Lambda],d^{C}]=
[\delta^{C},d^{C}]=\Delta^{C}.\]
The proof of $\square=\bar{\square}$ is similar and it is left to 
the reader.\end{proof}

\begin{cor}\label{X.3.8} 
In the above notation, if $h$ is a K\"{a}hler metric, then
$\Delta$  commutes with all the operators $P_{a,b}$, $*$, $d$, $L$,  
$C$, $\Lambda$, 
$d^{C}$, $\de$, $\debar$, $\delta$, $\delta^{C}$, $\de^{*}$, 
$\debar^{*}$.\end{cor}

\begin{proof}
Since $\Delta$ is of type $(0,0)$ it is clear that commutes with 
the projections $P_{a,b}$. Recalling that $\delta=-*d*$ we get 
$d=*\delta*$ and then 
\[ *\Delta=*d\delta+*\delta d=-*d*d*+*\delta*\delta*=
\delta d*+d\delta *=\Delta *.\]
\[ [L,\Delta]=[L,[d,\delta]]=[[L,d],\delta]+[[L,\delta],d]=
[d^{C},d]=0.\]
\[ [d,\Delta]=[d,[d,\delta]]=\frac{1}{2}[[d,d],\delta]=0.\]
Now it is sufficient to observe that all the operators in the 
statement belong to the 
$\C$-algebra generated by $P_{a,b}$, $*$, $d$ 
and $L$.\end{proof}

\begin{defn}\label{X.3.6} A $p$-form $v$ 
is called \emph{harmonic} if $\Delta v=0$.\end{defn}

\begin{cor}\label{X.3.9} Let $h$ be a K\"{a}hler metric and let 
$v=\sum_{r}L^{r}v_{r}$ be the Lefschetz decomposition of a $p$-form.\\
Then $v$ is harmonic if and only if $v_{r}$ is harmonic for every 
$r$.\end{cor}

\begin{proof} Since $\Delta$ commutes with $L$, if $\Delta v_{r}=0$ for 
every $r$ then also $\Delta v=0$. Conversely, 
since $v_{r}=G^{n}_{p,r}(\Lambda,L)v$ for suitable 
noncommutative polynomials with rational coefficients $G^{n}_{p,r}$, 
and $\Delta$ commutes with $\Lambda,L$ then $v$ harmonic implies 
$\Delta v_{r}=0$ for every $r$.\end{proof}

\begin{cor}\label{X.3.10} 
In the above notation, if $h$ is a K\"{a}hler metric and $v$ is 
a closed primitive $(p,q)$-form then $v$ is harmonic. 
\end{cor}

Note that if either $p=0$ or $q=0$ then $v$ is always primitive.

\begin{proof} It is sufficient to prove that $\delta v=0$, we have 
\[ \delta v=C\delta^{C}C^{-1}v=i^{q-p}C\delta^{C}v=
i^{q-p}C[d,\Lambda]v=0.\]
\end{proof}

\bigskip

\section{K\"{a}hler metrics on compact manifolds}

In this section we assume $M$ compact complex manifold of dimension $n$. 
We denote by $L^{a,b}=\Gamma(M,\sA^{a,b})$, 
$L^{p}=\bigoplus_{a+b=p}L^{a,b}$, $L=\bigoplus_{p}L^{p}$.\\
Every Hermitian metric $h$ on $M$  induces a 
structure of pre-Hilbert space on $L^{a,b}$ for every 
$a,b$ (and then also on $L$) given by:
\[ (\phi,\psi)=\int_{M}\phi\wedge\bar{*\psi}.\]
We have already seen that the operator $*\colon L^{a,b}
\to L^{n-a,n-b}$ is an isometry commuting with the complex 
conjugation and then we also have: 
\[ (\phi,\psi)=\int_{M}\phi\wedge\bar{*\psi}=\int_{M}\phi\wedge*\bar{\psi}=
(-1)^{a+b}\int_{M}*\phi\wedge\bar{\psi}=
\int_{M}\bar{\psi}\wedge * \phi=\bar{(\psi,\phi)}.\]

\begin{prop}\label{X.4.1} With respect to the above pre-Hilbert structures we have 
the following pairs (written in columns) of formally adjoint operator:\\
\begin{center}\renewcommand\arraystretch{2}
\begin{tabular}{|l|l|l|l|l|l|}\hline
    ~operator~&~$d$~&~$d^{C}$~&~$\de$~&~$\debar$~&~$L$~\\ \hline
~formal adjoint~~& ~$\delta$~ &~$\delta^{C}$~&~$\de^{*}$~&
~$\debar^{*}$~&~$\Lambda$~\\ \hline
\end{tabular}
\renewcommand\arraystretch{1}
\end{center}
~\\
In particular, all the four Laplacians are formally self-adjoint 
operators.\end{prop}

\begin{proof} We show here only that  $\delta$ is the formal adjoint of 
$d$. The proof of the remaining assertions is essentially the same 
and it is left as exercise.\\
Let $\phi$ be a $p$-form and $\psi$ a $p+1$-form. By Stokes theorem 
\[ 0=\int_{M}d(\phi\wedge \bar{*\psi})=\int_{M}d\phi\wedge\bar{*\psi}+
(-1)^{p}\int_{M}\phi\wedge d\bar{*\psi}.\]
Since $d\bar{*\psi}=\bar{d*\psi}$ and 
$d*\psi=(-1)^{2n-p}*^{2}d*\psi=-(-1)^{p}*\delta\psi$ we get 
\[0=\int_{M}d\phi\wedge\bar{*\psi}-
\int_{M}\phi\wedge \bar{*\delta\psi}=(d\phi,\psi)-(\phi,\delta\psi).\]
\end{proof}

Let $D$ be any of the operator $d,d^{C},\de,\debar$; denote $D^{*}$ its 
formal adjoint and by 
$\Delta_{D}=DD^{*}+D^{*}D$ its Laplacian (i.e. $\Delta_{d}=\Delta$, 
$\Delta_{\debar}=\bar{\square}$ etc...).
The space of $D$-harmonic $p$-forms is denoted by 
$\sH^{p}_{D}=\ker\Delta_{D}\cap L^{p}$.

\begin{lem}\label{X.4.2} 
We have $\ker\Delta_{D}=\ker D\cap \ker D^{*}$.    
\end{lem}
\begin{proof} The inclusion $\supset$ is immediate from the 
definitions of the Laplacian. The inclusion $\subset$ 
comes from   
\[ (\Delta_{D} \phi,\phi)=(DD^{*}\phi,\phi)+(D^{*}D\phi,\phi)=
(D^{*}\phi,D^{*}\phi)+(D\phi,D\phi)=\|D^{*}\phi\|^{2}+\|D\phi\|^{2}.\]
\end{proof}

The theory of elliptic self-adjoint operators on  compact manifolds 
gives:
\begin{thm}\label{X.4.3} 
In the notation above the spaces of $D$-harmonic forms 
$\sH^{p}_{D}$ are finite dimensional and there exist
orthogonal decompositions
\[ L^{p}=\sH^{p}_{D}\somdir{}{\perp}\Image\Delta_{D}.\]
\end{thm}

\begin{proof} See e.g. \cite{Voisin}.\end{proof}

\begin{cor}\label{X.4.4} The natural projection maps 
\[ \sH^{p}_{d}\to H^{p}(M,\C),\qquad\sH^{p,q}_{\debar}\to 
H^{q}_{\debar}(M,\Omega^{p})\]
are isomorphism.\end{cor}

\begin{proof} We first note that, according to Lemma~\ref{X.4.2}, every 
harmonic form is closed and then the above projection maps makes sense. 
It is evident that $\Image\Delta\subset \Image d+\Image\delta$.
On the other hand, since $d,\delta$ are formally adjoint and 
$d^{2}=\delta^{2}=0$
we have   
$\ker d\perp\Image\delta$, $\ker\delta\perp \Image d$: this  implies 
that $\Image d$, $\Image\delta$ and $\sH^{p}_{d}$ are pairwise orthogonal. 
Therefore  $\Image\Delta=\Image d\oplus\Image\delta$ and $\ker 
d=\sH^{p}_{d}\oplus \Image d$; the conclusion follows by De Rham 
theorem.\\ 
The isomorphism $\sH^{p,q}_{\debar}\to H^{q}_{\debar}(M,\Omega^{p})$ is 
proved in the same way (with Dolbeault's theorem instead of De Rham)
and it is left as exercise.\end{proof}

\begin{cor}\label{X.4.7} The map $\Delta_{D}\colon\Image\Delta_{D}\to
\Image\Delta_{D}$ is  bijective.\end{cor}
\begin{proof} Trivial consequence of Theorem~\ref{X.4.3}.\end{proof}

We define the \emph{harmonic 
projection} $H_{D}\colon L^{p}\to\sH^{p}_{D}$ as the orthogonal 
projection and the \emph{Green operator} $G_{D}\colon 
L^{p}\to\Image\Delta_{D}$ as the composition of 
\[ G_{D}\colon L^{p}\xrightarrow{Id-H_D}
\Image\Delta_{D}
\xrightarrow{\Delta_{D}^{-1}}\Image\Delta_{D}.\]

Note that $\Delta_{D} G_{D}=G_{D}\Delta_{D}=Id-H_{D}$ 
and $G_{D}H_{D}=H_{D}G_{D}=0$.\\

\begin{lem}\label{X.4.8} 
If $K$ is an operator commuting with $\Delta_{D}$ then 
$K$ commutes with $G_{D}$.\end{lem}

\begin{proof} Exercise (Hint: $K$ preserves image and kernel of
$\Delta_{D}$).
\end{proof}

If $h$ is a K\"ahler metric, then the equality $\Delta=2\bar{\square}$
implies that 
\[ H_{d}=H_{d^{C}}=H_{\de}=H_{\debar},\qquad 
G_{d}=G_{d^{C}}=\frac{1}{2}G_{\de}=\frac{1}{2}G_{\debar}.\]
In particular, according to Lemma~\ref{X.4.8} and Corollary~\ref{X.3.8}, 
$G_{d}=G_{d^{C}}$ commutes with $d,d^{C}$.\\

\begin{cor}\label{X.4.holo} 
If $h$ is a  K\"ahler metric on a compact manifold then:
Every holomorphic $p$-form on $M$ is harmonic. 
\end{cor}

\begin{proof} 
According to  Corollary~\ref{X.4.4} the inclusion 
$\sH^{p,0}_{\debar}\subset \Gamma(M,\Omega^{p})$ is an isomorphism 
and then if $\eta$ is a holomorphic $p$-form we have 
$\Delta(\eta)=2\bar{\square}(\eta)=0$.
\end{proof}

\begin{exer} Let $v\not=0$ be a primitive $(p,q)$-form on a compact 
manifold $M$ with K\"ahler form $\omega$. Prove that 
\[\int_{M}v\wedge\bar{v}\wedge\omega^{\wedge n-p-q}\not=0.\]
\end{exer}

\bigskip
\section{Compact K\"{a}hler manifolds}

In this section we will prove that certain good properties concerning the 
topology and the complex structure of compact complex manifolds 
are true whenever we assume the existence of a K\"ahler metric. 
This motivates the following definition: 

\begin{defn}\label{X.3.2} A complex manifold $M$ s called a 
\emph{K\"ahler manifolds} if there exists a K\"ahler metric on 
$M$.\end{defn}

We note that, while every complex manifold admits a Hermitian 
metric (this is an easy application of partitions of unity, cf. 
\cite[Thm. 3.14]{Kobook}), not every complex manifold is 
K\"ahlerian. We recall the following

\begin{thm}\label{X.3.3} 
\begin{enumerate}
    \item  $\C^{n}$, $\proj^{n}$ and the complex tori
    are K\"ahler manifolds.

    \item  If $M$ is a K\"ahler manifold and $N\subset M$ is a 
    regular submanifold then also $N$ is a K\"ahler manifolds.
\end{enumerate}\end{thm}

For a proof of Theorem ~\ref{X.3.3} we refer to \cite{G-H}.\\

From now on $M$ is a fixed compact K\"ahler manifold on dimension 
$n$.\\
For every $m\le 2n$ we denote by 
$H^{m}(M,\C)=H^{m}(M,\R)\otimes_{\R}\C$ the De Rham cohomology 
$\C$-vector spaces. We note that a differential $m$-form $\eta$ is 
$d$-closed if and only if its conjugate $\bar{\eta}$ is. In 
particular the complex conjugation induce an isomorphism of vector 
spaces $H^{m}(M,\C)=\bar{H^{m}(M,\C)}$. 

If $p+q=m$ we denote by $F^{p,q}\subset H^{m}(M,\C)$ the subspace of 
cohomology classes represented by $d$-closed form of type $(p,q)$ 
(note that a $(p,q)$-form $\eta$ is $d$-closed if and only if 
it is $\de\eta=\debar\eta=0$).
It is clear that $\bar{F^{p,q}}=F^{q,p}$.\\

\begin{thm}[Hodge decomposition]
In the notation above we have
\[H^{m}(M,\C)=\somdir{p+q=m}{}F^{p,q}\]
and the natural morphisms 
$F^{p,q}\to H^{p,q}_{\de}(M)$, $F^{p,q}\to H^{p,q}_{\debar}(M)$
are isomorphisms.
\end{thm}

\begin{proof} Take a K\"ahler metric on $M$ and use it to define the 
four Laplacians,  the harmonic projectors and the Green operators.  
According to Corollary~\ref{X.3.7} 
the Laplacian $\Delta$ is bihomogeneous of bidegree $(0,0)$ and we have 
\[\ker\Delta\cap L^{q}=\somdir{a+b=q}{}\ker\Delta\cap L^{a,b}.\]
The isomorphism $\ker\Delta\cap L^{q}\to H^{q}(M,\C)$ induces
injective maps $\ker\Delta\cap L^{a,b}\to F^{a,b}$; 
this maps are also surjective because every closed form $\alpha$ is 
cohomologically equivalent to its harmonic projection $H\alpha$ and 
$H$ is bihomogeneous of bidegree $(0,0)$.\\ 
The last equalities follow from the isomorphisms 
\[ \ker\Delta\cap L^{a,b}=\ker{\square}\cap 
L^{a,b}=H^{a,b}_{\de}(M),\quad
\ker\Delta\cap L^{a,b}=\ker\bar{\square}\cap 
L^{a,b}=H^{a,b}_{\debar}(M).\]
\end{proof}

\begin{cor}\label{X.4.5} 
If $M$ is a compact K\"ahler manifold then:
\begin{enumerate}
\item $b_{i}=\sum_{a+b=i}h^{a,b}$.

\item  $h^{p,q}=h^{q,p}$, in particular 
$b_{i}$ is even if $i$ is odd.

\item  $h^{p,p}>0$, in particular $b_{i}>0$ if $i$ is even.

\item Every holomorphic $p$-form on $M$ is 
$d$-closed. 

\end{enumerate}
($b_{i}=\dim_{\C}H^{i}(M,\C)$ are the  Betti numbers, 
$h^{p,q}=\dim_{\C}H^{q}(M,\Omega^{p})$ the Hodge numbers.)
\end{cor}

\begin{proof} 
Items 1 and 2 are immediate consequence of the Hodge decomposition.
Take a K\"ahler metric on $M$ and use it to define the 
four Laplacians,  the harmonic projectors and the Green operators.      
Let $\omega$ be the associated form of the K\"ahler metric 
on $M$. According to Corollary~\ref{X.3.5}, $\omega^{\wedge p}$ is 
harmonic and then $\ker\bar{\square}\cap L^{p,p}=
\ker\Delta\cap L^{p,p}\not=0$.\\
Finally, by Corollary~\ref{X.4.holo} the holomorphic forms are 
$\Delta$-harmonic and therefore $d$-closed. 
\end{proof}

\begin{ex}\label{X.4.6} The Hopf surfaces (Example~\ref{I.1.6}) have 
$b_{1}=b_{3}=1$, $b_{2}=0$ and then are not K\"ahler.
\end{ex}

Finally we are in a position to prove the following 
\begin{thm}\label{X.4.9} \emph{($\de\debar$-Lemma)} Let $M$ be a
compact K\"ahler manifold.  Then
\begin{enumerate}
\item There exists a linear operator $\sigma\colon L\to L$ of bidegree $(0,-1)$ 
such that 
\[ [\de,\sigma]=0,\qquad [\debar,\sigma]\de=[\debar,\sigma\de]=\de.\]
\item  $\Image \de\debar=\ker \de\cap \Image \debar=\ker \debar\cap \Image 
\de$.
\end{enumerate}\end{thm}

\begin{proof} 
{[1]} Choose a K\"ahler metric and 
define $\sigma=G_{\debar}\debar^{*}$. According to \ref{X.3.7}, 
\ref{X.3.8} and \ref{X.4.8} we have $\sigma=\debar^{*}G_{\debar}$, 
$[\de,\sigma]=0$ and, denoting by $H$ the harmonic projection, 
\[
[\debar,\sigma]\de=G_{\debar}\Delta_{\debar}\de=(Id-H)\de=\de.\]
{[2]} (cf. Exercise \ref{esededebar}) We prove only 
$\Image \de\debar=\ker \debar\cap \Image 
\de$, being the other equality the conjugate of this one.
The inclusion $\subset$ is evident, conversely 
let $x=\de\alpha$ be a $\debar$-closed differential form; 
we can write 
\[ 
x=\de\alpha=[\debar,\sigma]\de\alpha=\debar\sigma\de\alpha+\sigma\debar\de\alpha=
-\debar\de\sigma\alpha-\sigma\debar x=\de\debar(\sigma\alpha).\]
\end{proof}

\begin{cor}\label{X.4.10} Let $M$ be a compact K\"ahler manifold. Then for
every 
$p,q$ the natural maps
\[\frac{\ker\de\cap\ker\debar\cap L^{p,q}}{\de\debar L^{p-1,q-1}}\to
\frac{\ker\de\cap\ker\debar\cap L^{p,q}}{\debar(\ker\de\cap L^{p,q-1})}\to
\frac{\ker\debar\cap L^{p,q}}{\debar L^{p,q-1}}=H^{q}(M,\Omega^{p})\]
\[\frac{\ker\de\cap\ker\debar\cap L^{p,q}}{\de\debar L^{p-1,q-1}}\to
\frac{\ker\de\cap\ker\debar\cap L^{p,q}}{\de(\ker\debar\cap L^{p-1,q})}
\to\frac{\ker\de\cap L^{p,q}}{\de L^{p-1,q}}\]
are isomophisms.\end{cor}

\begin{proof} 
The two lines are conjugates each other and then it is sufficient 
to prove that the maps on the first row are isomorphisms.\\
Choose a K\"ahler metric, every $\debar$-closed form $\phi$ 
can be written as 
$\phi=\alpha+\debar\psi$ with $\bar{\square}\alpha=0$. 
Since $\square=\bar{\square}$ we have $\de\alpha=0$ and then the 
above maps are surjective.\\
According to Theorem~\ref{X.4.9} we have 
\[\de\debar(L^{p-1,q-1})\subset \debar(\ker\de\cap L^{p,q-1})\subset
\ker\de\cap \debar(L^{p,q-1})\subset \de\debar(L^{p-1,q-1})\]
and then all the  maps  are injective.\end{proof}

\begin{exer}\label{esededebar}
Prove that for a double complex $(L^{*,*},d,\delta)$ of vector 
spaces (with  $d,\delta$ differentials of respective bidegrees 
$(1,0)$ and $(0,1)$) 
the following conditions are equivalent:
\begin{enumerate}
    \item  There exists a linear operator $\sigma\colon L^{*,*}\to 
    L^{*,*-1}$ of bidegree $(0,-1)$ such that 
\[ [d,\sigma]=0,\qquad [\delta,\sigma]d=[\delta,\sigma d]=d.\]

\item  $\Image d\delta=\ker \delta\cap \Image d$.

\end{enumerate}
(Hint: The implication \implica{1}{2} is the same as in Theorem~\ref{X.4.9}. 
In order to prove \implica{2}{1} write $L^{a,b}=F^{a,b}\oplus C^{a,b}$ 
with $F^{a,b}=dL^{a-1,b}$ and observe that the complexes 
$(F^{a,*},\delta)$ are acyclic. Define first $\sigma\colon F^{a,b}\to 
F^{a,b-1}$ such that $[\delta,\sigma]d=d$ and then 
$\sigma\colon C^{a,b}\to 
C^{a,b-1}$ such that $[d,\sigma]=0$.)
\end{exer}

\bigskip
\section{Historical survey,~\ref{CAP:KAEHLER}}

Most of the properties of K\"ahler manifolds are stable under 
deformation. For example:

\begin{thm} \label{X.6.4}
Let $f\colon M\to B$ be a family of compact complex manifolds and 
assume that $M_{b}$ is K\"ahlerian for some $b\in B$.\\
Then there exists an open neighbourhood $b\in U\subset B$ such 
the functions $h^{p,q}\colon U\to \N$, 
$h^{p,q}(u)=\dim_{\C}H^{p,q}(M_{u})$ are constant 
and $\sum_{p+q=m}h^{p,q}(u)=\dim_{\C}H^{m}(M_{u},\C)$ for every 
$u\in U$.\end{thm} 

\begin{proof} (Idea) Exercise~\ref{Hodgeinequality}
implies
$\sum_{p+q=m}h^{p,q}(u)\ge \dim_{\C}H^{m}(M_{u},\C)$ and the 
equality holds whenever $M_{u}$ is K\"ahlerian.
On the other side, by semicontinuity theorem \ref{I.4.5} 
the functions $h^{p,q}$ are semicontinuous and by Ehresmann's theorem 
the function $u\mapsto \dim_{\C}H^{m}(M_{u},\C)$ is locally 
constant.\end{proof}

Theorem \ref{X.6.4} is one of the main ingredients for the proof of 
the following theorem, proved by Kodaira (cf. \cite{Kobook}, 
\cite{Voisin}) 

\begin{thm} \label{X.6.5}
Let $f\colon M\to B$ be a family of compact complex manifolds. 
Then the subset $\{ b\in B\mid M_{b} \hbox{ is K\"ahlerian }\}$ is 
open in $B$.\end{thm}

The proof of \ref{X.6.5} requires hard functional  and harmonic analysis. 
\bigskip

It seems that the name \emph{K\"ahler manifolds} comes from the fact that they were 
defined in a note of Erich K\"ahler (1906-2000) of 1933 but all their (first)
good properties were estabilished by W.V.D. Hodge some years later.\\


\chapter[~Deformations of manifolds with trivial $K$]
{Deformations of manifolds with
trivial canonical bundle}
\label{CAP:TRIVIALK}

\piede
\setcounter{thm}{0}
In the first part of this chapter we prove, following \cite{GoMil2} 
and assuming Kuranishi theorem~\ref{VI.5.1},
the following

\begin{thm}[Bogomolov-Tian-Todorov]\label{XII.0.1}
Let $M$ be a compact K\"ahler manifold with trivial canonical bundle
$K_{M}=\Oh_{M}$. Then $M$ admits a semiuniversal deformation with
smooth base $(H^{1}(M,T_{M}),0)$.\end{thm}

According to Corollary~\ref{VI.5.2}, it is sufficient to to show
that the natural map
\[\Def_{M}\left(\frac{\C[t]}{(t^{n+1})}\right)\to
\Def_{M}\left(\frac{\C[t]}{(t^{2})}\right)\]
is surjective for every $n\ge 1$. This will be done using
Corollary~\ref{VIII.4.4} and the so called
Tian-Todorov's lemma.\\
A generalization of this theorem has been given recently by H. 
Clemens \cite{clemens}. We will prove of Clemens' theorem in 
Chapter~\ref{CAP:TOOLS}.\\

In the second part we introduce  some interesting  classes of 
dg-algebras which  arise naturally both in mathematics and in 
physics: in particular we introduce the notion of differential 
Gerstenhaber algebra and differential 
Gerstenhaber-Batalin-Vilkovisky algebra. Then we show 
(Example~\ref{XII.5.7}) that the algebra of 
polyvector fields on a  manifold 
with trivial canonical bundle 
carries the structure of differential 
Gerstenhaber-Batalin-Vilkovisky algebra.

\bigskip

\section{Contraction on exterior  algebras}
\label{sezioneVI.1a}

Let $\K$ be a fixed field and
$E$ a vector space over $\K$ of dimension $n$;
denote by
$E\dual$ its dual and by $\langle,\rangle\colon E\times E\dual\to\C$ the natural
pairing. Given $v\in E$, the (left) contraction by $v$ is the linear
operator $v\vdash\colon\bigwedge^{b}E\dual\to \bigwedge^{b-1}E\dual$
defined by the formula
\[v\vdash(z_{1}\wedge\ldots\wedge
z_{b})=\sum_{i=1}^{b}(-1)^{i-1}
\langle v,z_{i}\rangle
z_{1}\wedge\ldots\wedge\widehat{z_{i}}\wedge\ldots\wedge
z_{b}.\]
For every $a\le b$  the contraction
\[ \external{}{a}E\times\external{}{b}E\dual
\mapor{\vdash}\external{}{b-a}E\dual\]
is the bilinear extension of
\[\begin{split}
(v_{a}\wedge\ldots\wedge v_{1})\vdash(z_{1}\wedge\ldots\wedge
z_{b})&
=v_{a}\vdash((v_{a-1}\wedge\ldots\wedge
v_{1})\vdash(z_{1}\wedge\ldots\wedge  z_{b}))\\
&=\ds\sum_{\sigma\in G}(-1)^{\sigma}
\left(\prod_{i=1}^{a}\langle v_{i},z_{\sigma(i)}\rangle
\right)z_{\sigma(a+1)}\wedge\ldots\wedge
z_{\sigma(b)}
\end{split}\]
where $G\subset \Sigma_{b}$ is the subset of permutations $\sigma$
such that $\sigma(a+1)<\sigma(a+2)<\ldots<\sigma(b)$.
We note that if $a=b$ then the contraction is a nondegenerate pairing
giving a natural isomorphism
$(\external{}{a}E)\dual=\external{}{a}E\dual$. This isomorphism is,
up to sign, the same considered is
Section~\ref{CAP:KAEHLER}.\ref{sec:exteriorherm}.\\
If $a>b$ we use the
convention that $\vdash=0$.

\begin{lem}\label{XII.1.1}
\begin{enumerate}
\item For every $v\in E$ the operator $v\vdash$ is a derivation of
degree $-1$ of the graded algebra $\external{}{*}E\dual$. \item
For every $v\in \bigwedge^{a}E$, $w\in \bigwedge^{b}E$, $z\in
\bigwedge^{c}E\dual$, we have \[ (v\wedge w)\vdash
z=v\vdash(w\vdash z).\] In particular the operator ~$w\vdash\colon
\bigwedge^{c}E\dual\to \bigwedge^{c-b}E\dual$ is the adjoint of ~
$\wedge w\colon \bigwedge^{c-b}E\to \bigwedge^{c}E$. \item If
$v\in \bigwedge^{a}E\dual$, $w\in \bigwedge^{b}E$, $\Omega\in
\bigwedge^{n}E\dual$, where $\dim E=n$,  $a\le b$, then:
\[ v\wedge (w\vdash \Omega)=(v\vdash w)\vdash \Omega.\]
\end{enumerate}\end{lem}

\begin{proof} {[1]} Complete $v$ to a basis $v=e_{1},\ldots,e_{n}$ of $E$ and
let $z_{1},\ldots,z_{n}$ be its dual basis. Every $w\in
\bigwedge^{*}E\dual$ can be written in a unique way as $w=z_{1}\wedge
w_{1}+w_{2}$ with $w_{1},w_{2}\in \bigwedge^{*}v^{\perp}$. According to
the definition of $\vdash$ we have $v\vdash w=w_1$.\\
If $w=z_{1}\wedge w_{1}+w_{2}$, $u=z_{1}\wedge u_{1}+u_{2}$ are
decompositions as above then
\renewcommand\arraystretch{1.8}
\[ \begin{array}{rl}
(v\vdash w)\wedge u+(-1)^{\bar{w}}w\wedge(v\vdash u)&=
w_{1}\wedge(z_{1}\wedge u_{1}+u_{2})+(-1)^{\bar{w_{2}}}(z_{1}\wedge
w_{1}+w_{2})\wedge u_{1}\\
&=w_{1}\wedge u_{2}+(-1)^{\bar{w_{2}}}w_{2}\wedge u_{1}.\\
v\vdash(w\wedge u)&=v\vdash ((z_{1}\wedge w_{1}+w_{2})
\wedge (z_{1}\wedge u_{1}+u_{2}))\\
&=v\vdash (z_{1}\wedge w_{1}\wedge u_{2}+w_{2}\wedge z_{1}\wedge
u_{1}+w_{2}\wedge u_{2})\\
&=w_{1}\wedge u_{2}+(-1)^{\bar{w_{2}}}w_{2}\wedge u_{1}.
\end{array}\]
\renewcommand\arraystretch{1}
{[2]} Immediate from the definition.\\
{[3]} Induction on $a$; if
$a=1$ then complete $v$ to a basis $v=z_{1},\ldots,z_{n}$ of $E\dual$
and denote $e_{1},\ldots,e_{n}\in E$ its dual basis. Writing
\[ w=e_{1}\wedge w_{1}+w_{2},\qquad w_{i}\in \external{}{*}v^{\perp},
\qquad w_{i}\vdash \Omega
=v\wedge\eta_{i},\qquad \eta_{i}\in\external{}{*}e_{1}^{\perp},\]
we have by Item 2
\[ w\vdash \Omega=(e_{1}\wedge w_{1})\vdash \Omega+(w_{2}\vdash
\Omega)
=e_{1}\vdash(w_{1}\vdash \Omega)+(w_{2}\vdash
\Omega)=\eta_{1}+v\wedge\eta_{2},\]
and then
\[ v\wedge(w\vdash \Omega)=v\wedge \eta_{1}=w_{1}\vdash \Omega=(v\vdash
w)\vdash \Omega.\]
If $a>1$ and $v=v_{1}\wedge v_{2}$, with $v_{1}\in E\dual$,
$v_{2}\in\bigwedge^{a-1}E\dual$ then by item 2 and inductive assumption
\[ v_{1}\wedge v_{2}\wedge (w\vdash \Omega)=v_{1}\wedge(
(v_{2}\vdash w)\vdash \Omega)=(v_{1}\vdash(v_{2}\vdash w))
\vdash \Omega=
((v_{1}\wedge v_{2})\vdash w)\vdash \Omega.\]
\end{proof}

\begin{lem}\label{XII.contra}
For every vector space $E$ of
dimension $n$ and every integer
$a=0,\ldots,n$,
the contraction operator  defines a natural
isomorphism
\[\external{}{a}E\mapor{i}\external{}{n}E\otimes
\external{}{n-a}E\dual,\qquad i(v)=Z\otimes(v\vdash\Omega)\]
where $(Z,\Omega)\in \bigwedge^{n}E\times\bigwedge^{n}E\dual$
is any pair satisfying  $Z\vdash\Omega=1$.
\end{lem}

\begin{proof} Trivial.\end{proof}

\begin{exer} Let $0\mapor{}E\mapor{}F\mapor{}G\mapor{}0$ be an exact
sequence of vector spaces with $\dim G=n<\infty$. Use the contraction
operator to define,
for every $a\le \dim E$, a natural surjective linear map
$\external{}{a+n}F\to\external{}{a}E\otimes\external{}{n}G$.
\end{exer}

\bigskip

\section{The Tian-Todorov's lemma}
\label{sec:TianTodorov}

The isomorphism $i$ of Lemma~\ref{XII.contra} can be extended
fiberwise to vector bundles; in particular, if $M$ is a complex
manifold of dimension $n$ and $T_{M}$ is its holomorphic tangent
bundle, we have holomorphic isomorphisms
\[i\colon \external{}{a}T_{M}
\mapor{}\external{}{n}T_{M}\otimes\external{}{n-a}T_{M}\dual
=\Omega^{n-a}_{M}(K_{M}\dual)\]
which extend to isomorphisms between their Dolbeault's sheaf
resolutions
\[i\colon (\sA^{0,*}(\external{}{a}T_{M}),\debar)
\mapor{}(\sA^{0,*}(\external{}{n}T_{M}\otimes\external{}{n-a}T_{M}
\dual),\debar)=(\sA^{n-a,*}(K_{M}\dual),\debar).\]

If $z_{1},\ldots,z_{n}$ are local holomorphic coordinates then a local
set of generators of $\external{}{a}T_{M}$ is given by the polyvector fields
$\desude{~}{z_{I}}=\desude{~}{z_{i_{1}}}\wedge\ldots\wedge\desude{~}{z_{i_{a}}}$,
being $I=(i_{1},\ldots,i_{a})$ a multiindex.\\
If $\Omega$ is a local frame of $K_{M}$ and $Z$ a local  frame of
$K_{M}\dual$ such that $Z\vdash\Omega=1$, then
\[i\left(\desude{~}{z_{I}} d\bar{z}_{J}\right)=
Z\otimes\left(\desude{~}{z_{I}}d\bar{z}_{J}\vdash\Omega\right)=
Z\otimes\left(\desude{~}{z_{I}}\vdash\Omega\right)
d\bar{z}_{J}.\]

Given a fixed Hermitian metric $h$ on the line bundle $K_{M}\dual$
we denote
by $D=D'+\bar{\partial}$ the unique hermitian connection on $K_{M}\dual$
compatible with the complex structure.\\
We recall (cf. \cite{Kobabook}) that $D'\colon
\sA^{0,b}(K_{M}\dual\otimes
\Omega^{a}_{M})\to \sA^{0,b}(K_{M}\dual\otimes
\Omega^{a+1}_{M})$
is defined in local coordinates as
\[ D'(Z\otimes \phi)=Z\otimes (\theta\wedge \phi+\de\phi),
\qquad \phi\in \sA^{a,b},\]
where $\theta=\de\log(|Z|^{2})=\de\log(h(Z,Z))$
is the connection form of the frame $Z$.\\
We have moreover $(D')^{2}=0$ and
$D'\bar{\partial}+\bar{\partial}D'=\Theta$ is the curvature of the
metric.\\

We can now define a $\C$-linear operator (depending on
$h$)\footnote{don't confuse this $\Delta$ with the Laplacian}
\[\Delta\colon \sA^{0,b}(\external{}{a}T_{M})\to
\sA^{0,b}(\external{}{a-1}T_{M}),\qquad
\Delta(\phi)=i^{-1}D'(i(\phi)).\]

\begin{lem}\label{XII.1.2}
Locally on $M$, with $\Omega,Z$ and $\theta$ as above we have
\[ \Delta(\phi)\vdash \Omega=\theta\wedge (\phi\vdash
\Omega)+\de(\phi\vdash \Omega)\]
for every $\phi\in \sA^{0,b}(\bigwedge^{*}T_{M})$.
\end{lem}

\begin{proof} By definition
\[ i\Delta(\phi)=Z\otimes (\Delta(\phi)\vdash \Omega),\]
\[ i\Delta(\phi)=D'(i(\phi))=
D'(Z\otimes (\phi\vdash \Omega))=
Z\otimes (\theta\wedge(\phi\vdash\Omega)+\de(\phi\vdash\Omega)).\]
\end{proof}

\begin{lem}\label{XII.1.3}
In  local holomorphic coordinates $z_{1},\ldots,z_{n}$ we have
\[\Delta\left(f\desude{~}{z_{I}} d\bar{z}_{J}\right)=
\left((\theta f+\partial f)\vdash \desude{~}{z_{I}}\right)d\bar{z}_{J},
\qquad f\in \sA^{0,0},\]
where $\theta$ is the connection form of the frame
$Z=\desude{~}{z_{1}}\wedge\ldots\wedge\desude{~}{z_{n}}$ and
the right hand side is considered $=0$ when $I=\emptyset$.
\end{lem}

\begin{proof}
We first note that if $\phi\in\sA^{0,0}(\external{}{a}T_{M})$ then
$i(\phi d\bar{z}_{J})=i(\phi)d\bar{z}_{J}$ and
\[D'i(\phi d\bar{z}_{J})=
D'(Z\otimes (\phi\vdash\Omega)\otimes d\bar{z}_{J})=
D'(Z\otimes (\phi\vdash\Omega))\otimes d\bar{z}_{J}:\]
this implies
that $\Delta(\phi d\bar{z}_{J})=\Delta(\phi)
d\bar{z}_{J}$.
According to Lemma~\ref{XII.1.2}
\[\Delta\left(f\desude{~}{z_{I}}\right)\vdash\Omega=
\theta\wedge \left(f\desude{~}{z_{I}}\vdash
\Omega\right)+\de\left(f\desude{~}{z_{I}}\vdash
\Omega\right)\]
Since
$\Omega=dz_{n}\wedge\ldots\wedge dz_{1}$ we have
$\de\left(\desude{~}{z_{I}}\vdash
\Omega\right)=0$ and then, by Item 3 of Lemma~\ref{XII.1.1},
\[\Delta\left(f\desude{~}{z_{I}}\right)\vdash\Omega=
(\theta f+\de f)\wedge\left(\desude{~}{z_{I}}\vdash\Omega\right)
=\left((\theta f+\de f)\vdash \desude{~}{z_{I}}\right)\vdash\Omega.\]
\end{proof}

Setting $\sP^{a,b}=\sA^{0,b}(\bigwedge^{-a}T_{M})$ for every $a\le
0$, $b\ge 0$, the direct sum
$\sP=(\bigoplus_{a,b}\sP^{a,b},\debar)$ is a sheaf of dg-algebras,
where the sections of $\sA^{0,b}(\bigwedge^{a}T_{M})$ have total
degree $b-a$ and $\debar\colon \sA^{0,b}(\bigwedge^{a}T_{M})\to
\sA^{0,b+1}(\bigwedge^{a}T_{M})$ is the Dolbeault differential.
The product on $\sP$ is the `obvious' one:
\[ (\xi\otimes \phi)\wedge (\eta\otimes \psi)=(-1)^{\bar{\phi}\,\bar{\eta}}
(\xi\wedge \eta)\otimes (\phi\wedge\psi).\]

\begin{lem}\label{XII.1.4} The $\C$-linear operator $\Delta
\colon \sP\to\sP$ has degree $+1$; moreover
$\Delta^{2}=0$ and $[\Delta,\debar]=\Delta\debar+\debar\Delta=i^{-1}\Theta i$.
\end{lem}

\begin{proof} Evident.\end{proof}

Consider the bilinear symmetric map of degree 1, $Q\colon
\sP\times\sP\to\sP$
\[ Q(\alpha,\beta)=\Delta(\alpha\wedge\beta)-\Delta(\alpha)
\wedge\beta- (-1)^{\bar{\alpha}}\alpha\wedge\Delta(\beta).\] A
brutal computation in local coordinates shows that $Q$ is
independent of the metric.  In fact, for every pair of $C^{\infty}$
functions $f,g$
\[Q\left(f\desude{~}{z_{I}} d\bar{z}_{J},
g\desude{~}{z_{H}} d\bar{z}_{K}\right)=
(-1)^{|J|\,|H|}Q\left(f\desude{~}{z_{I}},
g\desude{~}{z_{H}}\right) d\bar{z}_{J}\wedge
d\bar{z}_{K}\]
and
\[Q\left(f\desude{~}{z_{I}},g\desude{~}{z_{H}}\right)=
(\theta fg+\partial(fg))\vdash
\left(\desude{~}{z_{I}}\wedge\desude{~}{z_{H}}\right)-\]
\[-g\left((\theta f+\partial f)\vdash
\desude{~}{z_{I}}\right)\wedge\desude{~}{z_{H}}
-(-1)^{|I|}f\desude{~}{z_{I}}
\wedge \left((\theta g+\partial g)\vdash \desude{~}{z_{H}}\right).\]
According to Lemma~\ref{XII.1.1}, Item 1:
\[Q\left(f\desude{~}{z_{I}},g\desude{~}{z_{H}}\right)=
f\left(\partial g\vdash
\desude{~}{z_{I}}\right)\wedge\desude{~}{z_{H}}
+(-1)^{|I|}g\desude{~}{z_{I}}
\wedge \left(\partial f\vdash \desude{~}{z_{H}}\right).\]
In particular if $|I|=0$, $|H|=1$ then
\[Q\left(fd\bar{z}_{J},
g\desude{~}{z_{h}}
d\bar{z}_{K}\right)=(-1)^{|J|}g\desude{f}{z_{h}}d\bar{z}_{J}\wedge
d\bar{z}_{K},\]
while, if $|I|=|H|=1$ then
\[Q\left(f\desude{~}{z_{i}} d\bar{z}_{J},
g\desude{~}{z_{h}} d\bar{z}_{K}\right)=
(-1)^{|J|}\left(f\frac{\partial g}{\partial z_{i}}\desude{~}{z_{h}}-
g\frac{\partial f}{\partial z_{h}}\desude{~}{z_{i}}\right)
d\bar{z}_{J}\wedge d\bar{z}_{K}.\]
Recalling the definition of the bracket $[~,~]$ in the Kodaira-Spencer
algebra $KS_{M}=\bigoplus_{b}\sA^{0,b}(T_{M})$ we have:
\begin{lem}[Tian-Todorov]\label{XII.1.5}
If $\alpha\in \sA^{0,a}(T_{M})$, $\beta\in \sA^{0,b}(T_{M})$ then
\[ (-1)^{a}[\alpha,\beta]=
\Delta(\alpha\wedge\beta)-\Delta(\alpha)\wedge
\beta-(-1)^{a-1}\alpha\wedge\Delta(\beta).\]
In particular the bracket of two $\Delta$-closed forms is
$\Delta$-exact.
\end{lem}

\begin{ex}\label{XII.1.6}
If $M$ is compact K\"{a}hler and $c_{1}(M)=0$ in $H^{2}(M,\C)$
then by \cite[2.23]{Kobabook} there exists a Hermitian metric on
$K_{M}\dual$ such that $\Theta=0$; in this case $[\Delta,\debar]=0$
and $\ker\Delta$ is a differential graded subalgebra of $KS_{M}$.
\end{ex}

\begin{ex}\label{XII.1.7} If $M$
has a nowhere vanishing holomorphic $n$-form $\Omega$ ($n=\dim M$)  we
can set on $K_{M}\dual$ the trivial Hermitian metric induced by the
isomorphism $\Omega\colon K_{M}\dual\to\Oh_{M}$.
In this case, according to Lemma~\ref{XII.1.2}, the operator
$\Delta$ is defined by the rule
\[(\Delta\alpha)\vdash\Omega=\partial(\alpha\vdash\Omega).\]
\end{ex}

\bigskip

\section{A formality theorem}
\label{sezioneVI.2}

\begin{thm}\label{XII.2.1}
Let $M$ be a compact K\"ahler manifold with trivial canonical bundle
$K_{M}=\Oh_{M}$. Then the Kodaira-Spencer DGLA
\[ KS_{M}=\somdir{p}{}\Gamma(M,\sA^{0,p}(T_{M}))\]
is quasiisomorphic to an abelian DGLA.\end{thm}

\begin{proof}
Let $\Omega\in \Gamma(M,K_{M})$ be a nowhere vanishing holomorphic
$n$-form ($n=\dim M$); via the isomorphism $\Omega\colon
K_{M}\dual\to\Oh_{M}$, the isomorphism of complexes
\[i\colon (\sA^{0,*}(T_{M}),\debar)\to (\sA^{n-1,*},\debar)\]
is given in local holomorphic coordinates by
\[ i\left(f\desude{~}{z_{i}}d\bar{z}_{I}\right)=
f\left(\desude{~}{z_{i}}\vdash\Omega\right)d\bar{z}_{I}\]
and induces a structure of DGLA, isomorphic to $KS_{M}$ on
\[ L^{n-1,*}=\somdir{p}{}\Gamma(M,\sA^{n-1,p}).\]
Taking on $K_{M}\dual$ the trivial metric induced by
$\Omega\colon K_{M}\dual\to\Oh_{M}$, the connection $D$
is equal to the De Rham differential and then
the Tian-Todorov's lemma implies that the bracket of two $\de$-closed
form of $L^{n-1,*}$ is $\de$-exact; in particular
\[ Q^{*}=\ker\de\cap L^{n-1,*}\]
is a DGL subalgebra of $L^{n-1,*}$.\\
Consider the complex $(R^{*},\debar)$, where
\[ R^{p}=\frac{\ker\de\cap L^{n-1,p}}{\de L^{n-2,p}}\]
endowed with the trivial bracket, again by Lemma~\ref{XII.1.5}
the projection $Q^{*}\to R^{*}$ is a morphism of DGLA.\\
It is therefore sufficient to prove that the DGLA morphisms
\[ \xymatrix{L^{n-1,*}&Q^{*}\ar[l]\ar[r]&R^{*}}\]
are quasiisomorphisms.\\
According to the $\de\debar$-lemma \ref{X.4.9},
$\debar(\ker\de)\subset\Image\de$ and then the operator $\debar$ is
trivial on $R^{*}$: therefore
\[ H^{p}(R^{*})=\frac{\ker\de\cap L^{n-1,p}}{\de L^{n-2,p}},\qquad
H^{p}(L^{n-1,*})=\frac{\ker\debar\cap L^{n-1,p}}{\debar L^{n-1,p-1}},\]
\[ H^{p}(Q^{*})=\frac{\ker\de\cap\ker\debar\cap L^{n-1,p}}
{\debar(\ker\de\cap L^{n-1,p-1})}.\]
The conclusion now follows immediately from Corollary~\ref{X.4.10}.
\end{proof}

\begin{cor}\label{XII.2.2}
Let $M$ be a compact K\"ahler manifold with trivial canonical bundle
$K_{M}=\Oh_{M}$. For every local Artinian $\C$-algebra $(A,\ide{m}_{A})$ we
have
\[ \Def_{M}(A)=H^{1}(M,T_{M})\otimes \ide{m}_{A}.\]
In particular
\[\Def_{M}\left(\frac{\C[t]}{(t^{n+1})}\right)\to
\Def_{M}\left(\frac{\C[t]}{(t^{2})}\right)\]
is surjective for every $n\ge 2$.\end{cor}

\begin{proof} According to Theorem~\ref{VIII.4.6} and Corollary~\ref{VIII.4.4}
we have $\Def_{M}=\Def_{R^{*}}$. Since $R^{*}$ is an abelian DGLA we have by
Proposition~\ref{VIII.4.2}
\[ \Def_{R^{*}}(A)=H^{1}(R^{*})\otimes \ide{m}_{A}=
H^{1}(KS_{M})\otimes \ide{m}_{A}=H^{1}(M,T_{M})\otimes \ide{m}_{A}.\]
\end{proof}

\bigskip

\section{Gerstenhaber algebras and Schouten brackets}
\label{sezioneVI.3}

\begin{lem}\label{XII.3.1}
Let $(G,\wedge)$ be a graded $\Z$-commutative algebra and let
$[,]\colon G[-1]\times G[-1]\to G[-1]$ be a skewsymmetric bilinear
map of degree 0  such that
\[ ad_{a}=[a,-]\in \Der^{\deg(a,G[-1])}(G,G),\qquad \forall a\in
G[-1].\]
(Note that this last condition is equivalent to the so-called \emph{Odd
Poisson identity}
\[ [a,b\wedge c]=[a,b]\wedge c+(-1)^{\bar{a}(\bar{b}-1)}b\wedge
[a,c],\]
\[ [a\wedge b,c]=a\wedge [b,c]+(-1)^{\bar{c}(\bar{b}-1)}[a,c]\wedge
b,\]
for every $a,b,c\in G[-1]$, $\bar{x}=\deg(x,G[-1])$.)\\
Let $\sG\subset G$ be a set of homogeneous generators of the algebra
$G$, then:\begin{enumerate}
\item $[,]$ is uniquely determined by the values $[a,b]$, $a,b\in \sG$.
\item A derivation $d\in \Der^{n}(G,G)$ satisfies
$[d,ad_{a}]=ad_{d(a)}$ for every $a\in G[-1]$ if and only if
\[ d[a,b]=[da,b]+(-1)^{n\bar{a}}[a,db]\]
for every $a,b\in\sG$.
\item $[,]$ satisfies the Jacobi condition
$ad_{[a,b]}=[ad_{a},ad_{b}]$ if and only if
\[ [[a,b],c]=[a,[b,c]]-(-1)^{\bar{a}\,\bar{b}}[b,[a,c]].\]
for every $a,b,c\in\sG$.\end{enumerate}\end{lem}

\begin{proof} 1) is clear.\\
If $a\in \sG$ then by 2) the derivations $[d,ad_{a}]$ and
$ad_{d(a)}$ take the same values in $\sG$ and then
$[d,ad_{a}]=ad_{d(a)}$. The skewsymmetry of $[,]$ implies that for
every $b\in G[-1]$ the derivations $[d,ad_{b}]$ and $ad_{d(b)}$
take the same values in $\sG$.\\
The proof of 3) is made by applying twice 2), first with $d=ad_{a}$,
$a\in \sG$, and then with $d=ad_{b}$, $b\in G[-1]$.\end{proof}

\begin{defn}\label{XII.3.2}  A
\emph{Gerstenhaber algebra} is the data of a graded
$\Z$-commutative algebra $(G,\wedge)$ and a morphism of graded vector
spaces $ad\colon G[-1]\to \Der^{*}(G,G)$ such that the bracket
\[ [,]\colon G[-1]_{i}\times G[-1]_{j}\to G[-1]_{i+j},\qquad
[a,b]=ad_{a}(b)\]
induce a structure of graded Lie algebra on $G[-1]$
(cf. \cite[p. 267]{Gerst}).\\
A morphism of Gerstenhaber algebras is a morphism of graded algebras
commuting with the bracket $[,]$.
\end{defn}

For every graded vector space $G$ there exists an isomorphism from
the space of bilinear skewsymmetric maps $[,]\colon G[-1]\times
G[-1]\to G[-1]$ of degree 0 and the space of bilinear symmetric
maps $Q\colon G\times G\to G$ of degree 1; this isomorphism,
called \emph{d\'ecalage}, is given by the formula\footnote{The
d\'ecalage isomorphism is natural up to sign; the choice of
$\deg(a,G[-1])$ instead of $\deg(a,G)$ is purely conventional.}
\[ Q(a,b)=(-1)^{\deg(a,G[-1])}[a,b].\]
Therefore a Gerstenhaber algebra can be equivalently defined as a
graded algebra $(G,\wedge)$ endowed with a bilinear symmetric map
$Q\colon G\times G\to G$ of degree 1 satisfying the identities

~\\
\begin{tabular}{lc}
\emph{Odd Poisson}&
$Q(a,b\wedge c)=Q(a,b)\wedge c+(-1)^{(\bar{a}+1)\bar{b}}b\wedge
Q(a,c),$\\
&\\
\emph{Jacobi}&
$Q(a,Q(b,c))=(-1)^{\bar{a}}Q(Q(a,b),c)+
(-1)^{\bar{a}\,\bar{b}}Q(b,Q(a,c)),$\\
&
\end{tabular}
~\\
where $\bar{a}=\deg(a,G)$, $\bar{b}=\deg(b,G)$.

\begin{ex}\label{XII.3.3} (Schouten algebras)
A particular class of Gerstenhaber algebras are the so called {\em
Schouten algebras}: here the bracket is usually called \emph{Schouten
bracket}.\\
Consider a commutative $\K$-algebra $A_{0}$ and let
$A_{-1}\subset \Der_{\K}(A_{0},A_{0})$ be an $A_{0}$-submodule such
that $[A_{-1},A_{-1}]\subset A_{-1}$. Define
\[ A=\somdir{i\ge 0}{}A_{-i},\qquad
A_{-i}=\external{A_{0}}{i}A_{-1}.\]
With the wedge product, $A$ is a graded algebra of nonpositive
degrees.\\
There exists
a unique structure of Gerstenhaber algebra $(A,\wedge,[,])$
such that for every $a,b\in
A[-1]_{1}=A_{0}$, $f,g\in A[-1]_{0}=A_{-1}$
\[ ad_{a}(b)=0,\quad ad_{f}(a)=f(a),\quad ad_{f}(g)=[f,g].\]
In fact $A$ is generated by $A_{0}\cup A_{-1}$ and, according to
Lemma~\ref{XII.3.1}, the skew-symmetric bilinear map
\[ [\xi_0\wedge\ldots\wedge\xi_n,h]=
\sum_{i=0}^n(-1)^{n-i}\xi_i(h)\xi_0\wedge\ldots\wedge
\widehat{\xi_i}\wedge\ldots\wedge\xi_n\]
\[ [\xi_0\wedge\ldots\wedge\xi_n,\zeta_0\wedge\ldots\wedge\zeta_m]=
\qquad\qquad\qquad\]
\[\qquad=\sum_{i=0}^n\sum_{j=0}^m
(-1)^{i+j}[\xi_i,\zeta_j]\wedge\xi_0\wedge\ldots\wedge
\widehat{\xi_i}\wedge\ldots\wedge\xi_n\wedge\zeta_0\wedge\ldots\wedge
\widehat{\zeta_j}\wedge\ldots\wedge\zeta_m \]
where $h\in A_{0}$,
$\xi_0,\ldots,\xi_n,\zeta_0,\ldots,\zeta_m\in A_{-1}$ is well defined
and it is the unique extension of the natural bracket such that
$ad(A[-1])\subset \Der^{*}(A,A)$.\\
We need to show that $[,]$ satisfies the Jacobi identity
\[ [[a,b],c]=[a,[b,c]]-(-1)^{\bar{a}\,\bar{b}}[b,[a,c]].\]
Again by Lemma~\ref{XII.3.1} we may
assume that $0\le \bar{a}\le \bar{b}\le \bar{c}$.
There are 5 possible cases, where
the Jacobi identity is satisfied for trivial reasons, as summarized
in the following table:

\vbox{
\begin{center}
{\renewcommand{\arraystretch}{1.5}
\begin{tabular}{|l|l|l|l|}\hline
$\bar{a}$&$\bar{b}$&
$\bar{c}$&Jacobi is true because..\\ \hline
$1$&$1$&
$1$& all terms are $=0$\\ \hline
$0$&$1$&
$1$& all terms are  $=0$\\ \hline
$0$&$0$&
$1$& definition of $[,]$ on $A_{-1}$\\ \hline
$0$&$0$&
$0$& Jacobi identity on $A_{-1}$\\ \hline
\end{tabular}}
\end{center}
}
\end{ex}

\begin{ex}\label{XII.3.4}
Let $M$ be a complex manifold of dimension $n$, the sheaf
of graded algebras $\sT=\oplus_{i\le 0}\sT_{i}$,
$\sT_{i}=\sA^{0,0}(\bigwedge^{-i}T_{M})$,
admits naturally a Schouten bracket.\\
In local holomorphic coordinates $z_{1},\ldots,z_{n}$, since
\[\left[\desude{~}{z_{i}},\desude{~}{z_{j}}\right]=0,\qquad
\left[\desude{~}{z_{I}},g\right]_{Sch}=
(-1)^{|I|-1}\left(\partial g\vdash
\desude{~}{z_{I}}\right),\]
the Odd Poisson identity implies that
the Schouten
bracket  takes the simple form
\[\left[f\desude{~}{z_{I}},g\desude{~}{z_{H}}\right]_{Sch}=
(-1)^{|I|-1}f\left(\partial g\vdash
\desude{~}{z_{I}}\right)\wedge\desude{~}{z_{H}}
-g\desude{~}{z_{I}}
\wedge \left(\partial f\vdash \desude{~}{z_{H}}\right).\]
\end{ex}

\begin{defn}\label{XII.3.5}
A \emph{differential Gertstenhaber algebra} is a
Gerstenhaber algebra $(G,\wedge,[,])$ endowed with a differential
$d\in \Der^{1}(G,G)$ making $(G,d,[,])$ a differential graded Lie
algebra.\end{defn}

\begin{ex}\label{XII.3.6}
Given any Gertstenhaber algebra $G$ and an element $a\in
G_{0}=G[-1]_{1}$ such that $[a,a]=0$ we have that
$d=ad_{a}$ gives a structure of differential Gerstenhaber
algebra.\end{ex}

\begin{exer} For every $f\in\K[x_{1},\ldots,x_{n}]$ the Koszul
complex of the sequence $\desude{f}{x_{1}},\ldots,\desude{f}{x_{n}}$
carries a structure of differential Gerstenhaber algebra.
\end{exer}

\bigskip

\section{$d$-Gerstenhaber structure on polyvector fields}
\label{sezioneVI.4}

Let $M$ be a fixed complex manifold, then
the sheaf of dg-algebras $\sP$ defined in Section~\ref{sec:TianTodorov},
endowed with the Schouten bracket
\[\left[f\desude{~}{z_{I}} d\bar{z}_{J},
g\desude{~}{z_{H}} d\bar{z}_{K}\right]_{Sch}=
(-1)^{|J|(|H|-1)}
\left[f\desude{~}{z_{I}},
g\desude{~}{z_{H}}\right]_{Sch}
d\bar{z}_{J}\wedge
d\bar{z}_{K}\]
is a sheaf of differential Gerstenhaber algebras.\\
We have only to verify that locally
$\debar$ is a derivation of the graded Lie algebra $(\sP,[,])$:
this follows immediately from Lemma~\ref{XII.3.1} and from
the fact that locally the Kodaira-Spencer DGLA generates $\sP$ as a
graded algebra.\\
Via the d\'ecalage isomorphism, the Schouten bracket corresponds to
the symmetric bilinear map of degree 1 $Q\colon \sP\times\sP\to\sP$
given in local holomorphic coordinates by the formulas
\[ Q\left(fd\bar{z}_{J}\desude{~}{z_{I}},
gd\bar{z}_{K}\desude{~}{z_{H}}\right)=
(-1)^{|K|(|I|-1)+|J|}
d\bar{z}_{J}\wedge
d\bar{z}_{K}
Q\left(f\desude{~}{z_{I}},
g\desude{~}{z_{H}}\right),
\]
where
\[Q\left(f\desude{~}{z_{I}},g\desude{~}{z_{H}}\right)=
f\left(\partial g\vdash
\desude{~}{z_{I}}\right)\wedge\desude{~}{z_{H}}
+(-1)^{|I|}g\desude{~}{z_{I}}
\wedge \left(\partial f\vdash \desude{~}{z_{H}}\right).\]

Notice that, in the notation of Section~\ref{sec:TianTodorov},
\[ Q(\alpha,\beta)=\Delta(\alpha\wedge\beta)-\Delta(\alpha)
\wedge\beta-
(-1)^{\bar{\alpha}}\alpha\wedge\Delta(\beta)\]
and therefore we also have the following

\begin{lem}[Tian-Todorov]
for every $\alpha,\beta\in \sP[-1]$,
\[ [\alpha,\beta]_{Sch}=\alpha\wedge\Delta\beta+
(-1)^{\deg(\alpha,\sP[-1])}
(\Delta(\alpha\wedge\beta)-\Delta\alpha\wedge\beta).\]
\end{lem}

There exists a natural morphism
$\widehat{~}\colon \sP\to
\HOM(\sA^{*,*},\sA^{*,*})$
of sheaves of bigraded vector spaces
on $M$ given in local coordinates by
\[ \widehat{\phi\desude{~}{z_{I}}}(\eta)=
\phi\wedge\left(\desude{~}{z_{I}}
\vdash\eta\right).\]
Since, for every $\phi\in \sP^{0,p}=\sA^{0,p}$, $\eta\in\sA^{*,*}$, we have
\[\desude{~}{z_{I}}
\vdash\left(\phi\wedge\eta\right)=
(-1)^{p|I|}\phi\wedge\left(\desude{~}{z_{I}}
\vdash\eta\right)\]
the hat morphism $\widehat{~}$ is a morphism of algebras, being the
product in $\HOM(\sA^{*,*},\sA^{*,*})$ the composition product.
We observe that the composition product
is associative and therefore
$\HOM(\sA^{*,*},\sA^{*,*})$ has also a natural structure of sheaf of
graded Lie algebras. Since $\sP$ is graded commutative,
$[\widehat{a},\widehat{b}]=0$ for every $a,b\in \sP$.

\begin{lem}\label{quadratic} For every $a,b\in\sP$ homogeneous,
\begin{enumerate}
\item  $\widehat{\debar a}=[\debar,\widehat{a}].$

\item $\widehat{Q(a,b)}=[[\de,\widehat{a}],\widehat{b}]=
-(-1)^{\bar{a}}\widehat{a}\de\h{b}-
(-1)^{\bar{a}\,\bar{b}+\bar{b}}\,\h{b}\de\h{a} \pm
\de\h{a}\h{b}\pm\h{b}\h{a}\de$
\end{enumerate}
\end{lem}

\begin{proof} The proof of the first identity is straightforward and
left to the reader.\\
By Jacobi identity,
\[ 0=[\de,[\h{a},\h{b}]]=[[\de,\h{a}],\h{b}]-(-1)^{\bar{a}\,\bar{b}}
[[\de,\h{b}],\h{a}]\]
and therefore both sides of the equality~\ref{quadratic} are graded
symmetric.\\
Moreover, since $\h{b\wedge c}=\h{b}\h{c}$ and
\[ Q(a,b\wedge c)=Q(a,b)\wedge c+(-1)^{(\bar{a}+1)\bar{b}}b\wedge
Q(a,c),\]
\[ [[\de,\h{a}],\h{b}\h{c}]=[[\de,\h{a}],\h{b}]\h{c}
+(-1)^{(\bar{a}+1)\bar{b}}\h{b} [[\de,\h{a}],\h{c}],\]
it is sufficient to check the equality only when
$a,b=f,d\bar{z}_{j},\desude{~}{z_{i}}$, $f\in \sP^{0,0}=\sA^{0,0}$.\\
i) If $\phi\in \sP^{0,*}$ then
\[ [\de,\h{\phi}]\eta=\de(\phi\wedge\eta)-
(-1)^{\bar{\phi}}\phi\wedge\de\eta=\de\phi\wedge\eta.\]
In particular $[\de,\h{d\bar{z}_{j}}]=0$, $Q(d\bar{z}_{j},b)=0$ for
every $b$.\\
ii) If $f,g\in\sP^{0,0}$ then $Q(f,g)\in\sP^{1,0}=0$ and
\[ [[\de,\h{f}],\h{g}]\eta=\de f\wedge g\eta-g(\de f\wedge \eta)=0.\]
If $f\in\sP^{0,0}$ then $\ds Q\left(f,\desude{~}{z_{i}}\right)=
\desude{~}{z_{i}}\vdash \de f=\desude{f}{z_{i}}$ and
\[ \left[[\de,\h{f}],\h{\desude{~}{z_{i}}}\right]\eta=\de f\wedge
\left(\desude{~}{z_{i}}\vdash\eta\right)
+\desude{~}{z_{i}}\vdash(\de f\wedge \eta)=
\left(\desude{~}{z_{i}}\vdash \de f\right)\wedge\eta\]
where the last equality follows from the Leibnitz rule
applied to the derivation  $\desude{~}{z_{i}}\vdash$.\\
Finally $\ds Q\left(\desude{~}{z_{i}},\desude{~}{z_{j}}\right)=0$;
since $\de, \h{\desude{~}{z_{i}}},\h{\desude{~}{z_{j}}}$ are
derivations of $\sA^{*,*}$,
also $\left[\left[\de, \h{\desude{~}{z_{i}}}\right],
\h{\desude{~}{z_{j}}}\right]$ is a derivation of bidegree $(-1,0)$
and then it is sufficient to check the equality
for $\eta=dz_{i}$. This last verification is completely
straightforward and it is left to the reader.\end{proof}

\begin{exer} Prove that $\bar{\Omega}^{*}=\{a\in\sP\mid
[\de,\h{a}]=0\,\}$.\end{exer}

\bigskip
\section{GBV-algebras}
\label{sezioneVI.5}

In this section $\K$ is a fixed field of characteristic 0.

\begin{defn}\label{XII.5.1}
A GBV \emph{(Gerstenhaber-Batalin-Vilkovisky)} algebra is the data of a
graded algebra $(G,\wedge)$ and a linear map $\Delta\colon G\to G$ of
degree 1 such that:
\begin{enumerate}
    \item  $\Delta^{2}=0$

    \item  The symmetric bilinear map of degree 1
    \[ Q(a,b)=\Delta(a\wedge b)-\Delta(a)
\wedge b-
(-1)^{\bar{a}}a\wedge\Delta(b)\]
satisfies the odd Poisson identity
\[ Q(a,b\wedge c)=Q(a,b)\wedge c+(-1)^{(\bar{a}+1)\bar{b}}b\wedge
Q(a,c).\]
\end{enumerate}
\end{defn}

Note that the second condition on the above definition means that for
every homogeneous $a\in G$, the linear map $Q(a,-)$ is a derivation of
degree $\bar{a}+1$.

The map $Q$ corresponds, via the d\'ecalage isomorphism,
to a skewsymmetric bilinear map of degree 0,
$[,]\colon G[-1]\times G[-1]\to G[-1]$;
the expression of $[,]$ in terms of $\Delta$ is
\[ [a,b]=a\wedge\Delta(b)+
(-1)^{\deg(a,G[-1])}
(\Delta(a\wedge b)-\Delta(a)\wedge b).\]

\begin{ex}\label{XII.5.2}
If $\Delta$ is a differential of a graded algebra
$(G,\wedge)$, then $Q=0$ and
$(G,\wedge,\Delta)$ is a GBV algebra called {\em
abelian}.\end{ex}

\begin{ex}\label{XII.5.3} The sheaf $\sP$ of polyvector fields on a
complex manifold, endowed with the operator $\Delta$ described in
Section~\ref{sec:TianTodorov} is a sheaf of GBV algebra.\end{ex}

\begin{exer} Let $(G,\wedge,\Delta)$ be a GBV algebra. If $G$ has a
unit 1, then $\Delta(1)=0$.\end{exer}

\begin{lem}\label{XII.5.4} For every $a,b\in G$ homogeneous
\[ \Delta Q(a,b)+Q(\Delta(a),b)+(-1)^{\bar{a}}Q(a,\Delta(b))=0.\]
\end{lem}

\begin{proof} It is sufficient to write $Q$ in terms of $\Delta$ and use
$\Delta^{2}=0$.\end{proof}

\begin{thm}\label{XII.5.5} If
$(G,\wedge,\Delta)$ is a GBV algebra then
$(G[-1],[,],\Delta)$ is a DGLA and therefore
$(G,\wedge,Q)$ is a Gerstenhaber algebra.
\end{thm}

\begin{proof}
Working in $G[-1]$ (i.e. $\bar{a}=\deg(a,G[-1])$) we have
from Lemma~\ref{XII.5.4}
\[ \Delta [a,b]=[\Delta(a),b]+(-1)^{\bar{a}}[a,\Delta(b)]\]
and then we only need to prove the Jacobi identity.\\
Replacing $a=\alpha$, $b=\beta\wedge\gamma$ in the above formula
we have
\[ [\alpha,\Delta(\beta\wedge\gamma)]=(-1)^{\bar{\alpha}}
(\Delta[\alpha,\beta\wedge\gamma]-[\Delta\alpha,\beta\wedge\gamma])\]
and then
$[\alpha,\Delta(\beta\wedge\gamma)]$ is equal to
\[(-1)^{\bar{\alpha}}
\Delta([\alpha,\beta]\wedge \gamma)+(-1)^{\bar{\alpha}\,\bar{\beta}}
\Delta(\beta\wedge[\alpha,\gamma])-(-1)^{\bar{\alpha}}
[\Delta\alpha,\beta]\!\wedge\!\gamma+\!(-1)^{(\bar{\alpha}+1)\bar{\beta}}
\beta\!\wedge\![\Delta\alpha,\gamma].\]
Writing
\[ [\alpha,[\beta,\gamma]]=[\alpha,\beta\wedge\Delta\gamma]+
(-1)^{\bar{\beta}}([\alpha,\Delta(\beta\wedge\gamma)]-
[\alpha,\Delta\beta\wedge\gamma]),\]
\[ [[\alpha,\beta],\gamma]=[\alpha,\beta]\wedge \Delta\gamma+
(-1)^{\bar{\alpha}+\bar{\beta}}
(\Delta([\alpha,\beta]\wedge \gamma)-\Delta[\alpha,\beta]\wedge
\gamma),\]
\[ [\beta,[\alpha,\gamma]]=\beta\wedge\Delta[\alpha,\gamma]+
(-1)^{\bar{\beta}}(\Delta(\beta\wedge[\alpha,\gamma])-
\Delta\beta\wedge[\alpha,\gamma])\]
we get
\[ [\alpha,[\beta,\gamma]]=[[\alpha,\beta],\gamma]+
(-1)^{\bar{\alpha}\,\bar{\beta}}
[\beta,[\alpha,\gamma]].\]
\end{proof}

\begin{defn}\label{XII.5.6} Let $(G,\wedge,\Delta)$ be a GBV-algebra
and $d$ a differential of degree 1 of $(G,\wedge)$.
If $d\Delta+\Delta d=0$ then
the gadget $(G,\wedge,\Delta, d)$ is called a \emph{differential GBV
algebra}.\end{defn}

\begin{ex}\label{XII.5.7}
Let $\sP$ be the algebra of polyvector fields on a complex manifold 
$M$. In the notation of Section~\ref{sec:TianTodorov},
$(\sP,\wedge,\Delta, \debar)$ is  a sheaf of differential GBV
algebras if and only if the connection $D$ is integrable.\\
This happen in particular when $M$ has trivial canonical bundle and 
$D$ is the trivial connection.  
\end{ex}

\begin{exer}
If $(G,\wedge,\Delta, d)$ is a differential GBV-algebra then
$(G[-1],[,], d+\hbar\Delta)$ is a DGLA for every $\hbar\in\K$.
\end{exer}

\bigskip

\section{Historical survey,~\ref{CAP:TRIVIALK}}

The Schouten bracket was introduced by Schouten in \cite{Schouten}
while the Jacobi identity was proved 15 years later by
Nijenhuis \cite{Nijenhuis}.\\
The  now called Gerstenhaber algebras have been first 
studied in \cite{Gerst} as a structure on the cohomology of an 
associative ring.\\
Concrete examples of GBV algebra arising from string theory were 
studied in 1981 by Batalin and Vilkovisky, while the abstract 
definition 
of GBV algebra given in this notes 
was proposed in \cite{LZ} (cf. also \cite{sta97}).\\


\chapter[~~Graded coalgebras]{Graded coalgebras}
\label{CAPITOLO13}
\piede

This \chaptername\ is a basic course on graded coalgebra, with 
particular emphasis on symmetric graded coalgebra. The aim is give the 
main definitions and to give 
all the preliminaries for a satisfactory theory of $L_{\infty}$-algebras.\\  
Through all the chapter we work over a fixed field $\K$ of characteristic
0. Unless otherwise specified all the tensor products are made over
$\K$.\\
The main references for this \chaptername\ are  \cite[Appendix B]{Qui}
\cite{grassi1}, \cite{Boualge}.

\bigskip

\section[~~Koszul sign and unshuffles]{Koszul sign and unshuffles}
\label{sec:koszul}

Let $V, W\in \mathbf{G}$
be  graded vector spaces over $\K$. We recall (Definition 
\ref{twistingmap}) 
that the \emph{twisting map}
$T\colon V\otimes W\to W\otimes V$ is defined by the
rule $T(v\otimes w)=(-1)^{\bar{v}\,\bar{w}}w\otimes v$, for every pair
of  homogeneous elements $v\in V$, $w\in W$.\\

The \emph{tensor algebra} generated by $V\in \mathbf{G}$
is by definition the graded vector space
\[ T(V)=\hbox{$\bigoplus$}_{n\ge 0}\hbox{$\bigotimes$}^{n}V\]
endowed with the associative product $(v_{1}\otimes\ldots\otimes
v_{p}) (v_{p+1}\otimes\ldots\otimes
v_{n})=v_{1}\otimes\ldots\otimes
v_{n}$.\\

Let $I\subset T(V)$ be the homogeneous
ideal generated by the elements
$x\otimes  y-T(x\otimes y)$, $x,y\in V$; the \emph{symmetric algebra}
generated by $V$  is defined as the quotient
\[ S(V)=T(V)/I=\hbox{$\bigoplus$}_{n\ge 0}\symm{}{n}V,\qquad
\symm{}{n}V=\tensor{}{n}V/(\tensor{}{n}V\cap I)\]
The product in $S(V)$ is denoted by $\odot$. In particular
if $\pi\colon T(V)\to S(V)$ is the projection to the quotient
then for every $v_{1},\ldots,v_{n}\in V$,
$v_{1}\odot\ldots\odot v_{n}=\pi(v_{1}\otimes\ldots\otimes v_{n})$.\\

The \emph{exterior algebra} generated by $V$  is the quotient of $T(V)$
by the  homogeneous ideal $J$ generated by the elements $x\otimes
y+T(x\otimes y)$.\\
\[ \external{}{}V=T(V)/J=\hbox{$\bigoplus$}_{n\ge 0}\external{}{n}V,\qquad
\external{}{n}V=\tensor{}{n}V/(\tensor{}{n}V\cap J).\]

Every morphism of graded vector spaces $f\colon V\to W$ induces
canonically three homomorphisms of graded algebras
\[ T(f)\colon T(V)\to T(W),\qquad S(f)\colon S(V)\to S(W),\qquad
\external{}{}(f)\colon \external{}{}V\to \external{}{}W.\]

The following convention is adopted in force: let $V,W$ be graded
vector spaces and $F\colon T(V)\to T(W)$ a linear map. We
denote by
\[ F^{i}\colon T(V)\to \hbox{$\bigotimes$}^{i}W,\quad
F_{j}\colon \hbox{$\bigotimes$}^{j}V\to T(W),\quad
F^{i}_{j}\colon \hbox{$\bigotimes$}^{j}V\to \hbox{$\bigotimes$}^{i}W\]
the compositions of $F$ with the inclusion
$\hbox{$\bigotimes$}^{j}V\to T(V)$ and/or the projection
$T(W)\to \hbox{$\bigotimes$}^{i}W$.\\
Similar terminology is adopted for linear maps $S(V)\to S(W)$.\\

If $v_{1},\ldots,v_{n}$ is an ordered tuple of homogeneous
elements of $V$ and
$\sigma\colon \{1,\ldots,s\}\to \{1,\ldots,n\}$ is any map,
we denote $v_{\sigma}=v_{\sigma{1}}\odot v_{\sigma{2}}\odot
\ldots\odot v_{\sigma{s}}\in \symm{}{s}V$.\\
If $I\subset \{1,\ldots,n\}$ is a subset of cardinality $s$ we define
$v_{I}$ as above, considering $I$ as a strictly increasing map
$I\colon \{1,\ldots,s\}\to \{1,\ldots,n\}$.\\
If $I_{1}\cup\ldots\cup I_{a}=J_{1}\cup\ldots\cup J_{b}=\{1,\ldots,n\}$
are decompositions of $\{1,\ldots,n\}$ into disjoint subsets, we define the
\emph{Koszul sign}
$\epsilon\left(V,
\genfrac{}{}{0pt}{0}{I_{1},\ldots,I_{a}}{J_{1},\ldots,J_{b}};
\{v_{h}\}\right)=\pm 1$ by the relation

\[
\epsilon\left(V,\genfrac{}{}{0pt}{0}{I_{1},\ldots,I_{a}}{
J_{1},\ldots,J_{b}}; \{v_{h}\}\right)
v_{I_{1}}\odot\ldots\odot v_{I_{a}}=
v_{J_{1}}\odot\ldots\odot v_{J_{b}}.
\]

Similarly, if $\sigma$ is a permutation of $\{1,\ldots,n\}$,
$\epsilon(V,\sigma;v_{1},\ldots,v_{n})=\pm 1$ is defined by
\[v_{1}\odot\ldots\odot v_{n}=
\epsilon(V,\sigma;v_{1},\ldots,v_{n})
(v_{\sigma(1)}\odot\ldots\odot v_{\sigma(n)}),\]
or more explicitly
\[ \epsilon(V,\sigma;v_{1},\ldots,v_{n})=
\prod_{i<j}\left(\frac{\sigma_{i}-\sigma_{j}}{|\sigma_{i}-\sigma_{j}|}
\right)^{\bar{v_{i}}\,\bar{v_{j}}},\qquad \bar{v}=\deg(v;V).\]
For notational simplicity we shall write 
$\epsilon(\sigma;v_{1},\ldots,v_{n})$ or   $\epsilon(\sigma)$ when 
there is no possible confusion about $V$ 
and $v_{1},\ldots,v_{n}$.\\ 

The action of the twisting map on $\tensor{}{2}V$ extends naturally,
for every $n\ge 0$, to an action of the
symmetric group $\Sigma_{n}$ on the graded vector space
$\tensor{}{n}V$. This action can be described by the use of Koszul
sign, more precisely
\[ \sigma(v_{1}\otimes\ldots\otimes v_{n})=
\epsilon(\sigma;v_{1},\ldots,v_{n})
(v_{\sigma(1)}\otimes\ldots\otimes v_{\sigma(n)})\]

Denote by
$N\colon S(V)\to T(V)$  the linear map
\[ \begin{split}
N(v_{1}\odot\ldots\odot v_{n})=&
\sum_{\sigma\in\Sigma_{n}}\epsilon(\sigma;v_{1},\ldots,v_{n})
(v_{\sigma(1)}\otimes\ldots\otimes v_{\sigma(n)})\\
=&\sum_{\sigma\in\Sigma_{n}}\sigma(v_{1}\otimes\ldots\otimes 
v_{n}),\quad v_{1},\ldots,v_{n}\in V.
\end{split}
\]
Since $\K$ has characteristic $0$, a left inverse of
$\pi\colon T(V)\to S(V)$ is given
by $\sum_{n}\dfrac{Id^{n}}{n!}N$, where, according to our convention,
$Id^{n}\colon T(V)\to
\tensor{}{n}V$ is the projection.\\
For every homomorphism of graded
vector spaces $f\colon V\to W$, we have
\[ N\circ S(f)=T(f)\circ N\colon S(V)\to T(W).\]

The image of $N\colon \symm{}{n}V\to \tensor{}{n}V$ is contained
in the subspace $(\tensor{}{n}V)^{\Sigma_{n}}$ of
$\Sigma_{n}$-invariant vectors.

\begin{lem}\label{XIII.1.2}
In the notation above, let
$W\subset \tensor{}{n}V$ be the subspace generated by all the vectors $v-\sigma(v)$,
$\sigma\in \Sigma_{n}$, $v\in\tensor{}{n}V$.\\
Then $\tensor{}{n}V=(\tensor{}{n}V)^{\Sigma_{n}}\oplus W$ and
$N\colon \symm{}{n}V\to (\tensor{}{n}V)^{\Sigma_{n}}$
is an isomorphism with inverse $\dfrac{\pi}{n!}$.
\end{lem}

\begin{proof} It is clear from the definition of $W$ that $\pi(W)=0$;
moreover $v-N\dfrac{\pi}{n!}v\in W$ for every
$v\in \tensor{}{n}V$, and therefore $\Image(N)+W=\tensor{}{n}V$.\\
On the other side
if $v$ is $\Sigma_{n}$-invariant then
\[v=\frac{1}{n!}\sum_{\sigma\in\Sigma_{n}}\sigma(v)=\frac{1}{n!}N\pi(v)\]
and therefore $\Image(N)=(\tensor{}{n}V)^{\Sigma_{n}}$,
$\Image(N)\cap W\subset \Image(N)\cap \ker(\pi)=0$.\end{proof}

For every $0\le a\le n$, the multiplication map
$V^{\otimes a}\bigotimes V^{\otimes n-a}\to V^{\otimes n}$ is an
isomorphism of graded vector spaces; we denote its inverse by
\[\copa_{a,n-a}\colon V^{\otimes n}\to
V^{\otimes a}\tensor{}{} V^{\otimes n-a},\]
\[
\copa_{a,n-a}(v_{1}\otimes\ldots\otimes v_{n})=
(v_{1}\otimes\ldots\otimes v_{a})\otimes
(v_{a+1}\otimes\ldots\otimes v_{n}).\]

The multiplication
$\mu\colon (\bigodot^{a} V)\otimes (\bigodot^{n-a} V)
\to \bigodot^{n}V$ is surjective but not injective; a left inverse
is given by
$\copl_{a,n-a}{\dbinom{n}{a}}^{-1}$, where
\[ \copl_{a,n-a}(v_{1}\odot\ldots\odot v_{n})=
\sum_{}
\epsilon\left(\genfrac{}{}{0pt}{0}{I,I^{c}}{
\{1,\ldots,n\}}; v_{1},\ldots,v_{n}\right)
v_{I}\otimes v_{I^{c}},\]
the sum is taken over all subsets $I\subset \{1,\ldots,n\}$ of
cardinality $a$ and $I^{c}$ is the complement of $I$ to 
$\{1,\ldots,n\}$.\\

\begin{defn}\label{XIII.1.4}
The
set of
\emph{unshuffles} of type $(p,q)$ is the subset
$S(p,q)\subset\Sigma_{p+q}$  of permutations
$\sigma$ such that $\sigma(i)<\sigma(i+1)$ for every $i\not=p$.
\end{defn}

Since $\sigma\in S(p,q)$ if and only if the restrictions
$\sigma\colon\{1,\ldots,p\}\to \{1,\ldots,p+q\}$,
$\sigma\colon\{p+1,\ldots,p+q\}\to \{1,\ldots,p+q\}$,
are increasing maps, it follows easily that
the unshuffles are
a set of representatives for the cosets of the canonical embedding of
$\Sigma_{p}\times \Sigma_{q}$ inside $\Sigma_{p+q}$. More precisely
for every 
$\sigma\in\Sigma_{p+q}$ there exists a unique decomposition
$\sigma=\tau\rho$ with $\tau\in S(p,q)$ and $\rho\in
\Sigma_{p}\times \Sigma_{q}$.\\

\begin{exer}\label{eseXIII.1}
Prove the formula
\[ \copl_{a,n-a}(v_{1}\odot\ldots\odot v_{n})=\!\!
\sum_{\sigma\in S(a,n-a)}\!\!\epsilon(\sigma)
(v_{\sigma(1)}\odot\ldots\odot v_{\sigma(a)})\otimes
(v_{\sigma(a+1)}\odot\ldots\odot v_{\sigma(n)})\]
\end{exer}

\begin{lem}\label{XIII.1.6}
In the above notation, for every $0\le a\le n$
\[ \copa_{a,n-a}N=(N\otimes N)\copl_{a,n-a}\colon
\symm{}{n} V\to \tensor{}{a}V\otimes\tensor{}{n-a}V.\]
\end{lem}
\begin{proof} Easy exercise.\end{proof}

Consider two graded vector spaces $V,M$ and a homogeneous linear
map $f\colon\tensor{}{m}V\to M$. The symmetrization $\widetilde{f}\colon
\symm{}{m}V\to M$ of $f$ is given by the formula
\[ \widetilde{f}(a_{1}\odot a_{2}\odot\ldots\odot a_{m})=
\sum_{\sigma\in\Sigma_{m}}\epsilon(V,\sigma;a_{1},\ldots,a_{m})
f(a_{\sigma_{1}}\otimes\ldots\otimes a_{\sigma_{m}}).\]

If $g\colon \tensor{}{l}V\to V$ is a homogeneous linear map of
degree $k$, the (non associative) Gerstenhaber composition product
$f\bullet g\colon \tensor{}{m+l-1}V\to M$ is defined as
\[\begin{split}
&f\bullet g(a_{1}\otimes\ldots\otimes a_{m+l-1})=\\
&=\sum_{i=0}^{m-1}(-1)^{k(\bar{a_{1}}+\ldots+\bar{a_{i}})}
f(a_{1}\otimes\ldots\otimes a_{i}\otimes
g(a_{i+1}\otimes\ldots\otimes a_{i+l})\otimes \ldots \otimes
a_{m+l-1}).
\end{split}\]
The behavior of $\bullet$ with respect to symmetrization is given
in the following lemma.

\begin{lem}\label{syle} \emph{(Symmetrization lemma)}
In the notation above
\[
\begin{split}
&\widetilde{f\bullet g}(a_{1}\odot\ldots\odot a_{m+l-1})=\\
&=\sum_{\sigma\in S(l,m-1)} \!\!\!\!\!\!\!\!\epsilon(V,\sigma;
a_{1},\ldots,a_{m})\widetilde{f}(\widetilde{g}(a_{\sigma_1} \odot\ldots\odot
a_{\sigma_l})\odot a_{\sigma_{l+1}}\odot\ldots\odot
a_{\sigma_{l+m-1}}).
\end{split}
\]
\end{lem}

\begin{proof} We give only some suggestion, leaving the details of
the proof as exercise. First, it is sufficient to prove the
formula in the "universal" graded vector space $U$ with
homogeneous basis $a_1,\ldots, a_{m+l-1}$ and $b_I$, where $I$
ranges over all injective maps $\{1,\ldots,l\}\to \{1,\ldots,
m+l-1\}$, $b_I$ is homogeneous of degree
$k+\bar{a_{I(1)}}+\ldots+\bar{a_{I(l)}}$ and
$g(a_I)=b_I$.\\
Second, by linearity we may assume that $M=\K$ and $f$ an element
of the dual basis of the standard basis of
$\tensor{}{m}U$.\\
With these assumption the calculation becomes easy.
\end{proof}

\bigskip
\section[~~Graded coalgebras]{Graded coalgebras}

\begin{defn}\label{XIII.2.2}
A coassociative $\Z$-graded coalgebra is the data of
a graded vector space $C=\oplus_{n\in\Z}C^{n}\in\mathbf{G}$ and of a
coproduct $\Delta\colon C\to C\otimes C$ such that:
\begin{itemize}
\item $\Delta$ is a morphism of graded vector spaces.
\item (coassociativity) $(\Delta\otimes Id_{C})\Delta=
(Id_{C}\otimes\Delta)\Delta\colon C\to C\otimes C\otimes C$.
\end{itemize}
The coalgebra is called \emph{cocommutative} if $T\Delta=\Delta$.\\
\end{defn}

For simplicity of notation, from now on with the term  \emph{graded
coalgebra} we intend a $\Z$-graded coassociative coalgebra.

\begin{defn}\label{XIII.2.4} Let $(C,\Delta)$ and $(B,\Gamma)$ be 
graded coalgebras.
A \emph{morphism} of graded coalgebras $f\colon C\to B$
is a morphism of graded vector spaces that commutes with coproducts, i.e. 
$\Gamma f=(f\otimes f)\Delta$.\\
The category of graded coalgebras is denoted by $\mathbf{GC}$.
\end{defn}

\begin{exer}\label{eseXIII.2}
A \emph{counity} of a graded coalgebra
is a morphism of graded vector spaces
$\epsilon\colon C\to\K$ such that
$(\epsilon\otimes Id_{C})\Delta=(Id_{C}\otimes
\epsilon)\Delta=Id_{C}$.\\
Prove that if a counity exists, then it is unique (Hint:
$(\epsilon\otimes\epsilon')\Delta=$?).
\end{exer}

\begin{ex}\label{XIII.2.6}
Let $C=\K[t]$ be the polynomial ring in one variable $t$
of even degree. A coalgebra structure is given by
\[ \Delta(t^{n})=\sum_{i=0}^{n}t^{i}\otimes t^{n-i}.\]
We left to the reader the verification of the coassociativity, of the
commutativity and the existence of the counity.\\
If the degree of $t$ is equal to $0$, then
for every sequence $\{f_{n}\}_{n>0}\subset\K$ it is associated a
morphism of coalgebras $f\colon C\to C$ defined as
\[ f(1)=1,\qquad f(t^{n})=\sum_{s=1}^{n}\sum_{
\genfrac{}{}{0pt}{1}{(i_1,\ldots,i_s)\in\N^s}{
i_{1}+\ldots+i_{s}=n}}
f_{i_{1}}f_{i_{2}}\ldots a_{i_{s}}t^{s}.\]
The verification that $\Delta f=(f\otimes f)\Delta$ can be done in the
following way: Let $\{ x^{n}\}\subset C\dual=\K[[x]]$ be the dual basis of
$\{t^{n}\}$.
Then for every $a,b,n\in N$ we have:
\[ \langle x^{a}\otimes x^{b},\Delta f(t^{n})\rangle=
\sum_{i_{1}+\ldots+i_{a}+j_{1}+\ldots+j_{b}=n}
f_{i_{1}}\ldots f_{i_{a}}f_{j_{1}}\ldots f_{j_{b}},\]
\[ \langle x^{a}\otimes x^{b},f\otimes f\Delta(t^{n})\rangle=
\sum_{s}\sum_{i_{1}+\ldots+i_{a}=s\,}\sum_{\, j_{1}+\ldots+j_{b}=n-s}
f_{i_{1}}\ldots f_{i_{a}}f_{j_{1}}\ldots f_{j_{b}}.\]
Note that the sequence $\{f_{n}\}$, $n\ge 1$, can be recovered from
$f$ by the formula $f_{n}=\langle x,f(t^{n})\rangle$.\\
We shall prove later that every coalgebra endomorphism of $\K[t]$ has
this form for some sequence $\{f_{n}\}$, $n\ge 1$.
\end{ex}

\begin{lemmadef}\label{XIII.2.8}
Let $(C,\Delta)$ be a graded coassociative coalgebra,
we define recursively $\Delta^{0}=Id_{C}$ and, for $n>0$,
$\Delta^{n}=(Id_{C}\otimes \Delta^{n-1})\Delta\colon C\to
\bigotimes^{n+1}C$. Then:
\begin{enumerate}
    \item For every $0\le a\le n-1$ we have
\[\Delta^{n}=(\Delta^{a}\otimes \Delta^{n-1-a})\Delta\colon C\to
\hbox{$\bigotimes$}^{n+1}C,\]
\[\copa_{a+1,n-a}\Delta^{n}=
(\Delta^{a}\otimes \Delta^{n-1-a})\Delta\] 

    \item  For every $s\ge 1$ and every $a_{0},\ldots,a_{s}\ge 0$ we 
    have
    \[ (\Delta^{a_{0}}\otimes\Delta^{a_{1}}\otimes\ldots\otimes 
    \Delta^{a_{s}})\Delta^{s}=\Delta^{s+\sum a_{i}}.\]
In particular, if $C$ is cocommutative then the image of $\Delta^{n-1}$ is contained 
in the set of $\Sigma_n$-invariant elements of $\bigotimes^nC$.

\item If $f\colon (C,\Delta)\to (B,\Gamma)$ is a morphism of graded 
coalgebras then, for every $n\ge 1$ we have 
\[ \Gamma^{n}f=(\otimes^{n+1}f)\Delta^{n}\colon C\to 
\hbox{$\bigotimes$}^{n+1}B.\]

\end{enumerate}
\end{lemmadef}
\begin{proof} \vale{1} If $a=0$ or $n=1$ there is nothing to prove, thus we can
assume $a>0$ and use induction on $n$.
we have:
\[(\Delta^{a}\otimes \Delta^{n-1-a})\Delta=
((Id_{C}\otimes\Delta^{a-1})\Delta\otimes \Delta^{n-1-a})\Delta=\]
\[=(Id_{C}\otimes\Delta^{a-1}\otimes\Delta^{n-1-a})(\Delta\otimes
Id_{C})\Delta=\]
\[=(Id_{C}\otimes\Delta^{a-1}\otimes\Delta^{n-1-a})(Id_{C}\otimes
\Delta)\Delta=
(Id_{C}\otimes(\Delta^{a-1}\otimes\Delta^{n-1-a})
\Delta)\Delta=
\Delta^{n}.\]
\vale{2} Induction on $s$, being the case $s=1$ proved in item 1. 
If $s\ge 2$ we can write
\[ (\Delta^{a_{0}}\otimes\Delta^{a_{1}}\otimes\ldots\otimes 
    \Delta^{a_{s}})\Delta^{s}=
    (\Delta^{a_{0}}\otimes\Delta^{a_{1}}\otimes\ldots\otimes 
    \Delta^{a_{s}})(Id\otimes\Delta^{s-1})\Delta=\]
    \[ (\Delta^{a_{0}}\otimes(\Delta^{a_{1}}\otimes\ldots\otimes 
    \Delta^{a_{s}})\Delta^{s-1})\Delta=(\Delta^{a_{0}}\otimes
    \Delta^{s-1+\sum_{i>0} a_{i}})\Delta=\Delta^{s+\sum a_{i}}.\]
The action of $\Sigma_n$ on $\bigotimes^nC$ is generated by the operators
$T_a=Id_{\bigotimes^aC}\otimes T\otimes Id_{\bigotimes^{n-a-2}C}$, $0\le a\le n-2$,
and, if $T\Delta=\Delta$ then 
\[T_a\Delta^{n-1}=T_a(Id_{\bigotimes^aC}\otimes\Delta\otimes 
Id_{\bigotimes^{n-a-2}C})\Delta^{n-2}=\]
\[=(Id_{\bigotimes^aC}\otimes\Delta\otimes 
Id_{\bigotimes^{n-a-2}C})\Delta^{n-2}=\Delta^{n-1}.\]
\vale{3} By induction on $n$, 
\[ \Gamma^{n}f=(Id_{B}\otimes\Gamma^{n-1})\Gamma f=
(f\otimes\Gamma^{n-1}f)\Delta=(f\otimes (\otimes^{n}f)\Delta^{n-1})\Delta=
(\otimes^{n+1}f)\Delta^{n}.\]
\end{proof}

\begin{ex}\label{XIII.2.10}
Let $A$ be a graded associative algebra with product
$\mu\colon A\otimes A\to A$ and $C$ a graded coassociative
coalgebra with coproduct $\Delta\colon C\to C\otimes C$.\\
Then $\Hom^{*}(C,A)$ is a graded associative algebra with product
\[ fg=\mu(f\otimes g)\Delta.\]
We left as an exercise the verification that
the product in $\Hom^{*}(C,A)$ is associative.\\
In particular $\Hom_{\mathbf{G}}(C,A)=\Hom^{0}(C,A)$ is an associative
algebra and $C\dual=\Hom^{*}(C,\K)$ is a graded associative algebra.
(Notice that  in general $A\dual$ is not a coalgebra.)
\end{ex}

\begin{ex}\label{XIII.2.12}
The dual of the coalgebra $C=\K[t]$ (Example~\ref{XIII.2.6}) is
exactly the algebra of formal power series $A=\K[[x]]=C\dual$. Every coalgebra
morphism $f\colon C\to C$ induces a local homomorphism of $\K$-algebras
$f^{t}\colon A\to A$. Clearly $f^{t}=0$ only if $f=0$,  $f^{t}$ is
uniquely determined by $f^{t}(x)=\sum_{n>0}f_{n}x^{n}$ and then
every morphism of coalgebras $f\colon C\to C$ is uniquely determined by the
sequence $f_{n}=\langle f^{t}(x), t^{n}\rangle=\langle x,
f(t^{n})\rangle$.\\
The map $f\mapsto f^{t}$ is functorial and then preserves the composition
laws.\end{ex}

\begin{defn}\label{XIII.2.14}
A graded coassociative coalgebra $(C,\Delta)$ is called
\emph{nilpotent} if $\Delta^{n}=0$ for $n>>0$.\\
It is called \emph{locally nilpotent} if it is the direct limit of
nilpotent graded coalgebras or equivalently if
$C=\cup_{n}\ker\Delta^{n}$.
\end{defn}

\begin{ex} The coalgebra $\K[t]$ of Example~\ref{XIII.2.6} is locally
nilpotent.\end{ex}

\begin{ex}\label{XIII.2.16}
Let $A=\oplus A_{i}$ be a finite dimensional graded associative
commutative $\K$-algebra and
let $C=A\dual=\Hom^{*}(A,\K)$ be its graded dual.\\
Since $A$ and $C$ are finite dimensional,
the pairing $\langle c_{1}\otimes c_{2},a_{1}\otimes
a_{2}\rangle=(-1)^{\bar{a_{1}}\,\bar{c_{2}}}
\langle c_{1},a_{1}\rangle\langle c_{2},a_{2}\rangle$
gives a natural isomorphism
$C\otimes C=(A\otimes A)\dual$ commuting with the twisting maps $T$;
we may define
$\Delta$ as the transpose of the multiplication map $\mu\colon
A\otimes A\to A$.\\
Then $(C,\Delta)$ is a coassociative cocommutative coalgebra. Note
that $C$ is nilpotent if and only if $A$ is nilpotent.\\
\end{ex}

\begin{exer} Let $(C,\Delta)$ be a graded coalgebra and $p\colon C\to V$ a 
morphism of graded vector spaces. We shall say that $p$ \emph{
cogenerates} $C$ if for every $c\in C$ there exists $n\ge 0$ such 
that $(\otimes^{n+1}p)\Delta^{n}(c)\not=0$ in $\bigotimes^{n+1}V$.\\
Prove that every morphism of graded coalgebras $B\to C$ is uniquely 
determined by its composition $B\to C\to V$ with a cogenerator $p$.
\end{exer}

\medskip
\subsection[~~The reduced tensor coalgebra]{
The reduced tensor coalgebra}

Given a graded vector space $V$, we denote
$\bar{T(V)}=\bigoplus_{n>0}\bigotimes^{n}V$.
When considered as a subset of $T(V)$ it becomes
an  ideal of the tensor algebra  generated by $V$.\\
The \emph{reduced tensor coalgebra} generated by $V$ is the graded
vector space $\bar{T(V)}$ endowed with
the coproduct $\copa\colon \bar{T(V)}\to \bar{T(V)}\otimes\bar{T(V)}$, 
\[ \copa=\sum_{n=1}^{\infty}\sum_{a=1}^{n-1}\copa_{a,n-a},\qquad
\copa(v_{1}\otimes\ldots\otimes v_{n})=
\sum_{r=1}^{n-1}(v_{1}\otimes\ldots\otimes v_{r})\otimes
(v_{r+1}\otimes\ldots\otimes v_{n})\]
The coalgebra $(\bar{T(V)},\copa)$ is coassociative (but not cocommutative)
and locally nilpotent; in fact,
for every $s>0$,
\[ \copa^{s-1}(v_{1}\otimes\ldots\otimes v_{n})=\sum_{1\le
i_{1}<i_{2}<\ldots<i_{s}=n}(v_{1}\otimes\ldots\otimes v_{i_{1}})\otimes
\ldots\otimes(v_{i_{s-1}+1}\otimes\ldots\otimes v_{i_{s}})\]
and then
$\ker\copa^{s-1}=\bigoplus_{n=1}^{s-1}\bigotimes^{n}V$.\\
If $\mu\colon \bigotimes^{s}\bar{T(V)}\to\bar{T(V)}$ denotes the
multiplication map then, for every $v_{1},\ldots,v_{n}\in V$, we have
\[\mu\copa^{s-1}(v_{1}\otimes\ldots\otimes
v_{n})={\binom{n-1}{s-1}}v_{1}\otimes\ldots\otimes v_{n}.\]

For every morphism of graded vector spaces $f\colon V\to W$
the induced morphism of graded algebras
\[T(f)\colon\bar{T(V)}\to
\bar{T(W)},\qquad  T(f)(v_{1}\otimes\ldots\otimes v_{n})=
f(v_{1})\otimes\ldots\otimes f(v_{n})\]
is also a morphism of graded coalgebras.

\begin{exer}\label{eseXIII.3}
Let $p\colon T(V)\to \bar{T(V)}$ be the projection
with kernel $\K=\bigotimes^{0}V$ and 
$\phi\colon T(V)\to T(V)\otimes T(V)$ the unique
homomorphism of graded algebras such that $\phi(v)=v\otimes 1+1\otimes v$ for
every $v\in V$. Prove that $p\phi=\copa p$.
\end{exer}

If $(C,\Delta)$ is locally nilpotent then, for every $c\in C$,
there exists $n>0$ such that $\Delta^{n}(c)=0$ and then it is defined a morphism of
graded vector spaces
\[ \frac{1}{1-\Delta}=\sum_{n=0}^{\infty}\Delta^{n}\colon C\to \bar{T(C)}.\]

\begin{prop}\label{XIII.2.18}
Let $(C,\Delta)$ be a  locally nilpotent
graded coalgebra, then:
\begin{enumerate}
\item  The map $\ds\frac{1}{1-\Delta}=\sum_{n\ge 0}\Delta^{n}
\colon C\to \bar{T(C)}$ is a
morphism of graded coalgebras.
\item  For every graded vector space $V$ and every morphism of
graded  coalgebras $\phi\colon C\to \bar{T(V)}$,
there exists a unique  morphism of graded vector spaces $f\colon C\to
V$ such that $\phi$ factors as
\[ \phi=T(f)\ds\frac{1}{1-\Delta}
=\sum_{n=1}^{\infty}(\otimes^{n}f)\Delta^{n-1}
\colon C\to \bar{T(C)}\to
\bar{T(V)}.\]
\end{enumerate}
\end{prop}
\begin{proof}
{[1]} We have
\[
\begin{split}
\left(\left(\sum_{n\ge 0}\Delta^{n}\right)\otimes \left(\sum_{n\ge 0}\Delta^{n}\right)\right)\Delta&=
\sum_{n\ge 0}\sum_{a=0}^{n}(\Delta^{a}\otimes \Delta^{n-a})\Delta\\
&=\sum_{n\ge 0}\sum_{a=0}^{n}\copa_{a+1,n+1-a}\Delta^{n+1}=
\copa \left(\sum_{n\ge 0}\Delta^{n}\right)
\end{split}
\]
where in the last equality we have used the relation
$\copa\Delta^{0}=0$.\\
{[2]}
The unicity of $f$ is clear, since by the formula
$\phi=T(f)(\sum_{n\ge 0} \Delta^{n})$ it follows that $f$ is the
composition of $\phi$ and the projection $\bar{T(V)}\to V$.\\
To prove the existence of the factorization,
take any morphism of graded coalgebras
$\phi\colon C\to \bar{T(V)}$ and denote by $\phi^{i}\colon C\to 
\tensor{}{i}V$ the composition of $\phi$ with the projection.
It is sufficient to show that for every $n>1$, $\phi^{n}$ is uniquely 
determined by $\phi^{1}$.
Now, the morphism condition
$\copa\phi=(\phi\otimes\phi)\Delta$ composed with the projection 
$\bar{T(V)}\otimes\bar{T(V)}\to 
\bigoplus_{i=1}^{n-1}(\bigotimes^{i}V\otimes\bigotimes^{n-1}V)$
gives the equality
\[\copa\phi^{n}=\sum_{i=1}^{n-1}(\phi^{i}\otimes\phi^{n-i})\Delta,
\qquad n\ge 2.\]
Using induction on $n$, it is enough
to observe that the restriction of $\copa$ to
$\bigotimes^{n}V$ is injective for every $n\ge 2$.\\
\end{proof}

It is useful to restate part of the
Proposition~\ref{XIII.2.18} in the following form

\begin{cor}\label{XIII.2.20}
Let $V$ be a fixed graded vector space; for every locally
nilpotent graded coalgebra $C$ the composition with the
projection $\bar{T(V)}\to V$ induces a bijection
\[ \Hom_{\mathbf{GC}}(C,\bar{T(V)})=\Hom_{\mathbf{G}}(C,V).\]
\end{cor}

When $C$ is a reduced tensor coalgebra, Proposition~\ref{XIII.2.18}
takes the following more explicit form

\begin{cor}\label{XIII.2.22}
Let $U,V$ be  graded vector spaces and $p\colon\bar{T(V)}\to V$
the projection. Given $f\colon \bar{T(U)}\to V$, the linear map
$F\colon \bar{T(U)}\to \bar{T(V)}$
\[ F(v_{1}\otimes\ldots\otimes v_{n})=\sum_{s=1}^{n}\sum_{\, 1\le
i_{1}<i_{2}<\ldots<i_{s}=n}f(v_{1}\otimes\ldots\otimes
v_{i_{1}})\otimes
\ldots\otimes f(v_{i_{s-1}+1}\otimes\ldots\otimes v_{i_{s}})\]
is the unique morphism of graded coalgebras such that
$pF=f$.\end{cor}

\begin{ex}\label{XIII.2.24}
Let $A$ be an associative graded algebra. Consider the
projection $p\colon \bar{T(A)}\to A$, the
multiplication map $\mu\colon \bar{T(A)}\to A$ and its conjugate
\[\mu^{*}=-\mu T(-1),\qquad \mu^{*}(a_{1}\otimes \ldots\otimes
a_{n})=(-1)^{n-1}\mu(a_{1}\otimes \ldots\otimes
a_{n})=(-1)^{n-1}a_{1}a_{2}\ldots a_{n}.\]
The two coalgebra
morphisms $\bar{T(A)}\to \bar{T(A)}$ induced by $\mu$ and $\mu^*$
are isomorphisms, the one
inverse of the other.\\
In fact, the coalgebra morphism
$ F\colon \bar{T(A)}\to \bar{T(A)}$
\[F(a_{1}\otimes\ldots\otimes  a_{n})=\sum_{s=1}^{n}\sum_{\, 1\le
i_{1}<i_{2}<\ldots<i_{s}=n}(a_{1}a_{2}\ldots a_{i_{1}})\otimes
\ldots\otimes (a_{i_{s-1}+1}\ldots a_{i_{s}})\]
is induced by $\mu$ (i.e. $pF=\mu$),
$\mu^{*}F(a)=a$ for every $a\in A$
and for every $n\ge 2$
\[ \mu^{*}F(a_{1}\otimes\ldots\otimes a_{n})=
\sum_{s=1}^{n}(-1)^{s-1}{\hskip-10pt}\sum_{1\le
i_{1}<i_{2}<\ldots<i_{s}=n}{\hskip-18pt}a_{1}a_{2}\ldots a_{n}=\]
\[=\sum_{s=1}^{n}(-1)^{s-1}{\binom{n-1}{s-1}}a_{1}a_{2}\ldots a_{n}=
\left(\sum_{s=0}^{n-1}(-1)^{s}{\binom{n-1}{s}}\right)
a_{1}a_{2}\ldots a_{n}=0.\]
This implies that $\mu^{*}F=p$ and therefore, if 
$F^{*}\colon \bar{T(A)}\to \bar{T(A)}$ is induced by $\mu^{*}$ then 
$pF^{*}F=\mu^{*}F=p$ and by Corollary \ref{XIII.2.20} $F^{*}F$ is the 
identity.\end{ex}

\begin{exer}\label{eseXIII.4}
Let $A$ be an associative graded algebra over the field $\K$,
for every local homomorphism of $\K$-algebras
$\gamma\colon \K[[x]]\to \K[[x]]$, $\gamma(x)=\sum\gamma_{n}x^{n}$,
we can associate a coalgebra
morphism $F_{\gamma}\colon \bar{T(A)}\to \bar{T(A)}$ induced by the
linear map
\[ f_{\gamma}\colon \bar{T(A)}\to A,\qquad f(a_{1}\otimes \ldots\otimes
a_{n})=\gamma_{n}a_{1}\ldots a_{n}.\]
Prove the composition formula $F_{\gamma\delta}=F_{\delta}F_{\gamma}$.
(Hint: Example~\ref{XIII.2.12}.)
\end{exer}

\begin{exer}\label{eseXIII.5}
A graded coalgebra morphism  $F\colon \bar{T(U)}\to \bar{T(V)}$ is surjective
(resp.: injective, bijective) if and only if $F^{1}_{1}\colon U\to V$
is surjective (resp.: injective, bijective).
\end{exer}

\medskip
\subsection[~~The reduced symmetric coalgebra]{
The reduced symmetric coalgebra}

\begin{defn}\label{XIII.2.30}
The reduced symmetric coalgebra is by definition
$\bar{S(V)}=\bigoplus_{n>0}\bigodot^{n}V$, with the coproduct
$\copl=\sum_{n}\sum_{i=0}^{n-1}\copl^{n+1}_{i+1}$,
\[\copl(v_{1}\odot\ldots\odot v_{n})=
\sum_{r=1}^{n-1}
\sum_{I\subset \{1,\ldots,n\}; |I|=r}
\epsilon\left(\genfrac{}{}{0pt}{0}{I,I^{c}}{
\{1,\ldots,n\}}; v_{1},\ldots,v_{n}\right)
v_{I}\otimes v_{I^{c}}.\]
\end{defn}

The verification that $\copl$ is a coproduct is an easy consequence of  
Lemma~\ref{XIII.1.6}. In fact, the injective map
$N\colon\bar{S(V)}\to\bar{T(V)}$ satisfies the equality
$\copa N=(N\otimes N)\copl$ and then $N$ is an isomorphism between 
$(\bar{S(V)},\copl)$ and the  subcoalgebra of symmetric tensors of  
$(\bar{T(V)},\copa)$.\\

\begin{rem} It is  often convenient to think the symmetric 
algebra as a quotient of the tensor algebra and the symmetric 
coalgebra as a subset of the tensor coalgebra.\end{rem}

The coalgebra $\bar{S(V)}$ is coassociative without
counity. It follows from the definition of $\copl$ that 
$V=\ker\copl$ and $T\copl=\copl$, where $T$ is the twisting map; in 
particular $(\bar{S(V)},\copl)$ is cocommutative. For every 
morphism of graded vector spaces $f\colon V\to W$, the morphism 
$S(f)\colon \bar{S(V)}\to \bar{S(W)}$ is a morphism of graded 
coalgebras.\\

If $(C,\Delta)$ is any cocommutative graded coalgebra, 
then the image of $\Delta^{n}$ is contained in
the subspace of symmetric tensors and therefore
\[\frac{1}{1-\Delta}=N\circ \frac{e^{\Delta}-1}{\Delta},\]
where
\[ \frac{e^{\Delta}-1}{\Delta}=
\sum_{n=1}^{\infty}\frac{\pi}{n!}\Delta^{n-1}\colon C\to
\bar{S(C)}.\]

\begin{prop}\label{XIII.2.32}
Let $(C,\Delta)$ be a  cocommutative
locally nilpotent  graded coalgebra, then:
\begin{enumerate}
\item  The map $\ds\frac{e^{\Delta}-1}{\Delta}\colon C\to \bar{S(C)}$ is a
morphism of graded coalgebras.

\item  For every graded vector space $V$ and every morphism of
graded coalgebras $\phi\colon C\to \bar{S(V)}$, there exists a unique
factorization
\[ \phi=S(\phi^{1})\ds\frac{e^{\Delta}-1}{\Delta}
=\sum_{n=1}^{\infty}\frac{\bigodot^{n}\phi^{1}}{n!}\Delta^{n-1}
\colon C\to
\bar{S(C)}\to \bar{S(V)},\]
where $\phi^{1}\colon C\to V$ is a morphism of
graded vector spaces $f\colon C\to V$. (Note that $\phi^{1}$ is
the composition of $\phi$ and the projection $\bar{S(V)}\to V$.)
\end{enumerate}
\end{prop}

\begin{proof}
Since $N$ is an injective morphism of coalgebras and
$\ds\frac{1}{1-\Delta}=N\circ \frac{e^{\Delta}-1}{\Delta}$, the proof
follows immediately from Proposition~\ref{XIII.2.18}.
\end{proof}

\begin{cor}\label{XIII.2.34}
Let $C$ be a locally nilpotent cocommutative
graded coalgebra, and $V$ a graded vector space.
A morphism
$\theta\in \Hom_{\mathbf{G}}(C,\bar{S(V)})$ is a morphism of graded coalgebras if
and only if there exists
$m\in \Hom_{\mathbf{G}}(C,V)\subset \Hom_{\mathbf{G}}(C,\bar{S(V)})$
such that
\[\theta=\exp(m)-1\:=\sum_{n=1}^{\infty}\frac{1}{n!}m^{n},\]
being the $n$-th power of $m$ is considered with respect to the algebra
structure on $\Hom_{\mathbf{G}}(C,\bar{S(V)})$ (Example~\ref{XIII.2.10}).
\end{cor}

\begin{proof} An easy computation gives the  formula
$m^{n}=(\symm{}{n}m)\Delta^{n-1}$ for the
product defined in Example~\ref{XIII.2.10}.
\end{proof}

\begin{exer}\label{eseXIII.6}
Let $V$ be a graded vector space. Prove that the formula
\[\copc(v_{1}\wedge\ldots\wedge v_{n})=
\sum_{r=1}^{n-1}\sum_{\sigma\in
S(r,n-r)}(-1)^{\sigma}\epsilon(\sigma)
(v_{\sigma(1)}\wedge\ldots\wedge v_{\sigma(r)})\otimes
(v_{\sigma(r+1)}\wedge\ldots\wedge v_{\sigma(n)}),\]
where $(-1)^{\sigma}$ is the signature of the permutation $\sigma$,
defines a coproduct on $\bar{\bigwedge(V)}=\bigoplus_{n>0}\bigwedge^{n}V$.
The resulting coalgebra is called
\emph{reduced exterior coalgebra} generated by $V$.
\end{exer}

\bigskip

\section[~~Coderivations]{Coderivations}

\begin{defn}\label{XIII.3.2}
Let $(C,\Delta)$ be a graded coalgebra.
A linear map  $d\in\Hom^{n}(C,C)$
is called a \emph{coderivation of degree $n$}
if it satisfies the \emph{coLeibnitz rule}
\[\Delta d=(d\otimes Id_{C}+Id_{C}\otimes d)\Delta.\]
A coderivation $d$ is called a \emph{codifferential} if $d^{2}=d\circ d=0$.\\
More generally, if $\theta\colon C\to D$ is a morphism of graded
coalgebras,
a morphism of graded vector spaces $d\in\Hom^{n}(C,D)$ is called a
coderivation of degree $n$ (with respect to $\theta$) if
\[\Delta_{D} d=(d\otimes \theta+\theta\otimes d)\Delta_{C}.\]
\end{defn}

In the above definition we have adopted the Koszul sign convention: 
i.e. if $x,y\in C$, $f,g\in \Hom^{*}(C,D)$,
$h,k\in \Hom^{*}(B,C)$ are homogeneous then
$(f\otimes g)(x\otimes y)=
(-1)^{\bar{g}\,\bar{x}}f(x)\otimes g(y)$ and
$(f\otimes g)(h\otimes k)=(-1)^{\bar{g}\,\bar{h}}
fh\otimes gk$.\\

The coderivations of degree $n$ with respect a coalgebra morphism
$\theta\colon C\to D$
form a vector space
denoted  $\Coder^{n}(C,D;\theta)$.\\
For simplicity of notation we denote $\Coder^{n}(C,C)=
\Coder^{n}(C,C;Id)$.

\begin{lem}\label{XIII.3.4}
Let $C\mapor{\theta}D\mapor{\rho}E$ be
morphisms of graded coalgebras.
The compositions with $\theta$ and $\rho$ induce linear
maps
\[ \rho_{*}\colon\Coder^{n}(C,D;\theta)\to \Coder^{n}(C,E;\rho\theta),
\qquad f\mapsto \rho f;\]
\[ \theta^{*}\colon \Coder^{n}(D,E;\rho)\to \Coder^{n}(C,E;\rho\theta),
\qquad f\mapsto f\theta.\]
\end{lem}

\begin{proof} Immediate consequence of the equalities
\[ \Delta_{E}\rho=(\rho\otimes\rho)\Delta_{D},\qquad
\Delta_{D}\theta=(\theta\otimes\theta)\Delta_{C}.\]
\end{proof}

\begin{exer}\label{eseXIII.7}
Let $C$ be a graded  coalgebra and
$d\in\Coder^{1}(C,C)$ a codifferential of degree 1.
Prove that the triple
$(L,\delta,[,])$, where:
\[ L=\somdir{n\in\Z}{}\Coder^{n}(C,C),\quad
[f,g]=fg-(-1)^{\bar{g}\,\bar{f}}gf,\quad \delta(f)=[d,f]\]
is a differential graded Lie algebra.
\end{exer}

\begin{lem}\label{XIII.3.6}
Let $V,W$ be graded vector spaces, $f\in\Hom_{\mathbf{G}}(V,W)$ and
$g\in \Hom^{m}(\bar{S(V)},W)$. Then the morphism
$d\in \Hom^{m}(\bar{S(V)},\bar{S(W)})$ defined by the rule
\[ d(v_{1}\odot\ldots\odot v_{n})=
\sum_{\emptyset\not=I\subset\{1,\ldots,n\}}
\epsilon\left(\genfrac{}{}{0pt}{0}{I,I^{c}}{
\{1,\ldots,n\}}; v_{1},\ldots,v_{n}\right)
g(v_{I})\odot S(f)(v_{I^{c}})\]
is a coderivation of degree $m$ with respect to the morphism of 
graded coalgebras
$S(f)\colon \bar{S(V)}\to\bar{S(W)}$.\end{lem}

\begin{proof} Let $v_{1},v_{2},\ldots,v_{n}$ be fixed homogeneous
elements of $V$, we need to prove that
\[\copl d(v_{1}\odot\ldots\odot v_{n})=
(d\otimes S(f)+S(f)\otimes d)\copl(v_{1}\odot\ldots\odot v_{n}).\]
If $A\subset W$ is the image of $f$ and $B\subset W$ is the image of
$g$, it is not restrictive to assume that
$W=A\oplus B$: in fact we
can always factorize
\[ \xymatrix{V\ar[dr]_{(f,0)}\ar[drr]^{f}&&\\
&A\oplus B\ar[r]^{+}&W\\
\bar{S(V)}\ar[ur]^{(0,g)}\ar[urr]^{g}&&}\]
and apply  Lemma~\ref{XIII.3.4} to the coalgebra morphism
$S(+)\colon \bar{S(A\oplus B)}\to \bar{S(W)}$.\\
Under this assumption we have 
$(S(A)B\otimes\bar{S(A)})\cap (\bar{S(A)}\otimes S(A)B)=\emptyset$
and the image of $d$ is contained in ${S(A)}B\subset\bar{S(A\oplus 
B)}$. Therefore
the images of $\copl d$ and $(d\otimes S(f)+S(f)\otimes d)\copl$
are both contained in
$(S(A)B\otimes\bar{S(A)})\oplus (\bar{S(A)}\otimes S(A)B)$.\\
Denoting by $p\colon \bar{S(W)}\otimes\bar{S(W)}\to
S(A)B\otimes\bar{S(A)}$ the natural projection induced by the 
decomposition $W=A\oplus B$,
since both the operators
$\copl d$ and $(d\otimes S(f)+S(f)\otimes d)\copl$ are
invariant under the twisting map,
it is sufficient to prove that
\[p\copl d(v_{1}\odot\ldots\odot v_{n})=
p(d\otimes S(f))\copl(v_{1}\odot\ldots\odot v_{n}).\]
We have (all Koszul signs are referred to $v_{1},\ldots,v_{n}$)
\[p\copl d(v_{1}\odot\ldots\odot v_{n})=
p\copl\left(
\sum_{\emptyset\not=J\subset\{1,\ldots,n\}}
\epsilon\left(\genfrac{}{}{0pt}{0}{J,J^{c}}{
\{1,\ldots,n\}}\right)
g(v_{J})\odot S(f)(v_{J^{c}})\right)=\]
\[ =\sum_{\emptyset\not=J\subset I\subset \{1,\ldots,n\}}
\epsilon\left(\genfrac{}{}{0pt}{0}{J,J^{c}}{
\{1,\ldots,n\}}\right)
\epsilon\left(\genfrac{}{}{0pt}{0}{J,I-J,I^{c}}{
J,J^{c}}\right)
g(v_{J})\odot S(f)(v_{I-J})\otimes S(f)(v_{I^{c}})=\]
\[=\sum_{\emptyset\not=J\subset I\subset \{1,\ldots,n\}}
\epsilon\left(\genfrac{}{}{0pt}{0}{J,I-J,I^{c}}{
\{1,\ldots,n\}}\right)
g(v_{J})\odot S(f)(v_{I-J})\otimes S(f)(v_{I^{c}}).\]
On the other hand
\[ p(d\otimes S(f))\copl (v_{1}\odot\ldots\odot v_{n})=
p(d\otimes S(f))\left(\sum_{I}\epsilon
\binom{I,I^{c}}{\{1,\ldots,n\}}
v_{I}\otimes v_{I^{c}}\right)=\]
\[=\sum_{J\subset I}\epsilon
\left(\genfrac{}{}{0pt}{0}{I,I^{c}}{
\{1,\ldots,n\}}\right)
\epsilon\left(\genfrac{}{}{0pt}{0}{J,I-J,I^{c}}{
I,I^{c}}\right) g(v_{J})\odot S(f)(v_{I-J})\otimes S(f)(v_{I^{c}})=\]
\[=\sum_{J\subset I}\epsilon\left(\genfrac{}{}{0pt}{0}{J,I-J,I^{c}}{
\{1,\ldots,n\}}\right)
g(v_{J})\odot S(f)(v_{I-J})\otimes S(f)(v_{I^{c}}).\]
\end{proof}

\begin{prop}\label{XIII.3.8}
Let $V$ be a graded vector space and $C$ a locally nilpotent
cocommutative coalgebra. Then
for every coalgebra morphism $\theta\colon C\to \bar{S(V)}$ and every
integer $n$,
the composition with the projection $\bar{S(V)}\to V$
gives an isomorphism
\[ \Coder^{n}(C,\bar{S(V)};\theta)\to \Hom^{n}(C,V)=\Hom_{\mathbf{G}}(C,V[n]).\]
\end{prop}
\begin{proof}
The injectivity is proved essentially in the same way as in
Proposition~\ref{XIII.2.18}: if $d\in \Coder^{n}(C,\bar{S(V)};\theta)$ we
denote by  $\theta^{i},d^{i}\colon C\to \symm{}{i}V$ the composition 
of $\theta$ and $d$ with the projection $\bar{S(V)}\to \symm{}{i}V$.
The coLeibnitz rule is equivalent to the countable set of equalities
\[\copl^{i}_{a}d^{i}=d^{a}\otimes\theta^{i-a}+\theta^{a}\otimes
d^{i-a},\qquad 0<a<i.\]
Induction on $i$ and the injectivity of
\[ \copl\colon\somdir{m=2}{n}\symm{}{m}V\to
\tensor{}{2}\left(\somdir{m=1}{n-1}\symm{}{m}V\right)\]
show that $d$ is uniquely determined by $d^{1}$.\\
For the surjectivity, consider
$g\in \Hom^{n}(C,V)$;  according to Proposition~\ref{XIII.2.32}
we can write
$\theta=S(\theta^{1})\ds\frac{e^{\Delta}-1}{\Delta}$ and, by
Lemma~\ref{XIII.3.6},
the map $d=\delta\ds\frac{e^{\Delta}-1}{\Delta}$, where
$\delta\colon \bar{S(C)}\to \bar{S(V)}$ is given by
\[
\delta(c_{1}\odot\ldots\odot c_{n})=
\sum_{i\in\{1,\ldots,n\}}
\epsilon\left(\genfrac{}{}{0pt}{0}{\{i\},\{i\}^{c}}{
\{1,\ldots,n\}}; c_{1},\ldots,c_{n}\right)
g(c_{i})\odot S(\theta^{1})(c_{\{i\}^{c}})\]
is a coderivation of degree $n$ with respect to $\theta$ that lifts
$g$.
\end{proof}

\begin{cor}\label{XIII.3.10}
Let $V$ be a graded vector space, $\bar{S(V)}$ its reduced symmetric
coalgebra. The application $Q\mapsto Q^{1}$ gives an isomorphism of
vector spaces
\[\Coder^{n}(\bar{S(V)},\bar{S(V)})
=\Hom^{n}(\bar{S(V)},V)\]
whose inverse is given by the formula
\[Q(v_{1}\odot\ldots\odot v_{n})=
\sum_{k=1}^{n}\,\sum_{\sigma\in S(k,n-k)}\!\!\!\!
\epsilon(\sigma)
Q^{1}_{k}(v_{\sigma(1)}\odot\ldots\odot v_{\sigma(k)})\odot
v_{\sigma(k+1)}\odot\ldots\odot v_{\sigma(n)}.\]
In particular for every coderivation $Q$ we have
$Q^{i}_{j}=0$ for every $i>j$ and then the subcoalgebras
$\bigoplus_{i=1}^{r}\bigodot^{i}V$ are preserved by $Q$.
\end{cor}
\begin{proof} The isomorphism follows from Proposition~\ref{XIII.3.8},
while the inverse formula comes from Lemma~\ref{XIII.3.6}.
\end{proof}


\chapter{$L_{\infty}$ and EDF tools}
\label{CAP:TOOLS}
\piede

In this chapter we introduce the category $\mathbf{L_{\infty}}$ of
$L_{\infty}$-algebras and we define a sequence of natural transformations
\[ \mathbf{DGLA}\to \mathbf{L_{\infty}}\to \mathbf{PreDef}\to \mathbf{Def}\]
whose composition is the functor  $L\mapsto\Def_{L}$ (cf. \ref{VIII.5.8}).\\
In all the  four categories there is a  notion of
quasi-isomorphism which is preserved by the above natural
transformations: we recall that
in the category $\mathbf{Def}$ quasi-isomorphism means
isomorphism in tangent spaces and then by Corollary~\ref{VIII.6.4}
every quasi-isomorphism is an isomorphism.\\
Through all the chapter we work over a fixed field $\K$ of characteristic
0. Unless otherwise specified all the tensor products are made over
$\K$.\\

\bigskip

\section{Displacing (D\'ecalage)}

For every $n$ and every graded vector space $V$,
the twisting map gives a natural isomorphism
\[\decala_n\colon \tensor{}{n}(V[1])\to (\tensor{}{n})V[n],\qquad 
V[a]=\K[a]\otimes V\]
\[\decala_n(v_1[1]\otimes\ldots\otimes v_n[1])=
(-1)^{\sum_{i=1}^{n}(n-i)\deg(v_i;V)} (v_1\otimes\ldots\otimes 
v_n)[n],\qquad v[a]=1[a]\otimes v.\]
It is easy to verify that  $\decala_n$, called the \emph{
displacing}\footnote{It is often used the french name 
\emph{d\'ecalage}.} 
isomorphism,
changes symmetric into skewsymmetric
tensors  and therefore it induces an
isomorphism
\[ \decala_{n}\colon \symm{}{n}(V[1])\to (\external{}{n}V)[n],\]
\[
\decala_n(v_1[1]\odot\ldots\odot v_n[1])=
(-1)^{\sum_{i=1}^{n}(n-i)\deg(v_i;V)} (v_1\wedge\ldots\wedge v_n)[n].\]

\bigskip
\section{DG-coalgebras and $L_{\infty}$-algebras}

\begin{defn}\label{XIV.3.1}
By a  dg-coalgebra we intend a triple $(C,\Delta,d)$, where
$(C,\Delta)$ is a
graded coassociative cocommutative
coalgebra  and $d\in\Coder^{1}(C,C)$ is a codifferential.
If $C$ has a
counit $\epsilon\colon C\to\K$, we assume that $\epsilon d=0$.
The category of dg-coalgebras, where morphisms are morphisms of
coalgebras commuting with codifferentials, is denoted by $\mathbf{DGC}$.
\end{defn}

\begin{ex}\label{XIV.3.2}
If $A$ is a finite dimensional dg-algebra with
differential $d\colon A\to A[1]$, then $A\dual$
(Example~\ref{XIII.2.16})
is a dg-coalgebra with
codifferential the transpose of $d$.
\end{ex}

\begin{lem}\label{XIV.3.3} Let $V$ be a graded vector space and
$Q\in\Coder^{1}(\bar{S(V)},\bar{S(V)})$.
Then $Q$ is a codifferential, i.e. $Q\circ Q=0$, if and only if
for every $n>0$ and every $v_{1},\ldots,v_{n}\in V$
\[\sum_{k+l=n+1}\,\sum_{\sigma\in S(k,n-k)}\!\!
\!\!\!\epsilon(\sigma;v_{1},\ldots,v_{n})
Q^{1}_{l}(Q^{1}_{k}(v_{\sigma(1)}\!\odot\ldots\odot v_{\sigma(k)})\!\odot
v_{\sigma(k+1)}\!\odot\ldots\odot v_{\sigma(n)})=0.\]
\end{lem}

\begin{proof} Denote $P=Q\circ Q=\frac{1}{2}[Q,Q]$: since $P$ is a
coderivation we have that $P=0$ if and only if $P^{1}=Q^{1}\circ Q=0$.
According to Corollary~\ref{XIII.3.10}
\[ Q(v_{1}\odot \ldots\odot v_{n})=
\sum_{I\subset\{1,\ldots,n\}}\epsilon
\binom{I,I^{c}}{\{1,\ldots,n\}}Q^{1}(v_{I})\odot v_{I^{c}}\]
and then
\[ P^{1}(v_{1}\odot \ldots\odot v_{n})=
\sum_{I\subset\{1,\ldots,n\}}\epsilon
\binom{I,I^{c}}{\{1,\ldots,n\}}Q^{1}(Q^{1}(v_{I})\odot v_{I^{c}}).\]
\end{proof}

Note that $P^{1}_{n}=0$ whenever $Q^{1}_{m}=0$ for every
$m\ge \ds\frac{n+1}{2}$ and,
if $Q$ is a codifferential in $\bar{S(V)}$ then
$Q^1_1$ is a differential in the graded vector space $V$.

\begin{defn}\label{XIV.3.4}
Let $V$ be a graded vector space; a codifferential of degree 1 on
the symmetric coalgebra $C(V)=\bar{S(V[1])}$ is called an
\emph{$L_{\infty}$-structure} on $V$. The dg-coalgebra $(C(V),Q)$ is called an
\emph{$L_{\infty}$-algebra}.\\
An $L_{\infty}$-algebra $(C(V),Q)$ is called \emph{minimal} if $Q^{1}_{1}=0$.
\end{defn}

\begin{defn}\label{XIV.3.5} A \emph{weak morphism}
$F\colon (C(V),Q)\to (C(W),R)$
of $L_{\infty}$-algebras is a morphism of dg-coalgebras.
By an $L_{\infty}$-morphism we always  intend a weak morphism of
$L_{\infty}$-algebras.\\
A weak morphism $F$ is called a \emph{strong morphism}
if there exists a
morphism of graded vector spaces $F^{1}_{1}\colon V\to W$ such that
$F=S(F^{1}_{1})$.\\
We denote by $\mathbf{L_{\infty}}$ the category having
$L_{\infty}$-algebras as objects and (weak)
$L_{\infty}$-morphisms as arrows.\end{defn}

Consider now two $L_{\infty}$-algebras
$(C(V),Q)$, $(C(W),R)$ and a morphism of graded coalgebras
$F\colon C(V)\to C(W)$. Since
$FQ-RF\in \Coder^{1}(C(V),C(W);F)$, we have that $F$ is an
$L_{\infty}$-morphism if and only if $F^{1}Q=R^{1}F$.

\begin{lem}\label{XIV.3.6}
Consider  two $L_{\infty}$-algebras
$(C(V),Q)$, $(C(W),R)$  and
a morphism of graded vector spaces
$F^{1}\colon C(V)\to W[1]$. Then
\[F=S(F^{1})\frac{e^{\copl}-1}{\copl}\colon (C(V),Q)\to(C(W),R)\]
is an $L_{\infty}$-morphism if and only if
\begin{equation}\label{equlinf}
\sum_{i=1}^{n}R^{1}_{i}F^{i}_{n}=\sum_{i=1}^{n}F^{1}_{i}Q^{i}_{n}
\end{equation}
for every $n>0$.\\
\end{lem}

\begin{proof}
According to Proposition~\ref{XIII.2.32} $F$ is a morphism
of coalgebras. Since
$FQ-RF\in \Coder^{1}(C(V),C(W);F)$, we have that $F$ is an
$L_{\infty}$-morphism if and only if $F^{1}Q=R^{1}F$.\\
\end{proof}

\begin{exer} An $L_{\infty}$-morphism $F$ is strong if and only if
$F^{1}_{n}=0$ for every $n\ge 2$.\end{exer}

If $F\colon (C(V),Q)\to(C(W),R)$ is an $L_{\infty}$-morphism, then by
Lemma~\ref{XIV.3.6} $R^{1}_{1}F^{1}_{1}=F^{1}_{1}Q^{1}_{1}$ and therefore
we have a  morphism in cohomology
$H(F^{1}_{1})\colon H^{*}(V[1],Q^{1}_{1})\to H^{*}(W[1],Q^{1}_{1})$.

\begin{defn}\label{XIV.3.7} An $L_{\infty}$-morphism
$F\colon (C(V),Q)\to(C(W),R)$ is a \emph{ quasiisomorphism} if
$H(F^{1}_{1})\colon H^{*}(V[1],Q^{1}_{1})\to H^{*}(W[1],Q^{1}_{1})$
is an isomorphism.\end{defn}

The following exercise shows that the above definition is not
ambiguous.

\begin{exer} An $L_{\infty}$-morphism
$F\colon (C(V),Q)\to (C(W),R)$ is a quasiisomorphism if and only if
$H(F)\colon H^{*}(C(V),Q)\to H^{*}(C(W),R)$
is an isomorphism.\end{exer}

Given a  coderivation $Q\colon \bar{S(V[1])}\to \bar{S(V[1])}[1]$, their
components $Q^1_j\colon \bigodot^{n}(V[1])\to V[2]$, composed with the
inverse of the displacement isomorphism, give linear maps
\[ l_j=(Q^1_j\circ \decala_n^{-1})[-n]\colon \external{}{n}V\to V[2-n].\]
More explicitly
\[ l_j(v_1\wedge\ldots\wedge v_n)=(-1)^{-n}
(-1)^{\sum_{i=1}^{n}(n-i)\deg(v_i;V)}
Q^1_j(v_1[1]\odot\ldots\odot v_n[1])\]
The conditions of Lemma~\ref{XIV.3.3}  become
\[\sum_{\genfrac{}{}{0pt}{1}{k+i=n+1}{\sigma\in S(k,n-k)}(-1)^{k(i-1)}}
(-1)^{\sigma}\epsilon(\sigma)
l_{i}(l_{k}(v_{\sigma(1)}\wedge\ldots\wedge v_{\sigma(k)})\wedge
v_{\sigma(k+1)}\wedge\ldots\wedge v_{\sigma(n)})=0.\]

Setting $l_1(v)=d(v)$ and $l_2(v_1\wedge v_2)=[v_1,v_2]$, the first three
conditions ($n=1,2,3$) becomes:

{\renewcommand{\arraystretch}{1.5}
\[\begin{array}{ll}
1:&d^2=0\\
2:&d[x,y]=[dx,y]+(-1)^{\bar{x}}[x,dy]\\
3:&(-1)^{\bar{x}\,\bar{z}}[[x,y],z]+(-1)^{\bar{y}\,\bar{z}}[[z,x],y]+
(-1)^{\bar{x}\,\bar{y}}[[y,z],x]=\\
&=(-1)^{\bar{x}\,\bar{z}+1}
(dl_3(x,y,z)+l_3(dx,y,z)+(-1)^{\bar{x}}l_3(x,dy,z)+
(-1)^{\bar{x}+\bar{y}}l_3(x,y,dz))
\end{array}
\]}

If $l_{3}=0$ we recognize, in the three formulas above, 
the axioms defining a differential graded 
Lie algebra structure on $V$.

\begin{exer}
Let $(C(V),Q)$ be an  $L_{\infty}$-algebra. Then the bracket
\[ [w_{1},w_{2}]=(-1)^{\deg(w_1;V)}Q^{1}_{2}(
w_{1}[1]\odot  w_{2}[1])=l_2(w_1\wedge w_2)\]
gives  a structure of graded Lie
algebra on the cohomology of the complex $(V,Q^{1}_{1})$.
\end{exer}

\bigskip
\section{From DGLA to  $L_{\infty}$-algebras}

In this section we show that to every DGLA structure on a graded 
vector space $V$ it is associated naturally a $L_{\infty}$ structure on 
the same space $V$, i.e. a codifferential $Q$ on 
$C(V)=\bar{S(V[1])}$.\\
The coderivation $Q$ is determined by its components $Q^{1}_{j}\colon 
\bigodot^{j}V[1]\to V[2]$.

\begin{prop}\label{XIV.3.8}
Let $(V,d,[\, ,\,])$ be a differential graded Lie algebra.
Then the coderivation $Q$ of components
\begin{enumerate}
    \item $Q^{1}_{1}(v[1])=-d(v)$. 

    \item  $Q^{1}_{2}(w_{1}[1]\odot  w_{2}[1])
=(-1)^{\deg(w_{1};V)}[w_{1},w_{2}]$

    \item  $Q^{1}_{j}=0$ for every $j\ge 3$.
\end{enumerate}
is a codifferential and then 
gives an $L_{\infty}$-structure on $V$.
\end{prop}

\begin{proof}
The conditions of Lemma~\ref{XIV.3.3} are trivially satisfied for 
every $n>3$. For $n\le 3$ they becomes (where
$\widehat{x}=x[1]$ and $\bar{x}=\deg(x;V)$): 
{\renewcommand{\arraystretch}{1.5}
\[\begin{array}{ll}
n=1:&Q^{1}_{1}Q^{1}_{1}(\widehat{v})=d^2(v)=0\\
n=2:& Q^{1}_{1}Q^{1}_{2}(\widehat{x}\odot 
\widehat{y})+Q^{1}_{2}(Q^{1}_{1}(\widehat{x})\odot
\widehat{y})+(-1)^{(\bar{x}-1)(\bar{y}-1)} Q^{1}_{2}(Q^{1}_{1}(\widehat{y})\odot
\widehat{x})=\\ 
&=-(-1)^{\bar{x}}(d[x,y]-[dx,y])+[x,dy]=0\\
n=3: &Q^{1}_{2}(Q^{1}_{2}(\widehat{x}\odot\widehat{y})\odot\widehat{z})+
(-1)^{\bar{x}-1}Q^{1}_{2}(\widehat{x}\odot Q^{1}_{2}(\widehat{y}\odot\widehat{z}))+\\
&+(-1)^{\bar{x}(\bar{y}-1)}Q^{1}_{2}(\widehat{y}\odot 
Q^{1}_{2}(\widehat{x}\odot\widehat{z}))=\\
&=(-1)^{\bar{y}}[[x,y],z]+(-1)^{\bar{y}-1}[x,[y,z]]+
(-1)^{(\bar{x}-1)\bar{y}}[x,[y,z]]=0\\
\end{array}
\]}
\end{proof}

It is also clear that every morphism of $DGLA$ $f\colon V\to W$
induces a strong morphism of the corresponding $L_{\infty}$-algebras
$S(f[1])\colon C(V)\to C(W)$. Therefore we get in this way a functor
\[ \mathbf{DGLA}\to \mathbf{L_{\infty}}\]
that preserves quasiisomorphisms.\\ 
This functor is faithful; the following
example, concerning differential graded Lie algebras arising from 
Gerstenhaber-Batalin-Vilkovisky algebras, shows that it is not fully
faithful.\\

Let $(A,\Delta)$ be a GBV-algebra 
(Section \ref{CAP:TRIVIALK}.\ref{sezioneVI.5}); we
have seen that  $(G[-1], [\, ,\, ],\Delta)$, where
\[ [a,b]=a\Delta(b)+
(-1)^{\deg(a,G[-1])} (\Delta(ab)-\Delta(a)b)\] is a differential
graded Lie algebra and then it makes sense to consider the
associated $L_{\infty}$-algebra
$(C(G[-1]),\delta)=(\bar{S(G)},\delta)$. The codifferential
$\delta$ is induced by the linear map of degree 1
$\delta^1=\Delta+Q\in\Hom_{\K}^1(\bar{S(G)},G)$, where $\delta^1_1=\Delta$ and 
  \[ \delta^1_2=Q\colon \symm{}{2}G\to G,\qquad
  Q(a\odot b)=\Delta(ab)-\Delta(a)b- (-1)^{\bar{a}}a\Delta(b)\]

\begin{lem}\label{contazzo}
In the notation above,
\[
\begin{split}
\Delta(a_{1}a_{2}\ldots a_{m})=&\sum_{\sigma\in S(1,m-1)}
\epsilon(\sigma; a_{1},\ldots,a_{m})\Delta(a_{\sigma_1})
a_{\sigma_2}\ldots
a_{\sigma_m}+\\
&+\!\!\!\!\!\!\sum_{\sigma\in S(2,m-2)}
\!\!\!\!\!\!\!\!\epsilon(\sigma; a_{1},\ldots,a_{m})
Q(a_{\sigma_1},a_{\sigma_2})a_{\sigma_3}\ldots a_{\sigma_m}
\end{split}
\]
for every $m\ge 2$ and every $a_{1},\ldots,a_{m}\in G$.\
\end{lem}
\begin{proof}
For $m=2$ the above equality becomes
\[\Delta(ab)=\Delta(a)b+(-1)^{\bar{a}}a\Delta(b)+Q(a\odot b)\]
which is exactly the definition of $Q$.\\
By induction on $m$ we may assume the Lemma true for all integers
$<m$ and then
\[ \Delta((a_1  a_2)a_3\ldots
a_m)=\sum_{i=1}^{m}(-1)^{\bar{a_1}+\ldots+
\bar{a_{i-1}}}a_1\ldots\Delta(a_i)a_{i+1}\ldots a_m+\]
\[ +\sum_{i\ge 3}\epsilon\,
Q(a_1a_2\odot a_i)a_3\ldots \widehat{a_i}\ldots a_m+ \sum_{2<i<j}
\epsilon\, Q(a_i\odot a_j)a_1a_2\ldots \widehat{a_i}\ldots
\widehat{a_j}\ldots a_m.\] 
Replacing the odd Poisson
identity 
\[ Q(a_1a_2\odot a_i)=(-1)^{\bar{a_1}}a_1Q(a_2\odot a_i)+
(-1)^{(\bar{a_1}+1)\bar{a_2}}a_2Q(a_1\odot a_i)\] 
in the above formula, we obtain the desired equality.
\end{proof}

As an immediate consequence we have 

\begin{thm}\label{GBVareabelian}
In the notation above,
let  $(C(G[-1]),\tau)$ be the (abelian) $L_{\infty}$-algebra
whose codifferential is induced by $\Delta\colon G\to G$.
Then the
morphism of graded vector spaces $f\colon \bar{S(G)}\to G$,
\[ f(a_1\odot\ldots\odot a_m)=a_1a_2\ldots a_m\]
induces an isomorphism of $L_{\infty}$-algebras $F\colon
(C(G[-1]),\delta)\to(C(G[-1]),\tau)$.\end{thm}

\begin{proof}
According to Lemmas \ref{XIV.3.6} and \ref{contazzo} the morphism
of graded coalgebras induced by $f$ is an $L_{\infty}$-morphism.\\
Moreover, according to Example \ref{XIII.2.24} $F$ is an
isomorphism of graded coalgebras whose inverse is induced  by
\[ g\colon \bar{S(G)}\to G,\qquad
g(a_1\odot\ldots\odot a_m)=(-1)^{m-1}a_1a_2\ldots a_m.\]
\end{proof}

\bigskip

\section{From $L_{\infty}$-algebras to predeformation functors}

Let $Q\in\Coder^{1}(C(V),C(V))$ be a $L_{\infty}$ structure on a
graded vector space $V$, we define the Maurer-Cartan  functor
$MC_{V}\colon\NA\to\mathbf{Set}$
by setting:
\[ MC_{V}(A)=\Hom_{\mathbf{DGC}}(A\dual, C(V)).\]
We first note that the natural isomorphism
\[(C(V)\otimes A)^{0}=\Hom_{\mathbf{G}}(A\dual, C(V)),\qquad
(v\otimes a)c=c(a)v\]
is an isomorphism of algebras and then, according to Corollary~\ref{XIII.2.34},
every coalgebra morphism $\theta\colon A\dual \to C(V)$ is written
uniquely as $\theta=\exp(m)-1$ for some
$m\in (V[1]\otimes A)^{0}=\Hom_{\mathbf{G}}(A\dual, V[1])$.
As in Lemma~\ref{XIV.3.6}, $\theta$ is a morphism of dg-coalgebras if and
only if $m d_{A\dual}=Q^{1}\theta$; considering $m$ as an element of
the algebra
$(C(V)\otimes A)^{0}$ this equality becomes the \emph{Maurer-Cartan
equation} of an $L_{\infty}$-structure:
\[(Id_{V[1]}\otimes d_{A})m=
\sum_{n=1}^{\infty}\frac{1}{n!}(Q^{1}_{n}\otimes Id_{A})m^{n},\qquad
m\in (V[1]\otimes A)^{0}.\]
Via the d\'ecalage isomorphism the Maurer-Cartan equation becomes
\[Id_{V}\otimes d_{A}(m)=
\sum_{n=1}^{\infty}\frac{1}{n!}(-1)^{\frac{n(n+1)}{2}}
(l_{n}\otimes Id_{A})m\wedge\ldots\wedge m,\qquad
m\in (V\otimes A)^{1}.\]
It is then clear that if the $L_{\infty}$ structure comes from a DGLA
$V$ (i.e.
$l_{n}=0$ for every $n\ge 3$) then the Maurer-Cartan equation reduces
to the classical one.\\

It is evident that $MC_{V}$ is a covariant functor and
$MC_{V}(0)=0$. Let $\alpha\colon A\to C$, $\beta\colon B\to C$ be
morphisms in $\NA$, then
\[ MC_{V}(A\times_{C}B)=MC_{V}(A)\times_{MC_{V}(C)}MC_{V}(B)\]
and therefore $MC_{V}$ satisfies condition 2) of Definition~\ref{VIII.5.2}; in
particular it makes sense the tangent space $TMC_{V}$.

\begin{prop}\label{XIV.4.1}
The functor $MC_{V}$ is a predeformation functor with
$T^{i}MC_{V}=H^{i-1}(V[1],Q^{1}_{1})$.
\end{prop}
\begin{proof} If $A\in\NA\cap\mathbf{DG}$ then
\[MC_{V}(A)=\{ m\in (V\otimes A)^{1}\,| \,
Id_{V}\otimes d_{A}(m)=-
l_{1}\otimes Id_{A}(m)\}=Z^{1}(V\otimes A)\]
the same computation of \ref{VIII.5.8} shows that there exists a
natural isomorphism
$T^{i}MC_{V}=H^{i}(V,l_{1})=H^{i-1}(V[1],Q^{1}_{1})$.\\
Let $0\mapor{}I\mapor{}A\mapor{}B\mapor{}0$ be a small acyclic
extension in $\NA$, we want to prove that $MC_{V}(A)\to MC_{V}(B)$ is
surjective.\\
We have a dual exact sequence
\[0\mapor{}B\dual\mapor{}A\dual\mapor{}I\dual\mapor{}0,\qquad
B\dual=I^{\perp}.\]
Since $IA=0$ we have $\Delta_{A\dual}(A\dual)\subset B\dual\otimes
B\dual$.\\
Let $\phi\in MC_{V}(B)$ be a fixed element and $\phi^{1}\colon B\dual\to
V[1]$; by Proposition~\ref{XIII.2.32} $\phi$ is uniquely
determined by $\phi^{1}$. Let $\psi^{1}\colon A\dual\to V[1]$ be an
extension of $\phi^{1}$, then, again by \ref{XIII.2.32},
$\psi^{1}$ is induced by a unique
morphism of coalgebras $\psi\colon A\dual\to C(V)$.\\
The map $\psi d_{A\dual}-Q\psi\colon A\dual\to C(V)[1]$ is a
coderivation and then,
setting $h=(\psi d_{I\dual}-Q\psi)^{1}\in \Hom_{\mathbf{G}}(I\dual,
V[2])$,  we have that $\psi$ is a
morphism of dg-coalgebras if and only if $h=0$.\\
Note that $\psi^{1}$ is defined up to elements of
$\Hom_{\mathbf{G}}(I\dual, V[1])=
(V[1]\otimes I)^{0}$ and, since
$\Delta_{A\dual}(A\dual)\subset B\dual\otimes
B\dual$, $\psi^{i}$ depends only by $\phi$ for every $i>1$.
Since $I$ is acyclic and $h d_{I\dual}+Q^{1}_{1}h=0$ there exists
$\xi\in \Hom_{\mathbf{G}}(I\dual, V[1])$ such that
$h=\xi d_{I\dual}-Q^{1}_{1}\xi$ and then $\theta^{1}=\psi^{1}-\xi$
induces a
dg-coalgebra morphism $\theta\colon A\dual\to C(V)$ extending $\phi$.
\end{proof}

Therefore the Maurer-Cartan functor can be considered as a functor
$\mathbf{L_{\infty}}\to \mathbf{PreDef}$ that preserves quasiisomorphisms.
We have already noted that the composition
$\mathbf{DGLA}\to\mathbf{L_{\infty}}\to \mathbf{PreDef}$
is the Maurer-Cartan functor of  DGLAs.

\bigskip

\section{From predeformation to deformation functors}

We first recall the basics of homotopy theory of dg-algebras.

We denote by $\K[t_1,\ldots ,t_n,dt_1,\ldots ,dt_n]$ the dg-algebra of polynomial
differential forms on the affine space $\A^n$ with the de Rham
differential.
We have $\K[t,dt]=\K[t]\oplus\K[t]dt$ and
\[
\K[t_1,\ldots ,t_n,dt_1,\ldots ,dt_n]=\tensor{i=1}{n}\K[t_i,dt_i].
\]
Since $\K$ has characteristic 0, it is immediate to see that
$H_{*}(\K[t,dt])=\K[0]$ and then by K\"{u}nneth formula
$H_*(\K[t_1,\ldots ,t_n,dt_1,\ldots ,dt_n])=\K[0]$.
Note that for every dg-algebras $A$ and every
$s=(s_1,\ldots ,s_n)\in\K^n$
we have an evaluation morphism
\[e_s\colon
A\otimes\K[t_1,\ldots ,t_n,dt_1,\ldots ,dt_n]\to A\]
defined by
\[
e_s(a\otimes p(t_1,\ldots ,t_n,dt_1,\ldots ,dt_n))=p(s_1,\ldots,s_n,0,\ldots,0)a
\]

For every dg-algebra $A$ we denote $A[t,dt]=A\otimes\K[t,dt]$; if $A$ is
nilpotent then $A[t,dt]$ is still nilpotent.
If $A\in\NA$, then $A[t,dt]$ is the direct limit of objects in $\NA$.
To see this it is sufficient to consider,  for every
positive real number $\epsilon>0$, the dg-subalgebra
\[
A[t,dt]_\epsilon=A\oplus \oplus_{n>0} (A^{\lceil n\epsilon\rceil}
t^n\oplus A^{\lceil n\epsilon\rceil}t^{n-1}dt)\subset A[t,dt],
\]
where $A^{\lceil n\epsilon\rceil}$ is the subalgebra generated by all
the products $a_{1}a_{2}\ldots a_{s}$, $s\ge n\epsilon$, $a_{i}\in A$.

It is clear that if $A\in \NA$ then $A[t,dt]_\epsilon\in \NA$ for every
$\epsilon>0$ and $A[t,dt]$ is the union of all $A[t,dt]_\epsilon
$, $\epsilon>0$.

\begin{lem}\label{XIV.5.1}
For every dg-algebra $A$ the evaluation map
$e_{h}\colon A[t,dt]\to A$ induces an isomorphism
$H(A[t,dt])\to H(A)$ independent from  $h\in\K$.\end{lem}

\begin{proof} Let $\imath\colon A\to A[t,dt]$ be the inclusion, since
$e_{h}\imath=Id_{A}$ it is sufficient to prove that $\imath\colon
H(A)\to H(A[t,dt])$ is bijective.\\
For every $n>0$ denote $B_{n}=At^{n}\oplus At^{n-1}dt$; since
$d(B_{n})\subset B_{n}$ and $A[t,dt]=\imath(A)\bigoplus_{n>0}B_{n}$
it is sufficient to prove that $H(B_{n})=0$ for every $n$.
Let $z\in Z_{i}(B_{n})$, $z=at^{n}+nbt^{n-1}dt$, then
$0=dz=dat^{n}+((-1)^{i}a+db)nt^{n-1}dt$ which implies $a=(-1)^{i-1}db$
and then $z=(-1)^{i-1}d(bt^{n})$.
\end{proof}

\begin{defn}\label{XIV.5.2}
Given two morphisms of dg-algebras $f,g\colon A\to B$, a homotopy
between $f$ and $g$ is a morphism $H\colon A\to B[t,dt]$ such
that $H_0:=e_0\circ H=f$, $H_1:=e_1\circ H=g$ (cf. \cite[p.
120]{GrMor}).\\
We denote by $[A,B]$ the quotient of $\Hom_{\mathbf{DGA}}(A,B)$ by the
equivalence relation $\sim$ generated by homotopies.
\end{defn}

According to Lemma~\ref{XIV.5.1}, homotopic morphisms induce the same
morphism in homology.

\begin{lem}\label{XIV.5.3}
Given morphisms of dg-algebras,
\[\xymatrix{A\ar@/^/[r]^{f}\ar@/_/[r]_{g}&B\ar@/^/[r]^{h}\ar@/_/[r]_{l}&C},\]
if $f\sim g$ and $h\sim l$ then $hf\sim lg$.\end{lem}

\begin{proof} It is obvious from the definitions that $hg\sim lg$. For every
$a\in\K$ there exists a commutative diagram
\[\xymatrix{B\otimes\K[t,dt]\ar[r]^{h\otimes Id}\ar[d]^{e_{a}}&
C\otimes\K[t,dt]\ar[d]^{e_{a}}\\ B\ar[r]^{h}&C}.\]
If $F\colon A\to B[t,dt]$ is a homotopy between $f$ and $g$, then,
considering the composition of $F$ with $h\otimes Id$, we get a
homotopy between $hf$ and $hg$.\end{proof}

Since composition respects homotopy equivalence we can also consider
the homotopy categories
$K(\mathbf{DGA})$ and $K(\NA)$.
By definition, the objects of
$K(\mathbf{DGA})$ (resp.: $K(\NA)$) are the same of $\mathbf{DGA}$
(resp.: $\NA$), while the morphisms are
$Mor(A,B)=[A,B]$.

If $A,B\in \mathbf{DG}\cap \NA$, then two morphisms $f,g\colon A\to B$ are
homotopic in the sense of \ref{XIV.5.2} if and only if $f$ is homotopic to
$g$ as morphism of complexes.
In particular every acyclic complex is contractible
as a dg-algebra.

\begin{lem}\label{XIV.5.4}
A predeformation functor $F\colon \NA\to \mathbf{Set}$
is a deformation functor if and only if
$F$ induces a functor $[F]\colon K(\NA)\to \mathbf{Set}$.
\end{lem}

\begin{proof}
One implication is trivial, since every acyclic $I\in\NA\cap
\mathbf{DG}$ is isomorphic to 0 in $K(\NA)$.\\
Conversely,  let $H\colon A\to B[t,dt]$ be a homotopy, we need to
prove that $H_0$ and $H_1$ induce the same morphism from $F(A)$ to
$F(B)$.  Since $A$ is finite-dimensional there exists $\epsilon>0$
sufficiently small such that $H\colon A\to B[t,dt]_{\epsilon}$; now
the evaluation map $e_0\colon B[t,dt]_{\epsilon}\to B$ is a finite
composition of acyclic small extensions and then, since $F$ is a
deformation functor $F( B[t,dt]_{\epsilon})=F(B)$. For every $a\in
F(A)$ we have $H(a)=iH_0(a)$, where $i\colon B\to B[t,dt]_{\epsilon}$
is the inclusion and then $H_1(a)=e_1H(a)=e_1iH_0(a)=H_0(a)$.\\
\end{proof}

\begin{thm}\label{XIV.5.5}

Let $F$ be a predeformation functor, then there exists a deformation
functor $F^+$ and a natural transformation $\eta\colon F\to F^+$ such
that:
\begin{enumerate}
    \item  $\eta$ is a quasiisomorphism.

    \item For every deformation functor $G$ and every natural transformation
$\phi\colon F\to G$ there exists  a unique natural transformation
$\psi\colon F^+\to G$ such that $\phi=\psi\eta$.
\end{enumerate}
\end{thm}

\begin{proof}
We first define a functorial relation $\sim$ on the sets $F(A)$, $A\in
\NA$; we set $a\sim b$ if and only if there exists $\epsilon>0$ and
$x\in F(A[t,dt]_\epsilon)$ such that $e_0(x)=a$, $e_1(x)=b$. By
\ref{XIV.5.4}
if $F$ is a deformation functor then $a\sim b$ if and only if $a=b$.
Therefore if we define $F^+$ as the quotient of $F$ by the equivalence
relation generated by $\sim$ and $\eta$ as the natural projection,
then there exists a unique $\psi$ as in the
statement of the theorem. We only need to prove that $F^+$ is a
deformation
functor.

\begin{step}{1} If $C\in\mathbf{DG}\cap\NA$ is acyclic  then $F^+(C)=\{0\}$.\\
Since $C$ is acyclic there exists a homotopy $H\colon C\to
C[t,dt]_\epsilon$, $\epsilon\le 1$,
such that $H_0=0$, $H_1=Id$; it is then clear that
for every $x\in F(C)$ we have $x=H_1(x)\sim H_0(x)=0$.
\end{step}

\begin{step}{2} $\sim$ is an equivalence relation on $F(A)$ for every
$A\in \NA$.\\
This is essentially standard  (see e.g. \cite{GrMor}).
In view of the inclusion $A\to A[t,dt]_\epsilon$ the relation $\sim$ is
reflexive. The symmetry is proved by remarking that the
automorphism of dg-algebras
\[ A[t,dt]\to A[t,dt];\qquad a\otimes p(t,dt)\mapsto a\otimes p(1-t,-dt)\]
preserves the subalgebras $A[t,dt]_\epsilon$ for every $\epsilon>0$.

Consider now $\epsilon>0$ and $x\in F(A[t,dt]_\epsilon)$, $y\in
F(A[s,ds]_\epsilon)$ such that $e_0(x)=e_0(y)$; we need to prove that
$e_1(x)\sim e_1(y)$.

Write $\K[t,s,dt,ds]=\oplus_{n\ge 0}S^n$, where $S^n$ is the $n$-th
symmetric power of the acyclic complex $\K t\oplus\K s\mapor{d}\K
dt\oplus\K ds$ and define $A[t,s,dt,ds]_\epsilon=A\oplus
\oplus_{n>0} (A^{\lceil n\epsilon\rceil}\otimes S^n)$. There exists a
commutative diagram

\[\begin{array}{ccc}
A[t,s,dt,ds]_\epsilon&\mapor{t\mapsto 0}&A[s,ds]_\epsilon\\
\mapver{s\mapsto 0}&&\mapver{s\mapsto 0}\\
A[t,dt]_\epsilon&\mapor{t\mapsto 0}&A\\
\end{array}\]

The kernel of the surjective morphism
\[ A[t,s,dt,ds]_\epsilon \mapor{\eta} A[t,dt]_\epsilon\times_A
A[t,dt]_\epsilon\]
is equal to $\oplus_{n>0}(A^{\lceil n\epsilon\rceil}\otimes
(S^n\cap I))$, where $I\subset \K[t,s,dt,ds]$ is the homogeneous
differential ideal generated by $st,sdt,tds, dtds$. Since $I\cap S^n$ is
acyclic for every $n>0$, the morphism $\eta$ is a finite
composition of acyclic small extensions.

Let $\xi\in F(A[t,s,dt,ds]_\epsilon)$ be a lifting of $(x,y)$ and let
$z\in F(A[u,du]_\epsilon)$ be the image of $\xi$ under the morphism
\[ A[t,s,dt,ds]_\epsilon\to A[u,du]_\epsilon,\qquad
t\mapsto 1-u,\quad s\mapsto u\]
The evaluation of $z$ gives $e_0(z)=e_1(x)$, $e_1(z)=e_1(y)$.
\end{step}

\begin{step}{3}
If $\alpha\colon A\to B$ is surjective then
\[ F(A[t,dt]_\epsilon)
\mapor{(e_0,\alpha)}F(A)\times_{F(B)}F(B[t,dt]_\epsilon)\]
is surjective.\\

It is not restrictive to assume $\alpha$ a small extension with kernel
$I$. The kernel of $(e_0,\alpha)$ is equal to
$\oplus_{n>0}(A^{\lceil n\epsilon\rceil}\cap I)\otimes (\K t^n\oplus \K
t^{n-1}dt)$ and therefore  $(e_0,\alpha)$ is an acyclic small extension.
\end{step}

\begin{step}{4}
The functor $F^+$ satisfies \ref{schlessi1} of \ref{VIII.5.2}.\\

Let $a\in F(A)$, $b\in F(B)$ be such that $\alpha(a)\sim\beta(b)$; by
Step 3 there exists $a'\sim a$, $a'\in F(A)$ such that
$\alpha(a')=\beta(b)$ and then the pair $(a',b)$ lifts to $F(A\times_C
B)$.
\end{step}

\begin{step}{5}
The functor $F^+$ satisfies \ref{schlessi2} of \ref{VIII.5.2}.\\

By  \ref{VIII.5.3} it is sufficient to verify the condition separately for
the cases $C=0$ and $B=0$.
When $C=0$ the situation is easy: in fact $(A\times B)[t,dt]_\epsilon=
A[t,dt]_\epsilon\times B[t,dt]_\epsilon$,
$F((A\times B)[t,dt]_\epsilon)=
F(A[t,dt]_\epsilon)\times F(B[t,dt]_\epsilon)$ and the relation $\sim$
over $F(A\times B)$ is the product of the relations $\sim$ over $F(A)$
and $F(B)$; this implies that $F^+(A\times B)=F^+(A)\times F^+(B)$.

Assume now $B=0$, then the fibred product $D:=A\times_C B$ is equal to
the kernel of $\alpha$. We need to prove that the map $F^+(D)\to F^+(A)$
is injective.
Let $a_0,a_1\in F(D)\subset F(A)$ and let $x\in F(A[t,dt]_\epsilon)$ be
an element such that $e_i(x)=a_i$, $i=0,1$. Denote by $\bar{x}\in
F(C[t,dt]_\epsilon)$ the image of $x$ by $\alpha$.

Since $C$ is acyclic there exists a morphism of graded vector spaces
$\sigma\colon C\to C[-1]$ such that $d\sigma+\sigma d=Id$ and we can
define a morphism of complexes
\[ h\colon C\to (\K s\oplus\K ds)\otimes C\subset C[s,ds]_1;\qquad
h(v)=s\otimes v+ds\otimes\sigma(v)\]
The morphism $h$ extends in a natural way to a morphism
\[ h\colon C[t,dt]_\epsilon\to (\K s\oplus\K ds)\otimes C[t,dt]_\epsilon\]
such that for every scalar $\zeta\in\K$  there exists a commutative diagram

\[\begin{array}{ccc}
C[t,dt]_\epsilon&\mapor{h}&(\K s\oplus\K ds)\otimes C[t,dt]_\epsilon\\
\mapver{e_{\zeta}}&&\mapver{Id\otimes e_{\zeta}}\\
C&\mapor{h}&(\K s\oplus\K ds)\otimes C
\end{array}\]

Setting $\bar{z}=h(\bar{x})$ we have $\bar{z}_{|s=1}=\bar{x}$,
$\bar{z}_{|s=0}=\bar{z}_{|t=0}=\bar{z}_{|t=1}=0$. By step 3 $\bar{z}$
lifts to an element $z\in F(A[t,dt]_\epsilon[s,ds]_1)$ such that
$z_{|s=1}=x$; Now the specializations ${z}_{|s=0}$, ${z}_{|t=0}$,
${z}_{|t=1}$ are annihilated by $\alpha$ and therefore give a chain of
equivalences in $F(D)$
\[ a_0=z_{|s=1, t=0}\sim z_{|s=0, t=0}\sim z_{|s=0, t=1}\sim z_{|s=1,
t=1}=a_1\]
proving that $a_0\sim a_1$ inside $F(D)$.
\end{step}

The combination of
Steps 1, 4 and 5 tell us that $F^+$ is a deformation functor.\\

\begin{step}{6}
The morphism $\eta\colon F\to F^{+}$ is a quasiisomorphism.\\
Let $\epsilon$ be of degree $1-i$, $\epsilon^{2}=0$, then
$\K\epsilon\oplus I_{i}$ is isomorphic to the dg-subalgebra
\[ \K\epsilon\oplus\K\epsilon t\oplus\K\epsilon dt\subset
\K\epsilon[t,dt]\]
and the map $p\colon F(I_{i})\to F(\K\epsilon)$  factors as
\[ p\colon F(I_{i})\hookrightarrow F(I_{i})\oplus F(\K\epsilon)=
F(\K\epsilon\oplus\K\epsilon t\oplus\K\epsilon dt)\xrightarrow{e_{1}-e_{0}}
F(\K\epsilon).\]
On the other hand the evaluation maps $e_{0}, e_{1}$ factor as
\[ e_{i}\colon \K\epsilon[t,dt]\mapor{h}
\K\epsilon\oplus\K\epsilon t\oplus\K\epsilon
dt\mapor{e_{i}}\K\epsilon,\qquad i=0,1\]
where $h$ is the morphism of dg-vector spaces
\[ h(\epsilon
t^{n+1})=\epsilon t,\quad h(\epsilon t^{n}dt)=\frac{\epsilon
dt}{n+1},\quad \forall n\ge 0. \]
\end{step}

\end{proof}

\begin{cor}\label{XIV.5.6}
Let $L$ be a differential graded Lie algebra, then there
exists a natural isomorphism $MC_{L}^{+}=\Def_{L}$.\end{cor}

\begin{proof} According to Theorem~\ref{XIV.5.5}
there exists a natural morphism of functors
$\psi\colon MC_{L}^{+}\to\Def_{L}$;
by  \ref{VIII.5.8} $\psi$ is a quasiisomorphism and then, by
Corollary~\ref{VIII.6.4} $\psi$ is an isomorphism.\end{proof}

\begin{defn}\label{XIV.5.7}
Let $(C(V),Q)$ be a $L_{\infty}$-algebra and let $\Def_{V}=MC_{V}^{+}$
be the deformation functor associated to the predeformation functor
$MC_{V}$ . We shall call $\Def_{V}$ the
\emph{deformation functor associated to the $L_{\infty}$-algebra
$(C(V),Q)$}.
\end{defn}

A morphism of $L_{\infty}$-algebras $C(V)\to C(W)$ induces in the obvious
way a natural transformation $MC_{V}\to MC_{W}$ and then,
according to \ref{XIV.5.5}, a morphism
$\Def_{V}\to \Def_{W}$. Finally, since $MC_{V}\to \Def_{V}$ is a
quasiisomorphism we have $T^{i}\Def_{V}=H^{i}(V,Q^{1}_{1})$.

The following result is clear.

\begin{cor}\label{XIV.5.8}
Let $\theta\colon C(V)\to C(W)$ be a morphism of $L_{\infty}$-algebras.
The induced morphism $\Def_{V}\to \Def_W$ is an isomorphism if and only
if $\theta^1_1\colon V\to W$ is a quasiisomorphism of complexes.
\end{cor}

\bigskip
\section{Cohomological constraint to deformations of 
K\"ahler manifolds}

Theorem \ref{GBVareabelian} shows that the category of 
$L_{\infty}$-algebras is more flexible than the category of differential 
graded Lie algebras. Another example in this direction is given by 
the main theorem of \cite{ManettiCCKM}.

Let $X$ be a fixed compact K\"ahler manifold of dimension $n$ and
consider the graded vector space
$M_{X}=\Hom_{\C}^{*}(H^{*}(X,\C),H^{*}(X,\C))$
of linear endomorphisms of the singular cohomology of $X$. The
Hodge decomposition gives natural isomorphisms
\[ M_{X}=\somdir{i}{}M^{i}_{X},\qquad
M^{i}_{X}=\somdir{r+s=p+q+i}{}\Hom_{\C}
(H^{p}(\Omega^{q}_{X}),H^{r}(\Omega^{s}_{X}))\]
and the composition of the cup product and the contraction
operator $T_{X}\otimes\Omega^{p}_{X}\mapor{\vdash}\Omega^{p-1}_{X}$
gives natural linear maps
\[ \theta_{p}\colon H^{p}(X,T_{X})\to
\somdir{r,s}{}\Hom_{\C}^{*}
(H^{r}(\Omega^{s}_{X}),H^{r+p}(\Omega^{s-1}_{X}))
\subset M[-1]_{X}^{p}=M^{p-1}_{X}.\]
By Dolbeault's theorem
$H^{*}(KS_{X})=H^{*}(X,T_{X})$ and then the maps $\theta_{p}$ give
a morphism of graded vector spaces $\theta\colon H^{*}(KS_{X})\to
M[-1]_{X}$. This morphism is  generally nontrivial: consider for
instance a Calabi-Yau manifold where the map $\theta_{p}$ induces
an isomorphism
$H^{p}(X,T_{X})=\Hom_{\C}(H^{0}(\Omega^{n}_{X}),H^{p}(\Omega^{n-1}_{X}))$.

\begin{thm}\label{thmA} In the above notation,
consider $M[-1]_{X}$ as a differential graded Lie algebra with trivial
differential and trivial bracket.\\
Every choice of a K\"ahler metric
on $X$ induces a canonical lifting of $\theta$ to an $L_{\infty}$-morphism
from $KS_{X}$ to $M[-1]_{X}$.
\end{thm}

The application of Theorem \ref{thmA} to 
deformation theory, see \cite{ManettiCCKM},  are based on the idea 
that $L_{\infty}$-morphisms 
induce natural transformations of (extended) deformation functors commuting with tangential 
actions and obstruction maps (cf. Theorem \ref{VIII.6.2}). Being the deformation functor of the 
DGLA $M[-1]$ essentially trivial, the lifting of $\theta$ impose several constraint on 
deformations of $X$.

Denote by:\begin{itemize}
\item $A^{*,*}=\bigoplus_{p,q}A^{p,q}$, where
$A^{p,q}=\Gamma(X,\sA^{p,q})$ the vector space of global $(p,q)$-forms.
\item $N^{*,*}=\Hom_{\C}^{*}(A^{*,*},A^{*,*})=\bigoplus_{p,q}N^{p,q}$,
where $N^{p,q}=\bigoplus_{i,j}\Hom_{\C}^{*}(A^{i,j},A^{i+p,j+q})$ is
the space of homogeneous endomorphisms of $A^{*,*}$ of bidegree $(p,q)$.
\end{itemize}

The space $N^{*,*}$, endowed with the composition product and total
degree $\deg(\phi)=p+q$ whenever $\phi\in N^{p,q}$, is a graded associative
algebra and therefore,  with the standard bracket
\[ [\phi,\psi]=\phi\psi-(-1)^{\deg(\phi)\deg(\psi)}\psi\phi\]
becomes a graded Lie algebra. We note that the adjoint operator
$[\debar,~]\colon N^{*,*}\to N^{*,*+1}$ is a differential
inducing a structure of DGLA.\\

\begin{lem}\label{corotau}
Let $X$ be a compact K\"ahler manifold, then there exists
$\tau\in N^{1,-1}$ such that:\begin{enumerate}
\item $\tau$ factors to a linear map $A^{*,*}/\ker{\de}\to \Image\de$.
\item $[\debar,\tau]=\de$.
\end{enumerate}
In particular $\de\in N^{1,0}$ is a coboundary  in the
DGLA $(N^{*,*}, [~,~], [\debar,~])$.
\end{lem}
\begin{proof} In the notation of Theorem \ref{X.4.9} it is sufficient to consider
$\tau=\sigma\de=-\de\sigma$.
Note that  the above $\tau$ is defined
canonically from the choice of the K\"ahler metric.\end{proof}

We fix a
K\"ahler metric on $X$ and denote by: $\sH\subset A^{*,*}$ the
graded vector space of harmonic forms, $i\colon \sH\to A^{*,*}$
the inclusion and $h \colon A^{*,*}\to \sH$ the harmonic
projector.

We identify the graded vector space $M_X$ with the space of
endomorphisms of harmonic forms $\Hom_{\C}^{*}(\sH,\sH)$. We also
we identify $\Der^{*}(\sA^{*,*},\sA^{*,*})$ with its image into
$N=\Hom_{\C}^{*}(A^{*,*},A^{*,*})$.\\

According to Lemma~\ref{corotau} there exists $\tau\in N^{0}$
such that
\[ h\de=\de h=\tau h=h\tau=\de\tau=\tau\de=0,\qquad
[\debar,\tau]=\de.\]

For
simplicity of notation we denote by $L=\oplus L^{p}$ the
$\Z$-graded vector space $KS[1]_{X}$, this means that
$L^{p}=\Gamma(X,\sA^{0,p+1}(T_{X}))$, $-1\le p\le n-1$. The local
description of the  two linear maps of degree $+1$, 
$d\colon L\to
L$, $Q\colon \odot^{2}L\to L$ 
introduced, up to d\'ecalage, in Proposition~\ref{XIV.3.8} is: if $z_{1},\ldots,z_{n}$ are
local holomorphic coordinates, then
\[ d\left(\phi\desude{~}{z_{i}}\right)=
(\debar\phi)\desude{~}{z_{i}},\qquad \phi\in \sA^{0,*}.\]
If $I,J$ are ordered subsets of $\{1,\ldots,n\}$,
$a=fd\bar{z}_{I}\desude{~}{z_{i}}$,
$b=gd\bar{z}_{J}\desude{~}{z_{j}}$, $f,g\in\sA^{0,0}$ then
\[ Q(a\odot b)=(-1)^{\bar{a}}d\bar{z}_{I}\wedge d\bar{z}_{J}
\left(f\desude{g}{z_{i}}\desude{~}{z_{j}}-
g\desude{f}{z_{j}}\desude{~}{z_{i}}\right),\qquad \bar{a}=\deg(a,L).\]

The formula 
\begin{equation}\label{codiff}
\begin{split}
\delta(a_{1}\odot\ldots\odot a_{m})=&\sum_{\sigma\in S(1,m-1)}
\epsilon(L,\sigma; a_{1},\ldots,a_{m})da_{\sigma_1}\odot
a_{\sigma_2}\odot\ldots\odot
a_{\sigma_m}+\\
&+\!\!\!\!\!\!\sum_{\sigma\in S(2,m-2)}
\!\!\!\!\!\!\!\!\epsilon(L,\sigma; a_{1},\ldots,a_{m})Q(a_{\sigma_1}
\odot a_{\sigma_2})\odot a_{\sigma_3}\odot\ldots\odot
a_{\sigma_m}
\end{split}
\end{equation}
gives a
codifferential $\delta$ of degree 1 on $\bar{S}(L)$ and the
differential graded coalgebra $(\bar{S}(L),\delta)$ is exactly the
$L_{\infty}$-algebra associated to the Kodaira-Spencer DGLA $KS_{X}$.\\

If $\Der^{p}(\sA^{*,*},\sA^{*,*})$ denotes the  vector space of
$\C$-derivations of degree $p$ of the sheaf of graded algebras
$(\sA^{*,*},\wedge)$, where the degree of a $(p,q)$-form is $p+q$
(note that $\de,\debar\in
\Der^{1}(\sA^{*,*},\sA^{*,*})$), then  we have a morphism of
graded vector spaces
\[ L\mapor{\hat{~}}\Der^{*}(\sA^{*,*},\sA^{*,*})=\somdir{p}{}
\Der^{p}(\sA^{*,*},\sA^{*,*})
,\qquad a\mapsto\h{a}\] 
given in local coordinates by
\[ \h{\phi\desude{~}{z_{i}}}(\eta)=\phi\wedge\left(\desude{~}{z_{i}}
\vdash\eta\right).\]

\begin{lem} If $[~,~]$ denotes the standard bracket on
$\Der^{*}(\sA^{*,*},\sA^{*,*})$, then for every pair of homogeneous
$a,b\in L$ we have:
\begin{enumerate}
    \item  $\h{da}=[\debar,\h{a}]=\debar\h{a}-(-1)^{\bar{a}}
    \h{a}\debar.$

    \item  $\h{Q(a\odot b)}=-[[\de,\h{a}],\h{b}]=
    (-1)^{\bar{a}}\h{a}\de\h{b}+(-1)^{\bar{a}\,\bar{b}+\bar{b}}\,\h{b}\de\h{a}
    \pm \de\h{a}\h{b}\pm\h{b}\h{a}\de$.
\end{enumerate}
\end{lem}
\begin{proof} This is a special case of Lemma 
\ref{quadratic}.\end{proof}

Consider  the morphism
\[ F_{1}\colon L\to M_X, \qquad F_{1}(a)=h\h{a}i.\]
We note that $F_{1}$ is a morphism of complexes, in fact
$F_{1}(da)=h\h{da}i=h(\debar\h{a}\pm\h{a}\debar)i=0$. By
construction $F_{1}$ induces the morphism $\theta$ in cohomology
and therefore the theorem is proved whenever we lift $F_1$ to a
morphism of graded vector spaces $F\colon \bar{S}(L)\to M_{X}$
such that $F\circ\delta=0$.\\

Define, for every $m\ge 2$, the following morphisms of graded
vector spaces
\[f_{m}\colon \tensor{}{m}L\to M_X,\qquad
F_{m}\colon \symm{}{m}L\to M_X,\qquad
F=\sum_{m=1}^{\infty}F_{m}\colon \bar{S}(L)\to M_X,\]
\[ f_{m}(a_{1}\otimes a_{2}\otimes\ldots\otimes a_{m})=
h\h{a_{1}}\tau\h{a_{2}}\tau\h{a_{3}}\ldots\tau\h{a_{m}}i.\]
\[ F_{m}(a_{1}\odot a_{2}\odot\ldots\odot a_{m})=
\sum_{\sigma\in\Sigma_{m}}\epsilon(L,\sigma;a_{1},\ldots,a_{m})
f_{m}(a_{\sigma_{1}}\otimes\ldots\otimes a_{\sigma_{m}}).\]

\begin{thm}\label{1.2} In the above notation $F\circ\delta=0$ and therefore
\[\Theta=\sum_{m=1}^{\infty}\frac{1}{m!}F^{\odot m}\circ\Delta^{m-1}_{C(KS_X)}
\colon (C(KS_X),\delta)\to (C(M[-1]_X),0)\]
is an $L_{\infty}$-morphism with linear term $F_{1}$.\end{thm}

\begin{proof} We need to prove that
for every $m\ge 2$ and $a_{1},\ldots,a_{m}\in L$ we
have
\[F_{m}\left(\sum_{\sigma\in S(1,m-1)}\epsilon(L,\sigma)
da_{\sigma_{1}}\odot a_{\sigma_{2}}\odot \ldots\odot
a_{\sigma_{m}}\right)=\]
\[=-F_{m-1}\left(\sum_{\sigma\in S(2,m-2)}\epsilon(L,\sigma)
Q(a_{\sigma_{1}}\odot a_{\sigma_{2}})\odot a_{\sigma_{3}}\odot \ldots\odot
a_{\sigma_{m}}\right),\]
where $\epsilon(L,\sigma)=\epsilon(L,\sigma;a_{1},\ldots,a_{m})$.\\
It is convenient to introduce the auxiliary operators
$q\colon \bigotimes^{2}L\to N[1]$, $q(a\otimes
b)=(-1)^{\bar{a}}\h{a}\de\h{b}$ and
$g_{m}\colon \bigotimes^{m}L\to M[1]_X$,
\[ g_{m}(a_{1}\otimes \ldots\otimes a_{m})=
-\sum_{i=0}^{m-2}(-1)^{\bar{a_{1}}+\bar{a_{2}}+\ldots+\bar{a_{i}}}
h\h{a_{1}}\tau\ldots\h{a_{i}}\tau q(a_{i+1}\otimes
a_{i+2})\tau\h{a_{i+3}}\ldots\tau\h{a_{m}}i.\]
Since for every choice of operators
$\alpha=h,\tau$ and $\beta=\tau,i$ and every
$a,b\in L$ we have
\[ \alpha\h{Q(a\odot b)}\beta=\alpha((-1)^{\bar{a}}\h{a}\de\h{b}+
(-1)^{\bar{a}\,\bar{b}+\bar{b}}\h{b}\de\h{a})\beta=
\alpha(q(a\otimes b)+(-1)^{\bar{a}\,\bar{b}}q(b\otimes a))\beta,
\]
the symmetrization lemma \ref{syle} gives
\[\sum_{\sigma\in\Sigma_{m}}\epsilon(L,\sigma)
g_{m}(a_{\sigma_{1}}\otimes\ldots\otimes a_{\sigma_{m}})=
-F_{m-1}\left(\sum_{\sigma\in S(2,m-2)}\epsilon(L,\sigma)
Q(a_{\sigma_{1}}\odot a_{\sigma_{2}})\odot a_{\sigma_{3}}\odot \ldots\odot
a_{\sigma_{m}}\right).\]
On the other hand
\renewcommand{\arraystretch}{2.7}
\[ \begin{array}{l}
\ds f_{m}\left(\sum_{i=0}^{m-1}(-1)^{\bar{a_{1}}+\ldots+\bar{a_{i}}}
a_{1}\otimes\ldots\otimes a_{i}\otimes
da_{i+1}\otimes\ldots\otimes a_{m}\right)=\\
\ds =\sum_{i=0}^{m-1}(-1)^{\bar{a_{1}}+\ldots+\bar{a_{i}}}
h\h{a_{1}}\ldots\h{a_{i}}\tau(\debar
\h{a_{i+1}}-(-1)^{\bar{a_{i+1}}}\h{a_{i+1}}\debar)\tau\ldots \tau
\h{a_{m}}i\\
\ds =\sum_{i=0}^{m-2}(-1)^{\bar{a_{1}}+\ldots+\bar{a_{i}}}
h\h{a_{1}}\ldots\h{a_{i}}\tau
(-(-1)^{\bar{a_{i+1}}}\h{a_{i+1}}\debar\tau\h{a_{i+2}}+
(-1)^{\bar{a_{i+1}}}\h{a_{i+1}}\tau\debar\h{a_{i+2}})
\tau\ldots \tau
\h{a_{m}}i\\
\ds =-\sum_{i=0}^{m-2}(-1)^{\bar{a_{1}}+\ldots+\bar{a_{i}}}
h\h{a_{1}}\ldots\h{a_{i}}\tau
((-1)^{\bar{a_{i+1}}}\h{a_{i+1}}[\debar,\tau]\h{a_{i+2}})
\tau\ldots \tau
\h{a_{m}}i\\
\ds =-\sum_{i=0}^{m-2}(-1)^{\bar{a_{1}}+\ldots+\bar{a_{i}}}
h\h{a_{1}}\ldots\h{a_{i}}\tau q(a_{i+1}\otimes a_{i+2}) \tau\ldots
\tau
\h{a_{m}}i\\
\ds =g_{m}(a_{1}\otimes\ldots\otimes a_{m}).\end{array}\]
\renewcommand{\arraystretch}{1}
Using again Lemma~\ref{syle} we have
\[ \sum_{\sigma\in\Sigma_{m}}\epsilon(L,\sigma)
g_{m}(a_{\sigma_{1}}\otimes\ldots\otimes a_{\sigma_{m}})=
F_{m}\left(\sum_{\sigma\in S(1,m-1)}\epsilon(L,\sigma)
da_{\sigma_{1}}\odot a_{\sigma_{2}}\odot \ldots\odot
a_{\sigma_{m}}\right).\]
\qed\end{proof}

\medskip
~\\
\emph{Remark.} If $X$ is a Calabi-Yau manifold with holomorphic volume
form $\Omega$, then the composition of $F$ with the evaluation at
$\Omega$ induces an $L_{\infty}$-morphism $C(KS_X)\to C(\sH[n-1])$.\\
For every $m\ge 2$, $\operatorname{ev}_{\Omega}\circ F_{m}\colon
\bigodot^{m}L\to \sH[n]$ vanishes on
$\bigodot^{m}\{a\in L\mid
\de(a\vdash\Omega)=0\}$.
\medskip
~\\

\bigskip

\section{Historical survey, \ref{CAP:TOOLS}}

$L_{\infty}$-algebras, also called \emph{strongly homotopy Lie algebras},
are  the Lie analogue of the $A_{\infty}$ ( strongly homotopy 
associative algebras), introduced by Stasheff \cite{Sta} in the 
context of algebraic topology.\\
The popularity of $L_{\infty}$-algebras has been increased recently by 
their application in deformation theory (after \cite{SS1}), in 
deformation quantization (after \cite{K}) and in string theory (after 
\cite{WZ}, cf. also \cite{LadaStas}).


\begin{thebibliography}{99}


\bibitem{Artin68} M.~Artin: \emph{On the solutions of analytic equations.}
Invent. Math. \textbf{5} (1968) 277-291.

\bibitem{Artinbook} M. Artin: \emph{Deformations of singularities.}
Tata Institute  of Fundamental Research, Bombay (1976).

\bibitem{A-M} M.F.~Atiyah, I.G.~Macdonald:
\emph{Introduction to commutative algebra.}
Addison-Wesley, Reading, Mass. (1969).

\bibitem{Ba-Sta} C. Banica, O. Stanasila: \emph{M\'ethodes
alg\'ebrique dans la th\'eorie globale des espaces complexes.}
Gau\-thi\-er-Villars (1977).

\bibitem{BPV} W.~Barth, C.~Peters, A.~van de Ven:
\emph{Compact complex surfaces.}
Springer-Verlag Ergebnisse Math. Grenz. \textbf{4} (1984).

\bibitem{Boualge} N. Bourbaki: \emph{Algebra 1}.

\bibitem{Canonaco} A. Canonaco: \emph{
$L_{\infty}$-algebras and Quasi-Isomorphisms.}
In: \emph{Seminari di Geometria Algebrica 1998-1999} Scuola
Normale Superiore (1999).


\bibitem{cartan} H. Cartan: \emph{Elementary theory of analytic 
functions of one or several complex variables.} 
Addison-Wesley (1963)


\bibitem{Catacime} F. Catanese:
\emph{Moduli of algebraic surfaces.} Springer L.N.M. \textbf{1337} (1988)
1-83.

\bibitem{clemens} H. Clemens: \emph{Cohomology and obstructions, I: on the
geometry of formal Kuranishi theory.} preprint math.AG/9901084.

\bibitem{DGMS} P.~Deligne, P.~Griffiths, J.~Morgan, D.~Sullivan:
\emph{Real homotopy theory of K\"ahler manifolds.} Invent. Math. 
\textbf{29} (1975) 245-274.

\bibitem{Douady} A.~Douady: {\em
Obstruction primaire \`a la d\'eformation.}
S\'em.~Cartan \textbf{13}  (1960/61)  Exp.~4.

\bibitem{FaMa1} B. Fantechi, M. Manetti:
\emph{Obstruction calculus for functors of Artin
rings, I} Journal of Algebra \textbf{202} (1998) 541-576.

\bibitem{FaMa2} B. Fantechi, M. Manetti: \emph{On the $T^1$-lifting theorem.}
J. Alg. Geom. \textbf{8} (1999) 31-39.

\bibitem{FHT} Y. F\'elix, S. Halperin, J.C. Thomas: \emph{
Rational homotopy theory.} Springer GTM \textbf{205} (2001).

\bibitem{Fischer} G.~Fischer: \emph{Complex analytic geometry.}
Springer-Verlag LNM \textbf{538} (1976).

\bibitem{Gerst} M. Gerstenhaber: \emph{The cohomology structure of an
associative ring.} Ann. of Math. \textbf{78} (1963) 267-288.

\bibitem{Ge-Sch} M. Gerstenhaber, S.D. Schack: \emph{On the
deformation of algebra morphisms and diagrams.} T.A.M.S. \textbf{279} (1983)
1-50.

\bibitem{Gode} R.~Godement:
\emph{Topologie alg\'ebrique et th\'eorie des faisceaux.}
Hermann, Paris (1958).

\bibitem{GoMil1} W.M. Goldman, J.J. Millson:
\emph{The deformation theory of
representations of fundamental groups of compact k\"{a}hler manifolds}
Publ. Math. I.H.E.S. \textbf{67} (1988) 43-96.

\bibitem{GoMil2} W.M. Goldman, J.J. Millson:
\emph{The homotopy invariance of the
Kuranishi space} Ill. J. Math. \textbf{34} (1990) 337-367.

\bibitem{grassi1} M. Grassi: \emph{$L_{\infty}$-algebras and
differential graded algebras, coalgebras and Lie algebras.}
In: \emph{Seminari di Geometria Algebrica 1998-1999} Scuola
Normale Superiore (1999).

\bibitem{Ana-Ste} H. Grauert, R. Remmert: \emph{Analytische
Stellenalgebren.} Spinger-Verlag Grundlehren \textbf{176} (1971).

\bibitem{CAS} H. Grauert, R. Remmert: \emph{Coherent Analytic sheaves.}
Spinger-Verlag Grundlehren \textbf{265} (1984).

\bibitem{green} M.L. Green: \emph{Infinitesimal methods in Hodge theory.}
In \emph{Algebraic cycles and Hodge theory (Torino, 1993)},
1-92, Lecture Notes in Math. \textbf{1594}, Springer (1994)

\bibitem{G-H} P.~Griffiths, J.~Harris: \emph{Principles of Algebraic Geometry.}
Wiley-Interscience publication (1978).

\bibitem{GrMor} P.H. Griffiths, J.W Morgan: \emph{Rational
Homotopy Theory and Differential Forms} Birkh\"{a}user Progress
in Mathematics \textbf{16} (1981).

\bibitem{Gu-Ro} R. Gunning, H. Rossi: \emph{Analytic functions of several complex
variables.} Prenctice-Hall (1965).

\bibitem{Ho} G. Hochschild: \emph{The structure of Lie groups.}
Holden-Day San  Francisco (1965).

\bibitem{Horidhm} E. Horikawa: \emph{Deformations of holomorphic maps I.}
J. Math. Soc. Japan \textbf{25} (1973) 372-396; II J. Math. Soc. Japan 
\textbf{26} (1974) 647-667; 
III Math. Ann. \textbf{222} (1976) 275-282.

\bibitem{Hum} J.E.~Humphreys: \emph{Introduction to Lie algebras and
representation theory.} Springer-Verlag (1972).

\bibitem{Illu} L. Illusie: \emph{Complexe cotangent et d\'eformations I 
et II.} Springer LNM \textbf{239} (1971) and \textbf{283} (1972).

\bibitem{Ja} N. Jacobson: \emph{Lie algebras.} Wiley \& Sons (1962).

\bibitem{Kawa} Y.~Kawamata: \emph{Unobstructed deformations - a remark
on a paper of Z.~Ran.} J.~Algebraic Geom. \textbf{1} (1992) 183-190.

\bibitem{Kobabook} S.~Kobayashi: \emph{Differential geometry of complex
vector bundles.} Princeton Univ. Press (1987).

\bibitem{Kod1} K. Kodaira: \emph{On stability of compact submanifolds of complex
manifolds.} Amer. J. Math. \textbf{85} (1963) 79-94.

\bibitem{Kobook} K. Kodaira:
\emph{Complex manifold and deformation of complex
structures.} Springer-Verlag (1986).

\bibitem{KNS} K. Kodaira, L. Nirenberg, D.C.~Spencer:
\emph{On the existence of deformations
of  complex analytic structures.} Annals of Math.
\textbf{68} (1958) 450-459.

\bibitem{KS1} K. Kodaira, D.C.~Spencer:
\emph{On the variation of almost complex structures.}
In \emph{Algebraic geometry and topology,} Princeton Univ. Press
(1957) 139-150.

\bibitem{KS2} K. Kodaira, D.C.~Spencer:
\emph{A theorem of completeness for complex analytic fibre spaces.}
Acta Math. \textbf{100} (1958) 281-294.

\bibitem{KS3} K. Kodaira, D.C.~Spencer:
\emph{On deformations of complex analytic structures, I-II, III.}
Annals of Math. \textbf{67} (1958) 328-466; \textbf{71} (1960) 43-76.

\bibitem{K3} M. Kontsevich: \emph{Enumeration of rational curves via 
torus action.} In \emph{Moduli Space of Curves} (R. Dijkgraaf, C. 
Faber, G. van der Geer Eds) Birkh\"{a}user (1995) 335-368.

\bibitem{Konts} M. Kontsevich: \emph{Topics in algebra-deformation
theory.} Notes (1994):

\bibitem{K} M. Kontsevich: \emph{Deformation quantization of
Poisson manifolds, I.} {\tt q-alg/9709040}.

\bibitem{Kura} M.~Kuranishi: \emph{New Proof for the existence of
locally complete families of complex structures.} In: Proc. Conf.
Complex Analysis (Minneapolis 1964) Springer-Verlag (1965) 142-154.

\bibitem{LadaMarkl} T.~Lada, M.~Markl: \emph{Strongly homotopy Lie
algebras.} Comm. Algebra \textbf{23} (1995) 2147-2161, {\tt hep-th/9406095}.

\bibitem{LadaStas} T.~Lada, J.~Stasheff: {\em
Introduction to sh Lie algebras for physicists.}  
Int. J. Theor. Phys. \textbf{32} (1993) 1087-1104, {\tt hep-th/9209099}.

\bibitem{LZ} B.H. Lian, G.J. Zuckerman: \emph{New perspectives on the
BRST-algebraic structure of string theory.} Commun.  Math.  Phys. 
\textbf{154} (1993), 613Ð646, {\tt hep-th/9211072}.


\bibitem{MacLane} S.~Mac Lane: \emph{Categories for the working
mathematician.} Springer-Verlag (1971).

\bibitem{malgrange} B. Malgrange: \emph{Ideals of differentiable 
functions.} Tata Institute Research Studies \textbf{3} (1967).

\bibitem{IntroGA} M.~Manetti:
\emph{Corso introduttivo alla Geometria Algebrica.}
Appunti dei corsi tenuti dai docenti della Scuola Normale Superiore
(1998) 1-247.

\bibitem{ManettiDGLA} M. Manetti: \emph{Deformation theory via differential graded
Lie algebras.} In \emph{Seminari di Geometria Algebrica 1998-1999} Scuola
Normale Superiore (1999).

\bibitem{ManettiEDF} M. Manetti: \emph{Extended deformation functors.}
Internat. Math. Res. Notices \textbf{14} (2002) 719-756. {\tt 
math.AG/9910071}

\bibitem{ManettiCCKM} M. Manetti:
\emph{Cohomological constraint to deformations of compact
K\"ahler manifolds.} Adv. Math. \textbf{186} (2004) 125-142; 
{\tt math.AG/0105175}.

\bibitem{Matsu} H. Matsumura: \emph{On algebraic groups of birational
transformations.}  Rend. Accad. Lincei Ser. 8 \textbf{34} (1963)
151-155.

\bibitem{Matsubook} H. Matsumura: \emph{Commutative Ring Theory.}
Cambridge University Press (1986).

\bibitem{JJMil} J.J. Millson: \emph{Rational Homotopy Theory and 
Deformation Problems from Algebraic Geometry.} Proceedings of ICM, 
Kyoto 1990.

\bibitem{Nijenhuis} A. Nijenhuis: \emph{Jacobi-type identities for
bilinear differential concomitants of
certain tensor fields I} Indagationes Math. \textbf{17} (1955) 390-403.

\bibitem{Pala1} V.P. Palamodov: \emph{Deformations of complex spaces.}
Uspekhi Mat. Nauk. \textbf{31:3} (1976) 129-194.
Transl. Russian Math. Surveys \textbf{31:3} (1976) 129-197.

\bibitem{Pala2} V.P. Palamodov: \emph{Deformations of complex spaces.}
In: \emph{Several complex variables IV.}
Encyclopaedia of Mathematical Sciences \textbf{10}, Springer-Verlag
(1986) 105-194.

\bibitem{Qui} D.~Quillen: \emph{Rational homotopy theory.}
Ann. of Math. \textbf{90} (1969) 205-295.

\bibitem{Qui1} D.~Quillen: \emph{On the (co)homology of commutative
rings.} Proc. Sympos. Pure Math. \textbf{17} (1970) 65-87.

\bibitem{Ran} Z.~Ran: \emph{Deformations of
manifolds with torsion or negative canonical bundle.} J.~Algebraic
Geom. \textbf{1} (1992) 279-291.

\bibitem{ran2} Z. Ran: \emph{ Hodge theory and deformations of maps.}
Compositio Math. \textbf{97} (1995) 309-328.

\bibitem{Ran3} Z.~Ran: \emph{ Universal variations of Hodge structure and
Calabi-Yau-Schottky  relations.} Invent. Math. \textbf{138} (1999) 425-449.
{\tt math.AG/9810048}


\bibitem{Rim} D.S.~Rim: \emph{Formal deformation theory.} In \emph{SGA 7 I},
Exp. VI. Lecture Notes in Mathematics, \textbf{288}
Springer-Verlag (1972) 32-132.

\bibitem{Sch} M. Schlessinger: \emph{Functors of Artin rings.}
Trans. Amer. Math. Soc. \textbf{130} (1968) 208-222.

\bibitem{SS1} M.~Schlessinger, J.~Stasheff:
\emph{Deformation theory and rational homotopy type.}
preprint (1979).

\bibitem{SS2} M.~Schlessinger, J.~Stasheff: \emph{The Lie algebra
structure of tangent cohomology and deformation theory.}
J. Pure Appl. Algebra \textbf{38} (1985) 313-322.

\bibitem{Schouten} J.A.~Schouten: \emph{\"{U}ber Differentialkonkomitanten
zweier kontravarianter Gr\"{o}\ss en.} Indagationes
Math. \textbf{2} (1940), 449-452.

\bibitem{segreB} B. Segre: \emph{On arithmetical properties of quartics.}
Proc. London Math. Soc. \textbf{49} (1944) 353-395.

\bibitem{Shafarevich} I.R.~Shafarevich:
\emph{Basic algebraic geometry.} Springer-Verlag (1972).

\bibitem{FAC} J.P. Serre: \emph{Faisceaux alg\'ebriques coh\'erents.}
Ann. of Math. \textbf{61} (1955) 197-278.


\bibitem{Sta} J.D. Stasheff: \emph{On the homotopy associativity of H-spaces 
I,II.} Trans. AMS \textbf{108} (1963),
275-292, 293-312.


\bibitem{sta97} J.D. Stasheff: \emph{The (secret?) homological algebra of 
the Batalin-Vilkovisky
approach.} Secondary calculus and cohomological physics (Moscow, 1997),  
Contemp. Math. \textbf{219},
Amer. Math. Soc. (1998) 195-210, {\tt hep-th/9712157}.

\bibitem{Tian} G. Tian: \emph{Smoothness of the universal deformation
space of compact Calabi-Yau manifolds and its Peterson-Weil metric.}
in: S.T. Yau ed. \emph{Math. aspects of String Theory.}
Singapore (1988) 629-646.

\bibitem{Todorov} A.N.~Todorov: \emph{The Weil-Peterson geometry of the
moduli space of $SU(n\ge 3)$ (Calabi-Yau) manifolds I.} Commun. Math.
Phys. \textbf{126} (1989) 325-346.

\bibitem{Voisin} C.~Voisin: \emph{Th\'eorie de Hodge et g\'eom\'etrie
alg\'ebrique complexe} Soci\'et\'e Math\'ematique de France, Paris (2002).

\bibitem{Wawrik1} J.J.~Wawrik: \emph{Obstruction to the existence of a
space of moduli.} In \emph{Global Analysis} Princeton Univ. Press
(1969) 403-414.

\bibitem{Weil} A.~Weil: \emph{Introduction \`a l'\'etude des
vari\'et\'es K\"{a}hl\'eriennes.} Hermann Paris (1958).

\bibitem{Wells} R. O. Wells: \emph{Differential analysis on complex manifolds.}
Springer-Verlag (1980).

\bibitem{WZ} E. Witten, B. Zwiebach: \emph{Algebraic Structures And Differential 
Geometry In $2D$ String Theory.} 
{\tt hep-th/9201056}, Nucl. Phys. B377 (1992), 55-112.





\end{thebibliography}
\end{document}